\documentclass{memo-l}
\usepackage{amssymb}

\newtheorem{thm}{Theorem}[chapter]

\newtheorem{prop}[thm]{Proposition}
\newtheorem{cor}[thm]{Corollary}

\theoremstyle{definition}
\newtheorem{defn}[thm]{Definition}

\theoremstyle{remark}
\newtheorem{rema}[thm]{Remark}

\numberwithin{section}{chapter}
\numberwithin{equation}{chapter}

\newcommand{\asqrt}[0]{a_{\scriptscriptstyle{\square}}}
\newcommand{\bsqrt}[0]{b_{\scriptscriptstyle{\square}}}

\newcommand{\dtheta}[0]{\frac{\partial}{\partial
\theta}}
\newcommand{\dz}[0]{\frac{\partial}{\partial z}}
\newcommand{\Dz}[0]{\frac{\partial}{\partial
\theta} + \theta \frac{\partial}{\partial z}} 
\newcommand{\nor}[0]{\Upsilon}
\newcommand{\sou}[0]{\Delta}

\newcommand{\Z}{\mathbb{Z}_+} 
\newcommand{\Dx}[0]{\frac{\partial}{\partial \varphi} + \varphi
\frac{\partial}{\partial x}}
\newcommand{\Lx}[0]{x^{j + 1}
\frac{\partial}{\partial x} + \Bigl(\frac{j +
1}{2}\Bigr)\varphi x^j \frac{\partial}{\partial
\varphi}}
\newcommand{\twoLo}[0]{2x \frac{\partial}{\partial x} + \varphi
\frac{\partial}{\partial \varphi}}   
\newcommand{\Gx}[0]{x^j \biggl( \frac{\partial}{\partial \varphi} -
\varphi \frac{\partial}{\partial x}\biggr)} 
\newcommand{\Lminus}[0]{x^{- j + 1} \frac{\partial}{\partial x} +
\Bigl(\frac{- j + 1}{2} \Bigr) \varphi x^{-j} \frac{\partial}{\partial \varphi}}
\newcommand{\Gminus}[0]{x^{- j + 1}\biggl( \frac{\partial}{\partial
\varphi} - \varphi \frac{\partial}{\partial x}\biggr)}
\newcommand{\twoLow}[0]{2w \frac{\partial}{\partial w} + \rho
\frac{\partial}{\partial \rho}} 
\newcommand{\Lw}[0]{w^{j + 1} \frac{\partial}{\partial w} + \Bigl(\frac{j +
1}{2} \Bigr) \rho w^j \frac{\partial}{\partial \rho}} 
\newcommand{\Gw}[0]{w^j \biggl( \frac{\partial}{\partial \rho} -
\rho \frac{\partial}{\partial w}\biggr)}

\newcommand{\dw}[0]{\frac{\partial}{\partial w}}
\newcommand{\drho}[0]{\frac{\partial}{\partial \rho}}

\begin{document}

\frontmatter
\title{The moduli space of $N=1$ superspheres with tubes
and the sewing operation}

\author{Katrina Barron}
\address{Department of Mathematics, University of Notre Dame, 255 Hurley
Hall, Notre Dame, IN 46556}
\email{kbarron@nd.edu}

\maketitle

\setcounter{page}{4}

\tableofcontents

\begin{Large}
\begin{center}
{\bf Abstract}
\end{center}
\end{Large}

\medskip

Within the framework of complex supergeometry and motivated by 
two-dimen-sional genus-zero holomorphic $N = 1$ superconformal 
field theory, we define the moduli space of $N=1$ genus-zero
super-Riemann surfaces with oriented and ordered half-infinite 
tubes, modulo superconformal equivalence.  We define a sewing 
operation on this moduli space which gives rise to the sewing 
equation and normalization and boundary conditions.  To solve this 
equation, we develop a formal theory of infinitesimal $N = 1$ 
superconformal transformations based on a representation of the 
$N=1$ Neveu-Schwarz algebra in terms of superderivations.  We solve 
a formal version of the sewing equation by proving an identity for  
certain exponentials of superderivations involving infinitely many  
formal variables.  We use these formal results to give a
reformulation of the moduli space, a more detailed description of 
the sewing operation, and an explicit formula for obtaining a 
canonical supersphere with tubes {}from the sewing together of two 
canonical superspheres with tubes.  We give some specific examples 
of sewings, two of which give geometric analogues of associativity 
for an $N=1$ Neveu-Schwarz vertex operator superalgebra.  We study 
a certain linear functional in the supermeromorphic tangent space 
at the identity of the moduli space of superspheres with $1 + 1$ 
tubes (one outgoing tube and one incoming tube) which is 
associated to the $N=1$ Neveu-Schwarz element in an $N=1$ 
Neveu-Schwarz vertex operator superalgebra.  We prove the 
analyticity and convergence of the infinite series arising {}from 
the sewing operation.  Finally, we define a bracket on the 
supermeromorphic tangent space at the identity of the moduli  
space of superspheres with $1+1$ tubes and show that this gives 
a representation of the $N=1$ Neveu-Schwarz algebra with central 
charge zero.

\vfill

Received by editor July 17, 2000.

\subjclass{AMS Subject Classification: Primary 17B65, 17B68, 17B81, 32A05,
32C11, 58A50, 81R10, 81T40, 81T60; Secondary  17B69, 30F10, 32G15.}

\thanks{The author was supported in part by an AAUW American Dissertation Fellowship and
by a University of California President's Postdoctoral Fellowship}

\mainmatter

\chapter{Introduction}

In this monograph, we give a detailed study of the geometry underlying
two-dimensional genus-zero holomorphic $N=1$ superconformal field theory.  
Conformal field theory (or more specifically, string theory) and related 
theories (cf. \cite{BPZ}, \cite{FS}, \cite{V}, and \cite{S}) are the 
most promising attempts at developing a physical theory that combines 
all fundamental particle interactions, including gravity.  The geometry 
of conformal field theory extends the use of Feynman diagrams, 
describing the interactions of point particles whose propagation in 
time sweeps out a line in space-time, to one-dimensional ``particles'' 
(strings) whose propagation in time sweeps out a two-dimensional surface.  
For two-dimensional genus-zero holomorphic conformal field theory, these 
surfaces are genus-zero Riemann surfaces with incoming and outgoing 
half-infinite tubes, and algebraically, the corresponding string 
interactions can be described by products of vertex operators or, more 
precisely, by means of vertex operator algebras (cf. \cite{Bo}, 
\cite{FLM}).  In \cite{H book}, Huang studies the moduli space of 
genus-zero Riemann surfaces with tubes and provides a rigorous  
correspondence between the geometric and algebraic aspects of
two-dimensional genus-zero holomorphic conformal field theory.

In \cite{Fd}, Friedan describes the extension of the physical model of
conformal field theory to that of $N=1$ superconformal field theory 
and the notion of a superstring whose propagation in time sweeps out a
supersurface.  Whereas conformal field theory attempts to describe the
interactions of bosons, superconformal field theory attempts to
describe the interactions of bosons paired with $N$ fermions, for
$N =1,2,...$.  Thus $N = 1$ superconformal field theory describes the
interactions of boson-fermion pairs, $N = 2$ superconformal field
theory describes the interactions of boson-fermion-fermion triplets,
and so on.  In this work, we will consider $N = 1$ superconformal
field theory.  This theory requires an operator $D$ such that $D^2 =
\frac{\partial}{\partial z}$.  Such an operator arises naturally in
$N = 1$ complex supergeometry.  This is the geometry of manifolds
over a complex Grassmann algebra (i.e., complex supermanifolds) where
there is one ``even'' dimension corresponding to the one boson and one
``odd'' dimension corresponding to the $N = 1$ fermion in the
boson-fermion pairs.  

Thus the geometric setting for genus-zero holomorphic $N=1$ 
superconformal field theory is the moduli space of $N=1$ genus-zero
super-Riemann surfaces with oriented and ordered half-infinite tubes
(which are superconformally equivalent to $N=1$ superspheres with 
oriented and ordered punctures, and local superconformal coordinates 
vanishing at the punctures), modulo superconformal equivalence.  In 
this paper, following \cite{S}, \cite{V}, \cite{H book} and \cite{Fd},
we give a rigorous and detailed study of this moduli space.  We
define a sewing operation on this moduli space which gives rise to the 
sewing equation and normalization and boundary conditions.  Physically, 
one can view each equivalence class of superspheres in the moduli space 
as representing some superstring interaction and the sewing operation
as corresponding to taking the resulting superstring output of a given 
interaction as one of the inputs of another superstring interaction and 
determining the resulting total interaction.  To determine this 
resulting total interaction, one must solve the sewing equation.  To 
solve this equation, we develop a formal theory of infinitesimal $N = 1$  
superconformal transformations based on a representation of the $N=1$  
Neveu-Schwarz algebra in terms of superderivations.  We solve a formal 
version of the sewing equation and normalization and boundary 
conditions by proving an identity for certain exponentials of 
superderivations involving infinitely many formal variables.  We use 
these formal results to give a reformulation of the moduli space and a 
more detailed description of the sewing operation.  This includes an 
explicit formula for obtaining a canonical supersphere with tubes {}from 
the sewing together of two canonical superspheres with tubes which we 
obtain {}from our formal results by interpreting the formal solution to 
the sewing equation in terms of the analytic and geometric structure of 
the moduli space and proving the necessary convergence conditions for 
this formal solution.  
 
Although this monograph is self-contained as a study of the
geometry underlying two-dimensional genus-zero holomorphic $N=1$
superconformal field theory, we next mention an important application 
of these results to the algebraic aspects of this theory and the
correspondence between the algebraic and geometric aspects. 

In \cite{B vosas}, we extended the notion of vertex operator 
superalgebra (cf. \cite{T}, \cite{G}, \cite{FFR}, \cite{DL}, and 
\cite{KW}) to be defined over a Grassmann algebra instead of just 
$\mathbb{C}$ and to include odd formal variables instead of just even 
formal variables by introducing the notion of {\it $N=1$ Neveu-Schwarz
vertex operator superalgebra over a Grassmann algebra and with odd 
formal variables}.  We also proved that the category of $N=1$ 
Neveu-Schwarz vertex operator superalgebras over a Grassmann algebra 
and with odd formal variables is isomorphic to the category of 
$N=1$ Neveu-Schwarz vertex operator superalgebras over a Grassmann 
algebra and without odd formal variables.  However, in a vertex 
operator superalgebra with odd formal variables, the fact that the 
endomorphism $G(-1/2)$ plays the role of the superconformal operator 
$D = \frac{\partial}{\partial \theta} + \theta \frac{\partial}{\partial 
z}$ (as defined in Chapter 2 below) is made explicit, and the 
correspondence with the supergeometry is more natural.

Most of the results contained in this paper and in \cite{B vosas} were 
first proved in \cite{B thesis} and a research announcement of these 
results appeared in \cite{B thesis announcement}.  The main result of 
\cite{B thesis} is the introduction of the notion of {\it $N=1$ 
supergeometric vertex operator superalgebra over a Grassmann algebra} 
and the following isomorphism theorem:

\medskip

{\bf The Isomorphism Theorem:}  Projectively, the category of $N=1$ 
Neveu-Schwarz vertex operator superalgebras over a Grassmann algebra 
with (or without) odd formal variables and with central charge $c \in 
\mathbb{C}$ is isomorphic to the category of $N=1$ supergeometric 
vertex operator superalgebras over a Grassmann algebra and with central 
charge $c \in \mathbb{C}$.    

\medskip

This isomorphism theorem provides a rigorous foundation for the 
correspondence between the algebraic and geometric aspects of 
two-dimensional genus-zero holomorphic $N=1$ superconformal field
theory.  The results contained in this paper give the geometric 
development necessary to establish this correspondence, in particular, 
allowing one to introduce the notion of $N=1$ supergeometric vertex 
operator superalgebra, to prove the isomorphism theorem, and to construct 
the isomorphism explicitly.  In addition, in this monograph, we establish 
certain identities for representations of the $N=1$ Neveu-schwarz algebra
which arise {}from the geometry but pertain to the general algebraic 
setting.

The moduli space of $N=1$ superspheres with tubes and the sewing operation
is an example of what is known mathematically as a ``partial operad''.
The notion of operad was introduced by May \cite{M}, and in \cite{HL2}, 
\cite{HL3} it was shown that the moduli space of spheres with tubes 
together with the sewing operation defined in \cite{H thesis} is a 
partial operad.  Analogously, the moduli space of $N=1$ superspheres 
with tubes, together with the sewing operation and action of the 
symmetric groups defined and studied in this paper, is also a partial 
operad.  The isomorphism theorem can be reformulated as stating that 
projectively the category of $N=1$ Neveu-Schwarz vertex operator 
superalgebras with central charge $c \in \mathbb{C}$ is isomorphic to 
the category of supermeromorphic graded algebras with central charge $c 
\in \mathbb{C}$ over this $N=1$ supersphere with tubes partial operad  
(cf. \cite{H book}).  In this sense, one can think of the isomorphism  
theorem as associating an $n$-ary product on an infinite-dimensional  
graded vector space $V$ for every equivalence class of $N=1$ 
superspheres with one outgoing tube and $n \in \mathbb{N}$ incoming 
tubes in this partial operad, and showing that with a couple of 
additional axioms, these $n$-ary products define an $N = 1$ 
Neveu-Schwarz vertex operator superalgebra structure on $V$. 

The results in this paper are an extension of Huang's \cite{H thesis}, 
\cite{H book} work on the geometry underlying genus-zero holomorphic
conformal field theory.  In particular, this paper contains the extension 
of Chapters 1-4 in \cite{H book} to the $N=1$ super case.  Throughout this 
paper, one may consider the nonsuper case by setting all odd variables, all 
``soul'' Grassmann terms (see Section 2.1 for the definitions), and all odd 
superalgebra elements equal to zero.  In other words, this work 
subsumes the analogous nonsuper case studied in \cite{H thesis} and 
\cite{H book}.  In \cite{BMS}, Beilinson, Manin and Schechtman
study some aspects of $N = 1$ superconformal symmetry, i.e., the
$N=1$ Neveu-Schwarz algebra, {}from the viewpoint of algebraic
geometry.  In this work our approach is decidedly differential
geometric.

The rigorous foundation for the correspondence between the algebraic 
and geometric aspects of two-dimensional genus-zero holomorphic conformal 
and $N=1$ superconformal field theory developed in \cite{H thesis}, 
\cite{H book}, and \cite{B thesis} has proved useful in furthering both the 
algebraic and geometric aspects of conformal and superconformal field theory.  
On the one hand, it is much easier to rigorously construct conformal field  
theories and study many of their properties using the algebraic  
formulation of vertex operator (super)algebra (or equivalently, chiral  
algebra) (e.g., \cite{Wi}, \cite{BPZ}, \cite{Za}, \cite{FLM}, \cite{KT}, 
\cite{FZ}, \cite{FFR}, \cite{DL}, \cite{FF}, \cite{DMZ}, \cite{Wa}, 
\cite{Zh1}, \cite{Zh2}).  In addition, algebraic settings such as 
orbifold conformal field theory \cite{DHVW1}, \cite{DHVW2}, \cite{FLM} 
and coset models \cite{GKO} allow one to construct new conformal field 
theories {}from existing ones.  On the other hand, the geometry of 
(super)conformal field theory can give insight and provide tools useful 
for studying the algebraic aspects of the theory, for example: giving 
rise to general results in Lie theory \cite{BHL}; giving the necessary 
foundation for developing a theory of tensor products for vertex operator 
algebras \cite{HL1}, \cite{HL4} -- \cite{HL7}, \cite{H tensor}; giving 
rise to constructions in orbifold conformal field theory \cite{BDM}, 
\cite{H orbifold}; and providing a setting for establishing change of
variables formulas and other general identities for vertex operator
(super)algebras (cf., \cite{H book}, \cite{B thesis}, and Sections 3.6 
and 3.7 below).  Furthermore, the geometric formulation is indispensable 
for understanding both the geometric and algebraic aspects of higher 
genus theory. 

This paper is organized as follows.  In Chapter 2, within the framework 
of complex supergeometry (cf. \cite{D}, \cite{Bc1}, \cite{Ro2} and 
\cite{CR}) and motivated by two-dimensional genus-zero holomorphic $N = 1$  
superconformal field theory, we define the moduli space of $N=1$
genus-zero super-Riemann surfaces with oriented and ordered 
punctures, and local superconformal coordinates vanishing at the 
punctures, modulo superconformal equivalence.  We define a sewing operation 
on this moduli space by taking two superspheres, cutting out closed discs 
in the DeWitt topology around a puncture in each supersphere and 
appropriately identifying DeWitt neighborhoods of the boundaries.  This 
then results in a new supersphere with punctures and local superconformal 
coordinates vanishing at the punctures and is well defined on 
superconformal equivalence classes.  This sewing operation on the moduli 
space gives rise to the sewing equation and the normalization and boundary 
conditions.

In Chapter 3, we develop a formal theory of infinitesimal $N=1$
superconformal transformations based on a representation of the $N=1$
Neveu-Schwarz algebra in terms of superderivations.  We show that any 
local superconformal coordinate can be expressed in terms of 
exponentials  of these superderivations.  We then study the sewing 
equation {}from a purely algebraic viewpoint using this new 
characterization of local superconformal coordinates and give a formal 
solution to the sewing equation satisfying the normalization and 
boundary conditions by proving an identity for certain exponentials of 
superderivations involving infinitely many formal variables.  We 
formulate some additional algebraic identities which arise {}from 
certain sewings; prove analogues to the sewing identities for a 
general representation of the $N=1$ Neveu-Schwarz algebra, and show 
that the infinite series occurring in these identities have certain 
nice properties when the representation is a positive-energy module,
for instance, an $N=1$ Neveu-Schwarz vertex operator superalgebra
\cite{B vosas}.

In Chapter 4, we give a reformulation of the moduli space of $N=1$
superspheres with tubes using the results of Chapters 2 and 3 and give  
a detailed analytic study of the sewing operation in terms of 
exponentials of representatives of Neveu-Schwarz algebra elements as 
formulated algebraically in Chapter 3.  We define an action of the 
symmetric groups on the moduli space.  This action is needed to 
establish the fact that the moduli space of superspheres with tubes 
and the sewing operation is a partial operad.  We define 
supermeromorphic superfunctions and supermeromorphic tangent spaces for 
the moduli space.  These supermeromorphic superfunctions on the moduli 
space include as examples correlation functions for vertex operators in 
the algebraic theory of $N=1$ Neveu-Schwarz vertex operator 
superalgebras with odd formal variables (see \cite{B thesis}, 
\cite{B vosas}).  We study the group-like structures of the moduli space 
of superspheres with $1+1$ tubes (i.e., with one incoming tube and one 
outgoing tube), and give some specific examples of sewing, two of which 
give geometric analogues of associativity for an $N=1$ Neveu-Schwarz 
vertex operator superalgebra (see \cite{B thesis}, \cite{B vosas}).  We 
study a certain linear functional in the supermeromorphic tangent space 
at the identity of the moduli space of superspheres with $1+1$ tubes 
which is associated to the $N=1$ Neveu-Schwarz element in an $N=1$ 
Neveu-Schwarz vertex operator superalgebra.  We then give an explicit 
formula for obtaining a canonical supersphere with tubes {}from the 
sewing together of two canonical superspheres with tubes using the 
results of Chapters 2 and 3.  We prove the analyticity and convergence 
of certain series resulting {}from sewing which were obtained 
algebraically in Chapter 3.  

Finally, in Chapter 4 we show that there is a representation of the 
$N=1$ Neveu-Schwarz algebra with central charge zero on a subspace of 
the supermeromorphic tangent space of the moduli space of superspheres 
with $1+1$ tubes at the identity.  Thus, projectively, one can think of 
this moduli space as the partial monoid associated to the $N=1$ 
Neveu-Schwarz Lie superalgebra in analogy to the super Lie group 
associated to a finite-dimensional Lie superalgebra, (cf. \cite{Ro1}).

We have written this paper in such a way as to parallel the work of 
Huang \cite{H thesis}, \cite{H book} in order to accentuate the 
similarities and differences between the super and nonsuper cases.
The results in this monograph which require significant work or new
methods in extending the analogous results of \cite{H thesis}, 
\cite{H book} include: Proposition \ref{canonicalcriteria}, 
Proposition \ref{superconformal}, Theorem \ref{uniformization}, 
Proposition \ref{first Theta prop}, Proposition \ref{second Theta 
prop}, Proposition \ref{the linear functional G}, Proposition 
\ref{use of FG}, and Proposition \ref{NS bracket}.  In this 
monograph, we correct several misprints which appeared in \cite{B 
thesis}, \cite{B thesis announcement}, \cite{H thesis}, and \cite{H 
book}, and present a new proof of Proposition \ref{use of FG} which 
we believe to be more concise and more indicative of the role of 
the super structure in the analyticity and convergence of the infinite 
series arising {}from the sewing operation.

We would like to thank James Lepowsky and Yi-Zhi Huang for their
advice, support and expert comments during the writing of the 
dissertation {}from which this paper is derived, and their advice
on the presentation of this paper. 

\bigskip 

\begin{center}
{\bf Notational conventions}
\end{center}

\noindent
$\mathbb{C}$: the complex numbers.

\medskip
\noindent
$\mathbb{F}$: a field of characteristic zero.

\medskip
\noindent
$\mathbb{N}$: the nonnegative integers.

\medskip
\noindent
$\mathbb{Q}$: the rationals.

\medskip
\noindent
$\mathbb{R}$: the real numbers.

\medskip
\noindent
$\mathbb{R}_+$: the positive real numbers.

\medskip
\noindent
$\mathbb{Z}$: the integers.

\medskip
\noindent
$\mathbb{Z}_+$: the positive integers.

\medskip
\noindent
$\mathbb{Z}_2$: the integers modulo 2.

\chapter[Introduction to the moduli space of superspheres]{An 
introduction to the moduli space of $N=1$ superspheres with tubes 
and the sewing operation}

In this chapter we introduce the moduli space of $N=1$ superspheres 
with tubes and the sewing operation.  We begin in Section 2.1 with 
some background material on superalgebras, Grassmann algebras and 
superanalytic superfunctions (cf. \cite{D}, \cite{Bc1}, \cite{Ro2},
\cite{CR}).  In Section 2.2, we give the definition of $N=1$ 
superconformal superfunction and make note of the power series 
expansions of such functions vanishing at zero or infinity.  In 
Section 2.3, we give the definitions of supermanifold and $N=1$ 
super-Riemann surface and state the uniformization theorem for 
$N=1$ super-Riemann surfaces {}from \cite{CR}.

In Section 2.4, we study $N=1$ superspheres with ordered and oriented 
tubes, showing that these are superconformally equivalent to 
superspheres with ordered and oriented punctures and local 
superconformal coordinates vanishing at the punctures, and we define a 
sewing operation on this space.  

In Section 2.5, we define the moduli space of superspheres with tubes
and introduce canonical superspheres with tubes.  We then show that
any supersphere with tubes is superconformally equivalent to a
canonical supersphere with tubes and that two different canonical
superspheres with tubes are not superconformally equivalent.  This
shows that there is a bijection between the set of canonical
superspheres with tubes and the moduli space of superspheres with
tubes.

Finally, in Section 2.6, we define the sewing operation on the moduli
space of superspheres with tubes by showing how to obtain a canonical
supersphere {}from the sewing together of two canonical superspheres.
This gives rise to the sewing equation and normalization and boundary
conditions on the uniformizing function mapping the two sewn canonical 
superspheres to the new canonical supersphere.  Solving this sewing 
equation along with the normalization and boundary conditions will be
the main objectives in Chapters 3 and 4.

\section{Grassmann algebras and superanalytic superfunctions}

Let $\mathbb{F}$ be  a field  of characteristic zero.
For a $\mathbb{Z}_2$-graded vector space $X = X^0 \oplus X^1$, define
the {\it sign function} $\eta$ on the homogeneous subspaces of $X$ by
$\eta(x) = i$, for $x \in X^i$ and $i = 0,1$.  If $\eta(x) = 0$, we say
that $x$ is {\it even}, and if $\eta(x) = 1$, we say that $x$ is {\it
odd}.  Note that $0 \in X$ is both even and odd, i.e., the sign function
is double valued for $x=0$.  However, in practice this is never a 
problem, e.g., in property (ii) below.  

A {\it superalgebra} is an (associative) algebra $A$ (with identity $1 
\in A$), such that

(i) $A$ is a $\mathbb{Z}_2$-graded algebra

(ii) $ab = (-1)^{\eta(a)\eta(b)} ba$ for $a,b$ homogeneous in $A$.
\hfill \mbox{(supercommutativity)}

We will often write a given element $u \in A$ in 
terms of its vector space grading, i.e., if $u = u^0 + u^1$ for 
$u^i \in A^i$, we write $u = (u^0,u^1) \in A^0 \oplus A^1$.  Algebraic 
operations will still, of course, be considered to take place in the 
algebra $A$, and we will merely use the vector space decomposition in 
order to emphasize the $\mathbb{Z}_2$-grading of a given element.  
Note that when working over a field of characteristic zero or of 
characteristic greater than two, property (ii), supercommutativity, 
implies that the square of any odd element is zero.

A $\mathbb{Z}_2$-graded vector space $\mathfrak{g} = \mathfrak{g}^0
\oplus \mathfrak{g}^1$ is said to be a {\it Lie superalgebra} if it
has a bilinear operation $[\cdot,\cdot]$ on $\mathfrak{g}$ such that for
$u,v$ homogeneous in $\mathfrak{g}$
 
(i) $\; [u,v] \in {\mathfrak g}^{(\eta(u) + \eta(v))\mathrm{mod} \; 2}$ 

(ii) $\; [u,v] = -(-1)^{\eta(u)\eta(v)}[v,u]$  \hfill (skew-symmetry)
 
(iii) $\; (-1)^{\eta(u)\eta(w)}[[u,v],w] + (-1)^{\eta(v)\eta(u)}[[v,w],u]$ 
\hfill (Jacobi identity) 

$\hspace{2in} + \; (-1)^{\eta(w)\eta(v)}[[w,u],v] = 0$. 

\begin{rema}\label{envelope}
Given a Lie superalgebra $\mathfrak{g}$ and a superalgebra 
$A$, the space $(A^0 \otimes \mathfrak{g}^0) \oplus (A^1 \otimes \mathfrak{g}^1)$ 
is a Lie algebra with bracket given by
\begin{equation}\label{super v. non} 
[au , bv] = (-1)^{\eta(b) \eta(u)} ab [u,v]
\end{equation}
for $a,b \in A$ and $u,v \in \mathfrak{g}$ homogeneous  
(with obvious notation), where in (\ref{super v. non}) we have suppressed 
the tensor product symbol.  Note that the bracket on the left-hand side 
of (\ref{super v. non}) is a Lie algebra bracket, and the bracket on the 
right-hand side is a Lie superalgebra bracket.  We will call $(A^0 
\otimes \mathfrak{g}^0) \oplus (A^1 \otimes \mathfrak{g}^1)$ the {\it 
$A$-envelope of $\mathfrak{g}$}.  Of course, given two 
superalgebras $A$ and $\hat{A}$, we can form the {\it $A$-envelope of 
$\hat{A}$} given by $(A^0 \otimes \hat{A}^0) \oplus (A^1 \otimes \hat{A}^1)$ 
which is naturally an algebra.  In fact, since $A \otimes \hat{A}$ is itself
a superalgebra, we see that the $A$-envelope of $\hat{A}$ is equal to the
even homogeneous subspace of $A \otimes \hat{A}$, i.e., $(A \otimes \hat{A})^0$.
\end{rema}

For any $\mathbb{Z}_2$-graded associative algebra $A$ and for $u,v \in A$ 
of homogeneous sign, we can define $[u,v] = u v - (-1)^{\eta(u)\eta(v)} v 
u$, making $A$ into a Lie superalgebra.  The algebra of endomorphisms of 
$A$, denoted $\mbox{End} \; A$, has a natural $\mathbb{Z}_2$-grading 
induced {}from that of $A$, and defining $[X,Y] = X Y - (-1)^{\eta(X)
\eta(Y)} Y X$ for $X,Y$ homogeneous in $\mbox{End} \; A$, this gives 
$\mbox{End} \; A$ a Lie superalgebra structure.  An element $D \in 
(\mbox{End} \; A)^i$, for $i \in \mathbb{Z}_2$, is called a {\it 
superderivation of sign $i$} (denoted $\eta(D) = i$) if $D$ satisfies the 
super-Leibniz rule
\begin{equation}\label{leibniz} 
D(uv) = (Du)v + (-1)^{\eta(D) \eta(u)} uDv  
\end{equation} 
for $u,v \in A$ homogeneous.

\begin{rema}\label{superphilosophy} The above definitions of
superalgebra and Lie superalgebra are the standard definitions
used. However, we would like to point out that in analogy to
associative algebra, commutative algebra, and Lie algebra, a more
appropriate definition would have been to define a {\it superalgebra}
to be a $\mathbb{Z}_2$-graded associative algebra, a {\it 
supercommutative superalgebra} to be what we called a superalgebra 
above, and retain our definition of Lie superalgebra.  Using the 
criterion that the ``super'' version of a structure should be one 
which introduces a minus sign when two odd elements are exchanged, a
$\mathbb{Z}_2$-graded associative algebra vacuously satisfies the
criterion for being called a superalgebra.  The zero-graded subalgebra
of such an algebra is an associative algebra, and the zero-graded
subalgebra of a supercommutative superalgebra is a commutative
algebra.  In addition, the introduction of a bracket operation on a
$\mathbb{Z}_2$-graded associative algebra called a superalgebra (and
thus with the understanding that a minus sign be introduced when two
odd elements are exchanged) would then naturally define a Lie
superalgebra in analogy to a Lie algebra arising {}from an associative
algebra.  However, to avoid confusion with the literature, we have
used the terminology ``superalgebra'' in the customary way above to
denote a ``supercommutative superalgebra" and will continue to do so 
throughout the remainder of this work.  
\end{rema}

Let $\mathcal{T}(V)$ be the tensor algebra over the vector space $V$, 
and let $\mathcal{J}$ be the ideal of $\mathcal{T}(V)$ generated by 
the elements $v \otimes w + w \otimes v$ for $v,w \in V$.  Then the 
exterior algebra generated by $V$ is given by $\bigwedge (V) = 
\mathcal{T}(V)/\mathcal{J}$, and $\bigwedge(V)$ has the structure of 
a superalgebra.  For $L \in \mathbb{N}$, fix $V_L$ to be an 
$L$-dimensional vector space over $\mathbb{C}$ with fixed basis 
$\{\zeta_1,\zeta_2, \ldots,\zeta_L\}$ such that $V_L \subset V_{L+1}$.  
We denote $\bigwedge(V_L)$ by $\bigwedge_L$ and call this the {\it 
Grassmann algebra on $L$  generators}.  In other words, {}from now on 
we will consider the Grassmann algebras to have a fixed sequence of 
generators.  Note that $\bigwedge_L \subset \bigwedge_{L+1}$, and 
taking the direct limit as $L \rightarrow \infty$, we have the {\it 
infinite Grassmann algebra} denoted by $\bigwedge_\infty$. 
Then $\bigwedge_L$ and $\bigwedge_\infty$ are the associative 
algebras over $\mathbb{C}$ with generators $\zeta_{i}$, for $i = 1, 
2, ..., L$ and $i = 1, 2, ...$, respectively, and with relations:
\[\zeta_{i} \zeta_{j} = - \zeta_{j} \zeta_{i}, \qquad \qquad \zeta_{i}^2 = 0 .\]    
Note that $\mathrm{dim}_{\mathbb{C}} \; \bigwedge_L = 2^L$, and if 
$L = 0$, then $\bigwedge_0 = \mathbb{C}$. We use the notation 
$\bigwedge_*$ to denote a Grassmann algebra, finite or infinite.  
The reason we take $\bigwedge_*$ to be over $\mathbb{C}$ is that 
we will mainly be interested in complex supergeometry and 
transition functions which are superanalytic as described 
below.  However, formally, we could just as well have taken 
$\mathbb{C}$ to be any field of characteristic zero, and in fact, in
Chapter 3 we will consider formal superanalytic functions over a
general superalgebra. 
 
Let 
\begin{eqnarray*}
I_L \! \!  &=&  \! \! \bigl\{ (i) = (i_1, i_2, \ldots, i_{2n}) \; | \; i_1 <
i_2 < \cdots < i_{2n}, \; i_l \in \{1, 2, ..., L\}, \; n \in
\mathbb{N} \bigr\}, \\ 
J_L  \! \! &=&  \! \! \bigl\{(j) = (j_1, j_2, \ldots, j_{2n + 1}) \; | \;
j_1 < j_2 < \cdots < j_{2n + 1}, \; j_l \in \{1, 2, ..., L\}, \;  n \in \mathbb{N} \bigr\},
\end{eqnarray*}
and $K_L = I_L \cup J_L$.  Let
\begin{eqnarray*}
I_\infty \! \! &=& \! \! \bigl\{(i) = (i_1, i_2, \ldots, i_{2n})\; | \; i_1 < i_2
< \cdots < i_{2n}, \; i_l \in \mathbb{Z}_+, \; n \in \mathbb{N} \bigr\}, \\
J_\infty \! \! &=& \! \! \bigl\{(j) = (j_1, j_2, \ldots, j_{2n + 1})\; | \; j_1 <
j_2 < \cdots < j_{2n + 1}, \; j_l \in \mathbb{Z}_+, \; n \in \mathbb{N} \bigr\},
\end{eqnarray*}
and $K_\infty = I_\infty \cup J_\infty$.  We use $I_*$, $J_*$, and
$K_*$ to denote $I_L$ or $I_\infty$, $J_L$ or $J_\infty$, and $K_L$ 
or $K_\infty$, respectively.  Note that $(i) = (i_1,...,i_{2n})$ for 
$n = 0$ is in $I_*$, and we denote this element by $(\emptyset)$.  
The $\mathbb{Z}_2$-grading of $\bigwedge_*$ is given explicitly by
\begin{eqnarray*}
\mbox{$\bigwedge_*^0$} \! &=& \! \Bigl\{a \in \mbox{$\bigwedge_*$} \; 
\big\vert \; a = \sum_{(i) \in I_*} a_{(i)}\zeta_{i_{1}}\zeta_{i_{2}} 
\cdots \zeta_{i_{2n}}, \; a_{(i)} \in \mathbb{C}, \; n \in \mathbb{N}
\Bigr\}\\ 
\mbox{$\bigwedge_*^1$} \! &=& \! \Bigl\{a \in \mbox{$\bigwedge_*$} \; 
\big\vert \; a = \sum_{(j) \in J_*} a_{(j)}\zeta_{j_{1}}\zeta_{j_{2}} 
\cdots \zeta_{j_{2n + 1}}, \; a_{(j)} \in \mathbb{C}, \; n \in \mathbb{N}
\Bigr\} . 
\end{eqnarray*} 
Note that $a^2 = 0$ for all $a \in \bigwedge_*^1$.

We can also decompose $\bigwedge_*$ into {\it body}, $(\bigwedge_*)_B = 
\{ a_{(\emptyset)} \in \mathbb{C} \}$,  and {\it soul} 
\[(\mbox{$\bigwedge_*$})_S \; = \; \Bigl\{a \in \mbox{$\bigwedge_*$} \; \big\vert \; 
a = \! \! \! \sum_{ \begin{scriptsize} \begin{array}{c}
(k) \in K_*\\
k \neq (\emptyset)
\end{array} \end{scriptsize}} \! \! \!
a_{(k)} \zeta_{k_1} \zeta_{k_2} \cdots \zeta_{k_n}, \; a_{(k)} \in \mathbb{C} 
\Bigr\}\] 
subspaces such that $\bigwedge_* = (\bigwedge_*)_B \oplus 
(\bigwedge_*)_S$.  For $a \in \bigwedge_*$, we write $a = a_B + 
a_S$ for its body and soul decomposition.  Note that for all $a 
\in \bigwedge_L$, we have $a_S^{L + 1} = 0$.  However, no such 
general nilpotency condition holds for the  soul of $a \in 
\bigwedge_\infty$, i.e., there exist $a \in \bigwedge_\infty$  
such that $a_S^n \neq 0$ for all $n \in \mathbb{N}$.

Let $U$ be a subset of $\bigwedge_*$, and write $U = U^0 \oplus U^1$
for the decomposition of $U$ into even and odd subspaces.  Let $z$ 
be an even variable in $U^0$, i.e., an indeterminate element of 
$U^0$, and let $\theta$ be an odd variable in $U^1$, i.e., an 
indeterminate element of $U^1$.  Here we would like to stress that 
in this context ``variable'' does not mean formal variable.  We will 
call a map
\begin{eqnarray*}
H: U &\longrightarrow& \mbox{$\bigwedge_*$}\\
(z,\theta) &\mapsto& H(z,\theta) 
\end{eqnarray*}  
a {\it $\bigwedge_*$-superfunction in (1,1)-variables on $U$}.  Any such
superfunction $H(z,\theta)$ can be decomposed as 
\begin{equation}\label{Hdecomp}
H(z,\theta) = \biggl( \sum_{(i) \in I_*} H_{(i)}(z,\theta) \zeta_{i_1}
\zeta_{i_2} \cdots \zeta_{i_{2n}}, \sum_{(j) \in J_*}
H_{(j)}(z,\theta) \zeta_{j_1} \zeta_{j_2} \cdots \zeta_{j_{2n + 1}}
\biggr)    
\end{equation}
where each $H_{(k)}(z,\theta) : U \longrightarrow \mathbb{C}$ is a 
complex function on $U$ in terms of the $z_{(i)}$'s and 
$\theta _{(j)}$'s, for $(k) \in K_*$, $(i) \in I_*$, and $(j) \in 
J_*$.  We will often use the notation $H(z,\theta) = (\tilde{z}, 
\tilde{\theta})$ for the decomposition in (\ref{Hdecomp}) and will 
call $\tilde{z}$ (resp., $\tilde{\theta}$) an {\it even} (resp., {\it 
odd}) {\it superfunction} on $U$.  We will also occasionally use the 
notation $H(z,\theta) = (H^0(z,\theta), H^1(z,\theta))$ for the 
decomposition of $H$ into even and odd superfunctions.  The set of 
$\bigwedge_*$-superfunctions in (1,1)-variables on a fixed set $U$ is 
a superalgebra in the obvious way.  We call the complex valued 
function $H_B (z, \theta) = H_{(\emptyset)} (z, \theta)$ defined on 
$U$ the {\it body} of $H$ and $H_S (z, \theta) = H(z, \theta) - 
H_B (z, \theta)$ the {\it soul} of $H$.  Note that if $H(z,\theta)$ 
is a $\bigwedge_*$-superfunction in (1,1)-variables on $U$, then 
$H_B (z, \theta)$ and $H_S (z,\theta)$ are also 
$\bigwedge_*$-superfunctions in (1,1)-variables on $U$.
 
This notion of superfunction in $(1,1)$-variables can be extended in
the obvious way to the notion of a $\bigwedge_*$-superfunction in
$(m,n)$-variables, i.e., in $m \in \mathbb{N}$ even variables and $n
\in \mathbb{N}$ odd variables.

\begin{rema} The second ``1'' in the expression
``$(1,1)$-variables'' refers to the number of odd variables and is
exactly the number ``$N = 1$'' in the expression ``$N = 1$
superconformal field theory''.  We also remark that we are following
the conventions in the literature (cf. \cite{Fd}, \cite{Bc2}) in
calling an even variable $z$ and an odd variable $\theta$.  We caution
that this may lead the reader to think of $z$ as being merely a
complex variable.  This is far {}from the case since in fact $z \in
\bigwedge_*^0$.  Recall that we use the notation $z_B$ or
$z_{(\emptyset)}$ to denote the complex portion of the variable $z$,
but $z$ itself represents either a complex $2^{L - 1}$-tuple if
$\bigwedge_* = \bigwedge_L$, or an infinite number of complex
variables, $z_{(i)}$, for $(i) \in I_\infty$, if $\bigwedge_* =
\bigwedge_\infty$.  If one wishes to restrict to the nonsuper case,
one should think of setting all odd variables equal to zero and the
soul portion of all even variables equal to zero, thus leaving the
``body'' portion of the theory giving the nonsuper case.
\end{rema}

In developing a notion of superanalyticity, we would like to say that
a $\bigwedge_*$-superfunction $H$ is ``superanalytic" in an open set
$U \subseteq \bigwedge_*$ in an appropriately defined topology on
$\bigwedge_*$ if the appropriately defined partial derivatives
$\dtheta H$ and $\dz H$ exist and are well defined on $U$.  For 
$\bigwedge_* = \bigwedge_\infty$ the conditions needed on $H$ for these 
partials to exist in a subset of $\bigwedge_\infty$ which is open in an 
appropriately defined topology are very straightforward. However, for a 
$\bigwedge_L$-superfunction, there can be unwanted cancellation for the 
odd variables.  For example, if $\theta \in \bigwedge_L^1$ is an odd 
variable, and $a = \zeta_1 \zeta_2 \cdots \zeta_L \in \bigwedge_L$ is a 
fixed supernumber, then the superfunction $\theta a$ is identically 
zero, and thus the partial with respect to $\theta$ of this 
superfunction on the one hand should be $a$, and on the otherhand should 
be zero. This illustrates the fact that for $\dtheta H$ to be well 
defined for $H$ a $\bigwedge_L$-superfunction, some conditions on $H$ 
must be imposed.  One of these conditions is that the Grassmann algebra 
$\bigwedge_L$ be large enough and the coefficients of the supervariables 
in $H$ restricted to a suitable subspace so that this unwanted  
cancellation does not occur (cf. \cite{D}, \cite{Bc1}, \cite{Ro2}, 
\cite{CR}, \cite{B thesis}).  In fact, one needs the underlying 
Grassmann algebra to be generated by a vector space of dimension equal 
to or greater than the number of odd variables so that one can restrict 
the coefficients of the power series expansion of a 
$\bigwedge_L$-superfunction so as to not contain one of the underlying 
basis elements $\zeta_i$, $i=1,...L$, for each of the odd  variables.  
That is, if the number of odd variables is $N \in \mathbb{N}$ and 
well-defined partial derivatives as well as multiple partial derivatives 
with respect to these odd variables are desired, then the underlying 
Grassmann algebra must be defined over a vector space of dimension $L$ 
where $L$ is greater than or equal to $N$ and the coefficients in the 
power series expansion of the superfunction must be restricted to a 
subspace isomorphic to $\bigwedge_{L - N}$.

Since the motivation for this work is $N=1$ superconformal field
theory in which the fields are $N=1$ superanalytic superfunctions and
the symmetries of the theory are $N=1$ superconformal transformations,
one might expect that we need only require $L \geq 1$.  However, in
this paper, we will be interested in having well-defined partials for
an infinite number of odd variables.  This is because the moduli space 
of $(1,1)$-dimensional superspheres with tubes has an infinite number 
of odd coordinates (see Remark \ref{infinite variables1}).  In Chapter 
4, we will consider the supermeromorphic tangent space of this moduli 
space and partials with respect to the moduli space coordinates.  In 
order to differentiate with respect to all these variables with  
abandon, we must work over an infinite Grassmann algebra, (although in  
practice, we never take more than two partials in a row).   Thus  
starting in Section 2.5, for simplicity, we will begin working over 
$\bigwedge_\infty$, and until then, we continue to use a finite or 
infinite Grassmann algebra $\bigwedge_*$ noting the restrictions as 
needed.

Before defining superanalyticity, we define what it means for a
complex analytic function to be defined for an even element of
$\bigwedge_*$.  Let $z_B$ be a complex variable and $h(z_B)$ a complex
analytic function in some open set $U_B \subset \mathbb{C}$.  For $z$ a
variable in $\bigwedge_*^0$, we define 
\begin{equation}\label{Taylor expansion}
h(z) = \sum_{n \in \mathbb{N}} \frac{z_S^n}{n!} h^{(n)}(z_B), 
\end{equation}
i.e., $h(z)$ is the Taylor expansion about the body of $z = z_B + z_S$.  
Then $h(z)$ is well defined (i.e., convergent) in the open neighborhood 
$\{z = z_B + z_S \in \bigwedge_*^0 \; | \; z_B \in U_B \} \subseteq 
\bigwedge_*^0$.  This is because for $n \in \mathbb{N}$, $z_S^n$ has 
terms involving at least $2n$ of the $\zeta_i$'s.  Thus any 
$\zeta_{i_1} \cdots \zeta_{i_{2m}}$ coefficient in (\ref{Taylor 
expansion}) will be a finite sum.  Since $h(z)$ is algebraic in each 
$z_{(i)}$, for $(i) \in I_*$, it follows that $h(z)$ is complex 
analytic in each of the complex variables $z_{(i)}$.

For $n \in \mathbb{N}$, we introduce the notation $\bigwedge_{*>n}$ to
denote a finite Grassmann algebra $\bigwedge_L$ with $L > n$ or an 
infinite Grassmann algebra.  We will use the corresponding index 
notations for the corresponding indexing sets $I_{*>n}, J_{*>n}$ and 
$K_{*>n}$.

\begin{defn}\label{define superanalytic}
A {\em superanalytic $\bigwedge_{*>0}$-superfunction in
$(1,1)$-variables} $H$ is a $\bigwedge_{*>0}$-superfunction in
$(1,1)$-variables of the form
\begin{eqnarray*}
H(z, \theta) \! \! &=& \! \! (f(z) + \theta \xi(z), \psi(z) + \theta
g(z)) \\ 
&=& \! \! \!  \Biggl(\sum_{(i) \in I_{ * - 1}} \!
\! \! f_{(i)} (z) \zeta_{i_1} \zeta_{i_2} \cdots \zeta_{i_n} + \;
\theta \! \! \! \sum_{(j) \in J_{ * - 1}} \! \! \!
\xi_{(j)} (z) \zeta_{j_1} \zeta_{j_2} \cdots \zeta_{j_{2n + 1}},
\Biggr. \\ 
& & \qquad \quad \Biggl. \sum_{(j) \in J_{ * - 1}} \!  \! \!
\psi_{(j)} (z) \zeta_{j_1} \zeta_{j_2} \cdots \zeta_{j_{2n + 1}} + \;
\theta \! \! \! \sum_{(i) \in I_{ * - 1}} \! \! \!
g_{(i)} (z) \zeta_{i_1} \zeta_{i_2} \cdots \zeta_{i_n} \Biggr)
\end{eqnarray*}
where $f_{(i)}(z_B)$, $g_{(i)}(z_B)$, $\xi_{(j)}(z_B)$, and
$\psi_{(j)}(z_B)$ are all complex analytic in some non-empty open
subset $U_B \subseteq \mathbb{C}$.  
\end{defn}

The restriction of the coefficients of the $f_{(i)}(z)$'s,
$g_{(i)}(z)$'s, $\xi_{(j)}(z)$'s, and $\psi_{(j)}(z)$'s to
$\bigwedge_{* - 1} \subseteq \bigwedge_{*>0}$ is so that the partial 
with respect to $\theta$ of $H$ (as defined below) is well defined. 
Actually, we could restrict the coefficients to be in any
Grassmann subalgebra of $\bigwedge_{*>0}$ on $*-1$ generators, but 
for simplicity we will follow the convention of restricting to
$\bigwedge_{*-1}$ (cf. \cite{Ro2}).  If $\bigwedge_{*>0} = 
\bigwedge_\infty$, then $\bigwedge_{*-n} = \bigwedge_\infty$. 

Note that if each $f_{(i)}(z_B)$, $g_{(i)}(z_B)$, $\xi_{(j)}(z_B)$, 
and $\psi_{(j)}(z_B)$ is complex analytic in $U_B \subseteq \mathbb{C}$,
then $f_{(i)}(z)$, $g_{(i)}(z)$, $\xi_{(j)}(z)$, and $\psi_{(j)}(z)$
are all well defined in $\{z = z_B + z_S \in \bigwedge_{*>0}^0 \; | \; 
z_B \in U_B \}$.  Thus $H(z,\theta)$ is well defined (i.e., convergent)
for $\{(z,\theta) \in \bigwedge_{*>0} \; | \; z_B \in U_B \} = U_B 
\times (\bigwedge_{*>0})_S = U$.  Consider the topology on 
$\bigwedge_{*>0}$ given by the product of the usual topology on 
$(\bigwedge_{*>0})_B = \mathbb{C}$ and the trivial topology on 
$(\bigwedge_{*>0})_S$.  This topology on $\bigwedge_{*>0}$ is called 
the {\it DeWitt topology}.  The natural domain of any superanalytic
$\bigwedge_{*>0}$-superfunction is an open set in the DeWitt topology 
on $\bigwedge_{*>0}$.

We define the (left) partial derivatives $\dz$ and $\dtheta$ acting on
superfunctions which are superanalytic in some DeWitt open
neighborhood $U$ of $(z,\theta) \in \bigwedge_{*>0}$ by
\begin{eqnarray*}
\Delta z \left(\dz H(z,\theta) \right) + O((\Delta z)^2) &=& H(z +
\Delta z, \theta) - H(z,\theta) \\
\Delta \theta \left(\dtheta H(z,\theta) \right) &=& H(z, \theta +
\Delta \theta) - H(z,\theta) 
\end{eqnarray*}
for all $\Delta z \in \bigwedge_{*>0}^0$ and $\Delta \theta \in
\bigwedge_{*>0}^1$ such that $z + \Delta z \in U^0 = U_B \times
(\bigwedge_{*>0}^0)_S$ and $\theta + \Delta \theta \in U^1 =
\bigwedge_{*>0}^1$.  See for example \cite{B thesis} for a proof 
and discussion of the fact that these partials are in fact 
well defined.  Note that $\dz$ and $\dtheta$ are endomorphisms of the 
superalgebra of superanalytic $(1,1)$-superfunctions, and in fact, are 
even and odd superderivations, respectively.

Let $(\bigwedge_*)^\times$ denote the set of invertible elements in
$\bigwedge_*$.  Then
\[(\mbox{$\bigwedge_*$})^\times = \{a \in \mbox{$\bigwedge_*$} \; | 
\; a_B \neq 0 \} \]
since 
\[\frac{1}{a} = \frac{1}{a_B + a_S} = \sum_{n \in \mathbb{N}}
\frac{(-1)^n a_S^n}{a_B^{n + 1}} \]
is well defined if and only if $a_B \neq 0$.  In light of this fact,
note that the DeWitt topology on $\bigwedge_*$ is non-Hausdorff.  Two 
points $a, b \in \bigwedge_*$ can be separated by disjoint open sets in 
the DeWitt topology if and only if $a_B \neq b_B$, i.e., if and only if 
their difference is an invertible element of $\bigwedge_*$.  In other 
words, the DeWitt topology fails to be Hausdorff exactly to the extent 
that the nonzero elements of $\bigwedge_*$ fail in general to be 
invertible.

If $h(z_B)$ is complex analytic in an open neighborhood of the complex
plane, then $h(z_B)$ has a Laurent series expansion in $z_B$, given by 
$h(z_B) = \sum_{l \in \mathbb{Z}} c_l z_B^l$, for $c_l \in \mathbb{C}$, 
and we have
\begin{eqnarray}
h(z) &=& \sum_{n \in \mathbb{N}} \frac{z_S^n}{n!} h^{(n)}(z_B) =
\sum_{n \in \mathbb{N}} \sum_{l \in \mathbb{Z}} c_l \frac{l!}{n!(l - n)!}
z_S^n z_B^{l - n} \nonumber\\
&=& \sum_{l \in \mathbb{Z}} c_l (z_B + z_S)^l = \sum_{l 
\in \mathbb{Z}} c_l z^l \label{analytic definition}
\end{eqnarray} 
where $(z_B + z_S)^l$, for $l \in \mathbb{Z}$, is always understood to 
mean expansion in positive powers of the second variable, in this
case $z_S$.  Thus if $H$ is a $\bigwedge_{*>0}$-superfunction in 
$(1,1)$-variables which is superanalytic in a (DeWitt) open 
neighborhood, $H$ can be expanded as
\begin{equation}\label{Laurent series}
H(z, \theta) = \biggl(\sum_{l \in \mathbb{Z}} a_l z^l + \theta \sum_{l
\in \mathbb{Z}} n_l z^l, \sum_{l \in \mathbb{Z}} m_l z^l + \theta
\sum_{l \in \mathbb{Z}} b_l z^l \biggr)
\end{equation}   
for $a_l, b_l \in \bigwedge_{ * - 1}^0$ and $m_l, n_l \in
\bigwedge_{ * - 1}^1$.  

In Chapter 4, we will encounter superanalytic superfunctions 
in more than one even and more than one odd variable.  Let $m,n \in 
\mathbb{N}$, and let $U$ be a subset of $(\bigwedge_*^0)^m \oplus 
(\bigwedge_*^1)^n$.  A $\bigwedge_*$-superfunction $H$ on $U$ in 
$(m,n)$-variables is given by
\begin{eqnarray*}
H: U &\longrightarrow& \mbox{$\bigwedge_*$} \\ 
(z_1,z_2,...,z_m ,\theta_1,\theta_2,...,\theta_n) &\mapsto&
H(z_1,z_2,...,z_m ,\theta_1,\theta_2,...,\theta_n)
\end{eqnarray*}
where $z_k$, for $k = 1,...,m$, are even variables in $\bigwedge_*^0$ 
and $\theta_k$, for $k = 1,...,n$, are odd variables in $\bigwedge_*^1$.  
Let $f((z_1)_B,(z_2)_B,...,(z_m)_B)$ be a complex analytic function in
$(z_k)_B$, for $k = 1,...,m$.  For $z_k \in \bigwedge_*^0$, and $k
= 1,..,m$, define
\begin{multline}\label{more than one variable}
f(z_1,z_2,...,z_m) = \!
\sum_{l_1,...,l_m \in \mathbb{N}} \! \! \frac{(z_1)_S^{l_1} (z_2)_S^{l_2} \cdots
(z_m)_S^{l_m}}{l_1 ! l_2 ! \cdots l_m !} 
\biggl(\frac{\partial \; \;}{\partial (z_1)_B}\biggr)^{l_1} \cdot \\
\biggl(\frac{\partial \; \;}{\partial (z_2)_B}\biggr)^{l_2} \! \! 
\cdots \biggl(\frac{\partial \;
\;}{\partial (z_m)_B}\biggr)^{l_m} \! \! \cdot f((z_1)_B,(z_2)_B,...,(z_m)_B) .
\end{multline}

\begin{defn}  Let $m, n \in \mathbb{N}$.  Let $U \subset
(\bigwedge_{*>n-1}^0)^m \oplus (\bigwedge_{*>n-1}^1)^n$, and let $H$
be a $\bigwedge_{*>n-1}$-superfunction in $(m,n)$-variables defined on
$U$. Let
\begin{eqnarray*}
A^0_n \! &=& \! \left\{(k) = (k_1,k_2, ..., k_{2l}) \; | \; l \in \mathbb{N} , 
\; k_j \in \{1,...,n\}, \; k_1<k_2<\cdots k_{2l}\right\}\\ 
A^1_n \! &=& \! \left\{(k) = (k_1,k_2, ..., k_{2l + 1}) \; | \; l \in \mathbb{N} , 
\; k_j \in \{1,...,n\}, \; k_1<k_2<\cdots k_{2l+1} \right\} .
\end{eqnarray*}
Then $H$ is said to be {\em superanalytic} if
$H$ is of the form
\begin{multline*}
H(z_1,z_2,...,z_m ,\theta_1,\theta_2,...,\theta_n) \\  
= \; \Biggl( \sum_{ (k) \in A_n^0} \! \theta_{k_1} \cdots \theta_{k_{2l}} 
f_{(k)}(z_1,z_2,...,z_m) + \! \! \sum_{(k) \in A_n^1} \! \theta_{k_1} 
\cdots \theta_{k_{2l + 1}} \xi_{(k)}(z_1,z_2,...,z_m),\Biggr.\\  
\Biggl. \sum_{(k) \in A_n^0} \! \theta_{k_1} \cdots \theta_{k_{2l}}
\psi_{(k)}(z_1,z_2,...,z_m) + \! \! \sum_{(k) \in A_n^1} \! \theta_{k_1} 
\cdots \theta_{k_{2l + 1}} g_{(k)}(z_1,z_2,...,z_m) \! \Biggr), 
\end{multline*} 
where each $f_{(k)}$ and $g_{(k)}$ is of the form 
\[f_{(k)}(z_1,z_2,...,z_m) = \sum_{(i) \in I_{ * - n}}
f_{(k),(i)}(z_1,z_2,...,z_m) \zeta_{i_1}\zeta_{i_2} \cdots
\zeta_{i_{2s}}, \]  
each $\xi_{(k)}$ and $\psi_{(k)}$ is of the form
\[ \xi_{(k)}(z_1,z_2,...,z_m) = \sum_{(j) \in J_{ * - n}}
\xi_{(k),(j)} (z_1,z_2,...,z_m) \zeta_{j_1}\zeta_{j_2} \cdots
\zeta_{j_{2s + 1}},\] 
and each $f_{(k),(i)}((z_1)_B,(z_2)_B,...,(z_m)_B)$, $\xi_{(k),(j)}
((z_1)_B,(z_2)_B,...,(z_m)_B)$, \\
$\psi_{(k),(j)}((z_1)_B,(z_2)_B,...,(z_m)_B)$, and $g_{(k),(i)}((z_1)_B,
(z_2)_B,...,(z_m)_B)$ is analytic  in $(z_l)_B$, for $l = 1,...,m$, and
$((z_1)_B,(z_2)_B,...,(z_m)_B) \in  U_B \subset \mathbb{C}^m$.
\end{defn}

We require the even and odd variables to be in
$\bigwedge_{* > n-1}$, and we restrict the coefficients of the
$f_{(k),(i)}$'s, $\xi_{(k),(j)}$'s, $\psi_{(k),(j)}$'s, and
$g_{(k),(i)}$'s to be in $\bigwedge_{* - n} \subseteq
\bigwedge_{*>n-1}$ in order for the partial derivatives with respect
to each of the $n$ odd variables to be well defined and for multiple 
partials to be well defined.  As in the $(m,n) = (1,1)$ case, we could 
merely require the coefficients to be in any Grassmann  subalgebra of 
$\bigwedge_{*>n-1}$ on $*-n$ generators, but for simplicity we instead 
restrict to $\bigwedge_{*-n}$.

Consider the projection
\begin{eqnarray}
\pi^{(m,n)}_B : (\mbox{$\bigwedge_{* > n-1}^0$})^m \oplus 
(\mbox{$\bigwedge_{* > n-1}^1$})^n & \longrightarrow & \mathbb{C}^m  
\label{projection onto body}\\ 
(z_1,...,z_m, \theta_1,...,\theta_n) & \mapsto & ((z_1)_B,(z_2)_B,...,
(z_m)_B) . \nonumber
\end{eqnarray}
We define the {\it DeWitt topology on $(\bigwedge_{* > n-1}^0)^m
\oplus (\bigwedge_{* > n-1}^1)^n$} by letting 
\[U \subseteq ((\mbox{$\bigwedge_{*> n-1}^0$})^m \oplus 
(\mbox{$\bigwedge_{* > n-1}^1$})^n)\] 
be an open set in the DeWitt topology if and only if $U = 
(\pi^{(m,n)}_B)^{-1} (V)$ for some open set $V \subseteq \mathbb{C}^m$.  
Note that the natural domain of a superanalytic 
$\bigwedge_{* > n-1}$-superfunction in $(m,n)$-variables is an open set 
in the DeWitt topology.

A superconformal field theory based on ``superfields'' which are
superanalytic superfunctions in $(1,n)$-variables would be referred to
as an ``$N = n$ superconformal field theory''.

\begin{rema}\label{H_L remark}
Recall that $\bigwedge_L \subset \bigwedge_{L+1}$ for $L \in 
\mathbb{N}$, and note that {}from (\ref{more than one variable}), any 
superanalytic $\bigwedge_L$-superfunction, $H_L$, in $(m,n)$-variables for 
$L \geq n$ can naturally be extended to a superanalytic 
$\bigwedge_{L'}$-superfunction in $(m,n)$-variables for $L'>L$ and hence 
to a superanalytic $\bigwedge_\infty$-superfunction.  Conversely, if 
$H_{L'}$ is a superanalytic $\bigwedge_{L'}$-superfunction (or 
$\bigwedge_\infty$-superfunction) in $(m,n)$-variables for $L' >n$, then 
we can restrict $H_{L'}$ to a superanalytic $\bigwedge_L$-superfunction for 
$L'> L\geq n$ by restricting $(z_1,...z_m,,\theta_1,...,\theta_n) \in 
(\bigwedge_L^0)^m \oplus (\bigwedge_L^1)^n$ and setting 
$f_{(k),(i)} \equiv g_{(k),(i)} \equiv 0$ if $(i) \notin I_{L-n}$ and 
$\xi_{(k),(j)} \equiv \psi_{(k),(j)} \equiv 0$ if $(j) \notin J_{L-n}$. 
\end{rema}

\section{Superconformal $(1,1)$-superfunctions and power series}

Let $z$ be an even variable in $\bigwedge_{*>0}^0$ and $\theta$ an odd
variable in $\bigwedge_{*>0}^1$.  Following \cite{Fd}, define $D$ to 
be the odd superderivation $D = \Dz$ acting on 
$\bigwedge_{*>0}$-superfunctions in $(1,1)$-variables which are 
superanalytic in some DeWitt open subset $U \subseteq \bigwedge_{*>0}$.  
Then $D^2 = \dz$, and if $H(z,\theta) =  (\tilde{z},\tilde{\theta})$ is 
superanalytic in some  DeWitt open subset  $U \subseteq \bigwedge_{*>0}$, 
then $D$ transforms under $H(z,\theta)$ by
\begin{equation}\label{transform D}
D = (D\tilde{\theta})\tilde{D} + (D\tilde{z} - \tilde{\theta} D
\tilde{\theta})\tilde{D}^2
\end{equation}
where $\tilde{D} = \frac{\partial}{\partial \tilde{\theta}} +
\tilde{\theta} \frac{\partial}{\partial \tilde{z}}$ with 
$\frac{\partial}{\partial \tilde{z}}$ and $\frac{\partial}{\partial
\tilde{\theta}}$ defined by 
\[ \dz = \frac{\partial \tilde{z}}{\partial z} \frac{\partial}{\partial 
\tilde{z}} + \frac{\partial \tilde{\theta}}{\partial z} 
\frac{\partial}{\partial \tilde{\theta}} \qquad \mbox{and} \qquad 
\dtheta =  \frac{\partial \tilde{z}}{\partial \theta} 
\frac{\partial}{\partial \tilde{z}} + \frac{\partial
\tilde{\theta}}{\partial \theta} \frac{\partial}{\partial
\tilde{\theta}}.\] 

Recall that a complex function $f$ defined on an open set $U_B$ in
$\mathbb{C}$, of one complex variable $z_B$, is conformal in $U_B$ if
and only if $\frac{d \;}{d z_B} f(z_B)$ exists for $z_B \in U_B$ and 
is not identically zero in $U_B$, i.e., if and only if $f(z_B) = 
\tilde{z}_B$ transforms $\frac{d \;}{d z_B}$ by $\frac{d \;}{d z_B} = 
f'(z_B) \frac{d \;}{d \tilde{z}_B}$ for $f'$ not identically zero.  
Such a transformation of $\frac{d \;}{d z_B}$ is said to be {\it
homogeneous of degree one}, i.e., $f$ transforms $\frac{d \;}{d z_B}$ 
by a non-zero analytic function times $\frac{d \;}{d \tilde{z}_B}$ to 
the first power with no higher order terms in $\frac{d\;}{d\tilde{z}_B}$.   
Analogously we define a superconformal superfunction on a DeWitt open 
subset $U$ of $\bigwedge_{*>0}$ to be a superanalytic superfunction 
$H$ under which $D$ transforms homogeneously of degree one in $U$.  
That is, $H$ transforms $D$ by a non-zero superanalytic superfunction 
times $\tilde{D}$ to the first power, and no higher terms in 
$\tilde{D}$.  Since a superanalytic function $H(z,\theta) = 
(\tilde{z}, \tilde{\theta})$ transforms $D$ according to  
(\ref{transform D}), $H$ is superconformal if and only if, in 
addition to being superanalytic, $H$ satisfies
\begin{equation}\label{basic superconformal condition}
D\tilde{z} - \tilde{\theta} D\tilde{\theta} = 0,
\end{equation}
for $D \tilde{\theta}$ not identically zero, thus transforming $D$ by 
$D = (D\tilde{\theta})\tilde{D}$.  If $H(z,\theta) = \left( f(z) + 
\theta \xi(z), \psi(z) + \theta g(z) \right)$, then the condition 
(\ref{basic superconformal condition}) is equivalent to the conditions
\begin{equation}\label{superconformal condition}
\xi = g \psi, \quad \mbox{and} \quad g^2 = f' + \psi \psi',
\end{equation}  
with $D \tilde{\theta} = g(z) + \theta \psi ' (z)$ not identically 
zero.  Thus a superconformal function $H$ is uniquely determined by 
$f(z)$, $\psi(z)$ and a choice of a well-defined square root for the 
function $f' + \psi \psi'$, and $H$ must have the form
\begin{multline}\label{stupidconformal}
H(z,\theta) = \left(f(z) + \theta \psi(z) \sqrt{f'(z) + \psi(z)
\psi'(z)}, \right. \\
\left. \psi(z) + \theta \sqrt{f'(z) + \psi(z) \psi'(z)} 
\right) ,
\end{multline}
with $f'$ or $\psi'$ not identically zero.  Note that $\sqrt{f' + 
\psi \psi'}$ is a solution to (\ref{superconformal condition}), if and
only if $-\sqrt{f' + \psi \psi'}$ is also a solution.  Thus $H$ given 
by (\ref{stupidconformal}) is superconformal, if and only if 
$H(z,-\theta)$ is superconformal.  Hence for every even superfunction 
$f$ and odd superfunction $\psi$ with $f' + \psi \psi' \neq 0$ and with
a given well-defined square root for $f' + \psi \psi'$, there are 
exactly two distinct superconformal functions $H$ satisfying $\theta 
H (z, \theta) = \theta (f(z),\psi(z))$.  

By definition, if $f_B'(z_B) \neq 0$ then $\sqrt{f'(z) + \psi(z)
\psi'(z)}$ is the Taylor expansion of the square root about $f_B'
(z_B)$, and thus a well-defined square root for $f' + \psi \psi'$ only 
depends on a well-defined square root for $f_B'(z_B)$.  In addition, in 
the case $f'(z) \neq 0$, since the  square of any odd superfunction is 
zero, (\ref{stupidconformal}) can  be simplified to
\begin{equation}\label{supconf}
H(z,\theta) = \left(f(z) + \theta \psi(z) \sqrt{f'(z)}, \psi(z) +
\theta \sqrt{f'(z) + \psi(z) \psi'(z)} \right) .
\end{equation}

\begin{rema} By ``square root of $f' + \psi \psi'$", we mean 
a solution to (\ref{superconformal condition}) and not necessarily a
square root defined for some open subset of $\bigwedge_*$.  For instance,
$H(z,\theta) = (z^3/3, \theta z)$ and $H(z,\theta) = 
(z^3/3, -\theta z)$ are both superconformal for $(z,\theta) \in 
\bigwedge_{* > 0 }$, i.e., $z$ is a valid solution to $\sqrt{z^2}$, as 
is $-z$.  We could also choose a branch cut for the complex logarithm, 
thus defining a square root in the complex plane, and then let $H(z,
\theta) = (z^3/3, \theta \sqrt{z^2})$ where $\sqrt{z^2}$ is the 
Taylor expansion about $z_B$ using this complex square root.  However, 
then $H$ is superconformal only for supernumbers with body in the 
complement of the branch cut being used.  For instance, $H$ would
not be superconformal in a neighborhood of $0 \in \bigwedge_{*>0}$.
\end{rema}
 

Note that the space of superconformal functions on $\bigwedge_{*>0}$ 
is closed under composition when defined.  However, the sum of two 
superconformal functions is not in general superconformal.

\begin{rema}\label{Ds-superconformality 2} The above notion of
superconformal is based on the operator $D = \Dz$.  But the crucial
property of $D$ is that it satisfies $D^2 = \dz$.  Thus we could just 
as well have used the operator $D_s = s \frac{\partial}{\partial 
\theta} + \frac{1}{s} \theta \frac{\partial}{\partial z}$ for any $s 
\in (\bigwedge^0_{* - 1})^\times$ for $(z,\theta) \in \bigwedge_{*>0}$.  
This is equivalent to transforming $\theta$ by $\theta \leftrightarrow
\frac{1}{s} \theta$.  Then a $D_s$-superconformal superfunction $H_s$
is of the form
\begin{equation}\label{alternate superconformal}
H_s(z,\theta) = \left(f(z) + \frac{1}{s} \theta \psi(z) \sqrt{f'(z)},
\psi(z) + \frac{1}{s} \theta \sqrt{f'(z) + \psi(z) \psi'(z)} \right),
\end{equation}
for $f$ and $\psi$ superanalytic in $z$.  Note that given a 
$D_s$-superconformal function $H_s$, the transformation $s 
\leftrightarrow -s$ transforms the square root of $f' + \psi \psi'$ 
that $H_s$ defines to the corresponding negative square root and 
visa versa.  Thus (\ref{alternate superconformal}) can be thought of 
as defining a continuous deformation {}from a superconformal 
function $H_{s=1}(z,\theta)$ using a given well-defined square root of 
$f' + \psi \psi'$ to the superconformal function $H_{s= -1}(z,\theta) = 
H_{s=1}(z,-\theta)$ using the negative of this square root.  This 
deformation and the issue of alternate notions of superconformality 
will be studied in a subsequent paper.
\end{rema}

In Section 2.4, we will study ``super-Riemann spheres with punctures 
and local superconformal coordinates vanishing at the punctures''.  
These punctures can be thought of as being at $0 \in \bigwedge_{*>0}$, 
a non-zero point in $\bigwedge_{*>0}$, or at a distinguished point on
the supersphere we denote by ``$\infty$''.  As will be shown in
Section 2.4, we can always shift a non-zero point in $\bigwedge_{*>0}$
(or on the super-Riemann sphere) to zero via a global superconformal
transformation.  Thus all local superconformal coordinates vanishing
at the punctures can be expressed as power series vanishing at zero or
vanishing as $(z,\theta) = (z_B + z_S, \theta)  \longrightarrow (\infty 
+ 0, 0) = \infty$.  

If the puncture is at zero, we are interested in invertible 
superconformal functions $H(z,\theta)$ defined in a neighborhood of 
zero vanishing at zero.  Such an $H$ is of the form  (\ref{supconf}) 
where $f(z)$ and $\psi(z)$ are even and odd superanalytic functions, 
respectively, with $f(0) = 0$, $f'(0) \in \bigwedge_{ * - 
1}^\times$ and $\psi(0) = 0$.  Thus by (\ref{analytic definition}), 
$f$ and $\psi$ can be expanded as
\begin{eqnarray*}
f(z) &=& \sum_{j \in \mathbb{N}} a_j z^{j + 1}, \qquad \mbox{for $a_j 
\in\bigwedge_{ * - 1}^0$ and $a_0 \in (\bigwedge_{ * - 1}^0)^\times$}, \\ 
\psi(z) &=& \sum_{j \in \mathbb{N}} m_j z^{j + 1} , \qquad 
\mbox{for $m_j \in \bigwedge_{ * - 1}^1$} .
\end{eqnarray*}
Then by (\ref{superconformal condition}),
\[g^2(z) = \sum_{j \in \mathbb{N}} a_j (j + 1) z^j + \sum_{j,k \in
\mathbb{N}} m_j m_k (k + 1) z^{j + k + 1} . \] 
The condition that $a_0$ be invertible is necessary and sufficient for
$H$ to be invertible in a neighborhood of zero, and in this case, we
can factor $a_0$ out of the expression above so that
\begin{equation}\label{g squared}
g^2(z) = a_0 \biggl( 1 + \sum_{j \in \mathbb{Z}_+} \frac{a_j}{a_0} (j + 1) z^j + 
\sum_{j,k \in \mathbb{N}} \frac{m_j m_k}{a_0} (k + 1) z^{j + k + 1} \biggr) .
\end{equation}
Fix a branch of the complex logarithm, and let $\sqrt{(a_0)_B}$ be the
corresponding square root of $(a_0)_B$.  Then solving (\ref{g
squared}), we have
\[g(z) = \pm \sqrt{a_0} \Bigl(1 + \sum_{j \in \mathbb{N}} c_j z^{j + 1}\Bigr) \]
where $1 + \sum_{j \in \mathbb{N}} c_j z^{j + 1}$ is the expansion
about zero of the square root of the term in parenthesis on the
right-hand side of (\ref{g squared}) using the principal branch of log,
and $\sqrt{a_0}$ is by definition the expansion about
$\sqrt{(a_0)_B}$.

Thus if $H$ is superconformal and invertible in a neighborhood of
zero, and vanishing at zero, $H$ can be expanded as
\begin{multline}
H(z, \theta) = \Biggl( \sum_{j \in \mathbb{N}} a_j z^{j + 1} + \theta
\asqrt (\sum_{j \in \mathbb{N}} m_j z^{j + 1})( 1 + \sum_{j \in 
\mathbb{N}} c_j z^{j + 1}), \Biggr. \label{power series} \\
\Biggl. \sum_{j \in \mathbb{N}} m_j z^{j + 1} +
\theta \asqrt ( 1 + \sum_{j \in \mathbb{N}} c_j z^{j + 1})\Biggr) 
\end{multline} 
where $\asqrt \in (\bigwedge_{ * - 1}^0)^\times$ such that 
$\asqrt^2 = a_0$.  In other words, a superconformal function vanishing at 
zero and invertible in a neighborhood of zero is uniquely determined by
$\sum_{j \in \mathbb{N}} a_j z^{j + 1}$, $\sum_{j \in \mathbb{N}} m_j
z^{j + 1}$ and $\asqrt \in (\bigwedge_{ * -1}^0)^\times$
satisfying $\asqrt^2 = a_0$.  Note that by giving the data $\asqrt \in 
(\bigwedge_{ * -1}^0)^\times$, we need not specify a branch 
of the complex logarithm.

\begin{rema}\label{not using branch cut}
In \cite{B thesis} and \cite {B thesis announcement}, instead 
of specifying $\asqrt \in (\bigwedge_{ * - 1}^0)^\times$ 
such that $\asqrt^2 = a_0$ in the data describing a power series about 
zero, we fixed a branch cut for the complex logarithm and then 
specified whether we were using the positive or negative square root 
with respect to this branch cut for the square root of $a_0$.  However, 
it is more natural to not specify a square root and instead specify the
value of $\asqrt$ whose square is $a_0$.  We take this more natural 
approach throughout this paper.  
\end{rema}

Similarly, we would like to express a superconformal function vanishing as
$(z,\theta) \rightarrow (\infty,0) = \infty$ as a power series in $z$ and 
$\theta$.  The superfunction 
\begin{eqnarray*}
I: \mbox{$\bigwedge_{*>0}^\times$} &\longrightarrow& 
\mbox{$\bigwedge_{*>0}^\times$} \\ 
(z,\theta) & \mapsto & \Bigl(\frac{1}{z}, \frac{i \theta}{z} \Bigr)
\end{eqnarray*}
is superconformal, well defined and vanishing as $(z,\theta) 
\rightarrow \infty$. In fact, $H$ is superconformal, well defined and 
invertible in a neighborhood of $\infty$ and vanishing at $\infty$ if 
and only if $H(1/z, i\theta/z)$ is of the form 
(\ref{power series}), i.e.,
\begin{multline} \label{power series at infinity1}
H(z, \theta) = \Biggl( \sum_{j \in \mathbb{N}} a_j z^{-j - 1} +
\frac{i\theta}{z} \asqrt (\sum_{j \in \mathbb{N}} m_j z^{-j - 1}) ( 1 +
\sum_{j \in \mathbb{N}} c_j z^{-j - 1}), \Biggr.  \\
\Biggl. \sum_{j \in \mathbb{N}} m_j z^{-j - 1} +
\frac{i\theta}{z} \asqrt ( 1 + \sum_{j \in \mathbb{N}} c_j z^{-j - 1} )
\Biggr) 
\end{multline} 
where $1 + \sum_{j \in \mathbb{N}} c_j z^{-j - 1}$ is the power series
expansion of
\begin{equation}\label{square root at infinity}
\Bigl( 1 + \sum_{j \in \mathbb{Z}_+} \frac{a_j}{a_0} (j + 1) z^{-j} + \sum_{j,k \in
\mathbb{N}} \frac{m_j m_k}{a_0} (k + 1) z^{-j - k - 1}\Bigr)^{1/2}
\end{equation}
about $\infty$, and $\asqrt \in (\bigwedge_{ * -
1}^0)^\times$ such that $\asqrt^2 = a_0$.  Thus a superconformal 
function vanishing at $\infty$ and invertible in a neighborhood of
$\infty$ is uniquely determined by $\sum_{j \in \mathbb{N}} a_j
z^{-j - 1}$, $\sum_{j \in \mathbb{N}} m_j z^{-j - 1}$ and $\asqrt 
\in (\bigwedge_{ * - 1}^0)^\times$ satisfying $\asqrt^2 
= a_0$.

\section{Complex supermanifolds and super-Riemann surfaces}  

A {\em DeWitt $(m,n)$-dimensional topological superspace over
$\bigwedge_*$} is a topological space $X$ with a countable basis which
is locally homeomorphic to an open subset of $(\bigwedge_*^0)^m \oplus
(\bigwedge_*^1)^n$ in the DeWitt topology.  A {\em DeWitt
$(m,n)$-chart on $X$ over $\bigwedge_*$} is a pair $(U, \Omega)$ such
that $U$ is an open subset of $X$ and $\Omega$ is a homeomorphism of
$U$ onto an open subset of $(\bigwedge_*^0)^m \oplus
(\bigwedge_*^1)^n$ in the DeWitt topology.  A {\em superanalytic atlas
of DeWitt $(m,n)$-charts on $X$ over $\bigwedge_{* > n-1}$} is a
family of charts $\{(U_{\alpha}, \Omega_{\alpha})\}_{\alpha \in A}$
satisfying 

(i) Each $U_{\alpha}$ is open in $X$, and $\bigcup_{\alpha
\in A} U_{\alpha} = X$. 

(ii) Each $\Omega_{\alpha}$ is a homeomorphism {}from $U_{\alpha}$ to a
(DeWitt) open set in 
\[(\mbox{$\bigwedge_{* > n-1}^0$})^m \oplus (\mbox{$\bigwedge_{* > n-1}^1$})^n, \] 
such that
\[\Omega_{\alpha} \circ \Omega_{\beta}^{-1}: \Omega_{\beta}(U_\alpha
\cap U_\beta) \longrightarrow \Omega_{\alpha}(U_\alpha \cap U_\beta)\]
is superanalytic for all non-empty $U_{\alpha} \cap U_{\beta}$, i.e.,
$\Omega_{\alpha} \circ \Omega_{\beta}^{-1} = (\tilde{z}_1,...,
\tilde{z}_m, \tilde{\theta}_1,...,\tilde{\theta}_n)$ where $\tilde{z}_i$ 
is an even superanalytic $\bigwedge_{* > n-1}$-superfunction in 
$(m,n)$-variables for $i = 1,...,m$, and $\tilde{\theta}_j$ is an odd 
superanalytic $\bigwedge_{* >n-1}$-superfunction in $(m,n)$-variables 
for $j = 1,...,n$.  

Such an atlas is called {\em maximal} if, given any chart $(U, 
\Omega)$ such that
\[\Omega \circ \Omega_{\beta}^{-1} : \Omega_{\beta} (U \cap U_\beta)
\longrightarrow \Omega (U \cap U_\beta)\] 
is a superanalytic homeomorphism for all $\beta$, then $(U, \Omega)
\in \{(U_{\alpha}, \Omega_{\alpha})\}_{\alpha \in A}$.
 
A {\em DeWitt $(m,n)$-supermanifold over $\bigwedge_{* > n-1}$} is a
DeWitt $(m,n)$-dimensional topological space $M$ together with a
maximal superanalytic atlas of DeWitt $(m,n)$-charts over
$\bigwedge_{* > n-1}$.  

Given a DeWitt $(m,n)$-supermanifold $M$ over $\bigwedge_{* > n-1}$,
define an equivalence relation $\sim$ on M by letting $p \sim q$
if and only if there exists $\alpha \in A$ such that $p,q \in U_\alpha$
and $\pi_B^{(m,n)} (\Omega_\alpha (p)) = \pi_B^{(m,n)} (\Omega_\alpha (q))$
where $\pi_B^{(m,n)}$ is the projection given by (\ref{projection onto
body}).  Let $p_B$ denote the equivalence class of $p$ under this 
equivalence relation.  Define the {\it body} $M_B$ of $M$ to be the 
$m$-dimensional complex manifold with analytic structure given by the 
coordinate charts $\{((U_\alpha)_B, (\Omega_\alpha)_B) \}_{\alpha \in A}$ 
where $(U_\alpha)_B = \{ p_B \; | \; p \in U_\alpha \}$, and
$(\Omega_\alpha)_B : (U_\alpha)_B \longrightarrow \mathbb{C}^m$ is given
by $(\Omega_\alpha)_B (p_B) = \pi_B^{(m,n)} \circ \Omega_\alpha (p)$.

Note that $M$ is a complex fiber bundle over the complex manifold 
$M_{B}$.  The fiber is $(\bigwedge_{* > n-1}^0)_S^m \oplus 
(\bigwedge_{* > n-1}^1)^n$, a possibly infinite-dimensional vector space 
over $\mathbb{C}$.  This bundle is not in general a vector bundle
since the transition functions are not in general linear.

For any DeWitt $(1,n)$-supermanifold $M$, its body $M_{B}$ is a Riemann
surface. A {\em super-Riemann surface over $\bigwedge_{*>0}$} is a
DeWitt $(1,1)$-supermanifold over $\bigwedge_{*>0}$ with coordinate
atlas $\{(U_{\alpha}, \Omega_{\alpha})\}_{\alpha \in A}$ such that the
coordinate transition functions $\Omega_{\alpha} \circ
\Omega_{\beta}^{-1}$ in addition to being superanalytic are also
superconformal for all non-empty $U_{\alpha} \cap U_{\beta}$. Since 
the condition that the coordinate transition functions be superconformal 
instead of merely superanalytic is such a strong condition (unlike in 
the nonsuper case), we again stress the distinction between a 
supermanifold which has {\it superanalytic} transition functions versus 
a super-Riemann surface which has {\it superconformal} transition 
functions.  It would be perhaps more appropriate to refer to the later 
as a ``superconformal super-Riemann surface'' in order to avoid 
confusion.  In fact, in the literature one will find the term 
``super-Riemann surface" or ``Riemannian supermanifold" used for both  
merely superanalytic structures (cf. \cite{D}) and for superconformal  
structures (cf. \cite{Fd}, \cite{CR}). However, we will follow the  
terminology of \cite{Fd} and refer to a superconformal super-Riemann  
surface simply as a super-Riemann surface.
  
Next we show that if $M$ is a DeWitt $(1,1)$-supermanifold over 
$\bigwedge_\infty$, then for $L \in \mathbb{N}$ we can define a 
DeWitt $(1,1)$-supermanifold $M_L$ over $\bigwedge_L$ which can in 
some sense be thought of as a sub-supermanifold of $M$.  In 
addition, we can define a DeWitt $(1,1)$-dimensional topological 
superspace $M_L^{\mathbb{C}}$ over $\bigwedge_L$ for $L \in \mathbb{N}$ 
which as a $2^L$-dimensional complex manifold embeds in $M$.   

Recall that $\bigwedge_L \subset \bigwedge_\infty$ for $L \in 
\mathbb{N}$.  Let $\{\zeta_i\}_{i\in \mathbb{N}}$ be the fixed basis 
for $\bigwedge_\infty$, and let $i_{L,\infty} : \bigwedge_L 
\longrightarrow \bigwedge_\infty$ be the inclusion map. Denote the 
corresponding projection by
\begin{eqnarray*}
p_{\infty,L} : \mbox{$\bigwedge_\infty$} & \longrightarrow & \mbox{$\bigwedge_L$} \\
\zeta_i & \mapsto & \left\{ \begin{array}{ccl} \zeta_i & \mbox{for} &
                                                  i \leq L \\ 
                                                0 & \mbox{for} & 
                                                     i > L \\
                            \end{array} \right. .
\end{eqnarray*} 
Let $H$ be a superanalytic $\bigwedge_\infty$-superfunction in $(1,1)$-variables
defined on some DeWitt open set $U \subseteq \bigwedge_\infty$ with Laurent expansion
in $z \in \bigwedge_\infty^0$ given by
\[H (z,\theta) = \sum_{j \in \mathbb{Z}} a_j z^j +
\theta \sum_{j \in \mathbb{Z}} b_j z^j \]
for some $a_j, b_j \in \bigwedge_\infty$.
Define the superanalytic $\bigwedge_L$-superfunction $H_L$ in
$(1,1)$-variables on $U_L = p_{\infty,L}(U) \subset \bigwedge_L$ by
\[H_L (z,\theta) = \sum_{j \in \mathbb{Z}} p_{\infty, L - 1} (a_j) z^j +
\theta \sum_{j \in \mathbb{Z}} p_{\infty, L - 1} (b_j) z^j .\]

Let $M$ be a DeWitt $(1,1)$-supermanifold over $\bigwedge_\infty$ with
local coordinate atlas $(U_\alpha, \Omega_\alpha)_{\alpha \in A}$.
Define 
\[M_L = \Bigl( \bigsqcup_{\alpha \in A} p_{\infty,L} \circ \Omega_\alpha 
(U_\alpha) \Bigr)  / \approx \]
where $p \in p_{\infty,L} \circ \Omega_\alpha (U_\alpha)$ and 
$q \in p_{\infty,L} \circ \Omega_\beta (U_\beta)$ are equivalent
under the equivalence relation $\approx$ if and only if 
\[q = (\Omega_\beta \circ \Omega_\alpha^{-1})_L (p) .\]
Then $M_L$ is a DeWitt $(1,1)$-supermanifold over $\bigwedge_L$.
Note however that $M_L$ can not in general be embedded as a 
submanifold of $M$; the coordinate transition functions
of $M$ do not restrict to $M_L$ since $(\Omega_\beta \circ
\Omega_\alpha^{-1})_L$ is not a simple restriction of $\Omega_\beta 
\circ \Omega_\alpha^{-1}$ to $\bigwedge_L$.  If we instead form
\begin{equation}\label{super-subspace}
M_L^\mathbb{C} = \Bigl( \bigsqcup_{\alpha \in A} p_{\infty,L} \circ 
\Omega_\alpha (U_\alpha) \Bigr)  / \approx'
\end{equation}
where $p \in p_{\infty,L} \circ \Omega_\alpha (U_\alpha)$ and 
$q \in p_{\infty,L} \circ \Omega_\beta (U_\beta)$ are equivalent
under the equivalence relation $\approx'$ if and only if 
\[q = p_{\infty,L} \circ \Omega_\beta \circ \Omega_\alpha^{-1} 
\circ i_{L,\infty} (p) ,\]
then $M_L^\mathbb{C}$ is a DeWitt $(1,1)$-dimensional topological 
superspace over $\bigwedge_L$ which naturally embeds into $M$.  In 
addition, $M_L^\mathbb{C}$ is a $2^L$-dimensional complex analytic 
manifold but is not in general a DeWitt supermanifold due to the 
fact that the coordinate transition functions are not in general 
superanalytic.  

In addition to embedding into $M$, any $M_L^\mathbb{C}$ embeds into 
$M_{L'}^\mathbb{C}$ for $L,L' \in \mathbb{N}$ with $L' > L$.
Thus we have an infinite sequence of complex manifolds
\[M_B = M_0^\mathbb{C} \hookrightarrow M_1^\mathbb{C} \hookrightarrow 
M_2^\mathbb{C} \hookrightarrow \cdots \hookrightarrow M_{L-1}^\mathbb{C} 
\hookrightarrow M_L^\mathbb{C} \hookrightarrow M_{L+1}^\mathbb{C} 
\hookrightarrow \cdots ,\] 
and taking the direct limit of this sequence, we recover $M$, the DeWitt 
$(1,1)$-supermanifold over $\bigwedge_\infty$ (cf. \cite{Ro2}).
This construction can be extended to DeWitt $(m,n)$-supermanifolds
for general $m,n \in \mathbb{N}$ in the obvious way.

\section{Superspheres with tubes and the sewing operation} 

By {\em $N=1$ supersphere} we will mean a (superconformal) 
super-Riemann surface over $\bigwedge_{*>0}$ such that its body is a 
genus-zero one-dimensional connected compact complex manifold.  {}From
now on we will refer to such an object as a {\it supersphere}.

A {\em supersphere with $1+n$ tubes} for $n \in \mathbb{N}$, is a
supersphere $S$ with 1 negatively oriented point $p_0$ and $n$
positively oriented points $p_1,...,p_n$ (we call them {\em
punctures}) on $S$ which all have distinct bodies (i.e., $p_i$ is not
equivalent to $p_j$ for $i \neq j$ under the equivalence relation
$\sim$ defined in Section 2.3) and with local superconformal 
coordinates $(U_0,\Omega_0),...,(U_n, \Omega_n)$ vanishing at the 
punctures $p_0,...,p_n$, respectively.  We denote this structure by
\[(S;p_0,...,p_n;(U_0, \Omega_0),...,(U_n, \Omega_n)) . \]
We will always order the punctures so that the negatively oriented
puncture is $p_0$.  

\begin{rema} The reason we call a puncture with local 
superconformal coordinate vanishing at the puncture a ``tube" 
is that such a structure is indeed superconformally equivalent to a 
half-infinite superconformal tube representing an incoming (resp., 
outgoing) ``superparticle" or ``superstring" if the puncture is positively 
(resp., negatively) oriented.  For $r \in \mathbb{R}_+$, denote by  
\[ B_{z_B}^r = \{ w_B \in \mathbb{C} \; | \; |w_B - z_B| < r \} \qquad
(\mbox{resp.,} \; \; \bar{B}_{z_B}^r = \{ w_B \in \mathbb{C} \; | \; |w_B -
z_B| \leq r \}) \] 
an open (resp., closed) ball in the complex plane about the point
$z_B$ with radius $r$.  Denote a DeWitt open (resp., closed) ball 
in $\bigwedge_{*>0}$ about $(z, \theta)$ of radius $r$ by
\[ \mathcal{B}_z^r = B_{z_B}^r \times (\mbox{$\bigwedge_{*>0}$})_S \qquad
(\mbox{resp.,} \; \; \bar{\mathcal{B}}_z^r = \bar{B}_{z_B}^r \times 
(\mbox{$\bigwedge_{*>0}$})_S).\]
(Note that $\mathcal{B}_z^r$ depends only on $z_B$ and $r$.)  Let $p$ be
a positively oriented puncture with a local coordinate neighborhood $U$
and superconformal local coordinate map $\Omega: U \rightarrow
\bigwedge_{*>0}$ vanishing at the puncture.  Then for some $r \in 
\mathbb{R}_+$, we can find a DeWitt open disc $\mathcal{B}^r_0$ such that 
$\Omega^{-1}(\mathcal{B}^r_0) \subset U$.  Define the equivalence relation 
$\sim$ on $\bigwedge_{*>0}$ by $(z_1,\theta_1) \sim (z_2,\theta_2)$ if 
and only if $(z_1)_B = (z_2)_B + 2\pi i k$ for some integer $k$.  Then 
the set $\tau_r$ of all equivalence classes of elements of $(z,\theta) 
\in \bigwedge_{*>0}$ satisfying $\mathrm{Re}(z_B) < \log r$ (where
$\mathrm{Re}(z_B)$ is the real part of the complex number $z_B$)
together with the metric induced {}from the DeWitt metric on
$\bigwedge_{*>0}$ is a half-infinite tube in the body and is
topologically trivial in the soul.  Letting $H(z,\theta) = (\log z, 
\theta \sqrt{1/z})$, the map $q \mapsto H(\Omega(q))$ {}from
$\Omega^{-1}(\mathcal{B}^r_0)$ to $\tau_r$ is a well-defined invertible
superconformal map.  A closed curve on the supersphere shrinking to
$p$ corresponds to a closed loop or ``superstring" around this
half-infinite super-cylinder tending towards minus infinity in the body
coordinate.  We can perform a similar superconformal transformation
for the negative oriented puncture.
\end{rema}

\begin{rema}
In superconformal field theory, one generally wants to consider 
superspheres and higher genus super-Riemann surfaces with $m \in 
\mathbb{Z}_+$ negatively oriented (i.e., outgoing) tubes and
$n \in \mathbb{N}$ positively oriented (i.e., incoming) tubes.
However, for the purposes of this work, we restrict to genus zero
and $m=1$.
\end{rema}

Let
\[ (S_1;p_0,...,p_m;(U_0, \Omega_0),...,(U_m, \Omega_m)) \]
be a supersphere with $1+m$ tubes, for $m \in \mathbb{N}$, and let
\[ (S_2;q_0,...,q_n;(V_0, \Xi_0),...,(V_n, \Xi_n)) \]
be a superspheres with $1+n$ tubes, for $n \in \mathbb{N}$.  A map 
$F : S_1 \rightarrow S_2$ will be said to be superconformal if 
$\Xi_\beta \circ F \circ \Omega_\alpha^{-1}$ is superconformal for all 
charts $(U_\alpha, \Omega_\alpha)$ of $S_1$, for all charts $(V_\beta, 
\Xi_\beta)$ of $S_2$, and for all $(w, \rho) \in 
\Omega_\alpha (U_\alpha)$ such that $F \circ \Omega_\alpha^{-1} 
(w, \rho) \in V_\beta$.  If $m=n$ and there is a superconformal 
isomorphism $F : S_1 \rightarrow S_2$ such that for each $j = 0,..., 
n$, $F(p_j) = q_j$ and 
\[ \left. \Omega_j \right|_{W_j} = \left. \Xi_j \circ F \right|_{W_j} , \]
for $W_j$ some DeWitt neighborhood of $p_j$, then we say that these two 
superspheres with $1 + n$ tubes are {\em superconformally equivalent} and 
$F$ is a {\em superconformal equivalence} {}from   
\[ (S_1;p_0,...,p_n;(U_0, \Omega_0),...,(U_n, \Omega_n)) \]
to
\[ (S_2;q_0,...,q_n;(V_0, \Xi_0),...,(V_n, \Xi_n)) .\]
Thus the superconformal equivalence class of a supersphere with tubes
depends only on the supersphere, the punctures, and the germs of the local
coordinate maps vanishing at the punctures.

We now describe the sewing operation for two superspheres with tubes.
For the bodies of the two superspheres, the sewing will be identical
to the sewing operation defined in \cite{H thesis}, \cite{H book} for
spheres with tubes when the two superspheres can be sewn.  In general,
however, it is not true that if the bodies of two superspheres can be sewn 
then the two superspheres can be sewn since, as we shall see below, 
whether two superspheres (or spheres) can be sewn together depends on the
radius of convergence of the local coordinates of the punctures where the
sewing is taking place.  But in general, the radius of convergence of a 
superfunction may be smaller than the radius of convergence of its body
component.

Let
\[ (S_1;p_0,...,p_m;(U_0, \Omega_0),...,(U_m, \Omega_m)) \]
be a supersphere with $1+m$ tubes, for $m \in \mathbb{Z}_+$, and let 
\[ (S_2;q_0,...,q_n;(V_0, \Xi_0),...,(V_n, \Xi_n)) \]
be a superspheres with $1+n$ tubes, for $n \in \mathbb{N}$.  Assume 
that there exists a positive number $r$ such that the closed DeWitt ball 
about the origin of radius $r$ satisfies
\begin{equation}\label{open neighborhood at zero}
\bar{\mathcal{B}}_0^r = (\bar{B}_0^r \times (\mbox{$\bigwedge_{*>0}$})_S )
\subset \Omega_i(U_i), \quad \mbox{for some} \quad 0< i \leq m,
\end{equation} 
and the closed DeWitt ball about the origin of radius $1/r$ satisfies
\begin{equation}\label{open neighborhood at infinity} 
\bar{\mathcal{B}}_0^{1/r} = ( \bar{B}_0^{1/r} \times
(\mbox{$\bigwedge_{*>0}$})_S ) \subset \Xi_0 (V_0) ,
\end{equation}  
and also such that $p_i$ and $q_0$ are the only punctures in
$\Omega_i^{-1}(\bar{\mathcal{B}}_0^r)$ and $\Xi_0^{-1}
(\bar{\mathcal{B}}_0^{1/r})$, respectively.  In this case, we say 
that {\it the i-th puncture of the first supersphere with tubes, 
$S_1$, can be sewn with the 0-th puncture of the second supersphere 
with tubes, $S_2$}.  {}From these two superspheres with tubes, we 
obtain a supersphere with $1+(m + n - 1)$ tubes, i.e., with one 
negatively oriented tube and $m + n -1$ positively oriented tubes.  
We denote this supersphere by $S_1 \; _i\infty_0 \; S_2$, and it is 
obtained in the following way:

(i) By (\ref{open neighborhood at zero}) and (\ref{open neighborhood
at infinity}), there exist $r_1, r_2 \in \mathbb{R}_+$ satisfying
$0<r_2<r<r_1$ such that $\bar{\mathcal{B}}_0^{r_1} \subset \Omega_i(U_i)$
and $\bar{\mathcal{B}}_0^{1/r_2} \subset \Xi_0 (V_0)$.  Then the DeWitt
open subsets $S_1 \smallsetminus \Omega_i^{-1} (\bar{\mathcal{B}}_0^{r_2})$
and $S_2 \smallsetminus \Xi_0^{-1} (\bar{\mathcal{B}}_0^{1/r_1})$ of $S_1$
and $S_2$, respectively, are super-Riemann submanifolds of $S_1$ and $S_2$, 
respectively.  Let
\[\left(S_1 \smallsetminus \Omega_i^{-1} (\bar{\mathcal{B}}_0^{r_2})\right) 
\sqcup \bigl(S_2 \smallsetminus \Xi_0^{-1} (\bar{\mathcal{B}}_0^{1/r_1}) 
\bigr) \]
be the disjoint union of $S_1 \smallsetminus \Omega_i^{-1} 
(\bar{\mathcal{B}}_0^{r_2})$ and $S_2 \smallsetminus \Xi_0^{-1} 
(\bar{\mathcal{B}}_0^{1/r_1})$.  On this disjoint union we define the 
equivalence relation $\sim$ by setting $p \sim q$ for $p, q \in (S_1 
\smallsetminus \Omega_i^{-1} (\bar{\mathcal{B}}_0^{r_2})) 
\sqcup (S_2 \smallsetminus \Xi_0^{-1} (\bar{\mathcal{B}}_0^{1/r_1}))$ if and 
only if $p = q$ or $p \in \Omega_i^{-1} (\mathcal{B}_0^{r_1} \smallsetminus 
\bar{\mathcal{B}}_0^{r_2})$, $q \in \Xi_0^{-1} (\mathcal{B}_0^{1/r_2} 
\smallsetminus \bar{\mathcal{B}}_0^{1/r_1})$, and
\begin{equation}\label{equivalence condition}
\Xi_0^{-1} \circ I \circ \Omega_i (p) = q  
\end{equation} 
for $I(z,\theta) = (1/z,i\theta/z)$. Define the super-Riemann surface
$S_1 \; _i\infty_0 \; S_2$ to be the topological superspace
\[\Bigl( \bigl(S_1 \smallsetminus \Omega_i^{-1} (\bar{\mathcal{B}}_0^{r_2}) 
\bigr) \sqcup \bigl(S_2 \smallsetminus \Xi_0^{-1} (\bar{\mathcal{B}}_0^{1/r_1}) 
\bigr) \Bigr) / \sim \] 
with the superconformal structure determined by the superconformal
structures on $S_1 \smallsetminus \Omega_i^{-1} (\bar{\mathcal{B}}_0^{r_2})$ 
and $S_2 \smallsetminus \Xi_0^{-1} (\bar{\mathcal{B}}_0^{1/r_1})$ and the 
superconformal transition map
\[ \Xi_0^{-1} \circ I \circ \Omega_i : \Omega_i^{-1} (\mathcal{B}_0^{r_1} 
\smallsetminus \bar{\mathcal{B}}_0^{r_2}) \longrightarrow
\Xi_0^{-1} (\mathcal{B}_0^{1/r_2} \smallsetminus \bar{\mathcal{B}}_0^{1/r_1}) .\] 
In other words, $S_1 \; _i\infty_0 \; S_2$ is the union of $S_1 
\smallsetminus \Omega_i^{-1} (\bar{\mathcal{B}}_0^{r_2})$ and $S_2 
\smallsetminus \Xi_0^{-1} (\bar{\mathcal{B}}_0^{1/r_1})$ with the super-annulus 
$\Omega_i^{-1} (\mathcal{B}_0^{r_1} \smallsetminus 
\bar{\mathcal{B}}_0^{r_2}) \subset S_1 \smallsetminus \Omega_i^{-1} 
(\bar{\mathcal{B}}_0^{r_2})$ identified with the super-annulus $\Xi_0^{-1}
(\mathcal{B}_0^{1/r_2} \smallsetminus \bar{\mathcal{B}}_0^{1/r_1}) 
\subset S_2 \smallsetminus \Xi_0^{-1} (\bar{\mathcal{B}}_0^{1/r_1})$ via
the super-inversion map $I(z,\theta) = (1/z,i\theta/z)$. Note that $S_1
\smallsetminus \Omega_i^{-1} (\bar{\mathcal{B}}_0^{r_2})$ and $S_2 
\smallsetminus \Xi_0^{-1} (\bar{\mathcal{B}}_0^{1/r_1})$ are 
super-Riemann submanifolds of $S_1 \; _i\infty_0 \; S_2$ and that $S_1 \;  
_i\infty_0 \; S_2$ is a supersphere with $1 + (m+n-1)$ tubes.

(ii) The $1+ (m+n-1)$ ordered punctures of $S_1 \; _i\infty_0 \; S_2$ are
\[ p_0,p_1,...,p_{i - 1},q_1,...,q_n,p_{i + 1},...,p_m  \]
where $p_0$ is negatively oriented and the rest of the punctures are
positively oriented. 

(iii) The local coordinates vanishing at these punctures are 
\[\Bigl(U_k \smallsetminus \Omega_i^{-1} (\bar{\mathcal{B}}_0^r), \Omega_k |_{U_k 
\smallsetminus \Omega_i^{-1} (\bar{\mathcal{B}}_0^r)} \Bigr), \quad k \neq i,  \]
\[\Bigl(V_l \smallsetminus \Xi_0^{-1} (\bar{\mathcal{B}}_0^{1/r}), \Xi_l |_{V_l
\smallsetminus \Xi_0^{-1} (\bar{\mathcal{B}}_0^{1/r})} \Bigr), \quad l \neq 0 .\]

The above procedure to obtain a supersphere with $1+ (m + n -1)$ tubes 
{}from a supersphere $S_1$ with $1+m$ tubes and a supersphere $S_2$ with 
$1+n$ tubes is called the {\em sewing operation} for superspheres with 
tubes.  

\begin{rema} As in the non-super case of spheres with tubes, only
the local coordinate neighborhoods on the sewn supersphere $S_1 \;
_i\infty_0 \; S_2$ might depend on the positive number $r$.  Thus {}from
the definition of superconformal equivalence, we see that the
superconformal equivalence class of this supersphere with $1+(m + n -
1)$ tubes is independent of $r$.  It is also easy to see that this
superconformal equivalence class only depends on the superconformal
equivalence classes of
\[ (S_1;p_0,...,p_m;(U_0, \Omega_0),...,(U_m, \Omega_m)) \]
and
\[ (S_2;q_0,...,q_n;(V_0, \Xi_0),...,(V_n, \Xi_n)) .\]
\end{rema}

\section{The moduli space of superspheres with tubes}

The collection of all superconformal equivalence classes of
superspheres over $\bigwedge_{*>0}$ with $1+n$ tubes, for $n \in 
\mathbb{N}$, is called the {\it moduli space of superspheres over 
$\bigwedge_{*>0}$ with $1+n$ tubes}.  The collection of all 
superconformal equivalence classes of superspheres over 
$\bigwedge_{*>0}$ with tubes is called the {\it moduli space of 
superspheres over $\bigwedge_{*>0}$ with tubes}.

Let $S\hat{\mathbb{C}}$ be the supersphere with superconformal 
structure given by the covering of local coordinate neighborhoods 
$\{ U_\sou, U_\nor \}$ and the local coordinate maps
\begin{eqnarray*}
\sou : U_\sou  & \longrightarrow & \mbox{$\bigwedge_{*>0}$} \\
\nor: U_\nor & \longrightarrow & \mbox{$\bigwedge_{*>0}$} \\
\end{eqnarray*}
which are homeomorphisms of $U_\sou$ onto $\bigwedge_{*>0}$ and
$U_\nor$ onto $\bigwedge_{*>0}$, respectively, such that 
\begin{eqnarray*}
\sou \circ \nor^{-1} : \mbox{$\bigwedge_{*>0}^\times$} &\longrightarrow&
\mbox{$\bigwedge_{*>0}^\times$} \\ 
(w, \rho) & \mapsto & \Bigl(\frac{1}{w},\frac{i \rho}{w} \Bigr) =
I(w,\rho) .
\end{eqnarray*}
Thus the body of $S\hat{\mathbb{C}}$ is the Riemann sphere, 
$(S\hat{\mathbb{C}})_B = \hat{\mathbb{C}} = \mathbb{C} \cup \{\infty\}$,  
with coordinates $w_B$ near 0 and $1/w_B$ near $\infty$.   We will
call $S\hat{\mathbb{C}}$ the {\em super-Riemann sphere} and
will refer to $\nor^{-1} (0)$ as {\em the point at $(\infty, 
0)$} or just {\em the point at infinity} and to 
$\sou^{-1}(0)$ as {\em the point at $(0, 0)$} or just {\em the point at 
zero}.  

The Lie supergroup of superconformal isomorphisms of $S\hat{\mathbb{C}}$ 
is isomorphic to the connected component of $OSP(1|2)$ containing the 
identity (see \cite{D}, \cite{CR}).  This is the group of 
superprojective transformations and is given by
\[ \left\{ \left( \frac{aw + b}{cw + d} + \rho \frac{\gamma w +
\delta}{(cw + d)^2} , \frac{\gamma w + \delta}{cw + d} + \rho
\frac{1 + \frac{1}{2} \delta \gamma}{cw + d} \right)  \left| 
\begin{array}{c}
a,b,c,d \in \bigwedge_{ * - 1}^0, \\ \gamma, \delta \in \bigwedge_{ *
- 1}^1, \; ad - bc = 1
\end{array} \! \! \! \right. \! \right\}. \]
For such a superprojective transformation $T$, define
\[T_\sou : \mbox{$\bigwedge_{*>0}$} \smallsetminus \bigl( \{(- d/c)_B \} \times
(\mbox{$\bigwedge_{*>0}$})_S \bigr) \longrightarrow  \mbox{$\bigwedge_{*>0}$} 
\smallsetminus \bigl(\{(a/c)_B \} \times (\mbox{$\bigwedge_{*>0}$})_S \bigr) 
\hspace{.5in}\]  
\[\hspace{1.4in} (w, \rho) \mapsto \left( \frac{aw + b}{cw + d} + \rho 
\frac{\gamma w + \delta}{(cw + d)^2} , \frac{\gamma w + \delta}{cw + d} + 
\rho \frac{1 + \frac{1}{2} \delta \gamma }{cw + d} \right) , \]
and
\[T_\nor: \mbox{$\bigwedge_{*>0}$} \smallsetminus \bigl(\{(-a/b)_B \} \times
(\mbox{$\bigwedge_{*>0}$})_S \bigr) \longrightarrow  \mbox{$\bigwedge_{*>0}$} 
\smallsetminus (\{(d/b)_B \} \times \bigl(\mbox{$\bigwedge_{*>0}$})_S \bigr) 
\hspace{.5in}\]  
\[\hspace{1.3in} (w, \rho) \mapsto  \left( \frac{c + dw}{a + bw} - i\rho 
\frac{\gamma + \delta w}{(a + bw)^2} , - i\frac{\gamma  + \delta w}{a + bw} + 
\rho \frac{1 + \frac{1}{2} \delta \gamma }{a + bw} \right) ,\]
i.e., $T_\nor(w,\rho) =I^{-1} \circ T_\sou \circ I (w, \rho)$ for
$(w,\rho) \in \bigwedge_{*>0}^\times \smallsetminus (\{(-a/b)_B \} \times
(\mbox{$\bigwedge_{*>0}$})_S)$.

Define 
\begin{equation}\label{T1}
T(p) = \left\{
  \begin{array}{ll} 
      \sou^{-1} \circ T_\sou \circ \sou (p) & \mbox{if $p \in
           U_\sou \smallsetminus \sou^{-1}( \{(- d/c)_B \} \times
            (\bigwedge_{*>0})_S )$}, \\  
      \nor^{-1} \circ T_\nor \circ \nor (p) & \mbox{if $p \in
           U_\nor \smallsetminus \nor^{-1}( \{(-a/b)_B \} \times
             (\bigwedge_{*>0})_S )$ }.  
\end{array} \right.
\end{equation}
This defines $T$ for all $p \in S\hat{\mathbb{C}}$ unless \\
(i) $a_B = 0$ and $p \in \nor^{-1}( \{0 \} \times (\bigwedge_{*>0})_S)$.
In which case we define  \begin{equation}\label{T2}
T(p) = \sou^{-1} \left( \frac{a + bw}{c + dw} + i \rho \frac{\gamma + \delta w}{(c + dw)^2},
\frac{\gamma + \delta w}{c + dw} + i \rho \frac{1 + \frac{1}{2} \delta \gamma)}{c + dw} \right),
\end{equation} 
for $\nor (p) = (w,\rho) = (w_S, \rho)$; or \\
(ii) $d_B = 0$ and $p \in \sou^{-1}(\{0 \} \times (\bigwedge_{*>0})_S )$.
In which case we define  
\begin{equation}\label{T3}
T(p) = \nor^{-1} \left( \frac{cw + d}{aw + b} - \rho \frac{\gamma w + \delta}{(aw + b)^2}, 
-i \frac{\gamma w + \delta}{aw + b} - i \rho \frac{1 + \frac{1}{2} \delta \gamma}{aw + b}
\right) 
\end{equation} 
for $\sou (p) = (w,\rho) = (w_S, \rho)$.

Note that with this definition, $T$ is uniquely determined by
$T_\sou$, i.e., by its value on $\sou (U_\sou)$.

In \cite{CR}, Crane and Rabin prove the uniformization theorem for
super-Riemann surfaces.  We state the result for super-Riemann surfaces 
with genus-zero compact bodies.

\begin{thm}\label{uniformization theorem}(\cite{CR} Uniformization) 
Any super-Riemann surface with genus-zero compact body is
superconformally equivalent to the super-Riemann sphere
$S\hat{\mathbb{C}}$.
\end{thm}

Note that in our definition of the super-Riemann sphere
$S\hat{\mathbb{C}}$, we have chosen the transition function $\Delta
\circ \Upsilon^{-1} = I (z,\theta) = \left(1/z,i\theta/z \right)$, 
but we could just as well have chosen this transition to be 
$I(z, -\theta) = \left(1/z,-i\theta/z \right)$.  We will denote the
latter super-Riemann surface with genus-zero body by 
$(S\hat{\mathbb{C}})^-$ with coordinate charts $\{(U_\sou^-, \sou^-), 
(U_\nor^-, \nor^-)\}$ such that $\sou^- \circ (\nor^-)^{-1} (z, 
\theta) = I(z, -\theta)$.  These two superspheres are of course 
superconformally equivalent via $F: S\hat{\mathbb{C}} \rightarrow 
(S\hat{\mathbb{C}})^-$ defined by
\[F(p) = \left\{ \begin{array}{ll}
          (\sou^-)^{-1} \circ J \circ \sou(p) & \mbox{for $p \in U_\sou$} \\
          (\nor^-)^{-1} \circ \nor (p) & \mbox{for $p \in U_\nor$}
\end{array} \right. \]
where $J: \bigwedge_{*>0} \rightarrow \bigwedge_{*>0}$ is given 
by $J(z,\theta) = (z,-\theta)$.  In \cite{B thesis} and 
\cite{B thesis announcement}, we discuss the ``moduli space of 
superspheres with tubes with positive versus negative square root 
structure".  In this paper, we do not deal with the change of 
variables related to this symmetry, but rather leave that  
discussion for another paper.

\begin{rema}\label{infinite variables1}
In Chapter 4, we will want to consider functions on the moduli 
space of superspheres  with tubes which are superanalytic or 
supermeromorphic.  These superfunctions will in general involve 
an infinite number of odd variables --- not only the odd part of 
the coordinate for the finite number of punctures but also the 
possibly infinite amount of odd data involved in describing the 
local coordinates about the punctures.  (See Remarks 
\ref{remark giving power series} and \ref{infinite variables2} 
below.)  In this case, we need to work over $\bigwedge_\infty$ if we 
want all multiple partial derivatives with respect to these odd 
variables to be well defined.  However, in this work, it 
turns out that we only take up to two partial derivatives at a time 
with respect to odd variables before evaluating.  Thus we could work 
over a finite Grassmann algebra that allows us two extra degrees of 
freedom at any given time (or a few extra degrees of freedom just in 
case the need were to arise to take a few more partial  derivatives 
in a row before evaluating).  On the other hand, it is no harder to 
work over an infinite Grassmann algebra; the results we are 
interested in still hold, and then we may take partial derivatives 
without concern.  One may always later restrict to some
$\bigwedge_L$ for $L \in \mathbb{N}$ when substituting for these 
variables in the functional part of the theory or restrict to the 
supermanifold substructure defined in Section 2.3 for geometric 
aspects of the theory.  Thus for the remainder of this paper, we will   
mainly work over an infinite Grassmann algebra.
\end{rema}

Let 
\[\mathcal{B}^{r}_\infty = \Bigl\{ (w, \rho) \in \mbox{$\bigwedge_\infty$} 
\; \bigl| \; \Bigl|\frac{1}{w_B} \Bigr| < r  \bigr. \Bigr\} .\] 

\begin{prop}\label{canonicalcriteria}
Any supersphere over $\bigwedge_\infty$ with $1+n$ tubes for $n \in 
\mathbb{Z}_+$ is superconformally equivalent to a supersphere with 
$1+n$ tubes of the form
\begin{multline}\label{preliminarycanonical} 
\Bigl(S\hat{\mathbb{C}}; \nor^{-1}(0), \sou^{-1}(z_1,\theta_1),....., \sou^{-1}(z_{n
- 1}, \theta_{n - 1}), \sou^{-1}(0);  \Bigr. \\
\bigl(\nor^{-1}(\mathcal{B}^{1/r_0}_0), \Xi_0 \bigr), \bigl(\sou^{-1}(
\mathcal{B}^{r_1}_{z_1}), H_1 \circ \sou \bigr),...,\\
\Bigl. \bigl(\sou^{-1}(\mathcal{B}^{r_{n - 1}}_{z_{n - 1}}),H_{n -
1} \circ \sou \bigr), \bigl(\sou^{-1}(\mathcal{B}^{r_n}_0),H_n \circ \sou
\bigr) \Bigr) ,
\end{multline}
where 
\begin{equation}\label{at infinity for southern transform}
\Xi_0 |_{\sou^{-1}(\mathcal{B}^{r_0}_\infty)} = H_0 \circ \sou , 
\end{equation}
\begin{equation}\label{distinct punctures}
(z_1,\theta_1),....,(z_{n - 1}, \theta_{n - 1}) \in
\mbox{$\bigwedge_\infty^\times$}, \quad  \mbox{and} \quad (z_i)_B \neq 
(z_j)_B \quad \mbox{for} \quad i \neq j ,
\end{equation}   
\[r_0,..., r_n \in \mathbb{R}_+ = \{r \in \mathbb{R} \; | \; r > 0 \} , \]  
and $H_0, H_1,...,H_{n-1},H_n$ are superconformal functions
on $\mathcal{B}^{r_0}_\infty$, $\mathcal{B}^{r_1}_{z_1}$,...,
$\mathcal{B}^{r_{n - 1}}_{z_{n - 1}}$, $\mathcal{B}^{r_n}_0$, respectively, 
such that if we let $H_0 (w, \rho) = (\tilde{w_0}, \tilde{\rho}_0 )$, then
\begin{equation}\label{atinfinity}
\lim_{w \rightarrow \infty} H_0 (w, \rho) = 0, \quad \mbox{and} \quad
\lim_{w \rightarrow \infty} \frac{\partial}{\partial \rho} \frac{\partial}{\partial (\frac{1}{w})}
\tilde{\rho}_0 (w, \rho) = \lim_{w \rightarrow \infty} \frac{\partial}{\partial \rho} w
\tilde{\rho}_0 = i; 
\end{equation}
\begin{equation}\label{otherpunctures} 
H_j (z_j, \theta_j) = 0 , \quad \mbox{and} \quad
\Bigl. \frac{\partial}{\partial w} H_j (w, \theta_j) \Bigr|_{w = z_j} = 
\lim_{ w \rightarrow z_j} \frac{H_j(w, \theta_j)}{w - z_j} \in 
\mbox{$\bigwedge_\infty^\times$},  
\end{equation}
for $j = 1,..., n - 1$; and 
\begin{equation}\label{atzero}
H_n (0) = 0 , \quad \mbox{and} \quad \Bigl. \frac{\partial}{\partial w} 
H_n (w, 0) \Bigr|_{w = 0} = \lim_{ w \rightarrow 0} \frac{H_n(w,0)}{w} \in
\mbox{$\bigwedge_\infty^\times$}.
\end{equation}
\end{prop}

\begin{proof} By the uniformization theorem for
super-Riemann surfaces, Theorem \ref{uniformization theorem}, any
supersphere is superconformally isomorphic to the supersphere
$S\hat{\mathbb{C}}$.  Let
\begin{equation}\label{beginningsphere}
(S;p_0,...,p_n;(U_0, \Omega_0),...,(U_n, \Omega_n)) 
\end{equation} 
be a supersphere with $1+n$ tubes and $F : S \rightarrow
S\hat{\mathbb{C}}$ a superconformal isomorphism.  We have a
supersphere with $1+n$ tubes
\begin{equation}\label{nextsphere}
(S\hat{\mathbb{C}}; F(p_0),...,F(p_n);(F(U_0), \Omega_0 \circ
F^{-1}),...,(F(U_n), \Omega_n \circ F^{-1}))
\end{equation} 
which is superconformally equivalent to (\ref{beginningsphere}).  Let
\[\Omega_0 \circ F^{-1} \circ \sou^{-1} (w, \rho) = ((\Omega_0 \circ
F^{-1} \circ \sou^{-1})^0 (w, \rho), (\Omega_0 \circ F^{-1} \circ
\sou^{-1})^1 (w, \rho)) . \]
We have three cases:

(i) If $F(p_0), F(p_n) \in U_\sou$, let 
\[ \sou(F(p_0)) = (u_0,v_0) \in \mbox{$\bigwedge^0_\infty$} 
\oplus \mbox{$\bigwedge^1_\infty$} \]
\[ \sou(F(p_n)) = (u_n,v_n) \in \mbox{$\bigwedge^0_\infty$} 
\oplus \mbox{$\bigwedge^1_\infty$} . \] 
Then   
\[\lim_{v_0 \rightarrow 0} \frac{\partial}{\partial \rho} \Omega_0 \circ F^{-1} \circ \sou^{-1}
(u_0, \rho) = a \]
uniquely determines $a \in (\bigwedge^0_\infty)^\times$.  In this case, 
let $T : S\hat{\mathbb{C}} \rightarrow  S\hat{\mathbb{C}}$ such that 
\begin{multline*} 
\sou \circ T \circ \sou^{-1} (w, \rho) \; = \; T_\sou (w, \rho)  \\
= \; \left( \frac{-1}{a^2(u_n - u_0 - v_n v_0)} \cdot
\frac{w - u_n}{w - u_0} \; - \; \rho \; \frac{(v_0 - v_n)w + v_n u_0 - u_n v_0}{a^2(u_n -
u_0)(w - u_0)^2},\right. \\  
\left. - \; i \;
\frac{(v_0 - v_n)w + v_n u_0 - u_n v_0}{a (u_n - u_0)(w - u_0)} \; - \;
\frac{i\rho}{a (w - u_0)}
\right) .
\end{multline*}   

(ii) If $F(p_n) \notin U_\sou$, i.e., $\nor(F(p_n)) = (u_n', v_n') \in \bigwedge_\infty$
with $(u_n')_B = 0$, let $\sou(F(p_0)) = (u_0,v_0) \in \bigwedge_\infty$.  Then   
\[\lim_{v_0 \rightarrow 0} \frac{\partial}{\partial \rho} \Omega_0 \circ F^{-1} \circ \sou^{-1}
(u_0, \rho) = a \]
uniquely determines $a \in (\bigwedge^0_\infty)^\times$.  In this case, let
\begin{multline*}
T_\sou (w, \rho)  \\
= \; \left( \frac{-1}{a^2 (1 - u_n' u_0 - iv_n' v_0)} \cdot
\frac{u_n' w - 1}{w - u_0} \; - \; \rho \; \frac{(u_n' v_0 - iv_n')w + iv_n' u_0 -
v_0}{a^2(1 - u_n' u_0)(w - u_0)^2},\right. \\
\left. - \; i \;
\frac{(u_n' v_0 - iv_n')w + iv_n' u_0 - v_0}{a (1 -  u_n' u_0)(w - u_0)} \; - \;
\frac{i\rho}{a (w - u_0)}
\right) .
\end{multline*}  

(iii) If $F(p_0) \notin U_\sou$, i.e., $\nor(F(p_0)) = (u_0', v_0') \in \bigwedge_\infty$
with $(u_0')_B = 0$, let $\sou(F(p_n)) = (u_n,v_n)  \in \bigwedge_\infty$.  Then
\[\lim_{v_0' \rightarrow 0}\frac{\partial}{\partial \rho} \Omega_0 \circ
F^{-1} \circ \nor^{-1} (u_0', \rho) = a'\]  
uniquely determines $a' \in (\bigwedge_\infty^0)^\times$.  In this case, let
\begin{multline*}
T_\sou (w, \rho)  \\
= \; \left( \frac{1}{(a')^2 (1 - u_n u_0' + iv_n v_0')} \cdot
\frac{- w + u_n}{u_0'w - 1} \; + \; \rho \; \frac{(- iv_0' + u_0' v_n)w - 
v_n + i u_n v_0'}{(a')^2 (1 - u_n u_0')(u_0' w - 1)^2}, \right. \\
\left. \frac{(-i v_0' + u_0' v_n)w - v_n + i u_n v_0'}{a' (1 -  u_n u_0')
(u_0' w - 1)} \; + \; \frac{\rho}{a' ( u_0 'w - 1)} \right) .
\end{multline*} 

Then in each case (i) - (iii), $T_\sou$ uniquely defines $T: S\hat{\mathbb{C}}
\rightarrow S\hat{\mathbb{C}}$ by (\ref{T1}) - (\ref{T3}), and the supersphere 
with tubes (\ref{nextsphere}) is superconformally equivalent to
\begin{multline}\label{nexttolastsphere}
( S\hat{\mathbb{C}} ; \nor^{-1} (0), \sou^{-1}(z_1,\theta_1),.....,
\sou^{-1}(z_{n - 1}, \theta_{n - 1}), \sou^{-1}(0); \\
(T \circ F(U_0), \Omega_0 \circ F^{-1} \circ
T^{-1}),...,(T \circ F(U_n), \Omega_n \circ F^{-1}\circ T^{-1}))  
\end{multline} 
where
\[(z_j, \theta_j) = \sou \circ T \circ F(p_j),\]
for $j = 1,..., n - 1$.  Choose $r_0,...,r_n \in \mathbb{R}_+$ such that
\[\mathcal{B}^{r_0}_{\infty} \subset \sou \circ T \circ F(U_0) \]
\[\mathcal{B}^{r_j}_{z_j} \subset \sou \circ T \circ F(U_j) , \quad j = 
1,..., n - 1,  \]
\[\mathcal{B}^{r_n}_0 \subset \sou \circ T \circ F(U_n) . \]
Then the supersphere with tubes (\ref{nexttolastsphere}) is
superconformally equivalent to
\begin{multline*}
\Bigl( S\hat{\mathbb{C}} ; \nor^{-1} (0), \sou^{-1}(z_1,\theta_1),....., 
\sou^{-1}(z_{n - 1}, \theta_{n - 1}), \sou^{-1}(0); \\
\bigl(\sou^{-1} (\mathcal{B}^{r_0}_\infty) \cup \nor^{-1} (\{0 \} \times
(\mbox{$\bigwedge_\infty$})_S) , \Omega_0 \circ F^{-1}
\circ T^{-1}\bigr),\\
\bigl(\sou^{-1}(\mathcal{B}^{r_1}_{z_1}),\Omega_1 \circ F^{-1}\circ T^{-1}\bigr),
...,\bigl(\sou^{-1}(\mathcal{B}^{r_{n - 1}}_{z_{n - 1}}),
\Omega_{n - 2} \circ F^{-1}\circ T^{-1}\bigr), \\
\bigl(\sou^{-1}(\mathcal{B}^{r_n}_0), 
\Omega_n \circ F^{-1}\circ T^{-1}\bigr)\Bigr) 
\end{multline*} 
where
\begin{eqnarray*}
H_0 &=& \left. \Omega_0 \circ F^{-1} \circ T^{-1} \circ \sou^{-1} 
\right|_{\mathcal{B}^{r_0}_\infty} , \\
H_j &=& \left. \Omega_j \circ F^{-1} \circ T^{-1} \circ \sou^{-1} 
\right|_{\mathcal{B}^{r_j}_{z_j}} , \quad j = 1,..., n-1 , \\ 
H_n &=& \left. \Omega_n \circ F^{-1}\circ T^{-1} \circ \sou^{-1} 
\right|_{\mathcal{B}^{r_n}_0}  
\end{eqnarray*}
satisfy (\ref{atinfinity}), (\ref{otherpunctures}), and (\ref{atzero}),
respectively.
\end{proof}

A supersphere with $1+n$ tubes, for $n \in \mathbb{Z}_+$, of the form
(\ref{preliminarycanonical}) is called a {\em canonical supersphere
with $1+n$ tubes}.

\begin{rema}\label{remark giving power series}
A canonical supersphere with $1+n$ tubes, for $n \in \mathbb{Z}_+$, 
is determined by the punctures $(z_1,\theta_1),..., (z_{n - 1}, 
\theta_{n - 1}) \in \bigwedge_\infty^\times$ with $(z_i)_B \neq
(z_j)_B$ if $i \neq j$, the radii $r_0,...,r_n \in \mathbb{R}_+$ and 
the superconformal functions $H_0,...,H_n$ satisfying (\ref{atinfinity}), 
(\ref{otherpunctures}), and (\ref{atzero}), respectively. Consider the  
superconformal power series obtained by expanding the superconformal  
functions $H_0,...,H_n$ around $\nor^{-1} (0) = \infty$, $(z_1,\theta_1),
...,(z_{n-1},\theta_{n-1})$, and $0$, respectively.  We will denote by 
$H_i$ both the superconformal function and its power series expansion.  
By (\ref{power series}) and (\ref{power series at infinity1}), the
conditions (\ref{atinfinity}), (\ref{otherpunctures}), and
(\ref{atzero}) and the fact that the $H_j$'s are one-to-one as
superanalytic functions on their domains, we have
\begin{multline}\label{power series at infinity2}
H_0(w,\rho) =  \Biggl( \! \frac{1}{w} + \! \sum_{j \in \mathbb{Z}_+} 
a_j^{(0)} w^{-j-1} + \frac{i\rho}{w} \Bigl(\sum_{j \in \mathbb{N}} m_j^{(0)} 
w^{-j-1}\Bigr) \times \Biggr. \\
\Biggl.  \Bigl(1 + \! \sum_{j \in \mathbb{N}} c_j^{(0)} w^{-j-1} \Bigr), \;
\sum_{j \in \mathbb{N}} m_j^{(0)} w^{-j-1} +
\frac{i\rho}{w} \Bigl(1 + \! \sum_{j \in \mathbb{N}} c_j^{(0)} w^{-j-1}\Bigr) \! \Biggr) ,
\end{multline}
\begin{multline}\label{power series at ith puncture}
H_i(w,\rho) = \Biggl( \! (\asqrt^{(i)})^2 \biggl(w + \! \sum_{j \in 
\mathbb{Z}_+} a_j^{(i)} w^{j+1} + \rho \Bigl(\sum_{j \in \mathbb{N}} m_j^{(i)} 
w^{j+1} \Bigr) \times \biggr. \Biggr.\\
\Biggl. \left. \biggl. \Bigl(1 + \! \sum_{j \in \mathbb{N}} c_j^{(i)} w^{j+1}
\Bigr) \! \biggr) , \; \asqrt^{(i)} \biggl(\sum_{j \in \mathbb{N}} m_j^{(i)} 
w^{j+1} + \rho \Bigl(1 + \! \sum_{j \in \mathbb{N}} c_j^{(i)} w^{j+1} \Bigr) \! \biggr) \! \!
\Biggr) \right|_{ \! \! \! \begin{scriptsize} \begin{array}{l} 
(w,\rho) = (w - z_i \\
- \rho \theta_i, \rho -
\theta_i) \end{array} \end{scriptsize} } 
\end{multline}
for $i = 1,...,n$, where $(z_n, \theta_n) = 0$, $\asqrt^{(i)} \in
(\bigwedge_\infty^0)^\times$, $a_j^{(i)} \in \bigwedge_\infty^0$, for $j\in
\mathbb{Z}_+$, $m_j^{(i)} \in \bigwedge_\infty^1$, for $j\in \mathbb{N}$, and 
where $\left(1 + \sum_{j \in \mathbb{N}} c_j^{(0)} w^{-j-1} \right)$ is the 
power series expansion about infinity of
\[\Bigl(1 + \sum_{j \in \mathbb{Z}_+} (j + 1)a_j^{(0)} w^{-j} + \sum_{j,k \in \mathbb{N}} (k
+ 1) m_j^{(0)} m_k^{(0)} w^{- j - k - 1} \Bigr)^{1/2} , \] 
and $\left(1 + \sum_{j \in \mathbb{N}} c_j^{(i)} w^{j+1} \right)$ is
the power series expansion about zero of
\[\Bigl(1 + \sum_{j \in \mathbb{Z}_+} (j + 1) a_j^{(i)} w^j + \sum_{j,k \in \mathbb{N}}
(k + 1) m_j^{(i)} m_k^{(i)} w^{j + k + 1}\Bigr)^{1/2} . \]  
Thus a canonical supersphere with $1+n$ tubes, for $n \in \mathbb{Z}_+$, 
can be denoted by
\begin{equation}
((z_1,\theta_1),..., (z_{n - 1}, \theta_{n - 1});r_0,...,r_n;H_0,...,H_n)  
\end{equation} 
where $H_0,...,H_n$ are power series of the form (\ref{power series at
infinity2}) and (\ref{power series at ith puncture}).
\end{rema}

\begin{rema}\label{infinite variables2}
{}From the Remark \ref{remark giving power series} above, we can readily see 
that a point in the moduli space of superspheres with $1+n$ tubes, for
$n \in \mathbb{Z}_+$, will in general depend on an infinite number of odd 
variables --- the $\theta_1,...,\theta_{n-1}$ and the $m_j^{(i)} \in 
\bigwedge_\infty^1$, for $i=0,...,n$, and $j\in \mathbb{N}$.
\end{rema}

\begin{rema} \label{factoring out a in Chapter 2}
In (\ref{power series at ith puncture}), we have factored 
out $(\asqrt^{(i)})^2$ {}from the even part of $H_i$ and $\asqrt^{(i)}$ 
{}from the odd part (cf. equation (\ref{power series})).  This is due to 
the unique role that the $\asqrt^{(i)}$'s play, not only in determining 
which square root is involved in the superconformal structure of $H_i$ 
(see Remark \ref{not using branch cut}), but  also in determining the 
scaling of the local coordinate.  Factoring out this scaling operator is
necessary if one wants to achieve a certain symmetry in expressing the
infinitesimal local coordinate transformations at the punctures (see 
Remark \ref{factoring out a in Chapter 3}).
\end{rema}

\begin{prop}\label{equalspheres}
Two canonical superspheres with $1+n$ tubes, for $n \in \mathbb{Z}_+$,
\begin{equation}\label{stupidsphere1}
((z_1,\theta_1),..., (z_{n - 1}, \theta_{n - 1});r_0,...,r_n;H_0,...,H_n) 
\end{equation}
and 
\begin{equation}\label{stupidsphere2}
((\hat{z}_1,\hat{\theta}_1),..., (\hat{z}_{n - 1}, \hat{\theta}_{n -
1});\hat{r}_0,...,\hat{r}_n;\hat{H}_0,...,\hat{H}_n) 
\end{equation}
are superconformally equivalent if and only if $(z_j,\theta_j) = (\hat{z}_j,
\hat{\theta}_j)$ for $j = 1,..., n-1$, and $H_j = \hat{H_j}$, for $j = 0,..., 
n$, as superconformal power series.  
\end{prop}

\begin{proof} Let $F$ be a superconformal equivalence 
{}from (\ref{stupidsphere1}) to (\ref{stupidsphere2}).  The conclusion
of the proposition is equivalent to the assertion that $F$ must be the
identity map on $S\hat{\mathbb{C}}$.  By definition $F$ is a
superconformal automorphism of $S\hat{\mathbb{C}}$, i.e., a
superprojective transformation.  Also by definition we have
\begin{equation}\label{eq1}
F_\sou (0) = \sou \circ F \circ \sou^{-1} (0) = 0 , 
\end{equation}
\begin{equation}\label{eq2}
F_\nor (0) = \nor \circ F \circ \nor^{-1} (0) = 0 , 
\end{equation}
\begin{equation}\label{eq3}
\hat{H}_0 |_{\mathcal{B}^{\min (r_0, \hat{r}_0)}_\infty} = H_0
\circ F_\sou^{-1} |_{\mathcal{B}^{\min (r_0, \hat{r}_0)}_{\infty}} .  
\end{equation} {}From (\ref{eq1}) and (\ref{eq2}) and the fact that $F$ is a superprojective
transformation, we obtain
\begin{equation}\label{eq4}
F_\sou (w, \rho) = (a^2w, a \rho)
\end{equation}
for some $a \in (\bigwedge_\infty^0)^\times$.  Let $H_0
(w,\rho) = (\tilde{w}_0, \tilde{\rho}_0)$ and $\hat{H}_0 (w,\rho) =
(\hat{w}_0, \hat{\rho}_0)$. {}From (\ref{atinfinity}), we know that
\[\lim_{w \rightarrow \infty} \frac{\partial}{\partial \rho} w 
\tilde{\rho}_0 (w, \rho) = \lim_{w \rightarrow \infty} 
\frac{\partial}{\partial \rho} w \hat{\rho}_0 (w,\rho) = i . \] 
Thus by (\ref{eq3}) and (\ref{eq4})
\begin{eqnarray*}
i &=& \lim_{w \rightarrow \infty} \frac{\partial}{\partial \rho} w
\hat{\rho}_0 (w, \rho) =\lim_{w \rightarrow \infty}
\frac{\partial}{\partial \rho} w (H_0 \circ F_\sou^{-1})^1 (w, \rho)\\
&=& \lim_{w \rightarrow \infty} \frac{\partial}{\partial \rho} w
\tilde{\rho}_0 \Bigl(\frac{w}{a^2} , \frac{\rho}{a}\Bigr) \\ 
&=& \frac{i}{a},
\end{eqnarray*}
i.e., $a = 1$. Thus $F$ must be the identity map of
$S\hat{\mathbb{C}}$. \end{proof}

For superspheres with one tube, we have:

\begin{prop}\label{canonical criteria for one tube}
Any supersphere with one tube is superconformally equivalent to a
supersphere with one tube of the form
\begin{equation}\label{onetube}
\bigl(S\hat{\mathbb{C}}; \nor^{-1} (0) ; (\nor^{-1} (\mathcal{B}^{1/r_0}_0), 
\Xi_0) \bigr) 
\end{equation}
where $\Xi_0 |_{(\sou^{-1} (\mathcal{B}^{r_0}_{\infty})} = H_0 \circ
\Delta$, and $H_0$ can be expanded in a power series about infinity 
of the form (\ref{power series at infinity2}) with $a_1^{(0)} = 
m_0^{(0)} = 0$, i.e., such that the even coefficient of $w^{-2}$ and 
the odd coefficient of $w^{-1}$ in $H_0$ are zero.
\end{prop}

\begin{proof} Given a supersphere with one tube 
\begin{equation}\label{onetube1}
(S; p ; (U, \Omega)),
\end{equation}
by the uniformization theorem for super-Riemann surfaces, Theorem
\ref{uniformization theorem}, there exists a superconformal
isomorphism $F: S \rightarrow S\hat{\mathbb{C}}$ such that
\begin{equation}\label{onetube2}
(S\hat{\mathbb{C}}; F(p);(F(U), \Omega \circ F^{-1})) 
\end{equation}
is superconformally equivalent to (\ref{onetube1}).  By cases (i) and
(ii) in the proof of Proposition \ref{canonicalcriteria}, we know 
there exists a (non-unique) superprojective transformation $T$ {}from 
(\ref{onetube2}) to $(S\hat{\mathbb{C}}; T \circ F(p);(T \circ F(U), 
\Omega \circ F^{-1} \circ T^{-1}))$ such that $\nor \circ T \circ F(p) 
= 0$, and there exists $r \in \mathbb{R}_+$ satisfying $\sou^{-1} 
(\mathcal{B}_\infty^r) \subseteq T \circ F(U)$, and $\left. \Omega 
\circ F^{-1} \circ T^{-1} \circ \sou^{-1} \right|_{\mathcal{B}_\infty^r} 
= H$ where $H(w,\rho)$ satisfies (\ref{atinfinity}), i.e., has a power 
series expansion of the form (\ref{power series at infinity2}).

Let $\hat{T}_\sou (w,\rho) = (w - a_1 - i \rho m_0, -im_0 + \rho)$
where $a_1$ is the even coefficient of $w^{-2}$ and $m_0$ is the odd
coefficient of $w^{-1}$ in the power series expansion of $H$ about
infinity.  Then $H_0 = \Omega \circ F^{-1} \circ T^{-1} \circ
\hat{T}^{-1} \circ \sou^{-1}$ has a power series expansion of the
form (\ref{power series at infinity2}) with the even coefficient of 
$w^{-2}$ and the odd coefficient of $w^{-1}$ equal to zero, and there 
exists some $r_0 \in \mathbb{R}_+$ such that $H_0$ is convergent in 
$\mathcal{B}_\infty^{r_0}$. \end{proof}

A supersphere with one tube of the form (\ref{onetube}) is called a
{\em canonical supersphere with one tube}.  A canonical supersphere
with one tube is determined by $r_0 \in \mathbb{R}_+$ and a
superconformal power series $H_0$ satisfying (\ref{power series at
infinity2}) with $a_1^{(0)} = m_0^{(0)} = 0$, and can be denoted by
$(r_0;H_0)$. The following proposition can be proved similarly to Proposition
\ref{equalspheres}:

\begin{prop}\label{equal spheres with one tube}
Two canonical superspheres with one tube $(r_0;H_0)$ and $(\hat{r_0};
\hat{H_0})$ are superconformally equivalent if and only if $H_0 =
\hat{H_0}$.  
\end{prop} 

{}From Propositions \ref{canonicalcriteria}, \ref{equalspheres},
\ref{canonical criteria for one tube}, and \ref{equal spheres with one
tube} we have the following corollary.

\begin{cor}\label{bijection} 
There is a bijection between the set of canonical superspheres with 
tubes and the moduli space of superspheres with tubes.  In particular, 
the moduli space of superspheres with $1+n$ tubes, for $n \in \Z$, can 
be identified with all $2n$-tuples $((z_1,\theta_1),...,(z_{n-1},
\theta_{n-1}); H_0,...,H_n)$ satisfying $(z_i,\theta_i) \in 
\bigwedge_\infty^\times$, with $(z_i)_B \neq (z_j)_B$ if $i \neq j$, 
for $i,j = 1,...,n-1$, and such that $H_0,...,H_n$ are of the form 
(\ref{power series at infinity2}) and (\ref{power series at ith puncture}),
respectively, and are absolutely convergent in neighborhoods of
$\infty$, $(z_1,\theta_1),...,(z_{n-1},\theta_{n-1})$, and $0$,
respectively.  The moduli space of superspheres with one tube
can be identified with the set of all power series $H_0$ of the
form (\ref{power series at infinity2}) such that $a_1^{(0)} = 
m_0^{(0)} = 0$ and such that $H_0$ is absolutely convergent in 
a neighborhood of $\infty$. \footnote{There is a misprint in the nonsuper 
analogue of Corollary \ref{bijection} in \cite{H book}.  Proposition 
1.3.8 in \cite{H book} should state that the local coordinate $f_0(w)$ 
(where here $w$ is simply a complex variable) is absolutely convergent 
in a neighborhood of $w = \infty$, not $w = 0$ as stated.  In addition, 
it should state that $f_n(w)$ is absolutely convergent in a 
neighborhood of $w = 0$; it is stated as being absolutely convergent in
a neighborhood of $w=z_n$.  Corrected, the statements correspond to 
$(H_0)_B (w_B,0)$ being convergent in a neighborhood of $w_B = \infty$ 
and $(H_n)_B (w_B,0)$ being convergent in a neighborhood of $w_B = 0$ in
Corollary \ref{bijection} above.}
\end{cor}

\section{The sewing equation}

Let
\[C_1 = \bigl((z_1,\theta_1),...,(z_{m-1},\theta_{m-1});H_0,...,H_m\bigr)\]
and 
\[C_2 = \bigl((\hat{z}_1,\hat{\theta}_1),...,(\hat{z}_{n-1},\hat{\theta}_{n-1});
\hat{H}_0,...,\hat{H}_n\bigr)\] 
represent canonical superspheres with $1+m$ tubes and $1+n$ tubes,
respectively, for $m \in \mathbb{Z}_+$ and $n \in \mathbb{N}$.
Suppose that it is possible to sew the 0-th puncture of $C_2$ to the
$i$-th puncture of $C_1$ as described in Section 2.4.  Then the
resulting sewn supersphere $C_1 \; _i\infty_0 \; C_2$ is
superconformally equivalent to some canonical supersphere
\footnote{There is a misprint in \cite{H book}
on p.28 in the nonsuper analogue to this statement.  The 
resulting sewn sphere $C_1 \; _i\infty_0 \; C_2$ should have
local coordinates $g_0,...,g_{m+n-1}$ vanishing at the $1 + (m+n-1)$ 
punctures, not $g_1,...,g_{m+n-1}$ as stated.  Corrected, this
corresponds to the bodies of the local coordinates, denoted
$(H_0^*)_B,...,(H_{m+n-1}^*)_B$ in our notation, in the super case.} 
\[C_3= \bigl((z_1^*,\theta_1^*),...,(z_{m+n-2}^*,\theta_{m+n-2}^*);
H_0^*,...,H_{m+n-1}^* \bigr).\]

Since $C_1$ and $C_2$ correspond to canonical superspheres, $H_i$ is
a convergent power series in some neighborhood $U_i$ of $(z_i,
\theta_i)$, and $\hat{H}_0$ is a convergent power series in some
neighborhood $V_0$ of $\infty$.  According to the definition of 
sewing, since $C_2$ can be sewn to the $i$-th puncture of $C_1$, 
there exists $r \in \mathbb{R}_+$ such that $\mathcal{\bar{B}}_0^r 
\subset H_i(U_i)$ and $\mathcal{\bar{B}}_0^{1/r} \subset \hat{H}_0(V_0)$ 
and such that $\sou^{-1}
\circ  H_i^{-1}(\bar{\mathcal{B}}_0^r)$ contains only the puncture at
$(z_i,\theta_i)$ and $\sou^{-1} \circ \hat{H}_0^{-1}(\bar{\mathcal{B}}_0^{1/r} 
\smallsetminus (\{0\} \times (\bigwedge_\infty)_S))$
contains no punctures.  Furthermore, there exist $r_1, r_2 \in
\mathbb{R}_+$ such that $0< r_2< r < r_1$, and letting $U \subset
C_2$ be given by
\[U =  \sou^{-1} \circ \hat{H}_0^{-1} \bigl(\bar{\mathcal{B}}_0^{1/r_1} 
\smallsetminus \bigl(\{0\} \times (\mbox{$\bigwedge_\infty$})_S \bigr)
\bigr) \cup \nor^{-1} \bigl( \{0\} \times (\mbox{$\bigwedge_\infty$})_S 
\bigr) ,\]   
we have
\[C_1 \; _i\infty_0 \; C_2 = \Bigl(\bigl(C_1 \smallsetminus \sou^{-1} 
\circ H_i^{-1} (\bar{\mathcal{B}}_0^{r_2}) \bigr) \sqcup (C_2 
\smallsetminus U) \Bigr)/ \sim\] 
where $p,q \in (C_1 \smallsetminus \sou^{-1} \circ H_i^{-1} 
(\bar{\mathcal{B}}_0^{r_2})) \sqcup (C_2 \smallsetminus U)$ are equivalent 
if and only if $p=q$ or $p \in \sou^{-1} \circ H_i^{-1} 
(\mathcal{B}_0^{r_1} \smallsetminus \bar{\mathcal{B}}_0^{r_2})$, $q 
\in \sou^{-1} \circ \hat{H}_0^{-1} (\mathcal{B}_0^{1/r_2} 
\smallsetminus \bar{\mathcal{B}}_0^{1/r_1})$, and
\begin{equation}\label{preliminary sewing equation}
\sou^{-1} \circ \hat{H}_0^{-1} \circ I \circ H_i \circ \sou (p) = q  .
\end{equation}

Let $F$ be the unique superconformal equivalence taking $C_1 \;
_i\infty_0 \; C_2$ to $C_3$.  Then we can express $F$ as two
functions, one $F^{(1)}$ {}from $C_1 \smallsetminus \sou^{-1} \circ
H_i^{-1} (\bar{\mathcal{B}}_0^{r_2})$ to $C_3$ and the other $F^{(2)}$ 
{}from $C_2 \smallsetminus U$ to $C_3$ such that
\begin{eqnarray}\label{uniformizing function}
F: C_1 \; _i\infty_0 \; C_2 &\longrightarrow& C_3 \\
q &\mapsto& \left\{ \begin{array}{ll}
                              F^{(1)}(q) & \mbox{for $q \in C_1 \smallsetminus \sou^{-1} \circ
H_i^{-1} (\bar{\mathcal{B}}_0^{r_2})$}\\
                              F^{(2)}(q) & \mbox{for $q \in C_2 \smallsetminus U$}.
\end{array} \right. \nonumber
\end{eqnarray}
{}From the definitions of the sewing operation and canonical
supersphere with tubes, the maps $F^{(1)}$ and $F^{(2)}$ must satisfy
the following:

(i) For $p \in \sou^{-1} \circ H_i^{-1} (\mathcal{B}_0^{r_1}
\smallsetminus \bar{\mathcal{B}}_0^{r_2}) \subset C_1$ and $q \in
\sou^{-1} \circ \hat{H}^{-1}_0 ( \mathcal{B}_0^{1/r_2}
\smallsetminus \bar{\mathcal{B}}_0^{1/r_1}) \subset C_2$, if $p \sim q$
in $C_1 \; _i\infty_0 \; C_2$, then we must have $F(p) = F(q)$ in
$C_3$.  Thus by (\ref{preliminary sewing equation}), we must have
\begin{equation}\label{sewing equation}
F^{(1)}(p) = F^{(2)} \circ \sou^{-1} \circ \hat{H}_0^{-1} \circ I
\circ H_i \circ \sou (p) .
\end{equation}
Equation (\ref{sewing equation}) is called the {\it sewing equation},
and the function $F$ is called the {\it uniformizing function}.  In terms
of the local coordinate chart $(U_\sou,\sou)$ of $S\hat{\mathbb{C}}$,
equation (\ref{sewing equation}) is equivalent to 
\begin{equation}\label{local sewing equation}
F_\sou^{(1)} (w,\rho) =  F^{(2)}_\sou  \circ \hat{H}_0^{-1} \circ I
\circ H_i (w,\rho) 
\end{equation} 
for $F_\sou^{(1)} = \sou \circ F^{(1)} \circ \sou^{-1}$, $F^{(2)}_\sou
= \sou \circ F^{(2)} \circ \sou^{-1}$, and $(w,\rho) \in H_i^{-1} 
(\mathcal{B}_0^{r_1}\smallsetminus \bar{\mathcal{B}}_0^{r_2})$.

(ii) In order for $C_3$ to be canonical, the first puncture of $C_3$
must be at $\nor^{-1} (0)$, i.e., at infinity, the last puncture must
be at $\sou^{-1}(0)$, i.e., at zero, and the coordinates at the punctures
must satisfy the conditions (\ref{atinfinity}) - (\ref{atzero}).  If 
$i = m$, that is if the second supersphere is being sewn into the last 
puncture of the first supersphere which is the puncture at zero, and if 
$n \neq 0$, then requiring that the resulting supersphere be canonical 
is equivalent to the conditions
\begin{eqnarray}
F^{(1)}_\nor(0) &=& 0 \label{normalization 1}\\
\lim_{w\rightarrow \infty} \frac{\partial}{\partial \rho} (F_\sou^{(1)})^1
(w,\rho) &=& 1 \label{normalization 2} \\ 
F_\sou^{(2)}(0) &=& 0 \label{normalization 3}
\end{eqnarray}
where $F^{(1)}_\nor = \nor \circ F^{(1)} \circ \nor^{-1}$.  Conditions 
(\ref{normalization 1}) - (\ref{normalization 3}) are called the {\it 
normalization conditions} for $F$.

If $i = m$ and $n \neq 0$, then the resulting canonical supersphere
$F(C_1 \; _i\infty_0 \; C_2)$ $= C_3$ is represented by 
\begin{multline*}
C_3 = \bigl(F^{(1)}_\sou(z_1,\theta_1),...,F^{(1)}_\sou(z_{m-1},
\theta_{m-1}),F^{(2)}_\sou (\hat{z}_1, \hat{\theta}_1),...,F^{(2)}_\sou(\hat{z}_{n-1},
\hat{\theta}_{n-1});  \\ 
H_0 \circ (F^{(1)}_\sou)^{-1}, ..., 
H_{m-1} \circ (F^{(1)}_\sou)^{-1}, \hat{H}_1 \circ
(F^{(2)}_\sou)^{-1},...,\hat{H}_n \circ (F^{(2)}_\sou)^{-1}\bigr).
\end{multline*}

Conditions (i) and (ii) completely determine $F$, and if $i=m$, and $n
\neq 0$, then $F(C_1 \; _i\infty_0 \; C_2) = C_3$ is a canonical
supersphere with $1 + (m + n -1)$ tubes.  If $i \neq m$ or $n = 0$, we
can still use $F$ determined by conditions (i) and (ii) to map $C_1 \; 
_i\infty_0 \; C_2$ to a super-Riemann sphere $C_3$, but this supersphere 
might not be canonical.  If $n=0$ and $i = m$ with $m > 1$, or if $n 
\neq m$, the last puncture of the supersphere will not be at zero since 
$F$ sends the zero of the second supersphere to zero.  And if $n = 0$ 
and $m = 1$ the resulting sewn supersphere will have only one puncture
(the outgoing or negatively oriented puncture at $\nor^{-1}(0)=\infty$), 
and the coordinate at $\infty$ might not satisfy the extra condition 
given in Proposition \ref{canonical criteria for one tube} which 
specifies that the even coefficient of $w^{-2}$ and the odd coefficient 
of $w^{-1}$ be zero for the power series expansion of the local 
coordinate at infinity for a canonical supersphere with only one tube.  
In order to account for these discrepancies one must transform $C_3$ to 
a canonical supersphere via a superprojective transformation as 
specified in the proofs of Propositions \ref{canonicalcriteria} and 
\ref{canonical criteria for one tube}.  In other words, if $n$ and $i$ do 
not satisfy the conditions that $i = m$ and $n \neq 0$, we can still use 
the unique $F$ determined by the sewing equation and the normalization 
conditions above, but then we must compose $F$ with an easily determined 
superprojective transformation in order to obtain a canonical supersphere.

Thus the real work needed to obtain a canonical supersphere {}from the
sewing of two canonical superspheres is in solving for $F$ using the
sewing equation and the normalization conditions.  The existence and
uniqueness of a solution for $F$, i.e., a solution to the sewing
equation and the normalization conditions, are guaranteed by
Propositions \ref{canonicalcriteria} and \ref{canonical criteria for one tube}.

We give two special examples.  Define the superconformal shift
\begin{eqnarray*}
s_{(z_i,\theta_i)} : \mbox{$\bigwedge_\infty$} & \longrightarrow &
\mbox{$\bigwedge_\infty$}\\
(w,\rho) & \mapsto & (w - z_i - \rho \theta_i, \rho - \theta_i).
\end{eqnarray*}
If $\hat{H}_0 (w,\rho) = I(w,\rho) = 
(1/w,i\rho/w)$, then \footnote{There is a misprint in the analogous
nonsuper case to equations (\ref{boundary conditions trivial
infinity 1}) and (\ref{boundary conditions trivial infinity 2}) given
in \cite{H book}.  Equations (1.4.9) and (1.4.10) in \cite{H book} should be 
$F^{(1)}(w) = w - z_i$ and $F^{(2)}(w) = f_i^{-1}(w) - z_i$, respectively, not
$F^{(1)}(w) = w$ and $F^{(2)}(w) = f_i^{-1}(w)$ as stated.  The correction is
necessary if $F^{(2)}$ is to satisfy the normalization condition (1.4.3) in
the case that $z_i \neq 0$.}
\begin{eqnarray}
F^{(1)}_\sou (w, \rho) &=& s_{(z_i,\theta_i)}(w, \rho)  \label{boundary conditions trivial
infinity 1}\\ 
F^{(2)}_\sou (w, \rho) &=& s_{(z_i,\theta_i)}\circ  H_i^{-1}(w, \rho)
.\label{boundary conditions trivial infinity 2}
\end{eqnarray}

If $H_i (w, \rho) = (\asqrt^2(w - z_i - \rho\theta_i), \asqrt(\rho - \theta))$
for $\asqrt \in (\bigwedge_\infty^0)^\times$, then
\begin{eqnarray}
\hspace{.5in} F^{(1)}_\sou (w, \rho) \! &=& \! \left(\asqrt^{-2} \left(\hat{H}_0^{-1} \circ I
(\asqrt^2(w - z_i - \rho\theta_i), \asqrt(\rho - \theta)) \right)^0, \right. 
\label{boundary conditions trivial puncture 1} \\ 
& &\left. \hspace{.6in} \asqrt^{-1} \left(\hat{H}_0^{-1} \circ I (\asqrt^2(w - z_i -
\rho\theta_i), \asqrt(\rho - \theta)) \right)^1 \right) \nonumber \\ 
F^{(2)}_\sou(w, \rho) \! &=& \! \left(\asqrt^{-2}w, \asqrt^{-1}\rho \right) . \label{boundary
conditions trivial puncture 2}
\end{eqnarray}
We call the conditions (\ref{boundary conditions trivial infinity 1})
-- (\ref{boundary conditions trivial puncture 2}) the {\it boundary
conditions} for $F$.

As with the nonsuper case (see \cite{H book}), this naturally leads to
several questions:  How does $F$ depend on the local coordinates at
the punctures being sewn, i.e., how does $F$ depend on $\hat{H}_0$ and
$H_i$?  Is this dependence analytic or even algebraic in some sense?
In this case is it possible to find a solution for $F$ directly and
explicitly using the sewing equation, the normalization conditions,
and the boundary conditions?

The answers to these last two questions is yes as will be shown in
Chapters 3 and 4.  In fact in Chapter 3, in a certain formal algebraic 
setting, we will solve for $F$ directly and explicitly by showing that 
$F$ depends on $\hat{H}_0$ and $H_i$ algebraically in a certain sense, 
and that formally, $F$ can be obtained uniquely using the sewing 
equation, the normalization conditions and the boundary conditions. 
Then in Chapter 4, we will prove that this dependence is also analytic 
in a certain sense, and that if the local coordinates $\hat{H}_0$ and 
$H_i$ are convergent superconformal local coordinates, then $F$ is a 
convergent superconformal function on the appropriate domain.

\chapter[An algebraic study of the sewing operation]
{A formal algebraic study of the sewing operation}

In this chapter, we develop a formal theory of infinitesimal $N = 1$
superconformal transformations based on a representation of the $N =
1$ Neveu-Schwarz algebra of superconformal symmetries in terms of
superderivations, and we use these results to solve a formal version
of the sewing equation along with the normalization and boundary 
conditions introduced in Chapter 2, thus answering several of the 
questions posed in Chapter 2 regarding the sewing equation.  This 
solution to the formal sewing equation and the normalization and 
boundary conditions gives an identity for certain exponentials of 
superderivations involving infinitely many formal variables.  In 
addition, we prove two other identities which are related to
certain sewings and which also involve certain exponentials of 
superderivations.  These superderivations give a representation of 
the $N=1$ Neveu-Schwarz algebra with central charge zero, and we 
use these identities for this representation to prove similar 
identities for the Neveu-Schwarz algebra itself and hence for any 
representation of the Neveu-Schwarz algebra.  Thus we obtain a
correspondence between the supergeometric sewing operation of Chapter 
2 and certain identities on representations of the Neveu-Schwarz 
algebra. 

The material in this chapter is algebraic and independent of the 
supergeometry studied in Chapter 2.  However, the results 
of this chapter do of course have geometric motivation and meaning, 
and will be applied to the supergeometric setting in Chapter 4.  

This chapter is organized as follows.  In Section 3.1, we give some 
preliminary definitions and give two generalizations of the 
``automorphism property'' {}from \cite{FLM}.  In Section 3.2, we 
develop a formal supercalculus, define what is meant by a formal 
$N=1$ superconformal power series, and study the formal theory of 
superconformal local coordinate maps for a super-Riemann surface.  
Using formal exponentiation, we characterize certain formal 
superconformal local coordinate maps in terms of exponentials of 
superderivations with infinitely many formal variable coefficients.  
The proof of this characterization relies on the generalizations of 
the automorphism property proved in Section 3.1.

In Section 3.3, we introduce the formal sewing equation and solve 
this equation along with the formal normalization and boundary  
conditions in terms of exponentials of infinite series of certain  
superderivations.  This answers in the affirmative some of the  
questions posed in Chapter 2.  That is, we show that formally the 
uniformizing function giving a canonical supersphere with tubes 
{}from the sewing together of two canonical superspheres with tubes 
does depend on the local coordinates at the punctures where the two 
superspheres are being sewn, and that in a certain sense, this 
dependence is algebraic.  In addition, we show that formally, this 
solution to the sewing equation gives an explicit solution for the 
uniformizing function.   We do this by interpreting the sewing  
equation as a product of exponentials of certain infinite series of  
superderivations with formal variable coefficients and by showing 
that this product is equal to a different product of exponentials of 
superderivations with coefficients involving these formal variables.  
In addition, in Section 3.3, we prove two other identities expressing 
the exponentials of certain infinite series of superderivations with 
formal variable coefficients as a product of certain other exponentials.  
These identities relate to determining the resulting local coordinates
of two particular sewings and are needed in \cite{B thesis} to prove the 
isomorphism theorem.

In Section 3.4, we define the $N = 1$ Neveu-Schwarz algebra and point 
out that the superderivations we used in Sections 3.2 and 3.3 give a 
representation of the Neveu-Schwarz algebra with central charge zero.  
This shows that the Neveu-Schwarz algebra is the algebra of 
infinitesimal superconformal transformations.  We briefly discuss 
the subalgebra of the Neveu-Schwarz algebra consisting of infinitesimal 
global superconformal transformations (i.e., infinitesimal 
superprojective transformations).  In Section 3.5, we discuss modules 
for the Neveu-Schwarz algebra.

In Section 3.6, we generalize the identities obtained in Section 3.3
to general representations of the Neveu-Schwarz algebra.  We obtain 
an additional formal series in infinitely many formal variables which 
is related to the central charge of the Neveu-Schwarz algebra.  In 
Section 3.7, we give the corresponding identities for positive-energy 
representations of the Neveu-Schwarz algebra showing that the 
resulting series have certain nice properties.  Thus we obtain a 
correspondence between the sewing identities occurring in the 
supergeometric setting and the analogous identities occurring in the 
algebraic setting of a positive-energy representation of the 
Neveu-Schwarz algebra, for instance an $N=1$  Neveu-Schwarz vertex 
operator superalgebra (cf. \cite{B thesis}, \cite{B  vosas}).

\section{An extension of the automorphism property}

In this section we work over a field $\mathbb{F}$ of characteristic zero.  
Recall that the space $\mbox{Der} \; A$ of all superderivations of $A$ is a Lie
sub-superalgebra of $\mbox{End} \; A$.  Note that $(\mbox{Der} \;
A)^0$ consists of ordinary derivations.  Thus the following
``automorphism property'' (8.2.10) in \cite{FLM} holds.

\begin{prop}\label{autoprop} (\cite{FLM})
Let $A$ be a superalgebra, $u,v \in A$, $T \in 
(\mbox{\em Der} \; A)^0$, and $y$ a formal variable
commuting with $A$.  Then  
\begin{equation}\label{autoequation} 
e^{yT} \cdot (uv) = (e^{yT} \cdot u)(e^{yT} \cdot v).
\end{equation}
\end{prop}


An expression of the form $e^x$ denotes the formal exponential
series. In the proof of Proposition \ref{superconformal}, we will need
the following proposition which is a generalization of the
automorphism property (\ref{autoequation}).

\begin{prop}\label{conformalproof}
Let $A$ be a superalgebra, $a \in A^0$, $u,v \in A$, $T \in (\mbox{\em
Der} \; A)^0$, and $y$ a formal variable commuting with $A$.  Then
\begin{eqnarray}
e^{y(a + T)} \cdot (u v) &=& (e^{yT} \cdot u) (e^{y(a + T)} \cdot v) 
\label{conformal equation 2} \\
&=&(e^{y(a + T)} \cdot u) (e^{yT} \cdot v) . \label{conformal equation 1} 
\end{eqnarray}
\end{prop}

Proposition \ref{conformalproof} can be proven directly by expanding both 
sides of (\ref{conformal equation 1}) and (\ref{conformal equation 2}) as 
power series in $y$, comparing coefficients of $y^n$ for each $n \in 
\mathbb{N}$, and using induction on $n$.  This proof is given in 
\cite{B thesis}.  However, here we present the following alternate proof 
which is more Lie theoretic and which also appeared in \cite{B thesis}.  
This approach was suggested by J. Lepowsky.

\begin{proof} Let $X_1,X_2,..., X_k \in (\mbox{End} \; A[[y]])^0$.
Set
\[[X_1,X_2, ..., X_k] = [X_1, [X_2, \cdots, [X_{k-1}, X_k] \cdots ]] ,\] 
and let $[X_1^{(i_1)}X_2^{(i_2)} \cdots X_k^{(i_k)}]$ denote the 
iterated bracket starting with the sequence of $i_1$ of the $X_1$'s then
$i_2$ of the $X_2$'s, etc.  Then, if well defined, the Campbell-Baker-Hausdorff 
formula gives $e^{X_1} e^{X_2} = e^{C(X_1,X_2)}$ where  
\[C(X_1,X_2) = \! \! \sum_{k \in \Z} \sum_{ 
\begin{tiny} \begin{array}{c}
i_1,..., i_k, j_1,..., j_k \in \mathbb{N}\\
i_r + j_r \geq 1
\end{array} \end{tiny}
}  \! \frac{(-1)^{k-1}[X_1^{(i_1)}X_2^{(j_1)} \cdots
X_1^{(i_k)}X_2^{(j_k)}]}{k((i_1
+ j_1) + \cdots + (i_k + j_k)) (i_1! j_1! \cdots i_k! j_k!)}  ,\] 
(cf. \cite{Re}). For any $b \in A^0$, let $l_b : A \rightarrow A$ be left
multiplication by $b$.  Then for $b, c \in A^0$
\[ \left[ l_b, T \right] = - l_{Tb}, \qquad 
\mbox{and} \qquad \left[ l_b, l_c \right] = 0 .\] 
Thus for $i_k, j_k \in \{ 0,1\}$, we have 
\[ [(l_a + T)^{(i_1)} T^{(j_1)} \cdots (l_a + T)^{(i_k)} T^{(j_k)}] =
l_{T^{i_1 + j_1 + \cdots + i_k + j_k} a} ,\]   
and 
\[e^{y(l_a + T)} e^{-yT} = e^{yl_b} \]
where
\[b = \sum_{k \in \Z} d_k y^{k-1} (T^{k-1} a) \]
for some $d_k \in A^0$ with $d_1 = 1$. Using the automorphism property 
(\ref{autoequation}), we have
\begin{eqnarray*}
e^{y(a + T)} \cdot (uv) &=& e^{y(a + T)} e^{-yT} e^{yT} \cdot (uv) \\ 
&=& e^{yl_b} e^{yT} \cdot (uv) \\
&=& e^{yl_b} (e^{yT} \cdot u) (e^{yT} \cdot v) = (e^{y(a + T)} \cdot u) 
(e^{yT} \cdot v) \\ 
&=& (e^{yT} \cdot u) e^{yl_b} (e^{yT} \cdot v) = (e^{yT} \cdot u) 
(e^{y(a + T)} \cdot v),  
\end{eqnarray*}
as desired.
\end{proof}

\section[Formal supercalculus]
{Formal supercalculus and formal superconformal power series}

Let $R$ be a superalgebra over $\mathbb{Q}$ (with identity).  Let $x$ be 
a formal variable which commutes with all elements of $R$, and let
$\varphi$ be a formal variable which commutes with $x$ and elements of
$R^0$ and anti-commutes with elements of $R^1$ and itself.  In
general, we will use the term {\it even formal variable} to denote a
formal variable which commutes with all formal variables and with all
elements in any coefficient algebra.   We will use the term {\it odd
formal variable} to denote a formal variable which anti-commutes with
all odd elements and commutes with all even elements in any
coefficient algebra, and in addition, odd formal variables will all
anti-commute with each other.  Consequently, an odd formal variable
has the property that its square is zero.   

For a vector space $V$, and for even formal variables $x_1$,
$x_2$,..., and odd formal variables $\varphi_1$, $\varphi_2$,...,
consider the spaces 
\begin{multline*}
V[[x_1,...,x_m]][\varphi_1,...,\varphi_n]  \\
= \; \Bigl\{ \Bigl. \sum_{ \begin{tiny}
\begin{array}{c} k_1,..., k_m \in \mathbb{N}\\
l_1,...l_n \in \mathbb{Z}_2
\end{array} \end{tiny}} \! a_{k_1,...,k_m, l_1,..., l_n} x_1^{k_1}\cdots x_m^{k_m}
\varphi_1^{l_1} \cdots \varphi_n^{l_n} \; \Bigr| \; a_{k_1,...,k_m, l_1,..., l_n} 
\in V \Bigr\} 
\end{multline*}
and
\[V((x))[\varphi] = \Bigl\{ \Bigl. \sum_{n = N}^{\infty} a_n x^n + \varphi \! \sum_{n
= N}^{\infty} b_n x^n \; \Bigr| \; N \in \mathbb{Z} , \; a_n, b_n \in V \Bigr\} \subset
V[[x,x^{-1}]] [\varphi] .\]  
Then $R[[x_1,...,x_m]][\varphi_1,...,\varphi_n]$ is a $\mathbb{Z}_2$-graded
vector space with sign function given by
\[\eta(a x_1^{k_1}\cdots x_m^{k_m}\varphi_1^{l_1} \cdots \varphi_n^{l_n}) = 
(\eta(a) + l_1 + \cdots + l_n ) \; \mathrm{mod} \; 2 ,\]
and $R((x))[\varphi]$ is a superalgebra as is $R((x^{-1}))[\varphi]$.

Define 
\begin{equation}\label{define D}
D = \Dx .  
\end{equation}
Then $D$ is an odd derivation in both $\mbox{Der}(R((x))[\varphi])$
and $\mbox{Der}(R((x^{-1})) [\varphi])$.  Furthermore, $D$ satisfies
the super-Leibniz rule (\ref{leibniz}) for the product of any two 
elements in $R[[x,x^{-1}]] [\varphi]$ if that product is well defined.  
Note that
\[ D^2 = \frac{\partial}{\partial x} . \]

Recall that a superanalytic $(1,1)$-superfunction over $\bigwedge_{*>0}$, 
$H(z,\theta)$, has a Laurent expansion about $z$ and $\theta$ (\ref{Laurent
series}) which is an element of $\bigwedge_{*>0} [[z,z^{-1}]] 
[\theta]$.  Taking a general coefficient superalgebra $R$, we can write a 
formal superfunction in one even formal variable and one odd formal 
variable over $R$ as 
\[H(x,\varphi) = (f(x) + \varphi \xi(x), \psi(x) + \varphi g(x)) \in
(R[[x,x^{-1}]][\varphi])^0 \oplus (R[[x,x^{-1}]][\varphi])^1\]
where $f(x), g(x) \in R^0[[x,x^{-1}]]$, and $\xi(x), \psi(x) \in 
R^1[[x,x^{-1}]]$.  

In Chapter 2, the operator $D = \frac{\partial}{\partial \theta} + \theta
\frac{\partial}{\partial z}$ was used to define the notion of 
superconformal $(1,1)$-superfunction over $\bigwedge_{*>0}$ which is
a superanalytic $(1,1)$-function with the condition that it transform $D$ 
homogeneously of degree one.   This is equivalent to the conditions 
(\ref{superconformal condition}).  Thus formally, we define a series 
$H(x,\varphi) = (f(x) + \varphi \xi(x), \psi(x) + \varphi g(x))$ in 
$R[[x, x^{-1}]][\varphi]$ to be {\it formally superconformal} if  
\begin{equation}\label{grungy superconformal condition}
\xi(x) = g(x) \psi(x), \quad \mbox{and} \quad g(x)^2 =
\frac{\partial}{\partial x} f(x) + \psi(x) \frac{\partial}{\partial x}
\psi (x)
\end{equation} 
with $g(x) + \varphi \psi'(x)$ not identically zero.  Therefore a formal 
superconformal series is uniquely determined by an even formal series 
$f(x)$, an odd formal series $\psi(x)$ and a square root of the formal 
power series $f'(x) + \psi(x) \psi'(x)$.  Since a formal series 
$H(x,\varphi) = (f(x) + \varphi \xi(x), \psi(x) + \varphi g(x))$ in 
$R[[x,x^{-1}]][\varphi]$ has the property that $\varphi H(x,\varphi) = 
\varphi (f(x), \psi(x))$, any formal superconformal series $H(x,\varphi) 
\in R[[x,x^{-1}]][\varphi]$ can be uniquely expressed by the formal 
series $\varphi H(x,\varphi) \in \varphi R[[x,x^{-1}]]$ and a square root 
for $f'(x) + \psi(x) \psi'(x) \in  (R[[x,x^{-1}]][\varphi])^0$.

We introduce the notation $\tilde{x} =  H^0(x,\varphi) = f(x) + \varphi 
\xi(x)$ for the even part of $H$ and $\tilde{\varphi} = H^1(x,\varphi) = 
\psi(x) + \varphi g(x)$ for the odd part of $H$.  Then for
$H(x,\varphi) =  (\tilde{x},\tilde{\varphi})$, the condition (\ref{grungy
superconformal condition}) for $H$ to be superconformal is equivalent
to the condition 
\begin{equation}\label{nice superconformal condition}
D\tilde{x} = \tilde{\varphi} D \tilde{\varphi} .
\end{equation}

In Chapter 2, we began the study of the moduli space of 
superspheres with punctures and local superconformal coordinates 
vanishing at the punctures.  The punctures on a supersphere 
with tubes can be thought of as being at $0 \in \bigwedge_\infty$, 
a non-zero point in $\bigwedge_\infty$, or at a distinguished point on 
the supersphere we denote by $\infty$.  Since we can always shift a 
non-zero point in $\bigwedge_\infty$ to zero, all local superconformal
coordinates vanishing at the punctures can be expressed as power 
series vanishing at zero composed with a shift or at infinity.  Thus 
we want to study in more detail certain formal superconformal power 
series in $(xR[[x]] \oplus \varphi R[[x]]) \subset R[[x]][\varphi]$ or 
in $x^{-1} R[[x^{-1}]] [\varphi]$ since these will include the formal 
superconformal coordinates over a superalgebra $R$ vanishing at 
$(x,\varphi) = 0$ and $(x,\varphi) = (\infty,0)$, 
respectively.    

First restricting our attention to formal superconformal power series 
in $xR[[x]] \oplus \varphi R[[x]]$, we note that if $H(x,\varphi)
\in (xR[[x]] \oplus \varphi R[[x]])$ is superconformal and the even
part of the $x$ coefficient in $H(x,\varphi)$ is one, i.e., if
$\varphi H(x,\varphi) = \varphi(f(x), \psi(x))$ for $f(x) \in x
R^0[[x]]$ and $\psi(x) \in R^1[[x]]$ such that $f(x)$ is of the form
$f(x) = x + \sum_{j = 1}^{\infty} a_j x^{j + 1}$, then $f' + \psi
\psi' \in \{1 + \sum_{n = 1}^{\infty} c_n x^n \; | \; c_n \in R^0 \}$.   
We define 
\begin{eqnarray}
\sqrt{\quad}  : \Bigl\{ \Bigl. 1 + \sum_{n = 1}^{\infty} c_n x^n \; \Bigr| 
\; c_n \in R \Bigr\} &\longrightarrow& \Bigl\{ \Bigl. 1 + 
\sum_{n = 1}^{\infty} c_n x^n \; \Bigr| \; c_n \in R \Bigr\} \label{sqrt} \\
h(x) &\mapsto & \sqrt{h(x)} \nonumber
\end{eqnarray}
to be the Taylor expansion of $\sqrt{h(x)}$ about $x = 0$ such that 
$\sqrt{1} = 1$.   Note that $\sqrt{\quad}$ as defined here is the unique 
square root defined on $\{ 1 + \sum_{n = 1}^{\infty} c_n x^n \; | \; c_n 
\in R \}$ with range contained in $\{ 1 + \sum_{n = 1}^{\infty} c_n x^n \; 
| \; c_n \in R \}$, and $- \sqrt{\quad}$ is the unique square root defined 
on $\{ 1 + \sum_{n = 1}^{\infty} c_n x^n \; | \; c_n \in R \}$ mapping into 
$\{ -1 - \sum_{n = 1}^{\infty} c_n x^n \; | \; c_n \in R \}$.  

Any superconformal $H \in R[[x]][\varphi]$ with $\varphi H(x, \varphi)
= \varphi(f(x), \psi(x))$ (i.e., determined by $f(x)$ and $\psi(x)$)
for some $f(x), \psi(x) \in xR[[x]]$ and with the coefficient of $x$
in $f(x)$ equal to one must have the form     
\begin{equation}\label{form of a superconformal function} 
H(x,\varphi) = \left(f(x) \pm \varphi \psi(x) \sqrt{f'(x)}, \psi(x) \pm \varphi 
\sqrt{f'(x) + \psi(x) \psi'(x)} \right) .
\end{equation} 
We can distinguish between these two by specifying the even coefficient 
of $\varphi$.  That is if $H$ is superconformal, vanishing at zero and 
$\varphi H(x, \varphi) = \varphi(f(x), \psi(x))$ with the coefficient 
of $x$ in $f(x)$ equal to one and the even coefficient of $\varphi$ 
equal to one (resp., -1), then $H$ must be of the form (\ref{form of a 
superconformal function}) with the positive (resp., negative) sign.

Note that {}from (\ref{form of a superconformal function}) it is clear
that any superconformal function $H$ is also completely determined
by its even (or odd) part and a specification of square root.  That is, 
given $\tilde{x} = H^0(x,\varphi)$ (or $\tilde{\varphi} = H^1(x,\varphi)$) 
and a specification of square root, one can recover $f$ and $\psi$.  

We wish to express any formal superconformal series vanishing at zero
in terms of a formal exponential of an infinite sum of certain
superderivations.  For any even formal series $f(x) \in xR^0[[x]]$ with
$x$ coefficient one, and any odd formal series $\psi(x) \in xR^1[[x]]$, 
we first express $\varphi(f(x),\psi(x)) \in \varphi xR[[x]]$ in terms 
of $\varphi$ times the exponential of an infinite sum of 
superderivations in $\mbox{Der} (R((x))[\varphi])$ acting on 
$(x,\varphi)$.  Then we prove that this exponential of superderivations 
is in fact superconformal with even coefficient of $\varphi$ equal to 
one and that there is a one-to-one correspondence between such 
exponential expressions and formal superconformal power series in
$R[[x]][\varphi]$ vanishing at $(x,\varphi) = 0$ with even coefficient 
of the $x$ and $\varphi$ terms equal to one.  To do this, we will want
to introduce more formal variables.

Let $\mathcal{A}_j$, for $j \in \Z$, be even formal variables, and let 
$\mathcal{M}_{j - 1/2}$, for $j \in \Z$, be odd formal variables.  Let 
$\mathcal{A}  = \{\mathcal{A}_j\}_{j \in \Z}$ and $\mathcal{M} = 
\{\mathcal{M}_{j - 1/2}\}_{j \in \Z}$, and consider the 
$\mathbb{Q}$-superalgebra $\mathbb{Q}[\mathcal{A},\mathcal{M}]$ of 
polynomials in the formal variables $\mathcal{A}_1, \mathcal{A}_2,...$ 
and $\mathcal{M}_{1/2}, \mathcal{M}_{3/2},....$.  Consider the even 
superderivations 
\begin{equation}\label{L notation}
L_j(x,\varphi) = - \biggl( \Lx \biggr)
\end{equation}
and the odd superderivations
\begin{equation}\label{G notation}
G_{j -\frac{1}{2}} (x,\varphi) = - \Gx
\end{equation}
in $\mbox{Der} (R((x))[\varphi])$, for $j \in \mathbb{Z}_+$. We define the
sequences 
\[E^0(\mathcal{A}, \mathcal{M}) = \left\{E_j(\mathcal{A}, 
\mathcal{M})\right\}_{j \in \Z} \qquad \mbox{and} \qquad 
E^1(\mathcal{A}, \mathcal{M}) = \bigl\{E_{j - \frac{1}{2}}(\mathcal{A},
\mathcal{M})\bigr\}_{j \in \Z}\]
of even and odd elements, respectively, in $\mathbb{Q}[\mathcal{A}, \mathcal{M}]$ by  
\begin{multline}\label{exp1}
\varphi \Bigl( x + \sum_{j \in \Z} E_j(\mathcal{A}, \mathcal{M})x^{j + 1} ,
\sum_{j \in \Z} E_{j - \frac{1}{2}}(\mathcal{A}, \mathcal{M})x^j \Bigr) \\
= \; \varphi \exp\Biggl( \! - \! \sum_{j \in \Z} \Bigl( 
\mathcal{A}_j L_j(x,\varphi) + \mathcal{M}_{j - \frac{1}{2}} G_{j -\frac{1}{2}}
(x,\varphi) \Bigr) \! \Biggr) \cdot (x,\varphi) .
\end{multline}
As usual, ``exp'' denotes the formal exponential series, when it is
defined, as it is in the case of the above exponential of the 
derivation in $\mbox{Der} \; (\mathbb{Q}[\mathcal{A}, \mathcal{M}][[x]][\varphi])$.  
The reason for the $\varphi$ multiplier in
(\ref{exp1}), is that we are in fact uniquely defining a series  
\[(f(x),\psi(x)) = \biggl( x + \sum_{j \in \Z} E_j(\mathcal{A}, \mathcal{M})
x^{j + 1} , \sum_{j \in \Z} E_{j - \frac{1}{2}}(\mathcal{A}, \mathcal{M})x^j
\biggr)  \]  
in $x\mathbb{Q}[\mathcal{A}, \mathcal{M}][[x]]$ by means of $\varphi$ times 
a certain series in $\mathbb{Q}[\mathcal{A}, \mathcal{M}][[x]][\varphi]$. 
{}From (\ref{exp1}), we see that, for $j \in \Z$,
\begin{equation}\label{E even}
E_j(\mathcal{A}, \mathcal{M}) = \mathcal{A}_j + r_j^0(\mathcal{A}_1,...,\mathcal{A}_{j
- 1}, \mathcal{M}_{\frac{1}{2}}...,\mathcal{M}_{j - \frac{3}{2}}) ,
\end{equation}
and 
\begin{equation}\label{E odd}
E_{j - \frac{1}{2}}(\mathcal{A}, \mathcal{M}) = \mathcal{M}_{j - \frac{1}{2}} +
r_j^1(\mathcal{A}_1,...,\mathcal{A}_{j - 1}, \mathcal{M}_{\frac{1}{2}},..., 
\mathcal{M}_{j - \frac{3}{2}})   
\end{equation}
where $r_j^0(\mathcal{A}_1,..., \mathcal{A}_{j - 1}, \mathcal{M}_{1/2},..., 
\mathcal{M}_{j - 3/2})$, and $r_j^1( \mathcal{A}_1,..., \mathcal{A}_{j - 1},
\mathcal{M}_{1/2},..., \mathcal{M}_{j - 3/2})$ are in
$\mathbb{Q}[\mathcal{A}_1,...,\mathcal{A}_{j - 1}, \mathcal{M}_{1/2},..., 
\mathcal{M}_{j - 3/2} ]$, both with constant term zero. 

Let $R$ be a superalgebra over $\mathbb{Q}$.  Let $(R^0)^\infty$ be 
the set of all sequences $\{A_j\}_{j \in \Z}$ of even elements in $R$, 
let $(R^1)^\infty$ be the set of all sequences $\{M_{j - 1/2}\}_{j 
\in \Z}$ of odd elements in $R$, and let $R^\infty = (R^0)^\infty
\oplus (R^1)^\infty$.  Given any 
\[(A,M) = \bigl(\{A_j \}_{ j \in \Z}, \{M_{j - \frac{1}{2}} \}_{j \in \Z}\bigr)
= \bigl\{(A_j, M_{j - \frac{1}{2}}) \bigr\}_{j \in \Z} \in R^\infty , \]
we have a well-defined sequence $E(A,M) = (E^0(A,M),E^1(A,M))$ in
$R^\infty$ by substituting $A$ and $M$ into $E^0(\mathcal{A}, \mathcal{M})$
and $E^1(\mathcal{A}, \mathcal{M})$, respectively, since $E_j(\mathcal{A},
\mathcal{M})$ and $E_{j - \frac{1}{2}}(\mathcal{A}, \mathcal{M})$ are in
$\mathbb{Q}[\mathcal{A}, \mathcal{M}]$ for $j \in \Z$.  This defines a
map    
\begin{eqnarray*}
E : R^\infty &\longrightarrow& R^\infty =  \; (R^0)^\infty \oplus
(R^1)^\infty \\
(A,M) &\mapsto& E(A,M) \; = \; \bigl(E^0(A,M),E^1(A,M)\bigr).
\end{eqnarray*}

\begin{prop}\label{Ebijection}
The map $E$ is a bijection.  In particular, $E$ has an inverse
$E^{-1}$.
\end{prop}

\begin{proof} Given $(a,m) \in R^\infty$ consider the
infinite system of equations  
\begin{equation}\label{system of equations}
(E^0(A,M),E^1(A,M)) = (a,m)
\end{equation}
for the unknown sequence $(A,M) = \{(A_j, M_{j - 1/2})\}_{j
\in \Z}$.  Define $\mathcal{E}_{2l}(A,M) = E_{l}(A,M)$, and $\mathcal{E}_{2l
+ 1} (A,M) = E_{l - \frac{1}{2}}(A,M)$ for $l \in \Z$.  Then the infinite
system of equations (\ref{system of equations}) becomes a system of equations 
for $\mathcal{E}_k (A,M)$, and using (\ref{E even}) and (\ref{E odd}), it is 
easy to show by induction on $k = 1,2,...$, that this system of equations has 
a unique solution, i.e.,
\begin{eqnarray*}
(A,M) \; = \; E^{-1} (a,m) \! &=& \! \bigl((E^{-1})^0 (a,m),(E^{-1})^1 (a,m)\bigr) \\
&=& \! \Bigl(\bigl\{(E^{-1})_j (a,m)\bigr\}_{j \in \Z},\bigl\{(E^{-1})_{j - \frac{1}{2}}
(a,m)\bigr\}_{j \in \Z}\Bigr)  \; \in R^\infty.
\end{eqnarray*}
The proposition follows immediately. 
\end{proof}

\begin{cor}
For any formal power series of the form 
\begin{equation}\label{fpsi} 
(f(x), \psi(x)) = \biggl( x +
\sum_{j \in \Z} a_j x^{j + 1}, \sum_{j \in \Z} m_{j - \frac{1}{2}} x^j
\biggr) \in x\left(R^0[[x]] \oplus R^1[[x]] \right), 
\end{equation}
we have 
\[\varphi f(x) = \varphi \exp  \Biggl(\! - \! \sum_{j \in \Z} \left(
E^{-1}_j(a,m) L_j(x,\varphi)  + E^{-1}_{j - \frac{1}{2}}(a,m)
G_{j - \frac{1}{2}}(x, \varphi) \right) \! \Biggr) \! \cdot x , \] 
and 
\[\varphi \psi(x) = \varphi \exp  \Biggl(\! - \! \sum_{j \in \Z} \left(
E^{-1}_j(a,m) L_j(x,\varphi)  +  E^{-1}_{j - \frac{1}{2}}(a,m)
G_{j - \frac{1}{2}}(x, \varphi)  \right)  \! \Biggr) \! \cdot \varphi  \]
where $E^{-1}_j(a,m)$ and $E^{-1}_{j - 1/2}(a,m)$, for $j\in \Z$, denote the 
even and odd components of the series $E^{-1}(a,m) = \{E^{-1}_j(a,m), E^{-1}_{j -
1/2}(a,m)\}_{j \in \Z} \in R^\infty$, respectively.
\end{cor} 

\begin{proof} Using equation (\ref{exp1}) and the fact that
$E$ is a bijection, we have
\begin{eqnarray*}
& & \hspace{-.5in} (\varphi f(x), \varphi \psi(x)) \\
&=& \! \! \varphi \Bigl( x + \sum_{j \in \Z}
a_j x^{j + 1}, \sum_{j \in \Z} m_{j - \frac{1}{2}} x^j \Bigr)\\ 
&=&  \! \! \varphi \Bigl( x + \sum_{j \in \Z} E_j(E^{-1}(a,m)) x^{j +
1}, \sum_{j \in \Z} E_{j - \frac{1}{2}}(E^{-1}(a,m)) x^j \Bigr) \\
&=&  \! \! \varphi \exp  \Biggl( \! - \! \sum_{j \in \Z} \left( 
E^{-1}_j(a,m) L_j(x,\varphi)  + E^{-1}_{j - \frac{1}{2}} (a,m) 
G_{j - \frac{1}{2}}(x,\varphi) \right) \! \Biggr)  \! \cdot (x,\varphi) .  
\end{eqnarray*}
\end{proof}

\begin{prop}\label{superconformal}
Let $R$ be a superalgebra and 
\[(A,M) = \{(A_j, M_{j - 1/2}) \}_{j \in \Z} \in (R^0)^{\infty} \oplus 
(R^1)^{\infty} = R^\infty.  \]
Then  
\begin{eqnarray}
\hspace{.2in} H(x, \varphi) \! \! &=& \! \! \exp\Biggl( \! - \! \sum_{j \in \Z} \left( 
A_j L_j(x,\varphi) + M_{j - \frac{1}{2}} G_{j - \frac{1}{2}}(x,\varphi) 
\right) \Biggr) \cdot (x, \varphi) \label{H}\\ 
&=& \! \! \exp\Biggl( \sum_{j \in \Z} \biggl( A_j \biggl( \Lx \biggr) \biggr. 
\Biggr.\nonumber \\
& & \hspace{1.1in} \Biggl. \biggl. + \; M_{j - \frac{1}{2}} \Gx \biggr) \!
\Biggr) \! \cdot (x, \varphi) \nonumber
\end{eqnarray}
is superconformal and is the unique formal superconformal power
series in $R[[x]][\varphi]$ with even coefficient of $\varphi$ equal to one such that
\[\varphi H(x, \varphi) = \varphi \biggl(x + \sum_{j \in \Z}
E_j(A,M) x^{j + 1}, \sum_{j \in \Z} E_{j - \frac{1}{2}}(A,M) x^j \biggr)
. \] 
\end{prop}

\begin{proof}  
Let 
\begin{equation}\label{define T}
T = - \sum_{j \in \Z} \left( A_j L_j(x,\varphi) + M_{j -
\frac{1}{2}} G_{j - \frac{1}{2}}(x,\varphi) \right).
\end{equation} 
Then $T \in (\mbox{Der}(R[[x]][\varphi]))^0$, i.e., $T$ is even, and
thus 
\[e^{T} \cdot (x, \varphi) = (e^{T} \cdot x, e^{T} \cdot \varphi)
\in (R[[x]][\varphi])^0 \oplus (R[[x]][\varphi])^1.  \]
Let
\[h(x,\varphi) = \sum_{j \in \Z} \biggl( A_j \Bigl(\frac{j + 1}{2}
\Bigr) x^j + \varphi M_{j - \frac{1}{2}} j  x^{j - 1} \biggr) .\] 
Then $h \in (R[[x]][\varphi])^0$, i.e., $h$ is even, and we have
\[ [D,T] = h(x,\varphi)D . \]   
Thus 
\[D e^T \cdot x = e^{(h + T)} \cdot D \cdot x = e^{(h +
T)} \cdot \varphi . \] 
By Proposition \ref{conformalproof},
\begin{equation}\label{thisequation} 
e^{y(h + T)} \cdot \varphi = e^{y(h + T)} \cdot (\varphi 1) =
\left(e^{yT} \cdot \varphi \right)\bigl(e^{y(h + T)} \cdot 1\bigr) .
\end{equation}
But in this case, the coefficient of $y^n$ for a fixed $n \in \mathbb{N}$
has terms with powers of $x$ greater than or equal to $n - 1$.  Thus
we can set $y = 1$, and each power series in equation
(\ref{thisequation}) has only a finite number of $x^j$ terms for a
given $j \in \mathbb{N}$, i.e., each term is a well-defined power series
in $x$.  Therefore  
\[e^{(h + T)} \cdot \varphi = \left(e^{T} \cdot \varphi
\right)\bigl(e^{(h + T)} \cdot 1\bigr) .\] 
Thus writing $H(x, \varphi) = (e^{T} \cdot x, e^{T} \cdot \varphi) =
(\tilde{x}, \tilde{\varphi})$, we have  
\begin{eqnarray*}
D \tilde{x} &=& e^{(h + T)} \cdot \varphi \; = \; \left(e^{T} \cdot \varphi
\right)\bigl(e^{(h + T)} \cdot 1\bigr) \\
&=&  \left(e^{T} \cdot \varphi \right)\bigl(e^{(h + T)} \cdot D \cdot
\varphi\bigr) \; = \; \left(e^{T} \cdot \varphi \right) \left(D e^{T} \cdot
\varphi \right) \\
&=& \tilde{\varphi} D \tilde{\varphi}  
\end{eqnarray*}  
which proves that $H$ satisfies the superconformal condition
(\ref{nice superconformal condition}).  The uniqueness follows {}from
Proposition \ref{Ebijection}, the uniqueness of $\sqrt{\quad}$ on $\{1
+ \sum_{n \in \Z} c_n x^n \; | \; c_n \in R \}$, and the fact that if we
write $H = (f(x) + \varphi \xi(x), \psi(x) + \varphi g(x))$ then {}from
equation (\ref{H}), we can directly observe that $g(x)$ is of the form
$g(x) = 1 + \sum_{j \in \Z} A_j \left(\frac{j + 1}{2} \right) x^j + $
(terms in powers of $x$ strictly greater than 1), i.e., $g(x) \in \{1
+ \sum_{n \in \Z} c_n x^n \; | \; c_n \in R \}$.  
\end{proof}  

\begin{rema}
Here we note a difference between the super and nonsuper cases.  
In the nonsuper case (see \cite{H book}), the fact that 
\[\exp \Bigl(\sum_{j \in \Z} A_j  x^{j + 1} \frac{\partial}{\partial x}\Bigr) \cdot x , \] 
for $A_j \in \mathbb{C}$, gives a formal analytic function is trivial.
However, in the super case, as one can observe {}from the necessary
machinery involved in the proof of Proposition \ref{superconformal}, 
the proof that the analogous expression involving superderivations 
(\ref{H}) is superconformal is highly nontrivial.  
\end{rema}

\begin{rema}\label{envelope again} 
That $T$ as defined by (\ref{define T}) is an even superderivation is 
due to the fact that $T$ exists in the $R$-envelope of 
$\mbox{Der}(\mathbb{C}[[x]][\varphi])$; see Remark \ref{envelope}.   
Indeed in the work that follows, much of our ability to extend Huang's  
methods in \cite{H thesis} and \cite{H book} to the $N=1$ super case  
relies on this fact -- that we are working in the envelope of a Lie  
superalgebra which, by definition, is an ordinary Lie algebra.
\end{rema}

Now we would like to include formal superconformal power series vanishing at
zero with the even part of the coefficient of $x$ not necessarily one.  
For $b \in R^0$, we define the linear operators $b^{2x
\frac{\partial}{\partial x}}$ and $b^{\varphi \frac{\partial}{\partial \varphi}}$
{}from $R[x, x^{-1}, \varphi]$ to itself
by    
\begin{eqnarray*}
b^{2x \frac{\partial}{\partial x}} \cdot c \varphi^m x^n &=& c
\varphi^m b^{2n} x^n\\
b^{\varphi \frac{\partial}{\partial \varphi}} \cdot c
\varphi^m x^n &=& c b^m \varphi^m x^n
\end{eqnarray*} 
for $c \in R$, $m \in \mathbb{Z}_2$, and $n \in \mathbb{Z}$.  Then the operator
$b^{\left(\twoLo\right)} = b^{2x \frac{\partial}{\partial x}} b^{\varphi
\frac{\partial}{\partial \varphi}}$ is a well-defined linear operator
on $R [x,x^{-1},\varphi]$.  These operators can be extended to
operators on $R[[x,x^{-1}]][\varphi]$ in the obvious way.  We note
that for $H(x,\varphi) \in R[[x,x^{-1}]][\varphi]$, we have   
\begin{equation}\label{a0property}
b^{\left( \twoLo \right)} \cdot H(x, \varphi) \; = \;
H(b^{\left( \twoLo \right)} \cdot (x, \varphi)) \; = \; H(b^2x, b\varphi).  
\end{equation}

If $H$ is of the form (\ref{H}), in order for $H(b^2x,b\varphi)$ to 
correspond to an invertible local coordinate chart vanishing at
zero, we must have $b \in (R^0)^\times$, i.e., $b$ must be an 
invertible even element of the underlying superalgebra $R$.

\begin{rema}\label{square root in R}
The operation $b^{\left( \twoLo \right)} \cdot H(x,\varphi)$ for $H$ 
a power series of the form (\ref{H}), results in a formal power 
series vanishing at zero with the even part of the coefficient of 
$x$ equal to $b^2$ and the even part of the coefficient of $\varphi$ 
equal to $b$.  Thus, in keeping with the notation that the even 
coefficient of $x^{j+1}$ in $f(x)$ is denoted by $a_j$, the  
coefficient $b^2$ can be thought of as $a_0$ and $b$ can be thought 
of as a square root of $a_0$.  In \cite{B thesis} and \cite{B thesis 
announcement}, we assumed a well-defined square root on $(R^0)^\times$, 
which is equivalent to choosing a branch cut for the complex logarithm 
when $R = \bigwedge_*$ (cf. Remark \ref{not using branch cut}).  We 
used the notation $a_0$ and $\sqrt{a_0}$,  making it necessary to keep 
track of what square root was being used, or we used ``$\sqrt{a_0}$" as 
a composite symbol to denote a specified element of $(R^0)^{\times}$ 
such that $\sqrt{a_0}^2 = a_0$.  However, it is more natural to give 
$\sqrt{a_0} = b$, an invertible even element of $R$, as the basic data 
avoiding the need to keep track of a well-defined square root on 
$(R^0)^\times$.  In order to avoid confusion as to whether $\sqrt{a_0}$ 
is a specified element of $(R^0)^{\times}$ or some well-defined square 
root of an element in $(R^0)^{\times}$, we will use the notation 
$\asqrt \in (R^0)^{\times}$ as our basic data such that $\asqrt$ is the 
even coefficient of $\varphi$  and $\asqrt^2$ is the even coefficient 
of $x$ in $\asqrt^{\left(\twoLo \right)} \cdot H(x,\varphi)$ for $H$ a 
power series of the form (\ref{H}).
\end{rema}

Extending the notation (\ref{L notation}) to $j = 0$, let
\[L_0(x,\varphi) =  - \biggl( x \frac{\partial}{\partial x} +
\frac{1}{2}\varphi \frac{\partial}{\partial \varphi} \biggr),\]
which is an even superderivation in $\mbox{Der}(R[[x,x^{-1}]][\varphi])$.

\begin{prop}\label{azero}
Let $\asqrt \in (R^0)^{\times}$.  Then 
\begin{equation}\label{L0}
\asqrt^{-2L_0(x,\varphi)} \cdot (x,\varphi) = \asqrt^{\left(\twoLo
\right)} \cdot (x,\varphi) = 
\left(\asqrt^2 x, \asqrt \varphi \right) 
\end{equation}
is superconformal.
\end{prop}
 
\begin{proof} 
Writing $\asqrt^{-2L_0(x,\varphi)} \cdot (x,\varphi) =
(\tilde{x},\tilde{\varphi})$, we have   
\[D \tilde{x} = D \cdot (\asqrt^2 x) = \varphi \asqrt^2 = (\asqrt\varphi)
(\asqrt) = \tilde{\varphi}D \tilde{\varphi}. \] 
\end{proof}

Note that the operator $\asqrt^{-2L_0(x,\varphi)}$ should not be read
as $(\asqrt^{-2})^{L_0(x,\varphi)}$ but as $\asqrt^{\left( \twoLo \right)}$ 
thus retaining the basic data $\asqrt$ rather than just $\asqrt^2$; see Remark
\ref{square root in R}.

For any $(A,M) \in R^{\infty}$, we define a map $\tilde{E}$ {}from
$R^\infty$ to  the set of all formal superconformal power series in
$xR[[x]][\varphi]$ with leading even coefficient of $\varphi$ equal to
one, by defining
\begin{eqnarray}
\varphi \tilde{E}^0(A,M)(x,\varphi) &=& \varphi \biggl( x + \sum_{j \in
\Z} E_j(A,M) x^{j + 1} \biggr), \label{Etilde1}\\
\varphi \tilde{E}^1(A,M)(x,\varphi) &=& \varphi \sum_{j \in \Z} E_{j -
\frac{1}{2}}(A,M) x^j , \label{Etilde2}
\end{eqnarray}
and letting $\tilde{E}(A,M) (x, \varphi)$ be the unique formal 
superconformal power series with even coefficient of $\varphi$ equal to one 
such that
\[\varphi (\tilde{E}(A,M) (x, \varphi)) = \varphi(\tilde{E}^0(A,M), \tilde{E}^1(A,M)). \]
For $\asqrt \in (R^0)^\times$, we define a map $\hat{E}$ {}from
$(R^0)^\times
\times R^\infty$ to the set of all formal superconformal power series in 
$xR[[x]][\varphi]$ with invertible leading even coefficient of $\varphi$, by
defining   
\begin{eqnarray*}
\hat{E}^0(\asqrt,A,M)(x,\varphi) &=& \asqrt^2 \tilde{E}^0(A,M)(x,\varphi) , \\ 
\hat{E}^1(\asqrt,A,M)(x,\varphi) &=& \asqrt \tilde{E}^1(A,M)(x,\varphi) , 
\end{eqnarray*} 
and setting
\[\hat{E}(\asqrt,A,M)(x,\varphi) \; = \; (\hat{E}^0(\asqrt, A, M)
(x,\varphi), \hat{E}^1(\asqrt, A, M) (x,\varphi)) .\]
Then $\hat{E}(\asqrt,A,M) (x, \varphi)$ is the unique formal
superconformal power series satisfying 
\[\varphi (\hat{E}(\asqrt,A,M)(x,  \varphi)) \; = \;
\varphi(\hat{E}^0(\asqrt,A,M),\hat{E}^1(\asqrt,A,M))  \]
with even coefficient of $\varphi$ equal to $\asqrt$. 
The following proposition is an immediate consequence of Propositions
\ref{Ebijection}, \ref{superconformal}, and \ref{azero}.

\begin{prop}\label{above}
The map $\hat{E}$ {}from $(R^0)^{\times} \times R^{\infty}$ to the
set of all formal superconformal power series $H(x, \varphi)$ of the form
\begin{equation}\label{Ehat} 
\varphi H(x,\varphi) = \varphi \Biggl(\asqrt^2 \Bigl( x + \sum_{j \in \Z}
a_j x^{j + 1} \Bigr), \asqrt \sum_{j \in \Z} m_{j - \frac{1}{2}} x^j
\Biggr)  
\end{equation} 
and with even coefficient of $\varphi$ equal to $\asqrt$ 
for $(\asqrt,a,m) \in (R^0)^{\times} \times R^\infty$, is a
bijection.

The map $\tilde{E}$ {}from $R^{\infty}$ to the set of
formal superconformal power series of the form (\ref{Ehat}) with
$\asqrt = 1$ and even coefficient of $\varphi$ equal to 1 is also a 
bijection. 

In particular, we have inverses $\tilde{E}^{-1}$ and $\hat{E}^{-1}$.
\end{prop}

We will use the notation:
\[T_H (x, \varphi) = - \sum_{j \in \Z} \left(E^{-1}_j (a,m) 
L_j(x,\varphi) + E^{-1}_{j - \frac{1}{2}} (a,m) G_{j - \frac{1}{2}}
(x,\varphi) \right) \]  
for $H(x, \varphi)$ formally superconformal of the form (\ref{Ehat}) 
with even coefficient of $\varphi$ equal to $\asqrt$. Thus any 
superconformal power series $H(x,\varphi)$ of this form can be 
written uniquely as  
\begin{equation}\label{Hexpansion} 
H(x, \varphi) = e^{T_H (x, \varphi)} \cdot \asqrt^{-2L_0(x,\varphi)} 
\cdot (x, \varphi) . 
\end{equation}

Recalling (\ref{power series at ith puncture}), we know that
a local superconformal coordinate map vanishing at $0 \in 
\bigwedge_{*>0}$ with $\rho H(w, \rho) = \rho(f(w), \psi(w))$ is 
completely determined by $f(w)$, $\psi(w)$ and a choice of square 
root for $f'(w) + \psi(w) \psi'(w)$, where $f(w)$ can be expanded 
in a power series of the form $\asqrt^2 (w + \sum_{j \in \Z}  a_j
 w^{j + 1}) $ with $a_j \in \bigwedge_{*>0}^0$ and $\asqrt \in 
(\bigwedge_{*>0}^0)^{\times}$, and $\psi(w)$ can be expanded in a 
power series of the form $\asqrt \sum_{j \in \Z} m_{j - 1/2} w^j$ 
with $m_{j - 1/2} \in \bigwedge_*^1$, such that these power 
series are absolutely convergent to $f(w)$ and $\psi(w)$, 
respectively, in some neighborhood of zero.  Thus we see that 
formal superconformal power series $H(x,\varphi)$ of the form
(\ref{Hexpansion}) can be thought of as the ``local formal
superconformal coordinate maps vanishing at zero'' or the ``local 
formal superconformal coordinate transformations fixing the 
coordinates of a fixed point to be zero".  {}From 
(\ref{Hexpansion}), we see that the ``local formal superconformal  
transformations superanalytic and vanishing at zero'' are generated  
uniquely by the ``infinitesimal formal superconformal transformations''
of the form   
\[ (\log \asqrt) 2L_0(x,\varphi) + \sum_{j \in \Z} \left(
A_j L_j(x,\varphi) +  M_{j - \frac{1}{2}} G_{j - \frac{1}{2}}(x,\varphi) \right) \]  
(except for single-valuedness), for $\asqrt \in (R^0)^\times$, and
$\{(A_j,M_{j - 1/2})\}_{j \in \Z} \in R^\infty$. 
Proposition \ref{above} states that these ``infinitesimal
superconformal transformations'' can be identified with elements in
$(R^0)^\times \times R^\infty$.

\begin{prop}\label{auto2}
Let $u, v \in R((x))[\varphi]$; let $(A,M) \in R^\infty$; and let
\begin{eqnarray}
\hspace{.4in} T \! \! &=& \! \! - \sum_{j \in \Z} \left( A_j L_j(x,\varphi) + 
M_{j - \frac{1}{2}} G_{j - \frac{1}{2}}(x,\varphi) \right)  \label{T} \\
&=& \! \! \sum_{j \in \Z} \biggl( A_j \Bigl( \Lx \Bigr) + M_{j - \frac{1}{2}} \Gx \biggr) .
\nonumber
\end{eqnarray}
Then
\begin{equation}\label{specialized automorphism property}
e^T \cdot (uv) = \left( e^T \cdot u \right)\left( e^T \cdot v \right).
\end{equation}
\end{prop}

In other words, for $T$, $u$, and $v$ given above, the automorphism
property, Proposition \ref{autoprop}, holds if $y$ is set equal to 1.

\begin{proof}  For $i \in \mathbb{Z}$, the term 
$\frac{y^n T^n}{n!} \cdot x^i$ has powers in $x$ greater than or equal 
to $i + n - 1$ for $n \in \Z$, and $\frac{y^n T^n}{n!} \cdot \varphi 
x^i$ has powers in $x$ greater than or equal to $i + n$ for $n \geq 0$.  
Thus setting $y = 1$ in equation (\ref{autoequation}) of Proposition 
\ref{autoprop} applied to this case, each power series in equation 
(\ref{specialized automorphism property}) has only a finite number of
$x^j$ terms for a given $j \in \mathbb{Z}$, i.e., each term is a 
well-defined power series in $R((x))[\varphi]$. \end{proof}

\begin{prop}\label{Switch}
Let $\overline{H}(x, \varphi) \in R((x))[\varphi]$, and $H(x, \varphi)
= e^{T_H (x, \varphi)} \cdot (x, \varphi)$.  Then   
\begin{equation}\label{switch}
\overline{H}(H(x, \varphi)) = \overline{H}(e^{T_H (x, \varphi)}
\cdot (x, \varphi)) = e^{T_H (x, \varphi)} \cdot \overline{H}(x,
\varphi) .  
\end{equation}
\end{prop}

\begin{proof}  Write $H(x, \varphi) = e^{T_H} \cdot (x,
\varphi) = (e^{T_H} \cdot x, e^{T_H} \cdot \varphi) = (\tilde{x},
\tilde{\varphi})$. Equation (\ref{switch}) is trivial for
$\overline{H}(x, \varphi) = 1$, $\overline{H}(x, \varphi) = x$, and
$\overline{H}(x, \varphi) = \varphi$.  

(i) We prove the result for $\overline{H}(x, \varphi) = x^n$, with $n \in
\mathbb{N}$ and $n > 1$, by induction on $n$.  Assume $e^{T_H} \cdot x^k =
(e^{T_H} \cdot x)^k$ for $k \in \mathbb{N}$, $k < n$. Let 
$\overline{H}(x, \varphi) = x^n$.  Then by Proposition \ref{auto2}, 
\begin{eqnarray*}
e^{T_H} \cdot \overline{H}(x, \varphi) &=& e^{T_H} \cdot (x x^{n - 1}) =
\left( e^{T_H} \cdot x \right) \left( e^{T_H} \cdot x^{n - 1} \right) \\ 
&=& \left( e^{T_H} \cdot x \right) \left( e^{T_H} \cdot x \right)^{n - 1} =
\tilde{x} \tilde{x}^{n - 1} = \tilde{x}^n = \overline{H}(H(x,
\varphi)) . 
\end{eqnarray*}

(ii) We next prove the result for $\overline{H}(x, \varphi) = x^{-n}$,
$n \in \Z$.  Again by Proposition \ref{auto2}, 
\[ 1 = e^{T_H} \cdot (x x^{-1}) = (e^{T_H} \cdot x)(e^{T_H} \cdot x^{-1}) .\]
Thus \[ e^{T_H} \cdot x^{-1} = (e^{T_H} \cdot x)^{-1} .\]

Assume $e^{T_H} \cdot x^{-k} = (e^{T_H} \cdot x)^{-k}$ for $k \in
\mathbb{N}$, $k < n$. Let $\overline{H}(x, \varphi) = x^{-n}$.  Then by
Proposition \ref{auto2},   
\begin{eqnarray*}
e^{T_H} \cdot \overline{H}(x, \varphi) \! &=& \! e^{T_H} \cdot (x^{-1}
x^{-(n - 1)} ) = \left( e^{T_H} \cdot x^{-1} \right) \left( e^{T_H}
\cdot x^{-(n - 1)} \right) \\ 
&=& \! \left( e^{T_H} \cdot x \right)^{-1} \left( e^{T_H} \cdot x
\right)^{-(n - 1)} = \tilde{x}^{-1} \tilde{x}^{-(n - 1)} =
\tilde{x}^{-n} \\
&=&\! \overline{H}(H(x, \varphi)) . 
\end{eqnarray*}
Thus the result is true for $\overline{H}(x, \varphi) = x^n$, $n \in
\mathbb{Z}$.

(iii) For $\overline{H}(x, \varphi) = \varphi x^n$, $n \in \mathbb{Z}$,
we note that by Proposition \ref{auto2}, and the above cases for  
$\overline{H}(x, \varphi) = x^n$, $n \in \mathbb{Z}$, and
$\overline{H}(x, \varphi) = \varphi$,   
\begin{eqnarray*}
e^{T_H} \cdot \overline{H}(x, \varphi) &=& e^{T_H} \cdot (\varphi x^n) =
\left( e^{T_H} \cdot \varphi \right) \left( e^{T_H} \cdot x^n \right) \\ 
&=& \left( e^{T_H} \cdot \varphi \right) \left( e^{T_H} \cdot x \right)^n =
\tilde{\varphi} \tilde{x}^n = \overline{H}(H(x, \varphi)) . 
\end{eqnarray*} 

Since $e^{T_H} \cdot \varphi^i x^n \in R((x))[\varphi]$, for $i = 0, 1$
and $n \in \mathbb{Z}$, the result follows by linearity. 
\end{proof}

\begin{prop}\label{inverses}
Any formal superconformal power series $H(x, \varphi)$ of the form
(\ref{Hexpansion}) has a unique inverse with respect to composition of
formal power series, and this inverse is superconformal.  That is,
there exists a formal superconformal power series $H^{-1}(x, \varphi)$
of the form (\ref{Hexpansion}) such that 
\[H(H^{-1}(x, \varphi)) = (x, \varphi), \quad \mbox{and} \quad
H^{-1}(H(x, \varphi))  = (x, \varphi) .\]
If $\hat{E}^{-1} (H(x, \varphi)) = (\asqrt, A, M)$, i.e.,
\begin{equation}\label{explicit H at zero}
H(x,\varphi) = \exp \Biggl(   - \sum_{j \in \Z} \left( A_j L_j(x,\varphi) 
+ M_{j - \frac{1}{2}} G_{j - \frac{1}{2}}(x,\varphi) \right) \! \! \Biggr) 
\! \cdot  \asqrt^{-2L_0(x,\varphi)} \! \cdot (x,\varphi),
\end{equation}
then
\begin{eqnarray}
& &  \hspace{-.4in} H^{-1}(x, \varphi) \label{inverse}\\
\hspace{.5in} &=& \! \! \exp \Biggl(  \sum_{j \in \Z} \left( \asqrt^{-2j} A_j L_j(x,\varphi) 
+ \asqrt^{-2j +1} M_{j - \frac{1}{2}} G_{j - \frac{1}{2}}(x,\varphi) \right) 
\! \! \Biggr) \! \cdot \nonumber \\ 
& &  \hspace{2.8in} \cdot \asqrt^{2L_0(x,\varphi)} \! \cdot (x,\varphi) \qquad
\nonumber \\    
&=& \! \! \asqrt^{2L_0(x,\varphi)} \cdot \exp  \Biggl( \sum_{j \in \Z} 
\left( A_j L_j(x,\varphi) + M_{j - \frac{1}{2}} G_{j - \frac{1}{2}}
(x,\varphi) \right) \! \! \Biggr) \! \cdot  (x, \varphi) \nonumber \\
&=& \! \! \asqrt^{ 2L_0(x,\varphi)} \cdot \exp \left( -T_H (x,
\varphi) \right) \cdot (x, \varphi) \; \in R[[x]][\varphi]
. \hspace{.55in} \nonumber
\end{eqnarray}  
\end{prop}

\begin{proof} {}from (\ref{inverse}) and Proposition
\ref{Switch}, we have
\begin{eqnarray*}
& & \hspace{-.5in} H(H^{-1}(x, \varphi)) \\
&=& \! \! \left. H(H^{-1}(\asqrt^2 x_1, \asqrt \varphi_1))
\right|_{(x_1, \varphi_1) = ( \asqrt^{-2} x, \asqrt^{-1} \varphi)} \\
&=& \! \! H  \Biggl( \exp \Biggl( \sum_{j \in \Z} \left( A_j L_j(x_1,\varphi_1)  + M_{j -
\frac{1}{2}} G_{j - \frac{1}{2}}(x_1,\varphi_1) \right) \! \Biggr) \Bigg. \cdot \\
& & \hspace{2.4in} \cdot \Biggl. \left. (x_1, \varphi_1) \Biggr) \right|_{(x_1, 
\varphi_1) = (\asqrt^{-2} x, \asqrt^{-1} \varphi)} \\   
&=& \! \! \left. H(\exp\left( - T_H (x_1, \varphi_1) \right) \cdot (x_1,
\varphi_1) ) \right|_{(x_1, \varphi_1) = ( \asqrt^{-2} x,
\asqrt^{-1} \varphi)} \\
&=& \! \! \left. \exp\left( - T_H (x_1, \varphi_1) \right) \cdot H(x_1,
\varphi_1) \right|_{(x_1, \varphi_1) = ( \asqrt^{-2} x, \asqrt^{-1}
\varphi)} \\  
&=& \! \! \exp\left( - T_H (x_1, \varphi_1) \right) \cdot
\exp\left(T_H (x_1, \varphi_1) \right) \cdot \asqrt^{-2L_0(x_1,\varphi_1)} \cdot \\ 
& & \hspace{2.5in} \bigl. \cdot (x_1, \varphi_1) ) \bigr|_{(x_1, \varphi_1) =
( \asqrt^{-2} x, \asqrt^{-1} \varphi)} \\  
&=& \! \! \left. (\asqrt^2 x_1, \asqrt \varphi_1) \right|_{(x_1, \varphi_1)
= (\asqrt^{-2} x, \asqrt^{-1} \varphi)} \\
&=& \! \! (x, \varphi) . 
\end{eqnarray*}

Similarly,  
\begin{eqnarray*}
& & \hspace{-.35in} H^{-1}(H(x, \varphi))\\
&=& \! \! \left. H^{-1} (H(\asqrt^{-2} x_1, \asqrt^{-1}
\varphi_1)) \right|_{(x_1, \varphi_1) = (\asqrt^2 x, \asqrt \varphi)} \\ [.05in]
&=& \! \! H^{-1}  \Biggl( \exp \Biggl(\! - \! \sum_{j \in \Z} \left( \asqrt^{-2j} A_j
L_j(x_1,\varphi_1) +  \asqrt^{-2j + 1} M_{j - \frac{1}{2}} G_{j -
\frac{1}{2}}(x_1,\varphi_1)
\right) \!
\! \Biggr) \cdot \Biggr. \\
& & \left. \hspace{2.8in} \cdot (x_1, \varphi_1) \Biggr) \right|_{(x_1, \varphi_1) = (\asqrt^2 x,
\asqrt \varphi)} \\ [.05in]  
&=& \! \! H^{-1} \Bigl( \asqrt^{ 2L_0(x_1,\varphi_1)} \! \cdot \exp
\left( T_H (x_1, \varphi_1) \right) \cdot \asqrt^{ -2L_0(x_1,\varphi_1)} \! \cdot \left.  (x_1,
\varphi_1) \Bigr) \right|_{(x_1, \varphi_1) = (\asqrt^2 x, \asqrt\varphi)} \\ [.05in]
&=& \! \! \asqrt^{ 2L_0(x_1,\varphi_1)} \cdot \exp
\left( T_H (x_1, \varphi_1) \right) \cdot \asqrt^{ -2L_0(x_1,\varphi_1)} \cdot \left.
H^{-1}(x_1, \varphi_1) \right|_{(x_1,
\varphi_1) = (\asqrt^2 x, \asqrt \varphi)} \\  [.05in]
&=& \! \! \asqrt^{ 2L_0(x_1,\varphi_1)} \cdot \exp
\left( T_H (x_1, \varphi_1) \right) \cdot \asqrt^{ -2L_0(x_1,\varphi_1)} \cdot 
\asqrt^{ 2L_0(x_1,\varphi_1)} \cdot\\ [.05in]
& & \hspace{1.6in} \bigl. \cdot  \exp \left(- T_H (x_1,
\varphi_1) \right)
\cdot (x_1, \varphi_1) ) \bigr|_{(x_1, \varphi_1) = (\asqrt^2 x,
\asqrt \varphi)} \\  
&=& \! \! \bigl. \asqrt^{ 2L_0(x_1,\varphi_1)} \cdot (x_1,
\varphi_1) \bigr|_{(x_1, \varphi_1) = (\asqrt^2 x, \asqrt \varphi)} \\ 
&=& \! \! (x, \varphi) . 
\end{eqnarray*}

Since the formal composition of two formal superconformal power
series is again superconformal, by Propositions \ref{superconformal}
and \ref{azero}, $H^{-1}(x, \varphi)$ is superconformal. 
\end{proof} 

\begin{rema} {}from the proposition above, we see that the set of
all formal superconformal power series of the form (\ref{Hexpansion})
with $\asqrt = 1$ is a group with composition as the group
operation.  This is the group of ``formal superconformal local
coordinate transformations fixing the coordinates of a fixed point to
be zero with leading even coefficient of $\varphi$ equal to one'' or the 
group of ``formal superconformal transformations vanishing at zero with 
leading even coefficient of $\varphi$ equal to one''. 
\end{rema}

Given $(A,M), (B,N) \in R^\infty$, let $H(x,\varphi)$ and $\overline{H}
(x, \varphi)$ be two formal superconformal power series of the form 
(\ref{H}) such that  
\begin{equation}\label{for composition}
\tilde{E}^{-1}(H(x, \varphi)) = (A,M) \quad \mbox{and} \quad
\tilde{E}^{-1}(\overline{H}(x, \varphi)) = (B,N) .
\end{equation}  
We define the {\it composition} $(A,M) \circ (B,N)$ of $(A,M)$, and 
$(B,N)$ by
\begin{equation}\label{sequencecompositiondef}
(A,M) \circ (B,N) = \tilde{E}^{-1}((\overline{H} \circ H) (x, \varphi)) 
\end{equation}
where $(\overline{H} \circ H) (x, \varphi)$ is the formal composition
of $H(x, \varphi)$ and $\overline{H}(x, \varphi)$.

\begin{prop}\label{sequencecompositionprop} 
The set $R^\infty$ is a group with the operation $\circ$. Let 
$(A,M) = \{ (A_j , M_{j - 1/2}) \}_{j \in \Z} \in R^\infty$,
and for any $t \in R^0$, define 
\[t(A,M) = \{ (t A_j, t M_{j - 1/2}) \}_{j \in \Z}. \]  
Then for $s,t \in R^0$,   
\begin{equation}\label{composition}
(s(A,M)) \circ (t(A,M)) = (s + t)(A,M) .
\end{equation} 
That is the map $t \mapsto t(A,M)$ is a homomorphism {}from the additive
group of $R^0$ to $R^\infty$.  In addition, $(R^0)^\infty$ and 
$(R^1)^\infty$ are  subgroups of $R^\infty$.
\end{prop}

\begin{proof}  Since the set of all formal superconformal
power series of the form (\ref{H}) is a group with composition as its
group operation, it is obvious {}from the definition of $\circ$ that
$R^{\infty}$ is a group with this operation.  Let $H_t (x,\varphi) =
\tilde{E} (t(A,M))$.  By Proposition \ref{Switch} 
\begin{eqnarray*}
H_t (H_s (x,\varphi)) &=& e^{T_{H_s} (x,\varphi)} \cdot H_t
(x,\varphi) \\ 
&=& e^{T_{H_s} (x,\varphi)} \cdot e^{T_{H_t} (x,\varphi)} \cdot
(x,\varphi) \\
&=& e^{s T_{H_1} (x,\varphi)} \cdot e^{t T_{H_1} (x,\varphi)} \cdot
(x,\varphi) \\
&=& e^{(s + t) T_{H_1} (x,\varphi)} \cdot (x,\varphi) \\
&=& H_{(s + t)} (x,\varphi) \\
&=& \tilde{E} ((s + t) (A,M)).
\end{eqnarray*}
Or equivalently,
\[ \tilde{E}^{-1} (H_t (H_s (x,\varphi)) = (s + t) (A,M) .\]
But then {}from the definition of $s (A,M) \circ t (A,M)$,
\[ \tilde{E}^{-1} (H_t (H_s (x,\varphi)) = s (A,M) \circ t (A,M) .\]
Thus we obtain equation (\ref{composition}).  Letting $\mathbf{0} = 
(\mathbf{0}, \mathbf{0})$ be the sequence consisting of all zeros in 
$R^\infty$, it is clear that $(R^0)^\infty \oplus \{\mathbf{0}\}$ and 
$(R^1)^\infty \oplus \{\mathbf{0}\}$ are subgroups of $R^\infty$.
\footnote{There is a misprint in the analogous proof to 
Proposition \ref{sequencecompositionprop} for the nonsuper case given 
in \cite{H book}.  In the proof of Proposition 2.1.14 in \cite{H book}, 
factoring out the $t_1$ and $t_2$ in the expressions $e^{l_{f_{t_1}(x)}}$ 
and $e^{l_{f_{t_2}(x)}}$ on the bottom line of p.46, one should obtain 
$e^{t_1 l_{f_1(x)}}$ and $e^{t_2 l_{f_1(x)}}$, not 
$e^{t_1 l_{f_0(x)}}$ and $e^{t_2 l_{f_0(x)}}$ as stated.}
\end{proof}

We can extend the composition $\circ$ defined by (\ref{sequencecompositiondef})
to $(R^0)^\times \times R^\infty$.  Given $(\asqrt,A,M), (\bsqrt, B,N) 
\in (R^0)^\times \times R^\infty$, let $H(x,\varphi)$ and $\overline{H} 
(x, \varphi)$ be two formal superconformal power series of the form  
(\ref{explicit H at zero}) such that  
\begin{equation*}
\hat{E}^{-1}(H(x, \varphi)) = (\asqrt,A,M) \quad \mbox{and} \quad
\hat{E}^{-1}(\overline{H}(x, \varphi)) = (\bsqrt,B,N) .
\end{equation*}    
Define
\begin{equation}\label{formal composition for R}
(\asqrt, A, M) \circ (\bsqrt, B,N) \; = \; \hat{E}^{-1}((
\overline{H}\circ H) (x,\varphi)) . 
\end{equation}
Then in terms of the composition defined on $R^\infty$ by 
(\ref{sequencecompositiondef}), we have
\begin{equation}\label{composition with a}
(\asqrt, A, M) \circ (\bsqrt, B,N) \; = \; \Bigl(\asqrt \bsqrt, (A, M) \circ 
\bigl\{\asqrt^{2j} B_j, \asqrt^{2j - 1} N_{j-\frac{1}{2}} \bigr\}_{j \in \Z}\Bigr) .
\end{equation}

\begin{rema}\label{group remark}
With the composition operation defined above, $(R^0)^\times
\times R^\infty$ is a group naturally isomorphic to the group of
all formal superconformal power series of the form (\ref{explicit H at zero}).
The subset $R^\infty$ is a subgroup of $(R^0)^\times \times R^\infty$
isomorphic to the group of all formal superconformal power series of
the form (\ref{H}).  The fact that we can define a group action on
$(R^0)^\times \times R^\infty$ allows us to study the group
$(R^0)^\times \times R^\infty$ instead of the group of ``formal
superconformal local coordinate transformations fixing the coordinates
of a fixed point to be zero". 
\end{rema}

We now want to consider the ``formal superconformal coordinate maps
vanishing at infinity.''  Let $H(x, \varphi) \in x^{-1} R[[x^{-1}]]
[\varphi]$ be superconformal with
\begin{equation}\label{atinfty}
\varphi H(x, \varphi) = \varphi \biggl( \frac{1}{x} + \sum_{j \in \Z}
a_j x^{-j - 1}, \sum_{j \in \Z} m_{j - \frac{1}{2}} x^{-j} \biggr) = 
\varphi(f(x),\psi(x)).
\end{equation}
Then $H$ must define a square root for 
\[f'(x) + \psi(x) \psi'(x) \in \Bigl\{ - x^{-2} - \sum_{n = 3}^{\infty} c_n x^{-n} \; 
\Big| \; c_n \in R \Bigr\} . \] 
Let 
\[1 + \sum_{m \in \Z} d_m x^m = \biggl(1 + \sum_{n = 3}^{\infty} c_n
x^{n - 2} \biggr)^{1/2}\]
as defined by (\ref{sqrt}).  The two possibilities for the square 
root that $H$ must define for $f'(x) + \psi(x) \psi'(x) =  - x^{-2} - 
\sum_{n = 3}^{\infty} c_n
x^{-n}$ are given by 
\[ \biggl( - x^{-2} - \sum_{n = 3}^{\infty} c_n x^{-n} \biggr)^{1/2} = \; 
\pm \;\frac{i}{x} \biggl( 1 + \sum_{m \in \Z} d_m x^{-m} \biggr) . \]
That is, if $H \in x^{-1} R[[x^{-1}]][\varphi]$ is superconformal with 
leading even coefficient of $x^{-1}$ equal to one, then
\begin{multline}\label{form of H at infinity} 
H(x,\varphi) = \Biggl( \! \frac{1}{x} + \! \sum_{j \in \Z}
a_j x^{-j - 1} \pm \frac{i\varphi}{x} \biggl( \sum_{j \in \Z} m_{j - \frac{1}{2}} x^{-j} \!
\biggr) \! \biggl( 1 + \! \sum_{m \in \Z} d_m x^{-m} \! \biggr), \Biggr.  \\ 
\Biggl. \sum_{j \in \Z} m_{j - \frac{1}{2}} x^{-j}\pm \frac{i\varphi}{x}
\biggl( 1 + \! \sum_{m \in \Z} d_m x^{-m} \biggr) \! \Biggr) . 
\end{multline}
Thus we can specify $H$ of the form (\ref{form of H at infinity}) by 
specifying $f(x)$, $\psi(x)$ and whether the even coefficient of 
$\varphi x^{-1}$ is $i$ or $-i$.

Define 
\begin{equation}
I (x, \varphi) = \Bigl(\frac{1}{x}, \frac{i \varphi}{x} \Bigr) \in 
x^{-1}R[[x^{-1}]][\varphi]. 
\end{equation}
Then $I$ is superconformal of the form (\ref{form of H at infinity}) with 
leading even coefficient of $\varphi x^{-1}$ equal to $i$, and $I^{-1} = 
(1/x, - i \varphi/x)$ is superconformal of the form 
(\ref{form of H at infinity}) with even coefficient of $\varphi x^{-1}$ 
equal to $-i$.  


We now want to use the results we have developed about formal
superconformal series vanishing at zero to express any formal
superconformal series vanishing at infinity and with even 
coefficient of $\varphi x^{-1}$ equal to $i$ in terms of
superderivations in $\mbox{Der}(R((x^{-1}))[\varphi])$.  

Let $H(x, \varphi)$ be superconformal of the form (\ref{form of H 
at infinity}) with leading even coefficient of $\varphi x^{-1}$ equal 
to $i$, and let $H_{-1}(x, \varphi) = H \circ I^{-1} (x, \varphi)$.  
Then $H_{-1}$ is superconformal satisfying (\ref{Ehat}) with $\asqrt 
= 1$ and leading even coefficient of $\varphi$ equal to one.  Thus
$H_{-1}$ is of the form (\ref{Hexpansion}) with $\asqrt = 1$ and has 
a well-defined compositional inverse $H_{-1}^{-1}(x,\varphi)$.  

Note that $H \circ I^{-1} \circ H_{-1}^{-1}(x, \varphi)$, and $I^{-1} 
\circ H_{-1}^{-1} \circ H(x, \varphi)$ are well-defined formal 
superconformal series in $R[[x]][\varphi]$ and $xR[[x^{-1}]][\varphi]$, 
respectively.  Moreover, it is clear that the compositional inverse of
$H$ is $H^{-1}(x,\varphi) = I^{-1} \circ H_{-1}^{-1}(x, \varphi) \in 
x^{-1} R[[x]][\varphi]$.


Recall the even and odd superderivations introduced in (\ref{L notation}) 
and (\ref{G notation}).  Extending these definitions to include $L_{-j}
(x,\varphi)$ and $G_{-j + 1/2}(x,\varphi)$, for $j \in \Z$, we see
that these are superderivations in $\mbox{Der}(R((x^{-1}))[\varphi])$.

\begin{prop}\label{Infinity} 
Given $H(x, \varphi)$ superconformal of the form (\ref{form of H at infinity})
with even coefficient of $\varphi x^{-1}$ equal to $i$, we have
\begin{eqnarray}
\hspace{.5in} H(x, \varphi) \! \! &=&  \! \! \exp \Biggl( \sum_{j \in \Z} \left(
E^{-1}_j (a,m) L_{-j}(x,\varphi)  \right. \Biggr. \label{infty}\\
& & \hspace{1.2in} \Biggl. \left. + \; iE^{-1}_{j - \frac{1}{2}} (a,m) G_{-j +
\frac{1}{2}}(x,\varphi)
\right) \! \Biggr) \! \cdot \Bigl(\frac{1}{x}, \frac{i \varphi}{x} \Bigr) \nonumber \\
&=& \! \! \exp \Biggl(\! - \! \sum_{j \in \Z} \biggl(E^{-1}_j (a,m) \biggl( \Lminus \biggr) \biggr.
\Biggr. \nonumber \\
& & \hspace{.7in} \Biggl. \biggl. + \; iE^{-1}_{j - \frac{1}{2}} (a,m)
\Gminus \biggr) \! \Biggr) \! \cdot \Bigl(\frac{1}{x}, \frac{i \varphi}{x} \Bigr) \nonumber \\
&=& \! \! \exp \biggl( T_{H_{-1}} \Bigl(\frac{1}{x}, \frac{i \varphi}{x}
\Bigr) \biggr) \cdot
\Bigl(\frac{1}{x}, \frac{i \varphi}{x} \Bigr) . \nonumber
\end{eqnarray} 
In addition, $H^{-1}(x, \varphi) = I^{-1} \circ H^{-1}_{-1} (x,
\varphi)$ as defined above is the inverse of $H(x, \varphi)$ with
respect to composition.  That is,  
\[ H(H^{-1} (x, \varphi)) = H^{-1}(H(x, \varphi)) = (x, \varphi) .\]
Moreover,
\begin{eqnarray}
\hspace{.4in} H^{-1} \circ I (x, \varphi) \! \! &=& \! \! H^{-1} \Bigl(\frac{1}{x}, \frac{i
\varphi}{x} \Bigr) \label{forsewing} \\
&=& \! \! \exp \Biggl( \! - \! \sum_{j \in \Z} \Bigl(
E^{-1}_j  (a,m) L_{-j}(x,\varphi) \Bigr. \Biggr.   \nonumber \\
& & \hspace{1in} \Biggl. \left.  + \; iE^{-1}_{j - \frac{1}{2}} (a,m)
G_{-j + \frac{1}{2}}(x,\varphi) \right) \! \Biggr) \! \cdot (x, \varphi) \nonumber \\
&=& \! \! \exp \biggl(- T_{H_{-1}} \Bigl(\frac{1}{x}, \frac{i \varphi}{x}\Bigr) \biggr) \cdot (x,
\varphi) .\nonumber
\end{eqnarray} 
\end{prop}

\begin{proof} Since $H$ is superconformal satisfying
(\ref{form of H at infinity}) with even coefficient of $\varphi x^{-1}$ 
equal to $i$, the power series $H_{-1} (x, \varphi) = H \circ I^{-1} (x,\varphi)
\in R[[x]][\varphi]$ is superconformal with 
\[\varphi H_{-1} (x, \varphi) = \varphi \biggl( x + \sum_{j \in \Z}
a_j x^{j + 1}, \sum_{j \in \Z} m_{j - \frac{1}{2}} x^j \biggr)  \]
and with the even coefficient of $\varphi$ equal to one.
Thus by Proposition \ref{superconformal}, we have
\begin{eqnarray*}
H_{-1} (x, \varphi) \! \! &=& \! \!  \exp \Biggl(\! - \! \sum_{j \in \Z} \left( E^{-1}_j
(a,m) L_j(x,\varphi) \right. \Biggr.\\
& & \hspace{1.5in} \Biggl. \left. + \; E^{-1}_{j - \frac{1}{2}} (a,m) G_{j -
\frac{1}{2}}(x,\varphi) \right) \!  \Biggr) \! \cdot (x, \varphi) \\
&=& \! \!  \exp ( T_{H_{-1}} (x, \varphi)) \cdot (x, \varphi) .
\end{eqnarray*}
Write $I(x, \varphi) = (1/x, i \varphi/x) = (\tilde{x}, \tilde{\varphi})$.  
By the chain rule
\[ ix^{-1} \frac{\partial}{\partial \tilde{\varphi}} = \frac{\partial
\tilde{\varphi}}{\partial \varphi} \frac{\partial}{\partial
\tilde{\varphi}} = \frac{\partial}{\partial \varphi} -  \frac{\partial
\tilde{x}}{\partial \varphi} \frac{\partial}{\partial \tilde{x}} =
\frac{\partial}{\partial \varphi} \]
and 
\[ -x^{-2} \frac{\partial}{\partial \tilde{x}} = \frac{\partial
\tilde{x}}{\partial x} \frac{\partial}{\partial \tilde{x}} =
\frac{\partial}{\partial x} - (- \frac{\partial
\tilde{\varphi}}{\partial x} \frac{\partial}{\partial \tilde{\varphi}})
= \frac{\partial}{\partial x} + \varphi x^{-1}
\frac{\partial}{\partial \varphi} .\]  
Therefore 
\begin{eqnarray*}
& & \hspace{-.4in} H(x, \varphi) = H \circ I^{-1} \circ I (x, \varphi)
= H \circ I^{-1} (\tilde{x}, \tilde{\varphi}) = H_{-1} (\tilde{x},
\tilde{\varphi})  \\
&=& \! \! \exp ( T_{H_{-1}} (\tilde{x}, \tilde{\varphi})) \cdot (\tilde{x},
\tilde{\varphi}) \\ 
&=& \! \! \exp \Biggl( \sum_{j \in \Z} \left( E^{-1}_j (a,m) \biggl( x^{- j 
- 1} \Bigl( -x \varphi \frac{\partial}{\partial \varphi} - x^2
\frac{\partial}{\partial x} \Bigr) \biggr. \right. \Biggr. \\
& &  \biggl. + \; \Bigl(\frac{j + 1}{2} \Bigr)
\frac{i \varphi}{x} x^{-j} \Bigl(-ix \frac{\partial}{\partial \varphi} 
\Bigr) \! \biggr)  +  E^{-1}_{j - \frac{1}{2}}(a,m)
x^{-j} \Bigl( -ix \frac{\partial}{\partial \varphi}   \Bigr. \\
& & \hspace{2.2in} \Biggl. \Bigl. - \; \frac{i \varphi}{x} \Bigl( -x \varphi 
\frac{\partial}{\partial \varphi} - x^2 \frac{\partial}{\partial x} \Bigr) 
\Bigr) \! \biggr) \! \! \Biggr) \! \cdot \! \Bigl(\frac{1}{x}, \frac{i \varphi}{x} \Bigr) \\
&=& \! \! \exp \Biggl( \!  -  \! \sum_{j \in \Z} \biggl( E^{-1}_j (a,m) \Bigl(
\Lminus \Bigr) \biggr. \Biggr. \\
& & \hspace{1.7in} \Biggl. \biggl. + \; E^{-1}_{j - \frac{1}{2}} (a,m)
i \Gminus \! \biggr) \! \! \Biggr) \! \cdot \! \Bigl(\frac{1}{x}, \frac{i \varphi}{x} \Bigr) 
\end{eqnarray*}  
which gives (\ref{infty}).

By Proposition \ref{inverses}, we know that $H_{-1} (x, \varphi)$ has
a unique inverse $H_{-1}^{-1} (x, \varphi)$ with
\begin{eqnarray*}
H_{-1}^{-1} (x, \varphi) \! \! &=& \! \! \exp(T_{H_{-1}^{-1}} (x, \varphi) ) \cdot
(x, \varphi) \\
&=& \! \! \exp \Biggl( \! - \! \sum_{j \in \Z} \biggl( E^{-1}_j (a,m) \biggl( \Lx 
\biggr) \biggr. \Biggr. \\ 
& & \Biggr. \biggl. \hspace{1.5in} + \; E^{-1}_{j - \frac{1}{2}} (a,m)
\Gx \! \biggr) \! \! \Biggl) \cdot (x, \varphi) \\
&=& \! \! \exp(- T_{H_{-1}} (x, \varphi) ) \cdot (x, \varphi) .
\end{eqnarray*}

Setting $H^{-1} (x, \varphi) = I^{-1} \circ H_{-1}^{-1} (x,
\varphi)$, and since $H (x, \varphi) = H_{-1} \circ I (x, \varphi)$,
we have 
\begin{eqnarray*}
H \circ H^{-1} (x, \varphi) \! &=& \! H_{-1} \circ I \circ I^{-1} \circ
H_{-1}^{-1} (x, \varphi) \; = \; (x, \varphi) \\
H^{-1} \circ H (x, \varphi) \! &=& \! I^{-1} \circ H_{-1}^{-1} \circ H_{-1}
\circ I (x, \varphi) \;= \;(x, \varphi) .
\end{eqnarray*}
Moreover, by Proposition \ref{Switch}, with $\overline{H} (x, \varphi)
= I^{-1} (x, \varphi)$ and $H(x, \varphi)$ in the Proposition replaced
by $H_{-1}^{-1} (x, \varphi)$, we have
\begin{eqnarray*}
H^{-1} \circ I (x, \varphi) \! \! &=& \! \! I^{-1} \circ H_{-1}^{-1} \circ I (x,
\varphi) \; = \; I^{-1} \circ H_{-1}^{-1} \Bigl(\frac{1}{x},\frac{i \varphi}{x} 
\Bigr) \\
&=& \! \!  I^{-1} \circ \exp \Bigl(- T_{H_{-1}} \Bigl(\frac{1}{x},
\frac{i \varphi}{x} \Bigr) \Bigr) \cdot \Bigl(\frac{1}{x},\frac{i \varphi}{x} 
\Bigr) \\
&=& \exp \Bigl(- T_{H_{-1}} \Bigl(\frac{1}{x},\frac{i \varphi}{x} \Bigr) \Bigr) 
\cdot I^{-1} \Bigl(\frac{1}{x},\frac{i \varphi}{x} \Bigr)\\  
&=& \! \!  \exp \Bigl(- T_{H_{-1}} \Bigl(\frac{1}{x},\frac{i \varphi}{x}
\Bigr) \Bigr) \cdot (x, \varphi) 
\end{eqnarray*}
which gives (\ref{forsewing}).
\end{proof}

\begin{rema} The formal superconformal power series of the form
(\ref{infty}) can be thought of as the ``formal superconformal local
coordinate maps vanishing at $\infty = (\infty, 0)$''.
\end{rema}

The following two propositions are analogous to Proposition \ref{auto2}
and Proposition \ref{Switch}, respectively.

\begin{prop}\label{auto3}
Let $u, v \in R((x^{-1}))[\varphi]$; let $(B,N) \in R^\infty$; and let  
\begin{eqnarray}
\hspace{.3in} \bar{T} \! \! \! &=& \! \! - \! \sum_{j \in \Z} 
\left( B_j L_{-j}(x,\varphi) + N_{j -
\frac{1}{2}} G_{-j + \frac{1}{2}}(x,\varphi)\right) \label{T bar} \\
&=& \! \! \! \sum_{j \in \Z} \! \left( \! B_j \Bigl( \Lminus \Bigr) \right.
\nonumber \\
& & \hspace{1.9in} \left. + \; N_{j - \frac{1}{2}} \Gminus \! \right) .\nonumber
\end{eqnarray}
Then 
\begin{equation} \label{specialized automorphism property at infinity} 
e^{\bar{T}} \cdot (uv) = \left( e^{\bar{T}} \cdot u \right)\left(
e^{\bar{T}} \cdot v 
\right). 
\end{equation}
\end{prop}

In other words, for $\bar{T}$, $u$, and $v$ given above, the
automorphism property, Proposition \ref{autoprop}, holds if $y$ is set
equal to 1. 

\begin{proof}  For $i \in \mathbb{Z}$, and $j = 0,1$, each
$\frac{y^n \bar{T}^n}{n!} \cdot \varphi^j x^i$ has powers in $x$ less
than or equal to $i - n$ for $n \in \mathbb{N}$.  Thus setting $y = 1$ in
equation (\ref{autoequation}) of Proposition \ref{autoprop} applied to
this case, each power series in equation (\ref{specialized automorphism 
property at infinity}) has only a finite number of $x^j$ terms for a 
given $j \in \mathbb{Z}$, i.e., each term is a well-defined power series 
in $R((x^{-1}))[\varphi]$. \end{proof} 

\begin{prop}\label{Switch2}
Let $H(x, \varphi) = e^{\bar{T}} \cdot (x, \varphi)$ with $\bar{T}$
given by (\ref{T bar}), and let $\overline{H}(x, \varphi) \in
R((x^{-1}))[\varphi]$.  Then   
\begin{equation}\label{switch2}
\overline{H}(H(x, \varphi)) = \overline{H}(e^{\bar{T}} \cdot (x,
\varphi)) = e^{\bar{T}} \cdot \overline{H}(x,\varphi) .  
\end{equation}
\end{prop}

\begin{proof}  The proof is identical to steps (i), (ii),
and (iii) in the proof of Proposition \ref{Switch}.  To finish the
proof, we only need note that since $e^{\bar{T}} \cdot \varphi^i x^n
\in R((x^{-1}))[\varphi]$ for $i = 0, 1$ and $n \in \mathbb{Z}$, the
result follows by linearity.  \end{proof}

Generalizing Propositions \ref{Switch} and \ref{Switch2}, we have
the following proposition.

\begin{prop}  If $H, \overline{H} \in R[[x,x^{-1}]] [\varphi]$ with 
$H(x,\varphi) = e^T \cdot (x,\varphi)$ for some even superderivation 
$T \in \mathrm{Der} (R[[x,x^{-1}]]  [\varphi])$ such that $Ta = 0$
for $a \in R$, and if either $\overline{H} \circ H(x,\varphi)$ or 
$e^T \cdot \overline{H}(x,\varphi)$  exist in $R[[x,x^{-1}]][\varphi]$, 
then both exist in $R[[x,x^{-1}]]  [\varphi]$ and are equal.
\end{prop}

\begin{proof}  We need only define the existence of the necessary 
expressions since the proof of equality is the same as that for
Propositions \ref{Switch} and \ref{Switch2}.  

Suppose $H(x,\varphi) = e^T \cdot (x,\varphi)$, $\overline{H}(x,\varphi)$ 
and $e^T \cdot \overline{H}(x,\varphi)$ exist in $R[[x,x^{-1}]] [\varphi]$.
Writing
\[\overline{H}(x,\varphi) = \sum_{n \in \mathbb{Z}} a_n x^n + \varphi
\sum_{n \in \mathbb{Z}} b_n x^n\] 
for $a_n, b_n \in R$, then 
\[e^T \cdot \overline{H}(x,\varphi) = \sum_{n \in \mathbb{Z}} a_n e^T \cdot x^n + 
\sum_{n \in \mathbb{Z}} b_n e^T \cdot \varphi x^n .\]
Thus since $e^T \cdot \overline{H}(x,\varphi)$ exists in $R[[x,x^{-1}]] [\varphi]$,
and  
\begin{eqnarray*}
e^T \cdot x^n &=& \sum_{m \in \mathbb{N}} \frac{1}{m!} T^m x^n \\
e^T \cdot \varphi x^n &=& \sum_{m \in \mathbb{N}} \frac{1}{m!} T^m \varphi x^n 
= \sum_{m \in \mathbb{N}} \frac{1}{m!} \sum_{k \in \mathbb{Z}} \binom{m}{k}
(T^{m-k} \varphi) (T^k x^n),
\end{eqnarray*}
this implies that $T^m \cdot x^n$, for $a_n \neq 0$, and $T^m \cdot 
\varphi x^n$ for $b_n \neq 0$, exist in $R[[x,x^{-1}]] [\varphi]$.
By definition, for $n \in \mathbb{Z}_+$
\begin{eqnarray}
T^m x^n &=& \sum_{j_1 + \cdots + j_n = m} \frac{m!}{j_1!\cdots j_n!} 
\prod_{l=1}^{n} \left(T^{j_l} x \right)  \label{for switch 1} \\
T^m \cdot \varphi x^n &=& \sum_{k \in \mathbb{N}} \binom{m}{k}
(T^{m-k} \varphi) (T^k x^n) \label{for switch 2} \\
&=& \sum_{j_1 + \cdots + j_{n+1} = m}  
\frac{m!}{j_1!\cdots j_{n+1}!} (T^{j_{n+1}} \varphi) 
\prod_{l=1}^{n} \left(T^{j_l} x \right) , \nonumber  
\end{eqnarray}
and thus the right-hand sides of (\ref{for switch 1}) and (\ref{for switch 2}) 
must exist in $R[[x,x^{-1}]][\varphi]$ for $a_n \neq 0$ and $b_n \neq 0$,
respectively.  But then
\begin{eqnarray}\label{use auto}
(e^T \cdot x)^n &=& \sum_{j_1,...,j_n \in \mathbb{N}} \prod_{l = 1}^n \left( 
\frac{1}{j_l!} T^{j_l} x \right) \\
&=& \sum_{m \in \mathbb{N}} \frac{1}{m!} \sum_{j_1 + \cdots + j_n = m} 
\frac{m!}{j_1!\cdots j_n!} \prod_{l=1}^{n} \left(T^{j_l} x \right) \nonumber\\
&=& e^T \cdot x^n \nonumber
\end{eqnarray}
and
\begin{eqnarray}\label{use auto again}
\lefteqn{(e^T \cdot \varphi) (e^T \cdot x)^n}  \\
&=& \biggl( \sum_{m \in \mathbb{N}} \frac{1}{m!}
T^m \varphi \biggr) \Biggl(\sum_{m \in \mathbb{N}} \frac{1}{m!} \sum_{j_1 + \cdots + j_n = m} 
\frac{m!}{j_1!\cdots j_n!} \prod_{l=1}^{n} \left(T^{j_l} x \right) \Biggr) \nonumber \\
&=& \sum_{m \in \mathbb{N}} \frac{1}{m!} \sum_{j_1 + \cdots + j_{n +1}= m}  
\frac{m!}{j_1!\cdots j_{n+1}!} (T^{j_{n+1}} \varphi) 
\prod_{l=1}^{n} \left(T^{j_l} x \right) \nonumber \\
&=& e^T \cdot \varphi x^n \nonumber
\end{eqnarray} 
must also exist in $R[[x,x^{-1}]][\varphi]$ for $a_n \neq 0$ and $b_n 
\neq 0$, respectively.  Now note that equations (\ref{for switch 1}) and
(\ref{for switch 2}) hold if we replace $x$ by $x^{-1}$, and that these 
equations exist in $R[[x,x^{-1}]] [\varphi]$ for $a_{-n} \neq 0$ and 
$b_{-n} \neq 0$, respectively.  In this case (\ref{use auto}) and
(\ref{use auto again}) also exist in $R[[x,x^{-1}]] [\varphi]$ if we
replace $x$ by $x^{-1}$, i.e., if $a_{-n} \neq 0$ then $e^T \cdot x^{-n}
= (e^T \cdot x^{-1})^n$ exists in $R[[x,x^{-1}]] [\varphi]$, and if
$b_{-n} \neq 0$ then $e^T \cdot \varphi x^{-n} = (e^T \cdot \varphi) 
(e^T \cdot x^{-1})^n$ exists in $R[[x,x^{-1}]] [\varphi]$.

But since $x^n x^{-n} = 1$ and $e^T \cdot 1 = 1$, by the automorphism 
property Proposition \ref{autoprop}, we have that $e^T \cdot x^{-1}
= (e^T \cdot x)^{-1}$ which exists in $R[[x,x^{-1}]][\varphi]$, and
thus $(e^T \cdot x)^{-n}$ and $(e^T \cdot \varphi) (e^T \cdot x)^{-n}$ 
exist in $R[[x,x^{-1}]][\varphi]$ for $a_{-n} \neq 0$ and $b_{-n} \neq 0$,
respectively. Therefore
\[\overline{H}(e^T \cdot (x,\varphi)) = \sum_{n \in \mathbb{Z}} a_n (e^T \cdot x)^n + 
\sum_{n \in \mathbb{Z}} b_n (e^T \cdot \varphi)(e^T \cdot x)^n\]
exists in $R[[x,x^{-1}]][\varphi]$.  
  
In the case that $H(x,\varphi) = e^T \cdot (x,\varphi)$, 
$\overline{H}(x,\varphi)$ and $\overline{H}(e^T \cdot (x,\varphi))$ 
exist in $R[[x,x^{-1}]][\varphi]$, by  reversing the steps of the 
argument above, we conclude that $e^T \cdot \overline{H}(x,\varphi)$ 
exists in $R[[x,x^{-1}]] [\varphi]$.
\end{proof}

\section{The formal sewing equation and formal sewing identities} {}from the previous section, we see that in general the ``formal
infinitesimal superconformal transformations'' are of the form 
\begin{multline}\label{infinitesimal}
\sum_{j \in \Z} \left(B_j L_{-j}(x,\varphi) + N_{j - \frac{1}{2}}
G_{-j + \frac{1}{2}} \right)  + (\log \asqrt) 2L_0(x,\varphi) \\
+  \sum_{j \in \Z} \left( A_j L_j(x,\varphi)  + M_{j - \frac{1}{2}} G_{j -
\frac{1}{2}}(x,\varphi)\right)  
\end{multline}  
for $\asqrt \in (R^0)^\times$, $A_j, B_j \in R^0$, and $M_{j -
1/2}, N_{j - 1/2} \in R^1$.  In this section, we will use these
infinitesimal superconformal transformations to solve a formal 
version of the sewing equation along with the normalization and
boundary conditions defined in Chapter 2.  Recall that when
one supersphere is being sewn with another at the $i$-th puncture
of the first supersphere, the sewing equation is given by
\[F_\sou^{(1)}(w,\rho) = F_\sou^{(2)} \circ \hat{H}_0^{-1} \circ I 
\circ H_i (w,\rho)\]
where $\hat{H}_0$ is the local coordinate vanishing at $\infty$ of 
the second supersphere being sewn and $H_i$ is the local coordinate
vanishing at the $i$-th puncture, $(z_i,\theta_i) \in \bigwedge_\infty$,
of the first supersphere being sewn.  Thus formally we can write $H_i$
and $\hat{H}_0$ in terms of exponentials of infinitesimal superconformal
transformations and then try to solve for formal versions of 
$F_\sou^{(1)}$ and $F_\sou^{(2)}$.  To do this formally and 
algebraically without having to address issues of convergence, we will 
use additional formal variables such that certain formal infinite sums 
are well defined.  The geometric interpretation of this formal 
algebraic solution and questions of convergence will be addressed in
Chapter 4. 

Let $\mathcal{A} = \{ \mathcal{A}_{j} \}_{j \in \Z}$ and $\mathcal{B} = \{
\mathcal{B}_j \}_{j \in \Z}$ be two sequences of even formal variables,
and let $\alpha_0^{1/2}$ be another even formal variable. Let $\mathcal{M} 
= \{\mathcal{M}_{j - 1/2} \}_{j \in \Z}$ and $\mathcal{N} = \{ 
\mathcal{N}_{j - 1/2} \}_{j \in \Z}$ be two sequences of odd formal
variables.  Take the algebra $R$ of Section 3.2 to be 
\begin{multline*}
\mathbb{C}[\alpha_0^{\frac{1}{2}}, \alpha_0^{-\frac{1}{2}}] [[ \mathcal{A},
\mathcal{B}]] [\mathcal{M}, \mathcal{N}] \\
= \; \mathbb{C} [\alpha_0^{\frac{1}{2}},\alpha_0^{-\frac{1}{2}}] [[ \mathcal{A}_1, 
\mathcal{A}_2,..., \mathcal{B}_1, \mathcal{B}_2,...]] [\mathcal{M}_{\frac{1}{2}},
\mathcal{M}_{\frac{3}{2}},..., \mathcal{N}_{\frac{1}{2}}, \mathcal{N}_{\frac{3}{2}},
...] . 
\end{multline*}     

Let
\begin{multline}\label{formalzero}
H^{(1)}_{\alpha_0^{1/2}, \mathcal{A}, \mathcal{M}} (x, \varphi) \\
= \; \exp \Biggl( \! - \! \sum_{j \in \Z} \left( \mathcal{A}_j L_j(x,\varphi)
 + \mathcal{ M}_{j - \frac{1}{2}} G_{j - \frac{1}{2}} (x,\varphi)
\right) \! \! \Biggr) \! \cdot (\alpha_0^\frac{1}{2})^{-2L_0(x,\varphi)}
\cdot (x, \varphi) ,
\end{multline}
\begin{equation}\label{formalinfty}
H^{(2)}_{ \mathcal{B}, \mathcal{N}} (x, \varphi)  =  \exp \Biggl( 
\sum_{j \in \Z} \left( \mathcal{B}_j  L_{-j}(x,\varphi) + 
\mathcal{N}_{j - \frac{1}{2}}  G_{-j + \frac{1}{2}}(x,\varphi) 
\right) \! \! \Biggr) \! \cdot \! \Bigl(\frac{1}{x}, \frac{i \varphi}{x}\Bigr) . 
\end{equation} 

The formal power series $H^{(2)}_{ \mathcal{B}, \mathcal{N}} (x, \varphi)
\in x^{-1}R[[x^{-1}]][\varphi]$ is superconformal of the form
(\ref{infty}).  By Proposition \ref{Infinity}, its compositional 
inverse $(H^{(2)}_{ \mathcal{B}, \mathcal{N}})^{-1}  (x, \varphi)$ is a
well-defined element of $ x^{-1} R[[x]][\varphi]$, and by
(\ref{forsewing}) we have 
\begin{multline}\label{anotherinfty}
(H^{(2)}_{ \mathcal{B}, \mathcal{N}})^{-1} \Bigl(\frac{1}{x}, \frac{i
\varphi}{x}\Bigr) \\
= \; \exp \Biggl( \! - \! \sum_{j \in \Z} \left( \mathcal{B}_j 
L_{-j}(x,\varphi) + \mathcal{N}_{j - \frac{1}{2}} G_{-j + \frac{1}{2}}(x,\varphi) 
\right) \! \Biggr) \! \cdot (x, \varphi) . 
\end{multline}

\begin{prop}\label{Sewing1}
In $R[[x,x^{-1}]][\varphi]$, we have
\begin{multline}\label{eqnSewing1}
(H^{(2)}_{ \mathcal{B}, \mathcal{N}})^{-1} \circ I \circ
H^{(1)}_{\alpha_0^{1/2}, \mathcal{A}, \mathcal{M}} (x, \varphi) \\
= \; \exp \Biggl( \! - \! \sum_{j \in \Z} \Bigl( \mathcal{A}_j L_j(x,\varphi)  + 
\mathcal{M}_{j - \frac{1}{2}} G_{j - \frac{1}{2}}(x,\varphi)
\Bigr) \! \Biggr) \! \cdot (\alpha_0^\frac{1}{2})^{-2L_0(x,\varphi)} \cdot \\
\exp \Biggl( \! - \! \sum_{j \in \Z} \left( \mathcal{B}_j L_{-j}(x,\varphi) +
\mathcal{N}_{j - \frac{1}{2}} G_{-j + \frac{1}{2}}(x,\varphi) \right) \! \Biggr) \! 
\cdot (x, \varphi) 
\end{multline}  
where we always expand formal series in nonnegative powers of
$\mathcal{A}_j, \mathcal{B}_j, \mathcal{M}_{j - 1/2}$, and $\mathcal{N}_{j -
1/2}$, for $ j \in \Z$. 
\end{prop}

\begin{proof} Let 
\begin{multline*}
(H^{(2)}_{ \mathcal{B}, \mathcal{N}})^{-1} \circ I (x, \varphi) \\
= \; (x, \varphi) + \! \! \! \! \sum_{\begin{tiny} \begin{array}{c}
n,l \in \mathbb{N}\\
(n,l) \neq (0,0)
\end{array} \end{tiny}} 
\sum_{\begin{tiny} \begin{array}{c}
j_1 > ... > j_n > 0\\
k_1 > ... > k_l > 0
\end{array} \end{tiny}}
\sum_{m_1,...,m_n
\in \mathbb{Z}_{+}} \! \! \! H_{m_1,...,m_n}^{j_1, ... , j_n ;k_1, ... , k_l} (x,
\varphi) \cdot\\
\cdot \mathcal{B}_{j_1}^{m_1} \cdots  \mathcal{B}_{j_n}^{m_n} i\mathcal{N}_{k_1 - \frac{1}{2}}
\cdots i\mathcal{N}_{k_l - \frac{1}{2}} .
\end{multline*} {}from (\ref{anotherinfty}), we see that $H_{m_1,...,m_n}^{j_1, ... ,
j_n ;k_1, ... , k_l} (x, \varphi) \in \mathbb{Q} ((x))[\varphi]$.  Then by 
Proposition \ref{Switch}, we have  
\begin{multline*}
H_{m_1,...,m_n}^{j_1, ... , j_n ;k_1, ... , k_l}
(H^{(1)}_{\alpha_0^{1/2}, \mathcal{A}, \mathcal{M}} (x, \varphi) ) \\
= \; \exp \Biggl( \! - \! \sum_{j \in \Z} \left( \mathcal{A}_j L_j(x,\varphi)
+ \mathcal{M}_{j - \frac{1}{2}} G_{j - \frac{1}{2}}(x,\varphi) \right) \! \Biggr) \!
\cdot  (\alpha_0^\frac{1}{2})^{- 2L_0(x,\varphi)} \cdot \\
\cdot H_{m_1,...,m_n}^{j_1, ... , j_n ;k_1, ... , k_l} (x, \varphi) .
\end{multline*}
Therefore
\begin{multline*}
(H^{(2)}_{ \mathcal{B}, \mathcal{N}})^{-1} \circ I \circ
H^{(1)}_{\alpha_0^{1/2}, \mathcal{A}, \mathcal{M}} (x, \varphi) \\
= \; \exp \Biggl( \! - \! \sum_{j \in \Z} \left( \mathcal{A}_j L_j(x,\varphi)
+ \mathcal{M}_{j - \frac{1}{2}} G_{j - \frac{1}{2}}(x,\varphi) \right) \! \Biggr) \!
\cdot  (\alpha_0^\frac{1}{2})^{-2L_0(x,\varphi)} \cdot \\
\cdot (H^{(2)}_{ \mathcal{B}, \mathcal{N}})^{-1} \circ I (x, \varphi) .
\end{multline*}
Using (\ref{anotherinfty}), we obtain (\ref{eqnSewing1}). 
\end{proof}

\begin{rema}\label{explainsewing}
{}From Proposition \ref{Sewing1}, we see that the composition
\[(H^{(2)}_{ \mathcal{B}, \mathcal{N}})^{-1} \circ I \circ
H^{(1)}_{\alpha_0^{1/2}, \mathcal{A}, \mathcal{M}} (x, \varphi) \]
is generated by the formal infinitesimal superconformal 
transformations given by (\ref{infinitesimal}).  As shown in 
Chapter 2, geometrically this composition is the formal 
superconformal transition function of the sewn  neighborhoods of a 
supersphere with tubes sewn {}from two canonical superspheres with 
tubes (see equation (\ref{preliminary sewing equation})).  This is why 
we construct the formal superconformal transformations {}from the formal 
infinitesimal superconformal transformations in this way.  Of course,
in the formal version of the superconformal transition function
(\ref{eqnSewing1}), we have assumed that the puncture on the first 
supersphere into which the second supersphere is being sewn is at zero.  
In general this will be at some point $p \in U_\sou$ corresponding 
to $(z_i,\theta_i) = \sou(p)$.  Of course this discrepancy can be
rectified by appropriately incorporating the superconformal shift 
$s_{(z_i,\theta_i)} (x,\varphi) = (x - z_i - \varphi \theta_i, \varphi 
- \theta_i)$. 
\end{rema}

\begin{rema}\label{factoring out a in Chapter 3} 
Note the symmetry in the operators acting on $(x,
\varphi)$ on the right-hand side of equation (\ref{eqnSewing1}).  
Replacing $(x,\varphi)$, $(\mathcal{A},\mathcal{M})$, and 
$(\mathcal{B},\mathcal{N})$ with $I(x,\varphi) = (x^{-1},i
\varphi  x^{-1})$, $(\mathcal{B},-i\mathcal{N})$, and $(\mathcal{A},
-i\mathcal{M})$, respectively, for this operator, we obtain its 
inverse.  If we had not factored out $\asqrt$ as we did (see Remark
\ref{factoring out a in Chapter 2}), this symmetry would be broken.
\end{rema}

\begin{prop}\label{ultimate switch}
For $H(x,\varphi) \in R [x, x^{-1},\varphi]$, we have  
\begin{multline*}
H\Biggl( \! \exp \Biggl(\! - \! \sum_{j \in \Z} \biggl( \mathcal{A}_j L_j(x,\varphi)
+ \mathcal{M}_{j - \frac{1}{2}} G_{j - \frac{1}{2}}(x,\varphi) \biggr) \! \Biggr) \!
\cdot (\alpha_0^\frac{1}{2})^{- 2L_0(x,\varphi)} \cdot \Biggr. \\
\Biggl. \exp \Biggl( \! - \! \sum_{j \in \Z} \biggl( \mathcal{B}_j 
L_{-j}(x,\varphi) + \mathcal{N}_{j - \frac{1}{2}} G_{-j + \frac{1}{2}}(x,\varphi) 
\biggr)  \! \Biggr)\! \cdot (x,\varphi)\Biggr)
\end{multline*}
\begin{multline*}
= \exp \Biggl(\! - \! \sum_{j \in \Z} \biggl( \mathcal{A}_j L_j(x,\varphi)
+ \mathcal{M}_{j - \frac{1}{2}} G_{j - \frac{1}{2}}(x,\varphi) \biggr) \! \Biggr) \!
\cdot (\alpha_0^\frac{1}{2})^{- 2L_0(x,\varphi)} \cdot \\
\exp \Biggl( \! - \! \sum_{j \in \Z} \biggl( \mathcal{B}_j 
L_{-j}(x,\varphi) + \mathcal{N}_{j - \frac{1}{2}} G_{-j + \frac{1}{2}}(x,\varphi) 
\biggr)  \! \Biggr)\! \cdot H (x,\varphi)
\end{multline*}
\end{prop}

\begin{proof}  The result follows {}from Propositions
\ref{Switch} and \ref{Switch2}.  
\end{proof}

For the theorem below, it will be convenient to fix the following
notation.  Let  
\begin{eqnarray*}
g_\mathcal{A} = \sum_{j \in \Z} \mathcal{A}_j x^{j +1} , &\quad&
 g_\mathcal{M} = \sum_{j \in \Z} \mathcal{M}_{j - \frac{1}{2}} x^j , \\  
g_\mathcal{B} = \sum_{j \in \Z} \mathcal{B}_j x^{- j +1} ,
&\quad& g_\mathcal{N} = \sum_{j \in \Z} \mathcal{N}_{j - \frac{1}{2}} x^{- j +1} . 
\end{eqnarray*}
Let $(\mathbf{0},\mathbf{0}) = \mathbf{0} \in R^\infty$ be the sequence
consisting of all zeros.  Let $\tilde{R}$ be any superalgebra.
For $p(x) \in \tilde{R}[[x,x^{-1}]]$, we let $(p(x))^-$ and $(p(x))^+$
be the unique series such that $(p(x))^- \in \tilde{R}[[x^{-1}]]$,
$(p(x))^+ \in x \tilde{R} [[x]]$, and  
\[p(x) = (p(x))^- + (p(x))^+ .\]

\begin{thm}\label{uniformization} 
There exist formal series 
\[ \bar{F}^{(1)} \in x \mathbb{C} [\alpha_0^{\frac{1}{2}},
\alpha_0^{-\frac{1}{2}}] 
[x^{-1}, \varphi]  [[\mathcal{A}, \mathcal{B}]] [\mathcal{M}, \mathcal{N}] \]
and 
\[\bar{F}^{(2)} \in x \mathbb{C} [\alpha_0^\frac{1}{2},
\alpha_0^{-\frac{1}{2}}] [x, \varphi] [[\mathcal{A}, \mathcal{B}]] 
[\mathcal{M},  \mathcal{N}] \]  
which are superconformal in $(x,\varphi)$ satisfying the ``formal 
boundary conditions''
\begin{eqnarray}
\hspace{.4in} \left. \bar{F}^{(1)} (x, \varphi) \right|_{(\mathcal{A}, \mathcal{M}) = 
\mathbf{0}} &=& \! (\alpha_0^\frac{1}{2})^{2L_0(x,\varphi) }
\cdot \exp \Biggl( \! - \! \sum_{j \in \Z} \left( \mathcal{B}_j L_{-j}(x,\varphi) \right. \Biggr.   
\label{F1 boundary condition}  \\ 
& & \hspace{.1in} \Biggl. \left. + \; \mathcal{N}_{j - \frac{1}{2}} G_{-j + 
\frac{1}{2}}(x,\varphi) \right) \! \Biggr) \! \cdot (\alpha_0^\frac{1}{2})^{-2L_0(x,\varphi) }
\cdot (x,\varphi) \nonumber \\  
\left. \bar{F}^{(1)} (x, \varphi) \right|_{(\mathcal{B}, \mathcal{N}) =
\mathbf{0}} &=& \! (x,\varphi) \label{F1 other boundary condition} \\
\left. \bar{F}^{(2)} (x, \varphi) \right|_{(\mathcal{A}, \mathcal{M}) =
\mathbf{0}} &=& \! (\alpha_0^{-1} x, \alpha_0^{-\frac{1}{2}} \varphi)
\label{F2 other boundary condition} \\
\left. \bar{F}^{(2)} (x, \varphi) \right|_{(\mathcal{B}, \mathcal{N}) =
\mathbf{0}} &=& \! (\alpha_0^\frac{1}{2})^{ 2L_0(x,\varphi) }
\cdot \exp \Biggl( \sum_{j \in \Z} \left ( \mathcal{A}_j L_j(x,\varphi) \right. \Biggr.
\label{F2 boundary condition} \\ 
& & \hspace{1in} \Biggl. \left. + \; \mathcal{M}_{j - \frac{1}{2}}
G_{j - \frac{1}{2}}(x,\varphi) \right) \! \Biggr) \! \cdot (x,\varphi) , \nonumber 
\end{eqnarray}
and the conditions
\begin{multline}\label{specific F1}
\varphi \bar{F}^{(1)} (x, \varphi) = \varphi \left. \bar{F}^{(1)} (x, 
\varphi) \right|_{(\mathcal{A},\mathcal{M}) = \mathbf{0}} + \varphi 
\Biggl(\Bigl( \alpha_0^{-1} g_\mathcal{A} (x) \frac{\partial}{\partial x} 
g_\mathcal{B}(\alpha_0 x) \\
\Bigl. - \; \alpha_0^{-1} g_\mathcal{B} (\alpha_0 x) \frac{\partial}{\partial x} 
g_\mathcal{A} (x)\Bigr)^-   +  2  \alpha_0^{-\frac{1}{2}} \left(g_\mathcal{M} (x)
g_\mathcal{N}(\alpha_0 x)  \right)^- , \\
\alpha_0^{-\frac{1}{2}} \Bigl( g_\mathcal{A} (x) \frac{\partial}{\partial
x} g_\mathcal{N}(\alpha_0 x) - \frac{1}{2} g_\mathcal{N} (\alpha_0 x)
\frac{\partial}{\partial x} g_\mathcal{A} (x) \Bigr)^-  \\
+ \; \alpha_0^{-1} \Bigl( \frac{1}{2} g_\mathcal{M}
(x) \frac{\partial}{\partial x} g_\mathcal{B}(\alpha_0 x) - g_\mathcal{B}
(\alpha_0 x) \frac{\partial}{\partial x} g_\mathcal{M} (x)\Bigr)^- \Biggl. \Biggr) +
\varphi R^{(1)} 
\end{multline}  
\begin{multline}\label{specific F2}
\varphi \bar{F}^{(2)} (x, \varphi) = \varphi \left. \bar{F}^{(2)} (x, \varphi) 
\right|_{(\mathcal{B}, \mathcal{N}) = \mathbf{0}}  + \varphi \Biggl(\Bigl(g_\mathcal{B} (x)
\frac{\partial}{\partial x} g_\mathcal{A} (\alpha_0^{-1}x) \Bigr. \Biggr. \\
\Bigl. - \; g_\mathcal{A} (\alpha_0^{-1}x) \frac{\partial}{\partial x} 
g_\mathcal{B}(x) \Bigr)^+ - 2  \alpha_0^{-\frac{1}{2}}
\left(g_\mathcal{M} (x) g_\mathcal{N}(\alpha_0 x) \right)^+  ,\\
+ \; \alpha_0^{\frac{1}{2}} \Bigl( \frac{1}{2} g_\mathcal{N} (x)
\frac{\partial}{\partial x} g_\mathcal{A} (\alpha_0^{-1}x) - g_\mathcal{A} 
(\alpha_0^{-1}x) \frac{\partial}{\partial x} g_\mathcal{N}(x) \Bigr)^+ \\  
\Biggl. + \Bigl( g_\mathcal{B} (x) \frac{\partial}{\partial x} 
g_\mathcal{M} (\alpha_0^{-1} x) - \frac{1}{2} g_\mathcal{M} (\alpha_0^{-1} x) 
\frac{\partial}{\partial x} g_\mathcal{B}(x) \Bigr)^+ \Biggr) + \varphi R^{(2)} 
\end{multline}
where $R^{(1)}$ and $R^{(2)}$ are elements in 
\[\mathbb{C} [\alpha_0^{\frac{1}{2}}, \alpha_0^{-\frac{1}{2}}][[\mathcal{A}, 
\mathcal{B}]] [\mathcal{M}, \mathcal{N}] [[x^{-1}]] \]  
and 
\[\mathbb{C} [\alpha_0^{\frac{1}{2}}, \alpha_0^{-\frac{1}{2}}][[\mathcal{A}, 
\mathcal{B}]][\mathcal{M}, \mathcal{N}][[x]] ,  \] 
respectively, containing only terms with the total degree in the
$\mathcal{A}_j$'s and $\mathcal{M}_{j - 1/2}$'s at least one, total
degree in the $\mathcal{B}_j$'s and $\mathcal{N}_{j - 1/2}$'s at least
one, and total degree in the $\mathcal{A}_j$'s, $\mathcal{M}_{j -
1/2}$'s, $\mathcal{B}_j$'s and $\mathcal{N}_{j - 1/2}$'s at least three. 

Then $\bar{F}^{(2)} \circ (H_{\mathcal{B}, \mathcal{N}}^{(2)})^{-1}
\circ I \circ H_{ \alpha_0^{1/2}, \mathcal{A}, \mathcal{M}}^{(1)} (x,
\varphi)$ exists in 
\[x \mathbb{C} [\alpha_0^{\frac{1}{2}}, \alpha_0^{-\frac{1}{2}}] [x^{-1},
\varphi] [[\mathcal{A}, \mathcal{B}]] [\mathcal{M}, \mathcal{N}] , \] 
and there exist unique $\bar{F}^{(1)}$ and $\bar{F}^{(2)}$ satisfying
the above such that 
\begin{equation}\label{F bar equation}
\bar{F}^{(1)} (x, \varphi) = \bar{F}^{(2)} \circ (H_{\mathcal{B}, 
\mathcal{N}}^{(2)})^{-1} \circ I \circ H_{ \alpha_0^{1/2}, \mathcal{A}, 
\mathcal{M}}^{(1)} (x,
\varphi) .  
\end{equation}
We call equation (\ref{F bar equation}) the ``formal sewing equation".
\footnote{There is a misprint in the formulas (2.2.11) and (2.2.12) of 
Theorem 2.2.4 in \cite{H book} giving the analogous nonsuper case to our 
Theorem \ref{uniformization}.  The first two terms in the right-hand 
side of (2.2.11) should be replaced by $\alpha_0^{-1}(f_{\mathcal{B}}^{(2)}
)^{-1}(\frac{1}{\alpha_0 x})$  and the first two terms in the  
right-hand side of (2.2.12) should be replaced by $(f_{\mathcal{A},
\alpha_0}^{(1)})^{-1} (x)$.  Similar misprints occurred in \cite{H thesis}
and \cite{B thesis}.  These misprints were first corrected in \cite{BHL}.}
\end{thm}

\begin{proof} Write
\begin{eqnarray}
\varphi \bar{F}^{(1)} (x, \varphi) &=& \varphi \Bigl( x + \sum_{m,s,n,t
\in \mathbb{N}} h_{msnt}^0 (x) , \sum_{m,s,n,t \in \mathbb{N}} h_{msnt}^1
(x) \Bigr) \label{expand F1}\\
\varphi \bar{F}^{(2)} (x, \varphi) &=& \varphi \Bigl( \alpha_0^{-1} x
+ \sum_{m,s,n,t \in \mathbb{N}} k_{msnt}^0 (x) , \sum_{m,s,n,t \in
\mathbb{N}} k_{msnt}^1 (x) \Bigr) \label{expand F2}
\end{eqnarray}  
where $\bar{F}^{(2)}$ is superconformal, and 
\begin{eqnarray*}
(h_{msnt}^0 (x), h_{msnt}^1 (x)) \; = \; h_{msnt} (x) &\in& \! \mathbb{C}
[\alpha_0^{\frac{1}{2}}, \alpha_0^{-\frac{1}{2}}][[\mathcal{A}, \mathcal{B}]] 
[\mathcal{M}, \mathcal{N}] [[x^{-1}]] \\  
(k_{msnt}^0 (x), k_{msnt}^1 (x)) \; = \; k_{msnt} (x) &\in& \! x \mathbb{C}
[\alpha_0^{\frac{1}{2}}, \alpha_0^{-\frac{1}{2}}][[\mathcal{A}, \mathcal{B}]] 
[\mathcal{M}, \mathcal{N}][[x]] ,
\end{eqnarray*} 
are both homogeneous of degree $m$ in the $\mathcal{A}_j$'s, degree $s$ in the
$\mathcal{M}_{j - 1/2}$'s, degree $n$ in the $\mathcal{B}_j$'s, and degree
$t$ in the $\mathcal{N}_{j - 1/2}$'s,  for $j \in \Z$.  The fact
that $\bar{F}^{(2)}$ is superconformal implies that 
\begin{multline*}
\bar{F}^{(2)} (x, \varphi) = (\tilde{x}, \tilde{\varphi}) \\ 
= \biggl( \! \alpha_0^{-1} x + \! \! \sum_{m,s,n,t \in \mathbb{N}} \! \left( k_{msnt}^0 
(x) + \varphi q_{msnt}^1 (x) \right), \! \sum_{m,s,n,t \in \mathbb{N}} \! \left( 
k_{msnt}^1 (x)  + \varphi q_{msnt}^0 (x) \right) \! \biggr) \! \! 
\end{multline*}  
where each $q_{msnt} = (q_{msnt}^0, q_{msnt}^1)$ is homogeneous of
degree $m$ in the $\mathcal{A}_j$'s, degree $s$ in the $\mathcal{M}_{j -
1/2}$'s, degree $n$ in the $\mathcal{B}_j$'s, and degree $t$ in the 
$\mathcal{N}_{j - 1/2}$'s, for $j \in \Z$, and where $\bar{F}^{(2)} (x, 
\varphi) =  (\tilde{x},\tilde{\varphi})$ satisfies $D \tilde{x} = 
\tilde{\varphi} D \tilde{\varphi}$ for $D = \frac{\partial}{\partial 
\varphi} + \varphi \frac{\partial}{\partial x}$. (Of course $h_{msnt} = 
h_{msnt}^0$ for $s + t$ even, $h_{msnt} = h_{msnt}^1$ for $s + t$ odd, 
and similarly for $k_{msnt}$ and $q_{msnt}$.)

Thus letting $(H_{\mathcal{B}, \mathcal{N}}^{(2)})^{-1} \circ I \circ
H_{\alpha_0^{1/2}, \mathcal{A}, \mathcal{M}}^{(1)} (x, \varphi) = (H^0 (x,
\varphi), H^1 (x, \varphi))$, equation (\ref{F bar equation}) gives  
\begin{multline}\label{picking out terms 1}
\varphi \Biggl( x + \! \! \sum_{m,s,n,t \in \mathbb{N}} h_{msnt}^0 (x),
\sum_{m,s,n,t \in \mathbb{N}} h_{msnt}^1 (x) \Biggr) \\
= \; \varphi \Biggl( \alpha_0^{-1} H^0(x,\varphi) + \! \! \sum_{m,s,n,t \in
\mathbb{N}} \biggl(k_{msnt}^0( H^0(x, \varphi) ) + H^1(x, \varphi)
q_{msnt}^1 (H^0(x, \varphi) ) \biggr) ,\Biggr. \\ 
\Biggl. \sum_{m,s,n,t \in \mathbb{N}} \biggl(k_{msnt}^1 (
H^0(x, \varphi)) + H^1(x, \varphi) q_{msnt}^0 (H^0(x, \varphi) )
\biggr) \! \Biggr) .
\end{multline} {}from the boundary conditions (\ref{F1 boundary condition}) -
(\ref{F2 boundary condition}), we have
\begin{eqnarray*}
& & \hspace{-.4in} \varphi \biggl( \! x + \! \! \sum_{n,t \in \mathbb{N}} 
h_{00nt}^0 (x) ,\sum_{n,t \in \mathbb{N}} h_{00nt}^1 (x) \! \biggr) \\
&=& \varphi \left( \!\alpha_0^{-1} \biggl( (H_{\mathcal{B}, \mathcal{N}}^{(2)} )^{-1} \!
\circ I ( \alpha_0 x ,\alpha_0^{\frac{1}{2}} \varphi) \biggr)^0  , 
\alpha_0^{-\frac{1}{2}} \left( (H_{\mathcal{B}, \mathcal{N}}^{(2)} )^{-1} \circ I 
( \alpha_0 x , \alpha_0^{\frac{1}{2}} \varphi) \right)^1 \right) \\
&=& \varphi (\alpha_0^\frac{1}{2})^{2L_0(x,\varphi) } \cdot \exp \Biggl( \! - 
\! \sum_{j \in \Z} \left( \mathcal{B}_j L_{-j}(x,\varphi)  +  \mathcal{N}_{j -
\frac{1}{2}} G_{-j +  \frac{1}{2}}(x,\varphi) \right) \! \Biggr) \! \cdot \\
& & \hspace{3in} \cdot (\alpha_0^\frac{1}{2})^{-2L_0(x,\varphi) } \cdot (x,\varphi) ,
\end{eqnarray*}
\begin{eqnarray*}
\varphi \biggl( \! x + \! \! \sum_{m,s \in \mathbb{N}} h_{ms00}^0 (x) ,
\sum_{m,s \in \mathbb{N}} h_{ms00}^1 (x) \! \biggr) \! \! &=& \! \! \varphi (x, 0) , \\
\varphi \biggl(\! \alpha_0^{-1} x + \! \! \sum_{n,t \in \mathbb{N}} k_{00nt}^0
(x), \sum_{n,t \in \mathbb{N}} k_{00nt}^1 (x) \! \biggr) \! \! &=& \! \! \varphi
(\alpha_0^{-1} x, 0) , \\
\varphi \biggl( \! \alpha_0^{-1} x + \! \! \sum_{m,s \in \mathbb{N}} k_{ms00}^0 (x)
, \sum_{m,s \in \mathbb{N}} k_{ms00}^1 (x) \! \biggr) \! \! &=& \! \! \varphi (H_{\alpha_0^{1/2},
\mathcal{A}, \mathcal{M}}^{(1)})^{-1} (x, \varphi) \\  & & \hspace{-2.5in} =  \varphi
(\alpha_0^\frac{1}{2})^{ 2L_0(x,\varphi) }
\cdot \exp \Biggl( \sum_{j \in \Z} \left ( \mathcal{A}_j L_j(x,\varphi)  + \mathcal{M}_{j -
\frac{1}{2}} G_{j - \frac{1}{2}}(x,\varphi) \right) \! \Biggr) \! \cdot (x,\varphi) .
\end{eqnarray*}
These equations give $h_{ms00}, h_{00nt}, k_{ms00}$, and $k_{00nt}$ 
uniquely for all $m,s,n,t \in \mathbb{N}$.

By the superconformal condition $D\tilde{x} = \tilde{\varphi} D \tilde{\varphi}$ 
for $\bar{F}^{(2)} (x,\varphi) = (\tilde{x}, \tilde{\varphi})$ and the boundary
conditions, we see that  
\begin{eqnarray*} 
\sum_{m,s,n,t \in \mathbb{N}} q_{msnt}^0 (x) \! \! &=& \! \! \alpha_0^{-\frac{1}{2}}
\Bigl( 1 + \alpha_0 \biggl(\sum_{m,s,n,t \in \mathbb{N}} k_{msnt}^0 (x)
\biggl)' \Biggr. \hspace{1.1in}\\ 
& & \hspace{.4in} \Biggl. + \; \alpha_0 \biggl( \sum_{m,s,n,t \in
\mathbb{N}} k_{msnt}^1 (x) \biggr) \biggl( \sum_{m,s,n,t \in \mathbb{N}}
k_{msnt}^1 (x) \biggr)' \Biggr)^{1/2} \\  
\sum_{m,s,n,t \in \mathbb{N}} q_{msnt}^1 (x) \! \! &=& \! \! \alpha_0^{-\frac{1}{2}}
\biggl( \sum_{m,s,n,t \in \mathbb{N}} k_{msnt}^1 (x) \biggr) \biggl(1 +
\alpha_0 \biggr(\sum_{m,s,n,t \in \mathbb{N}} k_{msnt}^0 (x)\biggr)'
\biggr)^{1/2}  
\end{eqnarray*} 
where the square root is defined to be the Taylor series expansion about
$x = 0$ and $\mathcal{A} = \mathcal{M} = \mathcal{B} = \mathcal{N} = \mathbf{0}$ with
$\sqrt{1} = 1$.  Note that each $q_{msnt}$ is determined by $k_{ijlp}$ 
for $0 \leq i \leq m$, $0 \leq j \leq s$, $0 \leq l \leq n$, and $0 \leq p \leq t$.   
Thus the $h_{ms00}, h_{00nt}, k_{ms00}$, and $k_{00nt}$ that we have
determined {}from the boundary conditions for all $m,s,n,t \in \mathbb{N}$, uniquely
determine $q_{ms00}$ and $q_{00nt}$ for all $m,s,n,t \in \mathbb{N}$. {}from (\ref{picking out terms 1}), it is clear that the $h_{msnt}$ term
on the right-hand side only depends on $H^0$, $H^1$, and
$k_{ijlp}(H^0(x,\varphi))$ for $0 \leq i \leq m$, $0 \leq j \leq s$,
$0 \leq l \leq n$, and $0 \leq p \leq t$.  

Note that for all $k_{ms00}$, $k_{00nt}$, $q_{ms00}$, and $q_{00nt}$
the coefficient of a given term
$\mathcal{A}_{j_1}\cdots \mathcal{A}_{j_m}\mathcal{M}_{i_1 - 1/2}\cdots
\mathcal{M}_{i_s - 1/2}$ or $\mathcal{B}_{j_1}\cdots \mathcal{B}_{j_n}
\mathcal{N}_{i_1 - 1/2}\cdots \mathcal{N}_{i_t - 1/2}$ is in $\mathbb{C}[x]$, i.e., 
\[k_{ms00}, k_{00nt}, q_{ms00}, q_{00nt} \in \mathbb{C} [\alpha_0^{\frac{1}{2}},
\alpha_0^{-\frac{1}{2}}][x][[\mathcal{A}, \mathcal{B}]] [\mathcal{M}, \mathcal{N}]
. \] 
Thus by Propositions \ref{Sewing1} and \ref{ultimate switch}, we have
\begin{eqnarray*}
& & \hspace{-.4in} k_{ms00} (H^0(x,\varphi))  \\
&=& \! \! \! \sum_{\begin{tiny} \begin{array}{c}
j_1 \leq ... \leq j_m\\
i_1< ... < i_s \end{array} \end{tiny}} 
\! \! p_{i_1 \dots i_s}^{j_1 \dots j_m} ( H^0(x,\varphi))
\mathcal{A}_{j_1}\cdots \mathcal{A}_{j_m}\mathcal{M}_{i_1 - \frac{1}{2}} \cdots
\mathcal{M}_{i_s - \frac{1}{2}} \\
&=& \! \! \! \sum_{\begin{tiny} \begin{array}{c}
j_1 \leq ... \leq j_m\\
i_1 < ... < i_s \end{array} \end{tiny}} \! \! 
p_{i_1 \dots i_s}^{j_1 \dots j_m} \left( \! \exp \! \Biggl( \! - \! \sum_{j \in \Z}
\biggl( \mathcal{A}_j L_j(x,\varphi) + \mathcal{M}_{j - \frac{1}{2}} G_{j -
\frac{1}{2}}(x,\varphi) \biggr) \! \Biggr) \!  \cdot \right. \\
& & \! \left.  (\alpha_0^\frac{1}{2})^{- 2L_0(x,\varphi)}
\cdot \exp \Biggl( \! - \! \sum_{j \in \Z} \biggl( \mathcal{B}_j L_{-j}(x,\varphi) 
+  \mathcal{N}_{j - \frac{1}{2}} G_{-j + \frac{1}{2}}(x,\varphi) \biggr) \! \Biggr) \! 
\cdot x \! \right) \\ 
& & \hspace{2.8in} \mathcal{A}_{j_1}\cdots \mathcal{A}_{j_m}\mathcal{M}_{i_1 -
\frac{1}{2}}\cdots \mathcal{M}_{i_s - \frac{1}{2}} \\ 
&=& \! \! \! \sum_{\begin{tiny} \begin{array}{c}
j_1 \leq ... \leq j_m\\
i_1 < ... < i_s \end{array} \end{tiny}} 
\! \! \exp \! 
\Biggl( \! - \! \sum_{j \in \Z} \biggl( \mathcal{A}_j L_j(x,\varphi) + 
\mathcal{M}_{j - \frac{1}{2}} G_{j - \frac{1}{2}}(x,\varphi) \biggr) \! \Biggr) \!  
\cdot (\alpha_0^\frac{1}{2})^{-2L_0(x,\varphi)} \cdot \\
& &  \hspace{.9in} \exp \Biggl( \! - \! \sum_{j \in \Z} \biggl( \mathcal{B}_j L_{-j}
(x,\varphi) + \mathcal{N}_{j - \frac{1}{2}} G_{-j + \frac{1}{2}}(x,\varphi) \biggr) 
\! \Biggr) \cdot p_{i_1 \dots i_s}^{j_1 \dots j_m} (x)\\
&=& \! \! \! \sum_{\begin{tiny} \begin{array}{c}
j_1 \leq ... \leq j_m\\
i_1 < ... < i_s \end{array} \end{tiny}} 
\! \! \exp \! \Biggl( \! - 
\! \sum_{j \in \Z} \biggl( \mathcal{A}_j L_j(x,\varphi) + \mathcal{M}_{j - \frac{1}{2}} 
G_{j - \frac{1}{2}}(x,\varphi) \biggr) \! \Biggr) \!  \cdot
(\alpha_0^\frac{1}{2})^{-2L_0(x,\varphi)} \cdot\\
& & \hspace{.9in} \exp \Biggl(\! - \! \sum_{j \in \Z} \biggl( \mathcal{B}_j L_{-j}
(x,\varphi) +  \mathcal{N}_{j - \frac{1}{2}} G_{-j + \frac{1}{2}}(x,\varphi) \biggr) 
\! \Biggr) \cdot k_{ms00} (x) .
\end{eqnarray*}  
And similarly for $k_{00nt}$, $q_{ms00}$, $q_{00nt}$, and in fact for any 
\[p(x) \in \mathbb{C} [\alpha_0^{\frac{1}{2}}, \alpha_0^{-\frac{1}{2}}][x]
[[\mathcal{A}, \mathcal{B}]] [\mathcal{M}, \mathcal{N}] .\] 

Now we can solve equation (\ref{picking out terms 1}) for $h_{msnt}
(x)$ and $k_{msnt} (x)$ by induction on $m$, $s$, $n$, and $t$.  In
fact, if we compare terms which are homogeneous of degree $m$ in the
$\mathcal{A}_j$'s, degree $s$ in the $\mathcal{M}_{j - 1/2}$'s, degree 
$n$ in the $\mathcal{B}_j$'s, and degree $t$ in the $\mathcal{N}_{j - 
1/2}$'s, for $j \in \Z$, on both  sides of (\ref{picking out terms 1}), 
we have 
\begin{equation}\label{ktilde}
h_{msnt} (x) = k_{msnt} (\alpha_0 x) + \tilde{k}^{(msnt)} (x) +
\tilde{q}^{(msnt)} (x) 
\end{equation}
where $\tilde{k}^{(msnt)} (x)$ and $\tilde{q}^{(msnt)} (x)$ are
homogeneous of degree $m$ in the $\mathcal{A}_j$'s, degree $s$ in the
$\mathcal{M}_{j - 1/2}$'s, degree $n$ in the $\mathcal{B}_j$'s, and
degree $t$ in the $\mathcal{N}_{j - 1/2}$'s, for $j \in \Z$, and
depend only on $h_{ijlp}$ and $k_{ijlp}$ for $0 \leq i \leq m$, $0
\leq j \leq s$, $0 \leq l \leq n$, and $0 \leq p \leq t$, where $(i,j,l,p) 
\neq (m,s,n,t)$.  Assume $h_{ijlp}$ and $k_{ijlp}$ (and thus $q_{ijlp}$)
have already been obtained for $0 \leq i \leq m$, $0 \leq j \leq s$,
$0 \leq l \leq n$, and $0 \leq p \leq t$, where $(i,j,l,p) \neq (m,s,n,t)$, 
and assume that we have shown that each $k_{ijlp}$ and $q_{ijlp}$ is in
$\mathbb{C} [\alpha_0^{\frac{1}{2}}, \alpha_0^{-\frac{1}{2}}][x][[\mathcal{A}, 
\mathcal{B}]] [\mathcal{M}, \mathcal{N}]$. Then by using Proposition
\ref{ultimate switch} on each polynomial coefficient of each
$k_{ijlp}$ (as we did for $k_{ms00}(x)$), we can determine
$\tilde{k}^{(msnt)} (x)$ and $\tilde{q}^{(msnt)} (x)$.  Then {}from
(\ref{ktilde}), we have   
\[ h_{msnt} (x) = (\tilde{k}^{(msnt)} (x))^- + (\tilde{q}^{(msnt)}
(x))^- \]
\[ k_{msnt} (\alpha_0 x) = - (\tilde{k}^{(msnt)} (x))^+ -
(\tilde{q}^{(msnt)} (x))^+  \]
where 
\[(\tilde{k}^{(msnt)} (x))^- , (\tilde{q}^{(msnt)} (x))^- \in
\mathbb{C} [\alpha_0^{\frac{1}{2}}, \alpha_0^{-\frac{1}{2}}][[\mathcal{A},
\mathcal{B}]] [\mathcal{M}, \mathcal{N}] [[x^{-1}]] \] 
and 
\[(\tilde{k}^{(msnt)}(x))^+, (\tilde{q}^{(msnt)} (x))^+ \in \mathbb{C}
[\alpha_0^{\frac{1}{2}}, \alpha_0^{-\frac{1}{2}}][[\mathcal{A}, \mathcal{B}]] 
[\mathcal{M}, \mathcal{N}][[x]] \] 
such that 
\[(\tilde{k}^{(msnt)} (x))^- + (\tilde{k}^{(msnt)} (x))^+ =
\tilde{k}^{(msnt)} (x) , \]
and 
\[(\tilde{q}^{(msnt)} (x))^- + (\tilde{q}^{(msnt)} (x))^+ =
\tilde{q}^{(msnt)} (x) . \]
By the principle of induction, we obtain $h_{msnt} (x)$ and $k_{msnt}
(\alpha_0 x)$ for all $m,s,n,t \in \mathbb{N}$, and thus we obtain
$\varphi \bar{F}^{(1)} (x, \varphi)$ and $\varphi \bar{F}^{(2)} (x, 
\varphi)$.  It is clear {}from the procedure to solve (\ref{picking out
terms 1}) that the solutions $\varphi \bar{F}^{(1)} (x, \varphi)$ and
$\varphi \bar{F}^{(2)} (x, \varphi)$ are unique.  Furthermore note 
that the even coefficient of $\varphi$ in $\bar{F}^{(1)} (x, \varphi)$ 
is one and the even coefficient of $\varphi$ in $\left. \bar{F}^{(2)} 
(x, \varphi) \right|_{(\mathcal{A}, \mathcal{M}) = \mathbf{0}}$ is 
$\alpha_0^{-1/2}$.  Let $\bar{F}^{(1)} (x, \varphi)$ and 
$\bar{F}^{(2)} (x, \varphi)$ be the unique formal superconformal series 
satisfying $\varphi \bar{F}^{(1)} (x, \varphi)$  and $\varphi \bar{F}^{(2)} 
(x, \varphi)$ such that the even coefficient of $\varphi$ in $\bar{F}^{(1)} 
(x, \varphi)$ is one and the even coefficient of $\varphi$ in $\left. 
\bar{F}^{(2)} (x, \varphi) \right|_{ (\mathcal{A}, \mathcal{M}) = \mathbf{0}}$ is
$\alpha_0^{-1/2}$.   

To complete the proof, we note that 
\begin{eqnarray*}
h_{1010}(x) \! \! &=& \! \! k_{1010}(\alpha_0 x) + \alpha_0^{-1} g_\mathcal{B}
(\alpha_0 x)  \frac{\partial}{\partial x} k_{1000}(\alpha_0 x) +
g_\mathcal{A} (x)  \frac{\partial}{\partial x} k_{0010}(\alpha_0 x) \\
& & \hspace{2.7in} + \alpha_0^{-1} g_\mathcal{A} (x)
\frac{\partial}{\partial x} g_\mathcal{B}(\alpha_0 x) \\
&=& \! \! k_{1010}(\alpha_0 x) - \alpha_0^{-1} g_\mathcal{B} (\alpha_0 x)
\frac{\partial}{\partial x} g_\mathcal{A} (x) + \alpha_0^{-1} g_\mathcal{A}
(x)  \frac{\partial}{\partial x} g_\mathcal{B}(\alpha_0 x) .
\end{eqnarray*} 
Thus 
\[h_{1010}(x) = \left( \alpha_0^{-1} g_\mathcal{A} (x)
\frac{\partial}{\partial x} g_\mathcal{B}(\alpha_0 x) - \alpha_0^{-1}
g_\mathcal{B} (\alpha_0 x) \frac{\partial}{\partial x} g_\mathcal{A}
(x)\right)^- ,\]
and 
\[k_{1010}(x) = \left(g_\mathcal{B} (x) \frac{\partial}{\partial x}
g_\mathcal{A} (\alpha_0^{-1}x) - g_\mathcal{A} (\alpha_0^{-1}x)
\frac{\partial}{\partial x} g_\mathcal{B}(x) \right)^+ .\]
We note that
\begin{eqnarray*}
h_{1001}(x) \! \! &=& \! \! k_{1001}(\alpha_0 x) + g_\mathcal{N} (\alpha_0 x)
q_{1000}(\alpha_0 x) + g_\mathcal{A} (x)  \frac{\partial}{\partial x}
k_{0001}(\alpha_0 x) \\
& & \hspace{2.5in} + \alpha_0^{-\frac{1}{2}} g_\mathcal{A} (x)
\frac{\partial}{\partial x} g_\mathcal{N}(\alpha_0 x) \\
&=& \! \! k_{1001}(\alpha_0 x) - \frac{\alpha_0^{-\frac{1}{2}}}{2} g_\mathcal{N}
(\alpha_0 x) \frac{\partial}{\partial x} g_\mathcal{A}(x) +
\alpha_0^{-\frac{1}{2}} g_\mathcal{A} (x) \frac{\partial}{\partial x}
g_\mathcal{N}(\alpha_0 x) .\\ 
\end{eqnarray*} 
Thus 
\[h_{1001}(x) = \alpha_0^{-\frac{1}{2}} \left( g_\mathcal{A} (x)
\frac{\partial}{\partial x} g_\mathcal{N}(\alpha_0 x) - \frac{1}{2}
g_\mathcal{N} (\alpha_0 x) \frac{\partial}{\partial x} g_\mathcal{A} (x)
\right)^- ,\]   
and 
\[k_{1001}(x) = \alpha_0^{\frac{1}{2}} \left( \frac{1}{2} g_\mathcal{N} (x)
\frac{\partial}{\partial x} g_\mathcal{A} (\alpha_0^{-1}x) - g_\mathcal{A}
(\alpha_0^{-1}x) \frac{\partial}{\partial x} g_\mathcal{N}(x) \right)^+
.\] 
We note that
\begin{eqnarray*}
h_{0110}(x) \! \! &=& \! \! k_{0110}(\alpha_0 x) + \alpha_0^{-1} g_\mathcal{B}
(\alpha_0 x) \frac{\partial}{\partial x} k_{0100}(\alpha_0 x) +
\alpha_0^{\frac{1}{2}} g_\mathcal{M} (x) q_{0010}(\alpha_0 x) \\
& & \hspace{2.6in} + \frac{\alpha_0^{-1}}{2} g_\mathcal{M} (x)
\frac{\partial}{\partial x} g_\mathcal{B}(\alpha_0 x) \\ 
&=& \! \! k_{0110}(\alpha_0 x) - \alpha_0^{-1} g_\mathcal{B} (\alpha_0 x)
\frac{\partial}{\partial x} g_\mathcal{M} (x) + \frac{\alpha_0^{-1}}{2}
g_\mathcal{M} (x) \frac{\partial}{\partial x} g_\mathcal{B}(\alpha_0 x) .
\end{eqnarray*} 
Thus 
\[h_{0110}(x) = \alpha_0^{-1} \left( \frac{1}{2} g_\mathcal{M} (x)
\frac{\partial}{\partial x} g_\mathcal{B}(\alpha_0 x) - g_\mathcal{B}
(\alpha_0 x) \frac{\partial}{\partial x} g_\mathcal{M} (x)\right)^- ,\]   
and 
\[k_{0110}(x) =  \left( g_\mathcal{B} (x) \frac{\partial}{\partial x}
g_\mathcal{M} (\alpha_0^{-1} x) - \frac{1}{2} g_\mathcal{M} (\alpha_0^{-1}
x) \frac{\partial}{\partial x} g_\mathcal{B}(x)  \right)^+ .\]
And finally, we note that
\begin{eqnarray*}
h_{0101}(x) \! \! &=& \! \! k_{0101}(\alpha_0 x) + g_\mathcal{N} (\alpha_0 x)
q_{0100}(\alpha_0 x) + \alpha_0^{\frac{1}{2}} g_\mathcal{M} (x)
q_{0001}(\alpha_0 x) \\
& & \hspace{2.8in} + \alpha_0^{-\frac{1}{2}} g_\mathcal{M} (x) g_\mathcal{N}(\alpha_0 x) \\   
&=& \! \! k_{0101}(\alpha_0 x) - \alpha_0^{-\frac{1}{2}} g_\mathcal{N} (\alpha_0
x) g_\mathcal{M} (x) + \alpha_0^{-\frac{1}{2}} g_\mathcal{M} (x) g_\mathcal{N}(\alpha_0 x) .  
\end{eqnarray*} 
Thus 
\[h_{0101}(x) = 2 \alpha_0^{-\frac{1}{2}} \left(g_\mathcal{M} (x) g_\mathcal{N}(\alpha_0 x)  \right)^- ,\]   
and 
\[k_{0101}(x) =  - 2 \alpha_0^{-\frac{1}{2}} \left(g_\mathcal{M} (x)
g_\mathcal{N}(\alpha_0 x)  \right)^+ .\]
Substituting these into (\ref{expand F1}) and (\ref{expand F2}), we
obtain (\ref{specific F1}) and (\ref{specific F2}).  
\end{proof}

\begin{rema}\label{uniformization remark}
We noted in Remark \ref{explainsewing} that the geometric meaning
of the left-hand side of (\ref{eqnSewing1}) is the formal
superconformal coordinate transition function of the sewn 
neighborhoods of a supersphere with tubes sewn {}from two canonical
superspheres with tubes.  If we denote the two canonical 
superspheres being sewn by $S\hat{\mathbb{C}} = \bigwedge_\infty \cup
(\{\infty \} \times (\bigwedge_\infty)_S)$ with the $i$-th ($i \in \Z$) 
puncture $(z_i,\theta_i)$ of the first supersphere being sewn with the 
puncture at $\infty$ of the second supersphere, then geometrically 
$\bar{F}^{(1)} \circ s_{(z_i,\theta_i)}$ and $\bar{F}^{(2)}$ are formal 
versions of the two halves $F_\sou^{(1)}$ and $F_\sou^{(2)}$ of the 
uniformizing function $F$ which maps the genus-zero superconformal 
surface resulting {}from the sewn canonical superspheres with tubes to a 
super-Riemann sphere with tubes.  Then this super-Riemann sphere with 
tubes can be mapped to a canonical supersphere with tubes via a global 
superconformal transformation, i.e., a superprojective  transformation.  
Thus the meaning of $\bar{F}^{(1)} \circ s_{(z_i,\theta_i)} (x,\varphi)$ 
is the ``formal superconformal coordinate transition function of the first 
canonical supersphere taking $\infty$ of the first canonical supersphere 
to $\infty$ of the resulting super-Riemann sphere'', and the geometric 
meaning of $\bar{F}^{(2)} (x,\varphi)$ is the ``formal superconformal 
coordinate transition function of the second canonical supersphere 
taking $0$ of the second canonical supersphere to $0$ of the resulting 
super-Riemann sphere''.  Furthermore Theorem \ref{uniformization} shows 
that the uniformizing function (\ref{uniformizing function}) which is a 
solution to the sewing equation along with the normalization and boundary 
conditions (\ref{sewing equation}) -- (\ref{boundary conditions trivial 
puncture 2})  does in fact depend algebraically on the local coordinate 
charts at the $i$-th puncture for the first supersphere being sewn and at 
infinity for the second supersphere being sewn.   Moreover, this 
uniformizing function is uniquely determined by the formal sewing 
equation (\ref{F bar equation}), the formal boundary conditions 
(\ref{F1 boundary condition}) - (\ref{F2 boundary condition}), and  
(\ref{specific F1}) and (\ref{specific F2}) which contain the  
normalization conditions.  
\end{rema}

Since $\bar{F}^{(1)}$ and $\bar{F}^{(2)}$ are superconformal with 
\[ \bar{F}^{(1)} \in x \mathbb{C} [\alpha_0^{\frac{1}{2}}, \alpha_0^{-\frac{1}{2}}]
[[\mathcal{A}, \mathcal{B}]] [\mathcal{M}, \mathcal{N}] [[x^{-1}]] [\varphi] \]
where the even coefficient of $\varphi$ in $\bar{F}^{(1)}$ is equal to one,
and 
\[\bar{F}^{(2)} \in x\mathbb{C} [\alpha_0^{\frac{1}{2}},
\alpha_0^{-\frac{1}{2}}] [[\mathcal{A}, \mathcal{B}]] [\mathcal{M}, \mathcal{N}]
[[x]] [\varphi] ,\]
where the even coefficient of $\varphi$ in $\left. \bar{F}^{(2)} \right|_{
(\mathcal{A}, \mathcal{M})  = \mathbf{0}}$ is equal to $\alpha_0^{-1/2}$,  
by Propositions \ref{above} and \ref{Infinity}, there exist a unique
pair of sequences
\begin{equation}\label{where psi lives}
(\Psi_j, \Psi_{j - \frac{1}{2}}) = (\Psi_j, \Psi_{j - \frac{1}{2}})
(\alpha_0^\frac{1}{2}, \mathcal{A}, \mathcal{M}, \mathcal{B}, \mathcal{N}) 
\end{equation} 
in $\mathbb{C} [\alpha_0^{1/2}, \alpha_0^{-1/2}][[\mathcal{A},
\mathcal{B}]] [\mathcal{M}, \mathcal{N}]$ for $j \in \mathbb{Z}$, such that  
\begin{equation}\label{F1}
\bar{F}^{(1)} (x, \varphi) = \exp \Biggl( \sum_{j \in \Z} \Bigl(
\Psi_{-j} L_{-j}(x,\varphi)  + \Psi_{-j + \frac{1}{2}} G_{-j + \frac{1}{2}}(x,\varphi)
\Bigr) \! \Biggr) \! \cdot (x , \varphi) 
\end{equation}
and
\begin{multline}\label{F2}
\bar{F}^{(2)} (x, \varphi) = \exp \left(- \Psi_0 2L_0(x,\varphi)
\right) \cdot (\alpha_0^\frac{1}{2})^{2L_0(x,\varphi) } \cdot \\
\exp \Biggl(\! - \! \sum_{j \in \Z} \Bigl( \Psi_j L_j(x,\varphi) + \Psi_{j -
\frac{1}{2}} G_{j - \frac{1}{2}}(x,\varphi) \Bigr) \! \Biggr) \! \cdot (x ,
\varphi) ,
\end{multline}
i.e.,
\begin{eqnarray*}
I \circ \bar{F}^{(1)} \circ I^{-1} (x,\varphi) &=& \tilde{E}\Bigl(\bigl\{
\Psi_{-j},-i\Psi_{-j + \frac{1}{2}} \bigr\}_{j \in \Z}\Bigr) (x,\varphi) \\ 
&=& \hat{E}\Bigl(1,\bigl\{ \Psi_{-j},-i\Psi_{-j + \frac{1}{2}} \bigr\}_{j 
\in \Z}\Bigr) (x,\varphi)\\
\bar{F}^{(2)} (x,\varphi) &=& \hat{E}\Bigl(e^{\Psi_0} \alpha_0^{-\frac{1}{2}},
\bigl\{ \Psi_{j},\Psi_{j - \frac{1}{2}} \bigr\}_{j \in \Z}\Bigr) (x,\varphi).
\end{eqnarray*}

\begin{prop}\label{normal order}
For $j \in \Z$, we have \footnote{There is a misprint in equation 
(2.2.27) of Proposition 2.2.5 in \cite{H book} which gives the nonsuper 
version of equation (\ref{psi condition 2}) above.  The first term in
the right-hand side of (2.2.27) should be $-\alpha_0^{-j} \mathcal{B}_j$, 
not $-\alpha_0^{j} \mathcal{B}_j$ as stated.}
\begin{eqnarray}
(\Psi_j, \Psi_{j - \frac{1}{2}}) \! \! &=& \! \! (- \mathcal{A}_j , - \mathcal{M}_{j -
\frac{1}{2}}) + \mathcal{P}_j (\alpha_0^\frac{1}{2}, \mathcal{A}, \mathcal{M},
\mathcal{B}, \mathcal{N}) , \label{psi condition 1} \\  
\hspace{.2in} (\Psi_{-j}, \Psi_{- j + \frac{1}{2}}) \! \! &=& \! \! (- \alpha_0^{-j} 
\mathcal{B}_j, - \alpha_0^{-j + \frac{1}{2}}  \mathcal{N}_{j - \frac{1}{2}}) + 
\mathcal{P}_{-j} (\alpha_0^\frac{1}{2}, \mathcal{A}, \mathcal{M}, \mathcal{B}, 
\mathcal{N}),  \label{psi condition 2} \\  
\Psi_0 \! \! &=& \! \! 0 + \mathcal{P}_0 (\alpha_0^\frac{1}{2}, \mathcal{A}, \mathcal{M},
\mathcal{B}, \mathcal{N}) , \label{psi condition 3}  
\end{eqnarray}
where each $\mathcal{P}_j (\alpha_0^{1/2}, \mathcal{A}, \mathcal{M},
\mathcal{B}, \mathcal{N})$, for $j \in \mathbb{Z}$, contains only terms with 
total degree at least one in the $\mathcal{A}_k$'s and $\mathcal{M}_{k - 1/2}$'s, 
for $k \in \Z$, and with total degree at least one in the $\mathcal{B}_k$'s and 
$\mathcal{N}_{k - 1/2}$'s, for $k \in \Z$.  Both
\begin{equation}\label{endomorphism 1}
\exp \Biggl(\sum_{j \in \Z} \biggl( \Psi_{-j} L_{-j}(x,\varphi) +  \Psi_{-j + \frac{1}{2}}
G_{-j + \frac{1}{2}}(x,\varphi) \biggr) \! \Biggr) 
\end{equation}
and 
\begin{equation}\label{endomorphism 2}
\exp \Biggl(\sum_{j \in \Z} \biggl( \Psi_j L_j(x,\varphi) + \Psi_{j 
- \frac{1}{2}} G_{j - \frac{1}{2}}(x,\varphi) \biggr) \! \Biggr) 
\end{equation} 
are in the algebra $((\mbox{\em End} \; \mathbb{C} [x, x^{-1}, \varphi])
[\alpha_0^{1/2}, \alpha_0^{-1/2}][[ \mathcal{A}, 
\mathcal{B}]] [\mathcal{M}, \mathcal{N}])^0$, and in this algebra we have
\begin{multline}\label{normal order equation}
\exp \Biggl( \! - \! \sum_{j \in \Z} \biggl( \mathcal{A}_j L_j(x,\varphi) +
\mathcal{M}_{j - \frac{1}{2}} G_{j - \frac{1}{2}}(x,\varphi) \biggr) \! \Biggr) \cdot
(\alpha_0^\frac{1}{2})^{- 2L_0(x,\varphi) } \cdot \\
\exp \Biggl(\! - \! \sum_{j \in \Z} \biggl( \mathcal{B}_j L_{-j}(x,\varphi)
 + \mathcal{N}_{j - \frac{1}{2}} G_{-j + \frac{1}{2}}(x,\varphi)
\biggr) \! \Biggr) 
\end{multline}
\begin{multline*}
= \; \exp \Biggl(  \sum_{j \in \Z} \biggl( \Psi_{-j} L_{-j}(x,\varphi) + 
\Psi_{-j + \frac{1}{2}} G_{-j + \frac{1}{2}}(x,\varphi) \biggr) \! \Biggr) \cdot \\  
\exp \Biggl( \sum_{j \in \Z} \biggl( \Psi_j L_j(x,\varphi) + 
\Psi_{j - \frac{1}{2}} G_{j - \frac{1}{2}}(x,\varphi) \biggr) \! \Biggr) \cdot \\
(\alpha_0^\frac{1}{2})^{- 2L_0(x,\varphi) } \cdot
\exp \biggl(\Psi_0 2L_0(x,\varphi) \biggr)  .
\end{multline*}
\end{prop}

\begin{proof}
Equations (\ref{psi condition 1}), (\ref{psi condition 2}), and
(\ref{psi condition 3}) follow immediately {}from (\ref{specific F1}),
(\ref{specific F2}), (\ref{F1}) and (\ref{F2}).  By (\ref{where psi
lives}), we know that (\ref{endomorphism 1}) and (\ref{endomorphism
2}) are in the algebra $((\mbox{End} \; \mathbb{C}
[[x,x^{-1}]][\varphi])[\alpha_0^{1/2}, \alpha_0^{-1/2}]
[[ \mathcal{A}, \mathcal{B}]] [\mathcal{M}, \mathcal{N}])^0$.  By definition,
(\ref{endomorphism 1}) and (\ref{endomorphism 2}) applied to
$(x,\varphi)$ are in 
\[x \mathbb{C}[\alpha_0^{\frac{1}{2}},\alpha_0^{-\frac{1}{2}}]  [x^{-1}, \varphi] 
[[\mathcal{A}, \mathcal{B}]][\mathcal{M}, \mathcal{N}] \quad \mbox{and} 
\quad x \mathbb{C}[\alpha_0^{\frac{1}{2}},\alpha_0^{-\frac{1}{2}}] [x,\varphi] 
[[\mathcal{A}, \mathcal{B}]] [\mathcal{M}, \mathcal{N}], \]
respectively, and we have 
\begin{eqnarray*}
\bigl. \bar{F}^{(1)} \bigr|_{(\mathcal{A}, \mathcal{M}) = 
(\mathcal{B},\mathcal{N}) = \mathbf{0}} &=& (x,\varphi)\\ 
\bigl. \bar{F}^{(2)} \bigr|_{(\mathcal{A}, \mathcal{M}) = (\mathcal{B}, 
\mathcal{N}) = \mathbf{0}} &=& (\alpha_0^{-1}x, \alpha_0^{-\frac{1}{2}} \varphi). 
\end{eqnarray*}
Thus for any $H(x,\varphi) \in \mathbb{C} [\alpha_0^{1/2},\alpha_0^{-1/2}] 
[x, x^{-1},\varphi] [[ \mathcal{A}, \mathcal{B}]] [\mathcal{M}, \mathcal{N}]$, we
have 
\begin{equation}\label{first switch}
H \Biggl(\! \exp \Biggl( \sum_{j \in \Z} \biggl( \Psi_{-j} L_{-j}(x,\varphi) 
+ \Psi_{-j + \frac{1}{2}} G_{-j + \frac{1}{2}}(x,\varphi) \biggl) \! \Biggl) \! 
\cdot (x,\varphi) \Biggr) 
\end{equation}
and 
\begin{equation}\label{second switch} 
H \Biggl( \! \exp \Biggl(\sum_{j \in \Z} \biggl( \Psi_j L_j(x,\varphi)
+ \Psi_{j  - \frac{1}{2}} G_{j - \frac{1}{2}}(x,\varphi) \biggr) \! \Biggr) 
\! \cdot (x, \varphi) \Biggr) 
\end{equation} 
are in $\mathbb{C} [\alpha_0^{1/2}, \alpha_0^{-1/2}] [x, x^{-1},
\varphi] [[ \mathcal{A}, \mathcal{B}]] [\mathcal{M}, \mathcal{N}]$ .  Then by
Proposition \ref{Switch}, (\ref{first switch}) is equal to  
\[\exp \Biggl(\sum_{j \in \Z} \biggl( \Psi_{-j} L_{-j}(x,\varphi) 
+ \Psi_{-j + \frac{1}{2}} G_{-j + \frac{1}{2}}(x,\varphi) 
\Biggr) \! \Biggr) \! \cdot H (x,\varphi) ,\] 
and by Proposition \ref{Switch2},  (\ref{second switch}) is equal to   
\[\exp \Biggl( \sum_{j \in \Z} \biggl( \Psi_j L_j(x,\varphi)
 +  \Psi_{j  - \frac{1}{2}} G_{j - \frac{1}{2}}(x,\varphi) \biggr) \!
\Biggr) \! \cdot H (x, \varphi) .\]   
Thus the above expressions are also in 
\[x \mathbb{C} [\alpha_0^\frac{1}{2}, \alpha_0^{-\frac{1}{2}}]
[x^{-1},\varphi] [[ \mathcal{A}, \mathcal{B}]] [\mathcal{M}, \mathcal{N}] \quad
\mbox{and} \quad x\mathbb{C} [\alpha_0^\frac{1}{2},
\alpha_0^{-\frac{1}{2}}] [x,\varphi] [[ \mathcal{A}, \mathcal{B}]] [\mathcal{M}, \mathcal{N}] , \] 
respectively.  Since $H(x, \varphi)$ was arbitrary, (\ref{endomorphism
1}) and (\ref{endomorphism 2}) are in  
\[(\mbox{End} \;  \mathbb{C} [x, x^{-1}, \varphi]) [\alpha_0^{\frac{1}{2}},
\alpha_0^{-\frac{1}{2}}] [[ \mathcal{A}, \mathcal{B}]] [\mathcal{M}, \mathcal{N}] .\]
Furthermore, they are obviously even.  

To prove (\ref{normal order equation}), we note that using
Propositions \ref{Switch} and \ref{Switch2} repeatedly and by
(\ref{F2}), we have
\[\bar{F}^{(2)} \circ (H_{\mathcal{B}, \mathcal{N}}^{(2)})^{-1} \circ I
\circ H_{\alpha_0^{1/2}, \mathcal{A}, \mathcal{M}}^{(1)} (x,\varphi) 
\hspace{2.6in}\]
\[= \exp \Biggl( \! - \! \sum_{j \in \Z} \biggl( \mathcal{A}_j L_j(x,\varphi) + 
\mathcal{M}_{j - \frac{1}{2}} G_{j - \frac{1}{2}} \biggr) \!
\Biggr) \! \cdot   (\alpha_0^\frac{1}{2})^{- 2L_0(x,\varphi)} \cdot \hspace{1in} \]  
\[ \exp \Biggl( \! - \! \sum_{j \in \Z} \biggl( \mathcal{B}_j L_{-j}(x,\varphi) + 
\mathcal{N}_{j - \frac{1}{2}} G_{-j + \frac{1}{2}}(x,\varphi) \biggr) \! \Biggr) \! 
\cdot (\alpha_0^\frac{1}{2})^{2L_0(x,\varphi) } \cdot \]
\[ \hspace{.5in} \exp \left(\Psi_0 2L_0(x,\varphi) \right)  \exp \Biggl( \! - \! 
\sum_{j \in \Z} \biggl( \Psi_j L_j(x,\varphi) + \Psi_{j - \frac{1}{2}} G_{j -
\frac{1}{2}}(x,\varphi) \biggr) \! \Biggr) \! \cdot (x,\varphi) \] 
in $\mathbb{C} [x, x^{-1}, \varphi] [\alpha_0^{1/2},
\alpha_0^{-1/2}] [[ \mathcal{A}, \mathcal{B}]] [\mathcal{M}, \mathcal{N}]$.  
Using Propositions \ref{Switch} and \ref{Switch2} again, we see
that for $H(x,\varphi) \in \mathbb{C} [x, x^{-1}, \varphi]
[\alpha_0^{1/2}, \alpha_0^{-1/2}] [[ \mathcal{A}, \mathcal{B}]] 
[\mathcal{M}, \mathcal{N}]$, the expressions
\begin{multline*}
\exp \Biggl( \! - \! \sum_{j \in \Z} \biggl( \mathcal{A}_j L_j(x,\varphi)
+ \mathcal{M}_{j - \frac{1}{2}} G_{j- \frac{1}{2}}(x,\varphi) \biggr) 
\! \Biggr) \! \cdot (\alpha_0^\frac{1}{2})^{ - 2L_0(x,\varphi)} \cdot \\  
\exp \Biggl(\! - \! \sum_{j \in \Z} \biggl( \mathcal{B}_j L_{-j}(x,\varphi)
 + \mathcal{N}_{j - \frac{1}{2}} G_{-j + \frac{1}{2}}(x,\varphi) \biggr) \!
\Biggr) \! \cdot (\alpha_0^\frac{1}{2})^{2L_0(x,\varphi) } \cdot\\ 
\exp \left(\Psi_0 2L_0(x,\varphi)\right) \exp 
\Biggl( \! - \! \sum_{j \in \Z} \biggl( \Psi_j L_j(x,\varphi) + \Psi_{j - \frac{1}{2}} 
G_{j - \frac{1}{2}}(x,\varphi) \biggr) \! \Biggr) \! \cdot H(x,\varphi) 
\end{multline*}
and  
\begin{multline*}
H \left( \exp \Biggl(\! - \! \sum_{j \in \Z} \biggl( \mathcal{A}_j L_j(x,\varphi)
+ \mathcal{M}_{j - \frac{1}{2}} G_{j- \frac{1}{2}}(x,\varphi) \biggr) 
\! \Biggr) \! \cdot (\alpha_0^\frac{1}{2})^{-2L_0(x,\varphi)} \cdot \right.\\
\exp \Biggl(\! - \! \sum_{j \in \Z} \biggl( \mathcal{B}_j L_{-j}(x,\varphi)
 + \mathcal{N}_{j - \frac{1}{2}} G_{-j + \frac{1}{2}}(x,\varphi) \biggr) \!
\Biggr) \! \cdot (\alpha_0^\frac{1}{2})^{2L_0(x,\varphi) } \cdot\\ 
\left. \exp \left(\Psi_0 2L_0(x,\varphi)\right) \exp 
\Biggl(\! - \!  \sum_{j \in \Z} \biggl( \Psi_j L_j(x,\varphi) + \Psi_{j - \frac{1}{2}} 
G_{j - \frac{1}{2}}(x,\varphi) \biggr) \! \Biggr) \! \cdot (x,\varphi) \right) 
\end{multline*}
exist in $\mathbb{C} [x, x^{-1}, \varphi] [\alpha_0^{1/2},
\alpha_0^{-1/2}] [[ \mathcal{A}, \mathcal{B}]] [\mathcal{M}, \mathcal{N}]$
and are equal.  Thus by (\ref{F1}), (\ref{F2}) and (\ref{F bar
equation}), we have
\begin{equation}\label{to prove normal ordering}
\exp \Biggl(\sum_{j \in \Z} \biggl( \Psi_{-j} L_{-j}(x,\varphi) 
+ \Psi_{-j + \frac{1}{2}} G_{-j + \frac{1}{2}}(x,\varphi) \biggr) \!
\Biggr) \! \cdot H(x,\varphi) \hspace{.8in} 
\end{equation}
\begin{multline*}
= \; \exp \Biggl( \! - \! \sum_{j \in \Z} \biggl( \mathcal{A}_j L_j(x,\varphi)
+ \mathcal{M}_{j - \frac{1}{2}} G_{j - \frac{1}{2}}(x,\varphi) \biggr) \! 
\Biggr) \! \cdot (\alpha_0^\frac{1}{2})^{-2L_0(x,\varphi) } \cdot \\ 
\exp \Biggl(\! - \! \sum_{j \in \Z} \biggl( \mathcal{B}_j L{-j}(x,\varphi)
 + \mathcal{N}_{j - \frac{1}{2}} G_{-j + \frac{1}{2}}(x,\varphi) \biggr)
\! \Biggr) \! \cdot (\alpha_0^\frac{1}{2})^{2L_0(x,\varphi)} \cdot \\
\exp \left(\Psi_0 2L_0(x,\varphi) \right) \cdot \exp 
\Biggl( \! - \! \sum_{j \in \Z} \biggl( \Psi_j L_j(x,\varphi)  + \Psi_{j - 
\frac{1}{2}} G_{j - \frac{1}{2}}(x,\varphi) \biggr) \! \Biggr) \!
\cdot H(x,\varphi) .
\end{multline*} 
Taking $H(x,\varphi)$ to be
\begin{multline*}
\exp \Biggl(\sum_{j \in \Z} \biggl( \Psi_j L_j(x,\varphi) 
+ \Psi_{j - \frac{1}{2}} G_{j - \frac{1}{2}}(x,\varphi) \biggr) \!
\Biggr) \! \cdot  (\alpha_0^\frac{1}{2})^{-2L_0(x,\varphi) } \cdot \\
\exp \left(\Psi_0 2L_0(x,\varphi) \right) \cdot H_1(x,\varphi)
\end{multline*} 
where $H_1$ is any element of $\mathbb{C} [x, x^{-1}, \varphi]$, equation (\ref{to
prove normal ordering}) becomes
\begin{multline*}
\exp \Biggl( \! - \! \sum_{j \in \Z} \biggl( \mathcal{A}_j L_j(x,\varphi) +
\mathcal{M}_{j - \frac{1}{2}} G_{j - \frac{1}{2}}(x,\varphi) \biggr) \! \Biggr) \! 
\cdot (\alpha_0^\frac{1}{2})^{-2L_0(x,\varphi)} \cdot \\
\exp \Biggl(\! - \! \sum_{j \in \Z} \biggl( \mathcal{B}_j L_{-j}(x,\varphi) + 
\mathcal{N}_{j - \frac{1}{2}} G_{-j + \frac{1}{2}}(x,\varphi)
\biggr) \! \Biggr) \! \cdot H_1(x,\varphi) 
\end{multline*} 
\begin{multline*}
= \; \exp \Biggl( \sum_{j \in \Z} \biggl( \Psi_{-j} L_{-j}(x,\varphi)
 + \Psi_{-j + \frac{1}{2}} G_{-j + \frac{1}{2}}(x,\varphi) \biggr) \! \Biggr) \!
\cdot\\   
\exp \Biggl( \sum_{j \in \Z} \biggl( \Psi_j L_j(x,\varphi) + 
\Psi_{j - \frac{1}{2}} G_{j - \frac{1}{2}}(x,\varphi) \biggr) \! \Biggr) \!
\cdot (\alpha_0^\frac{1}{2})^{ -2L_0(x,\varphi) }
\cdot \\
\exp \left(\Psi_0 2L_0(x,\varphi) \right) \cdot
H_1(x,\varphi)  . 
\end{multline*}  
Since $H_1$ is an arbitrary element of $\mathbb{C} [x, x^{-1}, \varphi]$,
we obtain (\ref{normal order equation}).
\end{proof}

\begin{rema} The proposition above can be understood
algebraically as in some sense establishing the relation between a 
``non-normally ordered" product (the left-hand side of 
(\ref{normal order equation})) and a ``normally ordered" product 
(the right-hand side of (\ref{normal order equation})), where by
``normally ordered" we mean ordered so that one first acts by 
``lowering operators" and then by ``raising operators" (cf. \cite{FLM}).  
Thus if one thinks of the superderivations in (\ref{normal order equation})
as acting on an element in $\mathbb{C}[x^{-1}, x,\varphi]$, then the 
operators $L_j(x,\varphi)$ and $G_{j - 1/2}(x,\varphi)$
raise the degree of this element in terms of negative powers of $x$
for $j < 0$ and lower the degree for $j>0$.  The significance of such 
normal ordering is due to the fact that superconformal field theory
is mainly interested in ``positive energy representations" of the $N=1$ 
Neveu-Schwarz algebra (see Section 3.4 and 3.7), e.g., 
$N=1$ vertex operator superalgebras (cf. \cite{B vosas}).  For such 
representations, acting by an expression such as the left-hand side of 
(\ref{normal order equation}) is not well defined, but acting by an 
expression such as the right-hand side of (\ref{normal order equation}) is
well defined.  
\end{rema}

The identities proved in the two propositions below can be used to
determine explicitly the resulting canonical supersphere {}from the 
sewing together of two canonical superspheres in certain cases 
\cite{B thesis}.  In addition, in Section 3.6 and 3.7, we show that 
these identities give certain analogous identities for any  
representation of the $N=1$ Neveu-Schwarz algebra and have certain 
nice properties for positive-energy representations.

Let \footnote{There is a misprint in \cite{H book} in the analogous 
definition of the body of $\tilde{x}$ in the nonsuper case.  On p.58
in \cite{H book}, in defining $y_{\mathcal{A}, \alpha_0}$, the function
$f_{\mathcal{A}, \alpha_0}$ in the expression is $f^{(1)}_{\mathcal{A}, 
\alpha_0}$, i.e., $y_{\mathcal{A}, \alpha_0} = (f^{(1)}_{\mathcal{A}, 
\alpha_0})^{-1}(y)$.  }
\[(\tilde{x}, \tilde{\varphi}) = (H^{(1)}_{\alpha_0^{1/2},\mathcal{A}, 
\mathcal{M}})^{-1} (x,\varphi) \in (\alpha_0^{-1} x,\alpha_0^{-\frac{1}{2}}
\varphi) + x\mathbb{C}[x,\varphi][\alpha_0^{-\frac{1}{2}}][[\mathcal{A}]]
[\mathcal{M}] . \]
Let $w$ be another even formal variable and $\rho$ another odd formal
variable, and define 
\[s_{(x,\varphi)} (w,\rho) = (w -x - \rho \varphi, \rho - \varphi) .
\] 
Then $s_{(x,\varphi)} \circ  H^{(1)}_{\alpha_0^{1/2},\mathcal{A}, \mathcal{M}}
\circ s_{(\tilde{x}, \tilde{\varphi})}^{-1} (\alpha_0^{-1} w, 
\alpha_0^{-1/2}\rho)$ is in 
\[w\mathbb{C} [x, \varphi] [\alpha_0^{\frac{1}{2}}, \alpha_0^{-\frac{1}{2}}] 
[[\mathcal{A}]] [\mathcal{M}][[w]] \oplus \rho \mathbb{C}[x, \varphi] 
[\alpha_0^{\frac{1}{2}}, \alpha_0^{-\frac{1}{2}}][[\mathcal{A}]][\mathcal{M}][[w]] 
,\] 
is superconformal in $(w,\rho)$, (i.e., letting $D = \frac{\partial}{\partial
\rho} + \rho \frac{\partial}{\partial w}$, then $D\tilde{w} = \tilde{\rho}D 
\tilde{\rho}$ for $(\tilde{w},\tilde{\rho}) = s_{(x,\varphi)} 
\circ H^{(1)}_{\alpha_0^{1/2},\mathcal{A}, \mathcal{M}} \circ s_{(\tilde{x}, 
\tilde{\varphi})}^{-1} (\alpha_0^{-1} w, \alpha_0^{-1/2}\rho)$), and 
the even coefficient of the monomial $\rho$ is an element in 
\[1 + x\mathbb{C}[x][\alpha_0^{\frac{1}{2}}, \alpha_0^{-\frac{1}{2}}]
[[\mathcal{A}]][\mathcal{M}] \oplus \varphi \mathbb{C} [x][\alpha_0^{\frac{1}{2}},
\alpha_0^{-\frac{1}{2}}] [[\mathcal{A}]][\mathcal{M}] . \]

Let
\[ \Theta^{(1)}_j = \Theta^{(1)}_j(\alpha_0^{\frac{1}{2}},\mathcal{A}, 
\mathcal{M}, (x, \varphi)) \in \mathbb{C}[x, \varphi]
[\alpha_0^{\frac{1}{2}}, \alpha_0^{-\frac{1}{2}}][[\mathcal{A}]][\mathcal{M}] , \]
for $j \in \frac{1}{2} \mathbb{N}$, be defined by
\begin{multline}\label{define first Theta}
\Bigl(\exp(\Theta^{(1)}_0(\alpha_0^{\frac{1}{2}}, \mathcal{A}, \mathcal{M}, 
(x, \varphi)), \Bigl\{\Theta^{(1)}_j(\alpha_0^{\frac{1}{2}}, \mathcal{A}, 
\mathcal{M}, (x, \varphi)), \Bigr. \Bigr. \\
\Bigl. \Bigl. \Theta^{(1)}_{j - \frac{1}{2}} (\alpha_0^{\frac{1}{2}}, 
\mathcal{A}, \mathcal{M}, (x, \varphi)) \Bigr\}_{j \in \Z} \Bigr)
\end{multline}  
\[ = \;\hat{E}^{-1}(s_{(x,\varphi)} \circ H^{(1)}_{\alpha_0^{1/2},\mathcal{A}, 
\mathcal{M}} \circ s_{(\tilde{x},\tilde{\varphi})}^{-1}
(\alpha_0^{-1}w,\alpha_0^{-\frac{1}{2}}\rho)) .  \hspace{1.5in} \]

In other words, the $\Theta_j^{(1)}$'s are determined uniquely by
\begin{eqnarray*}
& & \hspace{-.3in} s_{(x,\varphi)} \circ H^{(1)}_{\alpha_0^{1/2},\mathcal{A}, \mathcal{M}}
\circ s_{(\tilde{x}, \tilde{\varphi})}^{-1} (\alpha_0^{-1} w, 
\alpha_0^{-\frac{1}{2}}\rho) \\
&=& \! \! \exp \Biggl( \! - \! \sum_{j \in  \Z} \biggl( \Theta^{(1)}_j L_j(w,\rho) +
\Theta^{(1)}_{j - \frac{1}{2}} G_{j - \frac{1}{2}}(w,\rho) \biggr) \! \Biggr) 
\! \cdot  \\
& & \hspace{2.7in} \exp \left(- \Theta^{(1)}_0 2L_0(w,\rho) \right) \cdot (w,\rho)  \\  
&=& \! \! \exp \Biggl( \sum_{j \in  \Z} \biggl( \Theta^{(1)}_j \! \left( \Lw \right) +
\Theta^{(1)}_{j - \frac{1}{2}} \Gw \biggr) \! \Biggr) \! \cdot \\
& & \hspace{2.3in} \exp \left(\! \Theta^{(1)}_0 \left( \twoLow \right) \! \right) \! \cdot
(w,\rho) . 
\end{eqnarray*}  

This formal power series in $(w,\rho)$ gives the formal local 
superconformal coordinate at a puncture of the canonical supersphere 
obtained {}from the sewing together of two particular canonical 
superspheres with punctures.  Specifically, this is the formal power
series giving the resulting local coordinate at the $n$-th puncture, 
given formally by $(H_{\alpha_0^{1/2}, \mathcal{A}, \mathcal{M}}^{(1)}
)^{-1}  (x,\varphi)$, of the supersphere $S_1 \; _n\infty_0 \;S_2$ with 
$1 + (n+1)$ tubes obtained by sewing the supersphere $S_2$ with $1+2$ 
punctures given by
\[S_2 = ((z,\theta); I(w,\rho), s_{(z,\theta)}(w,\rho), (w,\rho))\]
to the $n$-th puncture of a supersphere $S_1$ with $1+n$ tubes where 
the local coordinate vanishing at the $n$-th puncture of $S_1$ is 
given formally by $H_{\alpha_0^{1/2}, \mathcal{A}, \mathcal{M}}^{(1)} 
(w,\rho)$ and where $(x,\varphi) = (z,\theta)$.

\begin{prop} \label{first Theta prop}  In the 
$\mathbb{C}[x,\varphi][\alpha_0^{1/2}, \alpha_0^{-1/2}]
[[\mathcal{A}]][\mathcal{M}]$-envelope of the  superalgebra $\mbox{\em End} \; 
\mathbb{C}[w,w^{-1},\rho]$, i.e., in the algebra
\begin{multline*}
(( \mbox{\em End} \; \mathbb{C}[w,w^{-1},\rho])[x,\varphi][\alpha_0^{\frac{1}{2}}, 
\alpha_0^{-\frac{1}{2}}] [[\mathcal{A}]][\mathcal{M}])^0 \\
\subset ( \mbox{\em End} \; 
\mathbb{C} [w, w^{-1}, \rho] [x,\varphi] [\alpha_0^{\frac{1}{2}}, 
\alpha_0^{-\frac{1}{2}}][[\mathcal{A}]][\mathcal{M}])^0 ,
\end{multline*}
the following identity holds. \footnote{There is a misprint in the analogous
nonsuper version of equation (\ref{first Theta equation}) in \cite{H book}. 
Equation (2.2.31) in \cite{H book} should have no $\alpha_0$ in front of
the exponential expression on the left-hand side.  As usual, one can always
obtain the nonsuper case {}from the super case by setting all odd formal
variables equal to zero.} 
\begin{multline}\label{first Theta equation}
\exp \Biggl(\! - \! \! \sum_{m = -1}^{\infty} \sum_{j \in \Z}
\binom{j+1}{m+1} \alpha_0^{-j} x^{j - m}  \\
\biggl(\!  \Bigl(
\mathcal{A}_j + 2\left(\frac{j-m}{j+1} \right)
\alpha_0^{\frac{1}{2}} x^{-1} \varphi \mathcal{M}_{j - \frac{1}{2}}
\Bigr)  \cdot  L_m(w,\rho) \\
+ \; x^{-1} \Bigl( \left(\frac{j-m}{j+1} \right) \alpha_0^{\frac{1}{2}}
\mathcal{M}_{j - \frac{1}{2}} + \varphi \frac{(j-m)}{2} \mathcal{A}_j
\Bigr)  G_{m + \frac{1}{2}}(w,\rho) \biggl) \Biggr) 
\end{multline}
\begin{multline*}
= \exp \left((\tilde{x} - \alpha_0^{-1}x)  L_{-1}(w,\rho)  +
(\tilde{\varphi} - \alpha_0^{-\frac{1}{2}}\varphi) G_{-\frac{1}{2}}(w,\rho) \right) 
\cdot \\
\exp \Bigl( \! - \! \sum_{j \in \Z} \biggl( \Theta^{(1)}_j L_j(w,\rho) +
\Theta^{(1)}_{j - \frac{1}{2}} G_{j - \frac{1}{2}}(w,\rho) \biggr) \! 
\Biggr) \! \cdot \exp \left(-\Theta^{(1)}_0 2L_0(w,\rho) \right), 
\end{multline*}  
for $(\tilde{x},\tilde{\varphi}) = (H^{(1)}_{\alpha_0^{1/2}, \mathcal{A}, 
\mathcal{M}})^{-1} (x,\varphi)$.
\end{prop}

\begin{proof}  Taylor's theorem implies that for any
$H(w,\rho) \in R[[w,w^{-1}]][\rho]$
\[\exp \left( \! x \dw + \varphi \Bigl( \drho - \rho \dw \Bigr) \! \right)
\cdot H(w,\rho) = H(w + x + \rho \varphi, \rho + \varphi) ,\]
i.e.,
\[e^{-xL_{-1}(w,\rho) - \varphi G_{-\frac{1}{2}}(w,\rho)} H(w,\rho) = H
\circ s^{-1}_{(x,\varphi)}(w,\rho).\]
Thus 
\begin{eqnarray*}
\lefteqn{s_{(x,\varphi)} \circ H^{(1)}_{\alpha_0^{1/2}, \mathcal{A}, \mathcal{M}} \circ
s^{-1}_{(\alpha_0^{-1}x,\alpha_0^{-1/2}\varphi)} (\alpha_0^{-1}w,
\alpha_0^{-\frac{1}{2}}\rho) }\\ 
&=& \! \! s_{(x,\varphi)} \circ H^{(1)}_{\alpha_0^{1/2},\mathcal{A},
\mathcal{M}} \bigl(\alpha_0^{-1}(w + x + \rho \varphi), \alpha_0^{-\frac{1}{2}}(\rho + 
\varphi)\bigr) \\
&=& \! \! s_{(x,\varphi)} \circ H^{(1)}_{\alpha_0^{1/2},\mathcal{A}, \mathcal{M}} 
\bigl(\alpha_0^{-1}w + \tilde{x} + \alpha_0^{-\frac{1}{2}}\rho \tilde{\varphi} + 
(\alpha_0^{-1}x - \tilde{x}) + \alpha_0^{-\frac{1}{2}}\rho (\alpha_0^{-\frac{1}{2}}
\varphi - \tilde{\varphi}) , \bigr.\\
& & \hspace{3in} \bigl. \alpha_0^{-\frac{1}{2}}\rho + \tilde{\varphi} + 
(\alpha_0^{-\frac{1}{2}}\varphi - \tilde{\varphi})\bigr) \\  
&=& \! \! \exp \left( \! - (\tilde{x} - \alpha_0^{-1}x) \dw - (\tilde{\varphi} 
- \alpha_0^{-\frac{1}{2}}\varphi) \Bigl( \drho - \rho \dw \Bigr) \! \right) \cdot \\ 
& & \hspace{1.5in} s_{(x,\varphi)} \circ H^{(1)}_{\alpha_0^{1/2},\mathcal{A},
\mathcal{M}} \bigl(\alpha_0^{-1}w + \tilde{x} + \alpha_0^{-\frac{1}{2}}\rho \tilde{\varphi}, 
\alpha_0^{-\frac{1}{2}} \rho + \tilde{\varphi}\bigr) \\  
&=& \! \! \exp \left((\tilde{x} - \alpha_0^{-1} x) L_{-1}(w,\rho) + (\tilde{\varphi} 
- \alpha_0^{-\frac{1}{2}} \varphi) G_{-\frac{1}{2}}(w,\rho) \right) \cdot \\
& & \hspace{2.1in} s_{(x,\varphi)} \circ H^{(1)}_{\alpha_0^{1/2}, \mathcal{A},
\mathcal{M}} \circ s_{(\tilde{x},\tilde{\varphi})}^{-1} (\alpha_0^{-1}w , 
\alpha_0^{-\frac{1}{2}} \rho) \\
&=& \! \! \exp \left((\tilde{x} - \alpha_0^{-1}x) L_{-1}(w,\rho) + (\tilde{\varphi} 
- \alpha_0^{-\frac{1}{2}} \varphi) G_{-\frac{1}{2}}(w,\rho)  \right) \cdot \\
& & \hspace{.5in} \exp \Biggl( \! - \! \sum_{j \in  \Z} \biggl( \Theta^{(1)}_j L_j(w,\rho) +
\Theta^{(1)}_{j - \frac{1}{2}} G_{j - \frac{1}{2}}(w,\rho) \biggr) \! \Biggr) \! \cdot \\
& & \hspace{2.7in} \exp \left(- \Theta^{(1)}_0 2L_0(w,\rho) \right) \cdot (w,\rho) .
\end{eqnarray*}
On the other hand, since
\begin{multline*}
H^{(1)}_{\alpha_0^{1/2},\mathcal{A}, \mathcal{M}} (w,\rho) \\
= \; \exp \Biggl( \! - \! \!
\sum_{j \in \Z} \biggl( \mathcal{A}_j L_j(w,\rho) + \mathcal{M}_{j - \frac{1}{2}}
G_{j - \frac{1}{2}}(w,\rho) \biggr) \! \Biggr) \! \cdot
(\alpha_0^{\frac{1}{2}})^{-2L_0(w,\rho)}\cdot (w,\rho) ,
\end{multline*}   
by Proposition \ref{Switch2}, we have
\begin{equation*}
s_{(x,\varphi)} \circ H^{(1)}_{\alpha_0^{1/2},\mathcal{A}, \mathcal{M}} 
\circ s^{-1}_{(\alpha_0^{-1}x,\alpha_0^{-1/2}\varphi)} 
(\alpha_0^{-1}w,\alpha_0^{-\frac{1}{2}}\rho) \hspace{1.9in}  
\end{equation*}
\begin{multline*}
=  \; s_{(x,\varphi)} \Biggl( \exp \Biggl( \sum_{j \in \Z} 
\biggl( \mathcal{A}_j \Bigl( \alpha_0^{-j-1}(x + w + \rho \varphi)^{j + 1} 
\frac{\partial \qquad \qquad}{\partial \alpha_0^{-1}
(x + w + \rho \varphi)}  \\
+ \; \frac{(j + 1)}{2} \alpha_0^{-j - \frac{1}{2}}
(\rho + \varphi)(x + w + \rho \varphi)^{j} \frac{\partial \qquad}{\partial 
\alpha_0^{-\frac{1}{2}} (\rho + \varphi)} \Bigr) \\ 
 + \; \mathcal{M}_{j - \frac{1}{2}} \alpha_0^{-j}(x + w
+ \rho \varphi)^{j} \Bigl( \frac{\partial \qquad}{\partial 
\alpha_0^{-\frac{1}{2}} (\rho + \varphi)}  \\
- \; \alpha_0^{-\frac{1}{2}} (\rho + 
\varphi) \frac{\partial \qquad \qquad}{\partial \alpha_0^{-1}(x + w + \rho
\varphi)} \Bigr) \biggr) \! \Biggr) \! \cdot \\
(\alpha_0^{\frac{1}{2}})^{2\left( \alpha_0^{-1}(x + 
w + \rho \varphi)\frac{\partial \qquad}{\partial  \alpha_0^{-1} (x + w + \rho 
\varphi)}  + \frac{1}{2}  \alpha_0^{-1/2}(\rho + \varphi) 
\frac{\partial \quad}{\partial  \alpha_0^{-1/2} (\rho + \varphi)} 
\right)} \cdot\\
 \Biggl. ( \alpha_0^{-1}(x + w + \rho \varphi), 
\alpha_0^{-\frac{1}{2}}(\rho + \varphi)) \! \Biggr) 
\end{multline*}
\begin{multline*}
=  \; s_{(x,\varphi)} \Biggl( \exp \Biggl(\sum_{j \in \Z} \biggl( \mathcal{A}_j 
\Bigl( \alpha_0^{-j}(x + w + \rho \varphi)^{j + 1} \dw  \\
+ \; \frac{(j + 1)}{2} \alpha_0^{-j}(\rho + \varphi)(x +
w)^{j} \Bigl( \drho - \varphi \dw \Bigr) \Bigr) \\ 
+ \; \mathcal{M}_{j -\frac{1}{2}} \alpha_0^{-j +
\frac{1}{2}} (x + w + \rho \varphi)^{j} \Bigl( \Bigl( \drho - \varphi
\dw \Bigr) -  (\rho + \varphi) \dw \Bigr) \biggr) \! \Biggr) \! \cdot \\ 
\Biggl.(x + w + \rho \varphi, \rho + \varphi) \! \Biggr) 
\end{multline*}
\begin{multline*}
= \; \exp \Biggl( \sum_{j \in \Z} \biggl( \mathcal{A}_j 
\Bigl( \alpha_0^{-j}(x +  w + \rho \varphi)^{j + 1} \dw \\
+ \;  \frac{(j + 1)}{2} \alpha_0^{-j}(\rho + \varphi)(x
+  w)^{j} \Bigl( \drho - \varphi \dw \Bigr) \Bigr) \\   
+ \; \mathcal{M}_{j - \frac{1}{2}} \alpha_0^{-j + 
\frac{1}{2}}(x + w + \rho \varphi)^{j} \Bigl( \drho - (\rho + 2
\varphi) \dw \Bigr) \biggr) \! \Biggr) \! \cdot (w,\rho) 
\end{multline*}
\begin{multline*}
= \; \exp \Biggl( \sum_{j \in \Z} \sum_{m \in \mathbb{N}} \biggl( 
\mathcal{A}_j \alpha_0^{-j}\Bigl( \binom{j + 1}{m} x^{j-m + 1} w^m \dw\\
+ \; ( j + 1) \binom{j}{m} \rho \varphi x^{j-m} w^m \dw +
\frac{(j + 1)}{2} \binom{j}{m} (\rho + \varphi)  x^{j-m} w^m
\drho \\
 - \; \frac{(j + 1)}{2} \binom{j}{m}\rho \varphi x^{j-m} w^m \dw
\Bigr) \\ 
+ \; \mathcal{M}_{j - \frac{1}{2}} \alpha_0^{-j + \frac{1}{2}}
\Bigl( \binom{j}{m} x^{j-m} w^m \drho + j 
\binom{j-1}{m} \rho \varphi x^{j-m-1} w^m \drho \\
- \; (\rho + 2\varphi) 
\binom{j}{m} x^{j-m} w^m \dw \Bigr)
\biggr) \! \Biggr) \! \cdot (w,\rho) 
\end{multline*}
\begin{multline*}
= \; \exp \Biggl( \sum_{j \in \Z} \sum_{m = -1}^\infty \biggl( 
\binom{j + 1}{m + 1} \alpha_0^{-j}x^{j-m} (\mathcal{A}_j + 
2\left(\frac{j-m}{j+1} \right) x^{-1} \varphi \alpha_0^{\frac{1}{2}}
\mathcal{M}_{j - \frac{1}{2}} ) \cdot \\
\hspace{4in} w^{m + 1} \dw \\
+ \; \binom{j + 1}{m + 1} \alpha_0^{-j}x^{j-m} (\mathcal{A}_j +
2\left(\frac{j-m}{j+1} \right) x^{-1}\varphi
\alpha_0^{\frac{1}{2}}\mathcal{M}_{j - \frac{1}{2}} ) \frac{(m + 1)}{2}
\rho w^m \drho \\ 
+ \; \binom{j + 1}{m + 1} \alpha_0^{-j}x^{j-m-1}
\Bigl(\left(\frac{j-m}{j+1} \right)
\alpha_0^{\frac{1}{2}} \mathcal{M}_{j - \frac{1}{2}} + \varphi
\frac{(j-m)}{2} 
\mathcal{A}_j \Bigr) w^{m + 1} \drho \\
 - \; \binom{j + 1}{m + 1} \alpha_0^{-j}x^{j-m-1} \Bigl(
\left(\frac{j-m}{j+1} \right) x^{-1}\alpha_0^{\frac{1}{2}}\mathcal{M}_{j
- \frac{1}{2}} + \varphi \frac{(j-m)}{2} 
\mathcal{A}_j \Bigr) \cdot \quad \; \\
\rho w^{m + 1} \dw \biggr) \! \Biggr) \cdot (w,\rho) .
\end{multline*}
Thus
\begin{multline*}
\exp \Biggl(\! - \! \! \sum_{m = -1}^{\infty} \sum_{j \in \Z}
\binom{j+1}{m+1} \alpha_0^{-j} x^{j - m}  \\
\Biggl. \biggl(\!  \Bigl( \mathcal{A}_j + 2\left(\frac{j-m}{j+1} \right)
\alpha_0^{\frac{1}{2}} x^{-1} \varphi \mathcal{M}_{j - \frac{1}{2}}
\Bigr)  \cdot  L_m(w,\rho) \\
+ \; x^{-1} \Bigl( \left(\frac{j-m}{j+1} \right) \alpha_0^{\frac{1}{2}}
\mathcal{M}_{j - \frac{1}{2}} + \varphi \frac{(j-m)}{2} \mathcal{A}_j
\Bigr)  G_{m + \frac{1}{2}}(w,\rho) \biggl) \Biggr) \cdot (w,\rho) 
\end{multline*} 
\begin{multline*}
= \; \exp \left( (\tilde{x} - \alpha_0^{-1}x) L_{-1}(w,\rho) + (\tilde{\varphi}- 
\alpha_0^{-\frac{1}{2}} \varphi) G_{-\frac{1}{2}}(w,\rho)  \right) \cdot \\
\exp \Biggl(\! - \! \sum_{j \in \Z} \biggl( \Theta^{(1)}_j L_j(w,\rho) +
\Theta^{(1)}_{j - \frac{1}{2}} G_{j - \frac{1}{2}}(w,\rho) \biggr) \! 
\Biggr) \! \cdot \exp \left( -\Theta^{(1)}_0 2L_0(w,\rho) \right)
\cdot (w,\rho) . 
\end{multline*} 
By Proposition \ref{Switch}, this equality implies (\ref{first Theta
equation}).
\end{proof}

Let
\[(\tilde{x}, \tilde{\varphi}) = (H^{(2)}_{\mathcal{B}, \mathcal{N}})^{-1}
\circ I (x,\varphi) \in (x,\varphi) + \mathbb{C}[x^{-1}, \varphi]
[[\mathcal{B}]][\mathcal{N}] . \]
Let $w$ be another even formal variable and $\rho$ another odd formal
variable.  Now write $s_{(x,\varphi)}(w,\rho) = (-x + w - \rho 
\varphi, \rho - \varphi)$.  We will use the convention that we should expand
$(-x + w -\rho \varphi)^j = (-x + w)^j - j\rho \varphi (-x + w)^{j-1}$ 
in positive powers of the second even variable $w$, for $j \in \mathbb{Z}$
(cf. \cite{FLM}, \cite{FHL}). Then 
\[s_{(x,\varphi)} \circ I^{-1} \circ H^{(2)}_{\mathcal{B},\mathcal{N}}
\circ s_{(\tilde{x}, \tilde{\varphi})}^{-1} (w,\rho)  \in  
w\mathbb{C} [x^{-1} \!, \varphi][[\mathcal{B}]][\mathcal{N}][[w]] \oplus 
\rho \mathbb{C} [x^{-1} \!, \varphi][[\mathcal{B}]][\mathcal{N}][[w]] ,\]
is superconformal in $(w,\rho)$, (i.e., letting $D = \frac{\partial}{\partial
\rho} + \rho \frac{\partial}{\partial w}$, then $D\tilde{w} = \tilde{\rho}D 
\tilde{\rho}$ for $(\tilde{w},\tilde{\rho}) = s_{(x,\varphi)} \circ I^{-1} 
\circ H^{(2)}_{\mathcal{B}, \mathcal{N}} \circ s_{(\tilde{x}, \tilde{\varphi})}^{-1} 
(w,\rho)$), and the even coefficient of the monomial $\rho$ is an element in $1 + x^{-1}
\mathbb{C} [x^{-1}, \varphi][[\mathcal{B}]][\mathcal{N}]$.  

Let
\[ \Theta^{(2)}_j = \Theta^{(2)}_j(\mathcal{B}, \mathcal{N}, (x, \varphi)) \in
\mathbb{C}[x^{-1}, \varphi][[\mathcal{B}]][\mathcal{N}] , \]
for $j \in \frac{1}{2} \mathbb{N}$, be defined by
\begin{equation}\label{define second Theta}
\Bigl(\exp(\Theta^{(2)}_0(\mathcal{B}, \mathcal{N}, (x, \varphi)), \Bigl\{
\Theta^{(2)}_j(\mathcal{B},\mathcal{N}, (x, \varphi)),
\Theta^{(2)}_{j - \frac{1}{2}} (\mathcal{B},
\mathcal{N}, (x, \varphi)) \Bigr\}_{j \in \Z} \Bigr) 
\end{equation}
\[ = \; \hat{E}^{-1}(s_{(x,\varphi)} \circ I^{-1} \circ
H^{(2)}_{\mathcal{B}, \mathcal{N}} \circ s_{(\tilde{x},
\tilde{\varphi})}^{-1} (w,\rho)) . \hspace{1.6in}  \] 
In other words, the $\Theta_j^{(2)}$'s are determined uniquely by
\begin{eqnarray*}
& & \hspace{-.4in} s_{(x,\varphi)} \circ I^{-1} \circ H^{(2)}_{\mathcal{B}, 
\mathcal{N}} \circ s_{(\tilde{x}, \tilde{\varphi})}^{-1} (w,\rho) \\
&=& \! \! \exp \Biggl( \! - \! \sum_{j \in  \Z} \biggl( \Theta^{(2)}_j L_j(w,\rho) +
\Theta^{(2)}_{j - \frac{1}{2}} G_{j - \frac{1}{2}}(w,\rho) \biggr) \! \Biggr) 
\! \cdot  \\
& & \hspace{2.7in} \exp \left(- \Theta^{(2)}_0 2L_0(w,\rho) \right) \cdot (w,\rho)  \\  
&=& \! \! \exp \Biggl( \sum_{j \in  \Z} \biggl( \Theta^{(2)}_j \! \left( \Lw \right) +
\Theta^{(2)}_{j - \frac{1}{2}} \Gw \biggr) \! \Biggr) \! \cdot \\
& & \hspace{2.4in} \exp \left(\! \Theta^{(2)}_0 \left( \twoLow \right) \! \right) \! 
\cdot (w,\rho) . 
\end{eqnarray*}

This formal power series in $(w,\rho)$ gives the formal local 
superconformal coordinate at a puncture of the canonical supersphere 
obtained {}from the sewing together of two particular canonical 
superspheres with punctures.  Specifically, this is the formal power
series giving the resulting local coordinate at the $1$-st
puncture, given formally by $H_{\mathcal{B}, \mathcal{N}}^{(2)} \circ I 
(x,\varphi)$, of the supersphere $S_1 \; _2\infty_0 \; S_2$ with $1 + 
2$ tubes obtained by sewing a supersphere $S_2$ with $1+1$ tubes to 
the $2$-nd puncture of the supersphere $S_1$ with $1+2$ punctures given
by
\[S_1 = ((z,\theta); I(w,\rho), s_{(z,\theta)}(w,\rho), (w,\rho))\]
where the local coordinate vanishing at the puncture at $\infty$ of 
$S_2$ is given formally by $H_{\mathcal{B}, \mathcal{N}}^{(2)} 
(w,\rho)$ and where $(x,\varphi) = (z,\theta)$.

\begin{prop} \label{second Theta prop}  
In the $\mathbb{C}[x^{-1},\varphi][[\mathcal{B}]][\mathcal{N}]$-envelope of  
$\mbox{\em End} \; \mathbb{C}[w,w^{-1},\rho]$, i.e., in the algebra 
\[ (( \mbox{\em End} \; \mathbb{C}[w,w^{-1},\rho])[x^{-1},\varphi] [[\mathcal{B}]]
[\mathcal{N}])^0 \subset ( \mbox{\em End} \; \mathbb{C} [w, w^{-1},
\rho] [x^{-1},\varphi] [[\mathcal{B}]][\mathcal{N}])^0 ,\] 
the following identity holds.
\begin{multline}\label{second Theta equation}
\exp \Biggl( \sum_{m = -1}^{\infty} \sum_{j \in \Z}
\binom{-j + 1}{m + 1} x^{-j - m} \biggl( \Bigl( \mathcal{B}_j + 2\varphi 
\mathcal{N}_{j - \frac{1}{2}} \Bigr) L_m(w,\rho) \\
+ \; \Bigl( \mathcal{N}_{j - \frac{1}{2}} + \varphi x^{-1}
\frac{(-j-m)}{2} \mathcal{B}_j \Bigr) G_{m + \frac{1}{2}}(w,\rho)
\biggr) \Biggr) 
\end{multline}  
\begin{multline*}
= \exp \left((\tilde{x} - x) L_{-1}(w,\rho)  +
(\tilde{\varphi} - \varphi) G_{-\frac{1}{2}}(w,\rho) \right) 
\cdot \\
\exp \Bigl( \! - \! \sum_{j \in \Z} \biggl( \Theta^{(2)}_j L_j(w,\rho) +
\Theta^{(2)}_{j - \frac{1}{2}} G_{j - \frac{1}{2}}(w,\rho) \biggr) \! 
\Biggr) \! \cdot \exp \left(-\Theta^{(2)}_0 2L_0(w,\rho) \right), 
\end{multline*} 
for $(\tilde{x},\tilde{\varphi}) = (H_{\mathcal{B},\mathcal{N}}^{(2)})^{-1}
\circ I (x,\varphi)$. 
\end{prop}

\begin{proof}  By Taylor's theorem 
\begin{eqnarray*}
\lefteqn{s_{(x,\varphi)} \circ I^{-1} \circ H^{(2)}_{\mathcal{B}, \mathcal{N}} 
\circ s^{-1}_{(x,\varphi)} (w,\rho)}\\
&=& \! \! s_{(x,\varphi)} \circ I^{-1} \circ H^{(2)}_{\mathcal{B}, \mathcal{N}} (w
+ x + \rho \varphi, \rho + \varphi) \\
&=& \! \! s_{(x,\varphi)} \circ I^{-1} \circ H^{(2)}_{\mathcal{B}, \mathcal{N}} (w
+ \tilde{x} + \rho \tilde{\varphi} + (x - \tilde{x}) + \rho (\varphi -
\tilde{\varphi}) , \rho + \tilde{\varphi} + (\varphi -
\tilde{\varphi})) \\  
&=& \! \! \exp \left( \! - (\tilde{x} - x) \dw - (\tilde{\varphi} - \varphi)
\Bigl( \drho - \rho \dw \Bigr) \! \right) \cdot \\ 
& & \hspace{2in} s_{(x,\varphi)} \circ I^{-1} \circ H^{(2)}_{\mathcal{B},
\mathcal{N}} (w + \tilde{x} + \rho \tilde{\varphi}, \rho +
\tilde{\varphi}) \\  
&=& \! \! \exp \left((\tilde{x} - x) L_{-1}(w,\rho) + (\tilde{\varphi} - \varphi)
G_{-\frac{1}{2}}(w,\rho) \right) \cdot \\
& & \hspace{2.3in} s_{(x,\varphi)} \circ I^{-1} \circ H^{(2)}_{\mathcal{B},
\mathcal{N}} \circ s_{(\tilde{x},\tilde{\varphi})}^{-1} (w , \rho) \\
&=& \! \! \exp \left((\tilde{x} - x) L_{-1}(w,\rho) + (\tilde{\varphi} - \varphi)
G_{-\frac{1}{2}}(w,\rho)  \right) \cdot \\
& & \hspace{.5in} \exp \Biggl( \! - \! \sum_{j \in  \Z} \biggl( \Theta^{(2)}_j L_j(w,\rho) +
\Theta^{(2)}_{j - \frac{1}{2}} G_{j - \frac{1}{2}}(w,\rho) \biggr) \! \Biggr) \! \cdot \\
& & \hspace{2.5in} \exp \left(- \Theta^{(2)}_0 2L_0(w,\rho) \right) \cdot (w,\rho) .
\end{eqnarray*}
On the other hand, since
\[I^{-1} \circ H^{(2)}_{\mathcal{B}, \mathcal{N}} (w,\rho) = \exp \Biggl(
\sum_{j \in \Z} \biggl( \mathcal{B}_j L_{-j}(w,\rho) + \mathcal{N}_{j - \frac{1}{2}}
G_{-j + \frac{1}{2}}(w,\rho) \biggr) \! \Biggr) \! \cdot (w,\rho) ,\]  
by Proposition \ref{Switch2}, we have
\begin{equation*}
s_{(x,\varphi)} \circ I^{-1} \circ H^{(2)}_{\mathcal{B}, 
\mathcal{N}} \circ s^{-1}_{(x,\varphi)} (w,\rho) \hspace{2.9in}
\end{equation*}
\begin{multline*}
= \; s_{(x,\varphi)} \Biggl( \exp \Biggl( \! - \! \sum_{j \in \Z} 
\biggl( \mathcal{B}_j \Bigl( (x + w + \rho \varphi)^{-j + 1} \frac{\partial}{\partial
(x + w + \rho \varphi)}  \Bigr. \biggr. \Biggr.\Biggr. \\
\Bigl. + \; \frac{(-j + 1)}{2} (\rho + \varphi)(x + w
+ \rho \varphi)^{-j} \frac{\partial}{\partial (\rho + \varphi)}
\Bigr) \\ 
\Biggl. \biggl. + \; \mathcal{N}_{j - \frac{1}{2}} (x + w
+ \rho \varphi)^{-j + 1} \Bigl( \frac{\partial}{\partial (\rho +
\varphi)} - (\rho + \varphi) \frac{\partial}{\partial (x + w + \rho
\varphi)} \Bigr) \biggr) \! \Biggr) \! \cdot \\
\Biggl. (x + w + \rho \varphi, \rho + \varphi) \!
\Biggr) 
\end{multline*}
\begin{multline*}
= \; s_{(x,\varphi)} \Biggl( \exp \Biggl(\!  - \! \sum_{j \in \Z} \biggl( 
\mathcal{B}_j \Bigl( (x + w + \rho \varphi)^{-j + 1} \dw 
\Bigr. \biggr. \Biggr. \Biggr. \\
\Bigl. + \; \frac{(-j + 1)}{2} (\rho + \varphi)(x +
w)^{-j} \Bigl( \drho - \varphi \dw \Bigr) \Bigr) \\ 
\Biggl. \biggl. + \; \mathcal{N}_{j - \frac{1}{2}} (x +
w + \rho \varphi)^{-j + 1} \Bigl( \Bigl( \drho - \varphi \dw \Bigr) - 
(\rho + \varphi) \dw \Bigr) \biggr) \! \Biggr) \! \cdot \\ 
\Biggl. (x + w + \rho \varphi, \rho + \varphi) \!
\Biggr) 
\end{multline*}
\begin{multline*}
= \; \exp \Biggl( \! - \! \sum_{j \in \Z} \biggl( \mathcal{B}_j \Bigl( (x + w +
\rho \varphi)^{-j + 1} \dw \Bigr. \biggr. \Biggr. \\
+ \; \Bigl.  \frac{(-j + 1)}{2} (\rho + \varphi)(x
+  w)^{-j} \Bigl( \drho - \varphi \dw \Bigr) \Bigr) \\   
\Biggl. \biggl.  + \; \mathcal{N}_{j - \frac{1}{2}} (x +
w + \rho \varphi)^{-j + 1} \Bigl( \drho - (\rho + 2 \varphi) \dw
\Bigr) \biggr) \! \Biggr) \! \cdot (w,\rho) 
\end{multline*}
\begin{multline*}
= \; \exp \Biggl( \! - \! \sum_{j \in \Z} \sum_{m \in \mathbb{N}} \biggl( 
\mathcal{B}_j \Bigl(\binom{-j + 1}{m} x^{-j-m + 1} w^m \dw\Bigr. \biggr. \Biggr.  \\
+ \; ( -j + 1)\binom{-j}{m}\rho \varphi x^{-j-m} w^m \dw +
\frac{(-j + 1)}{2} \binom{-j}{m}(\rho + \varphi)  x^{-j-m} w^m\drho \\
\left. - \; \frac{(-j + 1)}{2} \binom{-j}{m}\rho \varphi x^{-j-m} w^m \dw\right) \\
+ \; \mathcal{N}_{j - \frac{1}{2}} \Bigl( \binom{-j + 1}{m} x^{-j-m + 1} w^m 
\drho + (-j + 1) \binom{-j}{m} \rho \varphi x^{-j-m} w^m \drho \Bigr. \\
\Biggl. \biggl. \Bigl. - \; (\rho + 2\varphi) \binom{-j + 1}{m}
x^{-j-m + 1} w^m \dw \Bigr) \biggr) \! \Biggr) \! \cdot (w,\rho) 
\end{multline*}
\begin{multline*}
= \; \exp \Biggl( \! - \! \sum_{j \in \Z} \sum_{m = -1}^\infty \biggl( 
\binom{-j + 1}{m + 1} x^{-j-m} (\mathcal{B}_j + 2\varphi
\mathcal{N}_{j - \frac{1}{2}} ) w^{m + 1} \dw \biggr. \Biggr. \\
+ \; \binom{-j + 1}{m + 1} x^{-j-m} (\mathcal{B}_j + 2\varphi
\mathcal{N}_{j - \frac{1}{2}} ) \frac{(m + 1)}{2} \rho w^m \drho \\
+ \; \binom{-j + 1}{m + 1} x^{-j-m} \Bigl(\mathcal{N}_{j -
\frac{1}{2}} + \varphi x^{-1} \frac{(-j-m)}{2} \mathcal{B}_j \Bigr) w^{m
+ 1} \drho \\
\Biggl. \biggl. - \; \binom{-j + 1}{m + 1} x^{-j-m} \Bigl(\mathcal{N}_{j -
\frac{1}{2}} + \varphi x^{-1} \frac{(-j-m)}{2} \mathcal{B}_j \Bigr) \rho
w^{m + 1} \dw \biggr) \! \Biggr) \! \cdot (w,\rho) .
\end{multline*}
Thus
\begin{multline*}
\exp \Biggl(\sum_{m = -1}^{\infty} \sum_{j \in \Z}
\binom{-j + 1}{m + 1} x^{-j - m} \biggl( \Bigl( \mathcal{B}_j + 2\varphi 
\mathcal{N}_{j - \frac{1}{2}} \Bigr) L_{m}(w,\rho)  \\
 + \; \Bigl( \mathcal{N}_{j - \frac{1}{2}} + \varphi x^{-1} 
\frac{(-j-m)}{2} \mathcal{B}_j \Bigr) G_{m + \frac{1}{2}}(w,\rho) 
\biggr) \! \Biggr) \!  \cdot (w,\rho) 
\end{multline*}
\begin{multline*} 
= \; \exp \left( (\tilde{x} - x) L_{-1}(w,\rho) + (\tilde{\varphi}- \varphi)
G_{-\frac{1}{2}}(w,\rho)  \right) \cdot \\
\exp \Biggl(\! - \! \sum_{j \in \Z} \biggl( \Theta^{(2)}_j L_j(w,\rho) +
\Theta^{(2)}_{j - \frac{1}{2}} G_{j - \frac{1}{2}}(w,\rho) \biggr) \! 
\Biggr) \! \cdot \exp \left( -\Theta^{(2)}_0 2L_0(w,\rho) \right)
\cdot (w,\rho) . 
\end{multline*}
By Proposition \ref{Switch}, this equality implies (\ref{first Theta
equation}).  
\end{proof}

\section[The $N=1$ Neveu-Schwarz algebra and superderivations]
{The $N=1$ Neveu-Schwarz algebra and a
representation in terms of superderivations} 

Let $\mathfrak{v}$ denote the Virasoro algebra with central charge $d$,
basis consisting of the central element $d$ and $L_n$, for $n \in 
\mathbb{Z}$, and commutation relations  
\begin{equation}\label{Virasoro relation}
[L_m ,L_n] = (m - n)L_{m + n} + \frac{1}{12} (m^3 - m) \delta_{m + n 
, 0} \; d ,
\end{equation}
for $m, n \in \mathbb{Z}$.  Consider the $N=1$ Neveu-Schwarz Lie 
superalgebra $\mathfrak{ns}$ which is a super-extension of 
$\mathfrak{v}$ by the odd elements $G_{n + 1/2}$, for $n \in \mathbb{Z}$, 
such that $\mathfrak{ns}$ has a basis consisting of the central element 
$d$, $L_n$, and $G_{n + 1/2}$, for $n \in \mathbb{Z}$, with supercommutation
relations 
\begin{eqnarray}
\bigl[ G_{m + \frac{1}{2}},L_n \bigr] &=& \Bigl(m - \frac{n - 1}{2} \Bigr) 
G_{m + n + \frac{1}{2}}  \label{Neveu-Schwarz relation1} \\  
\bigl[ G_{m + \frac{1}{2}} , G_{n - \frac{1}{2}} \bigr] &=& 2L_{m +
n} + \frac{1}{3} (m^2 + m) \delta_{m + n , 0} \; d \label{Neveu-Schwarz relation2}
\end{eqnarray}
in addition to (\ref{Virasoro relation}).

It is easy to check that the superderivations in $\mbox{Der}(\mathbb{C} 
[x, x^{-1},\varphi])$ given by (\ref{L notation}) and (\ref{G notation}), 
i.e., the superderivations  
\begin{eqnarray}
L_n(x,\varphi) &=& - \biggl(x^{n + 1} \frac{\partial}{\partial x} + \Bigl(
\frac{n + 1}{2} \Bigr) 
\varphi x^n \frac{\partial}{\partial \varphi} \biggr) \label{L(n)}\\ 
G_{n - \frac{1}{2}}(x,\varphi) &=& -x^n \biggl(\frac{\partial}{\partial \varphi} - 
\varphi \frac{\partial}{\partial x} \biggr)  \label{G(n - 1/2)}
\end{eqnarray}
for $n \in \mathbb{Z}$, satisfy the $N = 1$ Neveu-Schwarz relations 
(\ref{Virasoro relation}) - (\ref{Neveu-Schwarz relation2}) with central 
charge zero (cf. \cite{Bc2}).

Let $W$ be a $\mathbb{Z}_2$-graded vector space over $\mathbb{C}$ such 
that dim$W^1 = 1$ and dim$W^0 = 2$.  Recall the classical Lie superalgebra
$\mathfrak{osp}_{\mathbb{C}}(1|2)$ (cf. \cite{K}) the
orthogonal-symplectic superalgebra  
\begin{eqnarray*}
\mathfrak{osp}_{\mathbb{C}}(1|2) = \left\{ \biggl. \left(\begin{array}{ccc}   
                                         0 & p & q \\
                                         q & a & b \\
                                         -p & c & -a 
                                            \end{array} \right) \in
\mathfrak{gl}_{\mathbb{C}}(1|2) \; \biggr| \; a,b,c,d,p,q \in \mathbb{C} \right\} 
\end{eqnarray*}
which is the subalgebra of $\mathfrak{gl}_{\mathbb{C}}(1|2)$ leaving
the non-degenerate form  $\beta$ on $W$ given by 
$\left(\begin{array}{ccc}    
                1 & 0 & 0 \\
                0 & 0 & 1 \\
                0 & -1 & 0 \end{array}
\right)$ invariant, meaning $\beta (Xu,v) + (-1)^{\eta(X) \eta(u)} 
\beta(u,Xv) = 0$ for $X \in \mathfrak{gl}_{\mathbb{C}}(1|2)$, $u,v \in W$,
and $X$ and $u$ homogeneous. The subalgebra 
of $\mathfrak{ns}$ given by $\mbox{span}_{\mathbb{C}} \{L_{\pm1}, L_0, G_{\pm 1/2} \}$
is isomorphic to $\mathfrak{osp}_{\mathbb{C}}(1|2)$. The
correspondence 
\begin{eqnarray*}
\left(\begin{array}{ccc}   
        0 & 0 & 0 \\
        0 & 0 & 1 \\
        0 & 0 & 0 \end{array} \right) &\longleftrightarrow&  - 
\frac{\partial}{\partial x} ,  \\
\frac{1}{2} \left( \begin{array}{ccc}    
         0 & 0 & 0 \\
         0 & 1 & 0 \\
         0 & 0 & -1 \end{array} \right)  &\longleftrightarrow&  - \left( x
\frac{\partial}{\partial x} + \frac{1}{2} \varphi
\frac{\partial}{\partial \varphi} \right), \\ 
\left(\begin{array}{ccc}   
0 & 0 & 0 \\
0 & 0 & 0 \\
0 & -1 & 0 \end{array} \right) &\longleftrightarrow&   - \left( x^2
\frac{\partial}{\partial x} +  \varphi x \frac{\partial}{\partial
\varphi} \right), \\ 
\left(\begin{array}{ccc}   
0 & 0 & 1 \\
1 & 0 & 0 \\
0 & 0 & 0 \end{array} \right)  &\longleftrightarrow&  - \left(  
\frac{\partial}{\partial \varphi} -  \varphi
\frac{\partial}{\partial x} \right), \\
\left(\begin{array}{ccc}    
0 & 1 & 0 \\
0 & 0 & 0 \\
-1 & 0 & 0 \end{array} \right) &\longleftrightarrow&  - x\left(
\frac{\partial}{\partial \varphi} -  \varphi 
\frac{\partial}{\partial x} \right),  
\end{eqnarray*} 
defines a Lie superalgebra isomorphism between
$\mathfrak{osp}_{\mathbb{C}}(1|2)$ and the Lie superalgebra generated by
$L_{\pm1}(x,\varphi), L_{0}(x,\varphi), G_{\pm 1/2}(x,\varphi)$.  Let $y$ be an even 
formal variable and $\xi$ an odd formal variable.  Letting $X$ denote
each of the five matrices above, we observe that
\begin{eqnarray*}
e^{-yX} = \left(\begin{array}{ccc}   
                  1 & 0 & 0 \\
                  0 & 1 & -y \\
                  0 & 0 & 1 \end{array} \right) , \quad
\left(\begin{array}{ccc}    
1 & 0 & 0 \\
0 & e^{-\frac{y}{2}} & 0 \\
0 & 0 & e^{\frac{y}{2}}  \end{array} \right) , \quad
\left(\begin{array}{ccc}    
1 & 0 & 0 \\
0 & 1 & 0 \\
0 & y & 1 \end{array} \right) 
\end{eqnarray*}
\begin{eqnarray*}
e^{-\xi X} =  \left(\begin{array}{ccc}   
                      1 & 0 & -\xi \\
                     -\xi & 1 & 0 \\
                      0 & 0 & 1 \end{array} \right) , \quad 
\left(\begin{array}{ccc}   
1 & -\xi & 0 \\
0 & 1 & 0 \\
\xi & 0 & 1 \end{array}\right),
\end{eqnarray*}
respectively.  These are all elements in the connected component 
of the Lie supergroup $OSP(1|2)$ containing the identity with matrix 
elements in $\mathbb{C}[[y]][\xi]$, and hence have superdeterminant 1, 
where the superdeterminant is defined as  
\[\mbox{sdet}\left(\begin{array}{cc}  A & B \\
                                      C & D \end{array}\right) =
\mbox{det}(A - BD^{-1}C)(\mbox{det}D)^{-1} \]  
(cf. \cite{D}, \cite{CR}).  (In our case, $A$ is a one-by-one matrix,
$B$ is two-by-one, $C$ is one-by-two, and $D$ is two-by-two.)  In fact, 
for $R$ a superalgebra with $y \in R^0$ and $\xi \in R^1$, the five 
matrices above generate the connected component of $OSP_R(1|2)$ 
containing the identity. Denote this group by $G$.

$G$ acts on an even and an odd formal variable by the
superprojective transformations,
i.e,  for $g = e^{-yX}$ and $g = e^{-\xi X}$ above, we have
\begin{equation}\label{group action}
g \cdot (x,\varphi) =  (x + y, \varphi), \quad (e^yx,
e^{\frac{y}{2}} \varphi), \quad \Bigl(\frac{x}{1 - yx}, \varphi \frac{1}{1
- yx}\Bigr) 
\end{equation}
\begin{eqnarray*}
\qquad \qquad \qquad (x + \varphi \xi, \xi + \varphi), \quad (x +
\varphi \xi x, \xi 
x + \varphi ) ,
\end{eqnarray*}
respectively.  These generate the supergroup of superprojective
transformations which is the group of global superconformal automorphisms of
the super-Riemann sphere studied in Chapter 2.  

Thus, $\mathfrak{osp}_{\mathbb{C}}(1|2)$ is the superalgebra of infinitesimal
superprojective transformations.  Note that for the representative
elements $L(\pm 1), L(0), G(\pm 1/2)$, we have 
\begin{eqnarray*}
e^{-yL(-1)} \cdot (x,\varphi) &=& ( x + y, \varphi), \\
e^{-yL(0)} \cdot (x,\varphi) &=& ( e^yx, e^{\frac{y}{2}} \varphi), \\  
e^{-yL(1)} \cdot (x,\varphi) &=& \Bigl(\frac{x}{1 - yx} , \varphi
 \frac{1}{1 - yx} \Bigr), \\
e^{-\xi G(-\frac{1}{2})} \cdot (x,\varphi) &=& ( x + \varphi \xi, \xi + 
\varphi ), \\ 
e^{-\xi G(\frac{1}{2})} \cdot (x,\varphi) &=& ( x + \varphi \xi x, \xi 
x + \varphi) ,
\end{eqnarray*} 
as expected.

\section{Modules for the $N=1$ Neveu-Schwarz algebra}

Let $\mathfrak{ns}$ denote the $N = 1$ Neveu-Schwarz algebra defined above.
For any representation of $\mathfrak{ns}$, we shall use $L(m)$, $G(m -
1/2)$ and $c \in \mathbb{C}$ to denote the representation images
of $L_m$, $G_{m - 1/2}$ and $d$, respectively.  We can think of the 
identities proved in Section 3.3 as identities for the representation of 
$\mathfrak{ns}$ given by (\ref{L(n)}), (\ref{G(n - 1/2)}) and $c = 0$.  We want
to prove the corresponding identities for any representation of $\mathfrak{ns}$ 
and then we will want to prove additional properties related to these
identities for certain representations of $\mathfrak{ns}$.  We do this by first 
proving the identities in a certain extension of the universal enveloping 
algebra of $\mathfrak{ns}$.  This extension must be one such that terms such as 
$(t^{1/2})^{kL_0}$, for any formal variable $t^{1/2}$ and 
$k \in 2\mathbb{Z} \smallsetminus \{0\}$, can be well defined.

Consider the two subalgebras of $\mathfrak{ns}$
\begin{eqnarray*}
\mathfrak{ns}_+ &=& \bigoplus_{n \in \Z} \mathbb{C} L_n \oplus \bigoplus_{n
\in \Z} \mathbb{C} G_{n - \frac{1}{2}}, \\
\mathfrak{ns}_- &=& \bigoplus_{-n \in \Z} \mathbb{C} L_n \oplus \bigoplus_{-n
\in \Z} \mathbb{C} G_{n + \frac{1}{2}}. 
\end{eqnarray*}
Let $U(\mathfrak{ns}_-)$ be the universal enveloping algebra of the Lie
algebra $\mathfrak{ns}_-$.  For any $h, c \in \mathbb{C}$, the Verma module
$M(c,h)$ (cf. \cite{KW}) for $\mathfrak{ns}$ is a free
$U(\mathfrak{ns}_-)$-module generated by an element $\mathbf{1}_{c,h}$ such
that 
\begin{eqnarray*}
\mathfrak{ns}_+ \mathbf{1}_{c,h} &=& 0 ,\\
L(0) \mathbf{1}_{c,h} &=& h \mathbf{1}_{c,h} ,\\
d \mathbf{1}_{c,h} &=& c \mathbf{1}_{c,h} .
\end{eqnarray*}

Let  
\begin{equation}\label{module}
V = \coprod_{n \in \frac{1}{2} \mathbb{Z}} V_{(n)}
\end{equation}
be a module for $\mathfrak{ns}$ of central charge $c \in \mathbb{C}$ 
(i.e., $dv = cv$ for $v \in V$) such that for $v \in V_{(n)}$,  
\begin{equation}\label{graded module}
L(0)v = nv .
\end{equation}
Let $P(n)$ be the projection {}from $V$ to $V_{(n)}$.  For any formal 
variable $t^{1/2}$ and $k \in 2\mathbb{Z} \smallsetminus \{0\}$, we 
define $(t^{1/2})^{kL(0)} \in (\mathrm{End} V) [[t^{1/2}, t^{-1/2}]]$ 
by
\[(t^{\frac{1}{2}})^{kL(0)}v = (t^{\frac{1}{2}})^{kn}v\]
for $v \in V_{(n)}$, or equivalently
\[(t^{\frac{1}{2}})^{kL(0)} v= \sum_{n \in \frac{1}{2}\mathbb{Z}} P(n) 
(t^{\frac{1}{2}})^{kn}v\]
for any $v \in V$.

Let $V_P$ be a vector space over $\mathbb{C}$ with basis $\{P_n \; | \; 
n \in \frac{1}{2} \mathbb{Z}\}$.  Let $T(\mathfrak{ns} \oplus V_P)$ be the 
tensor algebra generated by the direct sum of $\mathfrak{ns}$ and $V_P$, 
and let $\mathcal{I}$ be the ideal of $T(\mathfrak{ns} \oplus V_P)$ 
generated by
\begin{multline*}
\Bigl\{L_n \otimes L_m - L_m \otimes L_n - [L_n, L_m], \;  
L_m \otimes G_{n - \frac{1}{2}} - G_{n - \frac{1}{2}} \otimes L_m -
\bigl[L_m , G_{n - \frac{1}{2}}\bigr], \\
G_{n + \frac{1}{2}} \otimes G_{m - \frac{1}{2}} + G_{m - \frac{1}{2}} 
\otimes G_{n + \frac{1}{2}}  - \bigl[G_{n + \frac{1}{2}}, G_{m - 
\frac{1}{2}}\bigr], \; L_n \otimes d - d \otimes L_n,  \\
\Bigl. G_{m - \frac{1}{2}} \otimes d - d \otimes G_{m - \frac{1}{2}}, 
\; P_i \otimes P_j - \delta_{i,j} P_i, \; P_i \otimes L_m - L_m \otimes 
P_{i + m}, \Bigr. \\
\bigl. P_i \otimes G_{m - \frac{1}{2}} - G_{m - \frac{1}{2}}
\otimes P_{i + m -\frac{1}{2}}, \; P_i \otimes d - d \otimes P_i \; \bigr| \; m,n 
\in \mathbb{Z}, \; i,j \in \mbox{$\frac{1}{2}\mathbb{Z}$} \Bigr\} .
\end{multline*}

Then $U_P(\mathfrak{ns}) = T(\mathfrak{ns} \oplus V_P) / \mathcal{I}$ is an
associative superalgebra and the universal enveloping algebra
$U(\mathfrak{ns})$ of the Neveu-Schwarz algebra is a subalgebra.  Linearly
$U_P(\mathfrak{ns}) = U(\mathfrak{ns}) \otimes V_P$.

For any even formal variable $t^{1/2}$ and $k \in 2 \mathbb{Z} 
\smallsetminus \{0\}$, we define
\[(t^\frac{1}{2})^{kL_0} = \sum_{n \in \frac{1}{2}\mathbb{Z}} P_n 
(t^\frac{1}{2})^{kn} \; \in \; U_P(\mathfrak{ns}) [[t^\frac{1}{2}, t^{-\frac{1}{2}}]] .\]
Then for $k,n \in 2\mathbb{Z} \smallsetminus \{0\}$, and $m \in
\mathbb{Z}$, 
\begin{eqnarray*}
(t^\frac{1}{2})^{kL_0} (t^\frac{1}{2})^{nL_0} &=& (t^\frac{1}{2})^{(k + n)L_0},  \\
(t^\frac{1}{2})^{kL_0} L_m &=& L_m (t^\frac{1}{2})^{-km} (t^\frac{1}{2})^{kL_0}, \\
(t^\frac{1}{2})^{kL_0} G_{m - \frac{1}{2}} &=& G_{m - \frac{1}{2}} 
(t^\frac{1}{2})^{-k(m - \frac{1}{2})} t^{kL_0},
\end{eqnarray*} 
and $(t^{1/2})^{kL_0}$ commutes with $d$.

The following proposition is proved similarly to that for the universal 
enveloping algebra of a Lie algebra.   

\begin{prop}\label{universal enveloping}
Let $V$ be a module for $\mathfrak{ns}$ of the form (\ref{module}) such
that (\ref{graded module}) holds.  Then there is a unique algebra
homomorphism {}from $U_P (\mathfrak{ns})$ to $\mbox{End} \; V$ such that
$L_j$, $G_{j - 1/2}$, $d$ and $P_n$ are mapped to $L(j)$, $G(j
- 1/2)$, $c$ and $P(n)$, respectively, for $j \in \mathbb{Z}$ and
$n \in \frac{1}{2} \mathbb{Z}$.
\end{prop}

If 
\[V_{(n)} = 0 \quad \mbox{for $n$ sufficiently small}, \]
then we say that $V$ is a {\it positive energy module} for the
Neveu-Schwarz algebra and that the corresponding representation is a
{\it positive energy representation} of the Neveu-Schwarz algebra. 

The Verma module $M(c,h)$ with central charge $c \in \mathbb{C}$
generated {}from a lowest weight vector $\mathbf{1}_{c,h}$ with weight
$h \in \mathbb{Z}$ is an example of a positive energy representation of
the Neveu-Schwarz algebra as is an $N=1$ vertex operator superalgebra
(see \cite{KW}, \cite{B vosas}).

\section[Realizations of the sewing identities]{Realizations of the
sewing identities for general representations of the $N=1$ Neveu-Schwarz 
algebra} 

As mentioned before, the identities proved in Section 3.3 can be 
thought of as identities for the representation of the Neveu-Schwarz 
algebra on $\mathbb{C} [x,x^{-1}, \varphi]$ given by (\ref{L(n)}), 
(\ref{G(n - 1/2)}) and $c = 0$.  We now want to prove the corresponding 
identities for any representation of the $N=1$ Neveu-Schwarz algebra.  
Since in general, the central charge $c$ will not be zero, we will have 
an extra term in these identities involving $c$.  We first prove the 
identities for $U_P({\mathfrak n}{\mathfrak s})$ defined in Section 3.5.  

\begin{prop}\label{central charge sewing}
Let $\mathcal{A}$ and $\mathcal{B}$ be two sequences of even formal
variables, $\mathcal{M}$ and $\mathcal{N}$ two sequences of odd formal
variables, $\alpha_0^{1/2}$ another even formal variable, and $(\Psi_j,
\Psi_{j - 1/2})$, for $j \in \mathbb{Z}$, the canonical sequence of
formal series given by Proposition \ref{normal order}.  There
exists a unique canonical formal series 
\[\Gamma = \Gamma(\alpha_0^\frac{1}{2}, \mathcal{A}, \mathcal{M}, \mathcal{B},
\mathcal{N}) \in (\mathbb{C}[\alpha_0^{\frac{1}{2}}, \alpha_0^{-\frac{1}{2}}]
[[\mathcal{A}, \mathcal{B}]] [\mathcal{M}, \mathcal{N}])^0\] 
such that
\begin{equation}\label{first terms of Gamma}
\Gamma =  \sum_{j \in \Z} \! \left( \! \biggl( \frac{j^3 - j}{12} \biggr)
\alpha_0^{-j} \mathcal{A}_j \mathcal{B}_j + \biggl( \frac{j^2 - j}{3} \biggr)
\alpha_0^{-j + \frac{1}{2}} \mathcal{N}_{j - \frac{1}{2}} \mathcal{M}_{j -
\frac{1}{2}} \right) + \Gamma_0
\end{equation}
where $\Gamma_0$ contains only terms with total degree at least three
in the $\mathcal{A}_j$'s, $\mathcal{M}_{j - 1/2}$'s, $\mathcal{B}_j$'s,
and $\mathcal{N}_{j - 1/2}$'s for $j \in \Z$ with each term
containing at least one of the $\mathcal{A}_j$'s or $\mathcal{M}_{j -
1/2}$'s and at least one of the $\mathcal{B}_j$'s or $\mathcal{N}_{j
- 1/2}$'s such that in
\footnote{There is a misprint in the analogous 
nonsuper case to this proposition given in \cite{H book}.  In Proposition 
4.2.1 in \cite{H book}, the equality (4.2.2) (which is the nonsuper
part of formula (\ref{switching identity in NS})) takes place in 
$\mathcal{R}_\Pi = U_\Pi(\mathfrak{L}) [\alpha_0, \alpha_0^{-1}]
[[\mathcal{A},\mathcal{B}]]$, not $\mathcal{R} = U (\mathfrak{L}) 
[\alpha_0,\alpha_0^{-1}][[\mathcal{A},\mathcal{B}]]$ as stated.  Huang's 
$\mathcal{R}_\Pi$ is equivalent to our $U_P(\mathfrak{v})[\alpha_0,
\alpha_0^{-1}][[\mathcal{A},\mathcal{B}]] \subset (U_P(\mathfrak{ns})
[\alpha_0^{1/2}, \alpha_0^{-1/2}][[\mathcal{A},\mathcal{B}]][\mathcal{M}, 
\mathcal{N}])^0$.}
\[(U_P(\mathfrak{ns}) [\alpha_0^{\frac{1}{2}}, \alpha_0^{-\frac{1}{2}}] [[\mathcal{A}, 
\mathcal{B}]][\mathcal{M}, \mathcal{N}])^0 , \]
i.e., in the $\mathbb{C}[\alpha_0^{1/2}, \alpha_0^{-1/2}] [[\mathcal{A}, 
\mathcal{B}]][\mathcal{M}, \mathcal{N}]$-envelope of $U_P(\mathfrak{ns})$, we have
\begin{multline}\label{switching identity in NS}
e^{- \sum_{j \in \Z} (\mathcal{A}_j L_j + \mathcal{M}_{j - \frac{1}{2}}
G_{j - \frac{1}{2}})} \; (\alpha_0^\frac{1}{2})^{-2L_0} \; e^{-
\sum_{j \in \Z} (\mathcal{B}_j L_{-j} + \mathcal{N}_{j - \frac{1}{2}} G_{- j
+ \frac{1}{2}})} \\
= \; e^{\sum_{j \in \Z} (\Psi_{-j} L_{-j} + \Psi_{- j +
\frac{1}{2}} G_{- j + \frac{1}{2}})} \; e^{\sum_{j \in \Z} (\Psi_j L_j
+ \Psi_{j - \frac{1}{2}} G_{j - \frac{1}{2}})} \; e^{2\Psi_0 L_0} \;
(\alpha_0^\frac{1}{2})^{-2L_0} \; e^{\Gamma c} . 
\end{multline}
\end{prop}

\begin{proof} We use Huang's proof of Proposition
4.2.1 in \cite{H book} generalized to the present situation.  The idea 
is to use the Campbell-Baker-Hausdorff formula and compare the result 
with (\ref{normal order equation}).  However one must first establish the 
appropriate algebraic setting in which the Campbell-Baker-Hausdorff 
formula holds rigorously.  Let $\mathcal{C}$ be a sequence of even formal 
variables, $\mathcal{O}$ a sequence of odd formal variables and 
$\mathcal{C}_0$ another even formal variable.  Let $\mathcal{W}$ be the 
$\mathbb{Z}_2$-graded vector space $\mathfrak{ns} [\alpha_0^{1/2}, 
\alpha_0^{-1/2}] [[\mathcal{A}, \mathcal{B}, \mathcal{C}, \mathcal{C}_0]] 
[\mathcal{M}, \mathcal{N}, \mathcal{O}]$.  Consider the subspace $\mathcal{V}$ 
spanned by the elements 
\[L_\mathcal{A} = - \! \sum_{j \in \Z} \mathcal{A}_j L_j, \quad L_\mathcal{B} = - \!
\sum_{j \in \Z} \alpha_0^{-j}\mathcal{B}_j L_{-j}, \quad L_\mathcal{C} = - \! 
\sum_{j \in \Z} \mathcal{C}_j L_j,\]
\[ L_{\mathcal{C}_0} = - \;\mathcal{C}_0 L_0, \quad L_\mathcal{M} = - \!
\sum_{j \in \Z} \mathcal{M}_{j-\frac{1}{2}} G_{j - \frac{1}{2}}, \quad
L_\mathcal{N} = - \! \sum_{j \in \Z} \alpha_0^{-j + \frac{1}{2}} 
\mathcal{N}_{j- \frac{1}{2}} G_{-j + \frac{1}{2}},\]
and 
\[L_\mathcal{O} = - \! \sum_{j \in \Z} \mathcal{O}_{j - \frac{1}{2}} G_{j-\frac{1}{2}}. \] 
Then $\mathcal{V} \subset \mathcal{W}^0$, i.e., $\mathcal{V}$ 
contains only even elements.  In other words, letting $R = \mathbb{Q}[\alpha_0^{1/2},
\alpha_0^{-1/2}][[\mathcal{A}, \mathcal{B}, \mathcal{C}, \mathcal{C}_0]]
[\mathcal{M},\mathcal{N},\mathcal{O}]$, then $\mathcal{V}$ is in the 
$R$-envelope of $\mathfrak{ns}$ which is a Lie algebra; see Remarks 
\ref{envelope} and \ref{envelope again}.  Thus the same procedure used 
to prove Proposition 4.2.1 in \cite{H book} is valid for  $\mathcal{V}$ 
where we simply use Proposition \ref{normal order} instead of Huang's 
Proposition 2.2.5 in \cite{H book}, and the Neveu-Schwarz algebra 
relations and  corresponding universal enveloping algebra instead of 
just the Virasoro algebra relations and corresponding universal 
enveloping algebra. 
\end{proof}

\begin{rema}  Note that we used Proposition \ref{normal order} to prove 
Proposition \ref{central charge sewing} above.  In \cite{BHL}, we give 
a more straightforward and Lie-theoretic proof of Proposition 
\ref{central charge sewing} by proving a certain bijectivity property 
for the Campbell-Baker-Hausdorff formula in the theory of Lie algebras.  
This allows us to prove Proposition \ref{central charge sewing} 
directly for the Neveu-Schwarz algebra rather than go through the
representation in terms of superderivations and then lift to the
Lie algebra as we have done above.  
\end{rema}

Let $V = \coprod_{n \in \frac{1}{2} \mathbb{Z}} V_{(n)}$ be a module for
the Neveu-Schwarz algebra, and let $L(j), G(j - 1/2) \in
\mbox{End} \; V$ and $c \in \mathbb{C}$ be the representation images of 
$L_j$, $G_{j - 1/2}$ and $d$, respectively, such that for $v
\in V_{(n)}$, $L(0)v = nv$.  Combining Propositions \ref{universal 
enveloping} and \ref{central charge sewing}, we have the following 
corollary.  

\begin{cor}\label{normal order in End}
In the $\mathbb{C}[\alpha_0^{1/2}, \alpha_0^{-1/2}] [[\mathcal{A}, 
\mathcal{B}]] [\mathcal{M}, \mathcal{N}]$-envelope of $\mbox{\em End} \; V$, 
i.e., in $((\mbox{\em End} \; V) [\alpha_0^{1/2}, \alpha_0^{-1/2}] [[\mathcal{A}, 
\mathcal{B}]] [\mathcal{M}, \mathcal{N}])^0$, we have 
\begin{multline}\label{psi gamma corollary}
e^{- \sum_{j \in \Z} (\mathcal{A}_j L(j) + \mathcal{M}_{j - \frac{1}{2}}
G(j - \frac{1}{2}))} \cdot (\alpha_0^\frac{1}{2})^{-2L(0)} \cdot e^{-
\sum_{j \in \Z} (\mathcal{B}_j L(-j) + \mathcal{N}_{j - \frac{1}{2}} G(- j +
\frac{1}{2}))} 
\end{multline}
\begin{multline*}
= \; e^{\sum_{j \in \Z} (\Psi_{-j} L(-j) + \Psi_{- j + \frac{1}{2}} G(- j
+ \frac{1}{2}))} \cdot e^{\sum_{j \in \Z} (\Psi_j L(j) + \Psi_{j -
\frac{1}{2}} G(j - \frac{1}{2}))} \cdot \\
\cdot e^{2\Psi_0 L(0)} \cdot (\alpha_0^\frac{1}{2})^{-2L(0)} \cdot e^{\Gamma c} .
\end{multline*}
\end{cor}

We also have the following two propositions and corollaries corresponding 
to the identities (\ref{first Theta equation}) and (\ref{second Theta equation}).

\begin{prop} \label{first Theta prop in universal enveloping algebra}
Let $w$ be another formal even variable and $\rho$ be another formal
odd variable, and for $j \in \frac{1}{2} \mathbb{N}$, let $\Theta^{(1)}_j = 
\Theta^{(1)}_j (\alpha_0^{1/2}, \mathcal{A}, \mathcal{M}, (x,\varphi))$ be 
the sequence of formal series in $\mathbb{C} [x, \varphi][\alpha_0^{1/2}, 
\alpha_0^{-1/2}][[\mathcal{A}]][\mathcal{M}]$ given by (\ref{define first Theta}).  
Then in $(U(\mathfrak{ns}) [x, \varphi][\alpha_0^{1/2}, 
\alpha_0^{-1/2}] [[\mathcal{A}]][\mathcal{M}])^0$, we have   
\begin{multline*}
\exp \Biggl(\! - \! \! \sum_{m = -1}^{\infty} \sum_{j \in \Z}
\binom{j+1}{m+1} \alpha_0^{-j} x^{j - m}  \biggl(\!  \Bigl(
\mathcal{A}_j + 2\left(\frac{j-m}{j+1} \right)
\alpha_0^{\frac{1}{2}} x^{-1} \varphi \mathcal{M}_{j - \frac{1}{2}}
\Bigr)  L_m \\
+ \; x^{-1} \Bigl( \left(\frac{j-m}{j+1} \right) \alpha_0^{\frac{1}{2}}
\mathcal{M}_{j - \frac{1}{2}} + \varphi \frac{(j-m)}{2} \mathcal{A}_j
\Bigr)  G_{m + \frac{1}{2}} \biggl) \Biggr) 
\end{multline*}
\[= \; e^{\left(\!(\tilde{x} - \alpha_0^{-1}x) L_{-1}  + (\tilde{\varphi} -
\alpha_0^{-1/2}\varphi) G_{-\frac{1}{2}} \right)} \cdot e^{\left( - \! \sum_{j \in 
\Z} \left( \Theta^{(1)}_j L_j + \Theta^{(1)}_{j - \frac{1}{2}} G_{j - 
\frac{1}{2}} \right) \! \right)} \cdot e^{\left( - 2 \Theta^{(1)}_0 L_0 \right)},  \]
where $(\tilde{x}, \tilde{\varphi}) = (H_{\alpha_0^{1/2},\mathcal{A}, 
\mathcal{M}}^{(1)})^{-1} (x,\varphi)$ with $H_{\alpha_0^{1/2},\mathcal{A}, 
\mathcal{M}}^{(1)} (x,\varphi)$ given by (\ref{formalzero}).  
\end{prop}

\begin{proof} The proof is the same as that for
Proposition \ref{central charge sewing} except that in this case we
only consider the subalgebra with basis $L_j$, $G_{j + 1/2}$,
for $j \geq -1$ and use Proposition \ref{first Theta prop} instead of
Proposition \ref{normal order}. 
\end{proof}

{}From Propositions \ref{universal enveloping} and \ref{first Theta 
prop in universal enveloping algebra}, we have the following corollary.  

\begin{cor} \label{first Theta identity in End}
In $((\mbox{\em End} \; V) [x, \varphi][\alpha_0^{1/2}, 
\alpha_0^{-1/2}] [[\mathcal{A}]][\mathcal{M}])^0$, we have 
\begin{multline}\label{first Theta corollary}
\exp \Biggl(\! - \! \! \sum_{m = -1}^{\infty} \sum_{j \in \Z}
\binom{j+1}{m+1} \alpha_0^{-j} x^{j - m}  \\
\biggl(\!  \Bigl( \mathcal{A}_j + 2\left(\frac{j-m}{j+1} \right)
\alpha_0^{\frac{1}{2}} x^{-1} \varphi \mathcal{M}_{j - \frac{1}{2}}
\Bigr)  L(m) \\
+ \; x^{-1} \Bigl( \left(\frac{j-m}{j+1} \right) \alpha_0^{\frac{1}{2}}
\mathcal{M}_{j - \frac{1}{2}} + \varphi \frac{(j-m)}{2} \mathcal{A}_j
\Bigr)  G(m + \frac{1}{2}) \biggl) \Biggr) 
\end{multline} 
\begin{multline*} 
= \;  e^{\left((\tilde{x} - \alpha_0^{-1}x) L(-1)  + (\tilde{\varphi} - 
\alpha_0^{-1/2}\varphi)
G(-\frac{1}{2}) \right)} \cdot e^{\left( - \sum_{j \in \Z} \left( \Theta^{(1)}_j
L(j) + \Theta^{(1)}_{j - \frac{1}{2}} G(j - \frac{1}{2}) \right) \right)} 
 \cdot \\
\cdot e^{\left( - 2\Theta^{(1)}_0 L(0) \right)} 
\end{multline*} 
for $(\tilde{x}, \tilde{\varphi}) = (H_{\alpha_0^{1/2},\mathcal{A}, 
\mathcal{M}}^{(1)})^{-1} (x,\varphi)$.
\end{cor}

\begin{prop} \label{second Theta prop in universal enveloping algebra}
For $j \in \frac{1}{2} \mathbb{N}$, let $\Theta^{(2)}_j = 
\Theta^{(2)}_j (\mathcal{B}, \mathcal{N}, (x,\varphi))$ be the sequence of formal
series in $\mathbb{C} [x^{-1}, \varphi][[\mathcal{B}]][\mathcal{N}]$ given by
(\ref{define second Theta}).  Then in 
\[(U(\mathfrak{ns}) [x^{-1}, \varphi][[\mathcal{B}]][\mathcal{N}])^0,\] 
we have   
\begin{multline*}
\exp \biggl( \sum_{m = -1}^{\infty} \sum_{j \in \Z}
\binom{-j + 1}{m + 1} x^{-j - m} \biggl( \Bigl( \mathcal{B}_j + 2\varphi 
\mathcal{N}_{j - \frac{1}{2}} \Bigr)  L_m \\   
+ \; \Bigl( \mathcal{N}_{j - \frac{1}{2}} + \varphi x^{-1}
\frac{(-j-m)}{2} \mathcal{B}_j \Bigr) G_{m + \frac{1}{2}}
\biggr)  \biggr) 
\end{multline*}   
\[= \; e^{\left((\tilde{x} - x) L_{-1}  + (\tilde{\varphi} - \varphi)
G_{-\frac{1}{2}} \right)} \cdot e^{\left( - \sum_{j \in \Z} \left( \Theta^{(2)}_j
L_j + \Theta^{(2)}_{j - \frac{1}{2}} G_{j - \frac{1}{2}} \right) \right)}
\cdot e^{\left( - 2 \Theta^{(2)}_0 L_0 \right)}, \hspace{.5in} \]
where $(\tilde{x}, \tilde{\varphi}) = (H_{\mathcal{B}, 
\mathcal{N}}^{(2)})^{-1} \circ I (x,\varphi)$ is given by (\ref{anotherinfty}).  
\end{prop}

\begin{proof} The proof is the same as that for Proposition 
\ref{central charge sewing} except that in this case we only consider 
the subalgebra with basis $L_j$, $G_{j + 1/2}$, for $j \geq -1$ and 
use Proposition \ref{second Theta prop} instead of Proposition \ref{normal order}.  
\end{proof}

{}From Propositions \ref{universal enveloping} and \ref{second Theta 
prop in universal enveloping algebra}, we have the following corollary.  

\begin{cor} \label{second Theta identity in End}
In $((\mbox{\em End} \; V) [x^{-1}, \varphi] [[\mathcal{B}]]
[\mathcal{N}])^0$, we have 
\begin{multline}\label{second Theta corollary}
\exp \biggl(\sum_{m = -1}^{\infty}  \sum_{j \in \Z}
\binom{-j + 1}{m + 1} x^{-j - m} \biggl( \Bigl( \mathcal{B}_j + 2\varphi 
\mathcal{N}_{j - \frac{1}{2}} \Bigr) L(m) \\
+ \; \Bigl( \mathcal{N}_{j - \frac{1}{2}} + \varphi x^{-1}
\frac{(-j-m)}{2} \mathcal{B}_j \Bigr) G(m + \frac{1}{2}) \biggr)  \biggr)
\end{multline} 
\[= \; e^{\left((\tilde{x} - x) L(-1)  + (\tilde{\varphi} - \varphi)
G(-\frac{1}{2}) \right)} \cdot e^{\left( - \sum_{j \in \Z} \left( \Theta^{(2)}_j
L(j) + \Theta^{(2)}_{j - \frac{1}{2}} G(j - \frac{1}{2}) \right) \right)}
\cdot e^{\left( -2\Theta^{(2)}_0 L(0) \right)} ,\]
for $(\tilde{x}, \tilde{\varphi}) = (H_{\mathcal{B}, \mathcal{N}}^{(2)})^{-1} 
\circ I (x,\varphi)$.
\end{cor}

\section[Positive-energy representations]{The corresponding
identities for positive-energy representations of the $N=1$ Neveu-Schwarz 
algebra}

Using positive-energy representations of the Neveu-Schwarz algebra, we
study the formal series $\Gamma$, $\Psi_j$, for $j \in \frac{1}{2} 
\mathbb{Z}$, and $\Theta^{(1)}_j$ and $\Theta^{(2)}_j$, for $j \in 
\frac{1}{2} \mathbb{N}$, and will see that for these representations the 
identities (\ref{psi gamma corollary}), (\ref{first Theta corollary}) 
and (\ref{second Theta corollary}) become identities containing only a 
finite number of formal variables.  

\begin{prop}\label{where Psi's actually live}
The formal series $\Psi_j$, for $j \in \frac{1}{2}\mathbb{Z}$, and $\Gamma$
are actually in 
\[\mathbb{C}[\mathcal{A}][\mathcal{M}][\mathcal{B}][\mathcal{N}]
[[\alpha_0^{-\frac{1}{2}} ]] . \]  
\end{prop}

\begin{proof} Using the Verma modules $M(c,h)$ defined in
Section 3.5 and Corollary \ref{normal order in End}, the proof is
completely analogous to the proof of Lemma 4.3.1 in \cite{H book}. 
\end{proof}
 
We have the following immediate corollary of Proposition 
\ref{where Psi's actually live}. 

\begin{cor}\label{Psi, Gamma well-defined}
Let $(A,M)$ and $(B,N)$ be two sequences in $\bigwedge_*^\infty$;
let $\asqrt \in (\bigwedge_*^0)^\times$; and let $t^{1/2}$ be an even 
formal variable.  Then $\Psi_j(t^{-1/2}\asqrt, A, M, B, N)$, for $j \in
\frac{1}{2}\mathbb{Z}$, and $\Gamma(t^{-1/2} \asqrt, A,M,B,N)$ are well
defined and belong to $\bigwedge_*[[t^{1/2}]]$.  
\end{cor}

\begin{prop}\label{positive energy gamma}
Let $(A,M)$, $(B,N)$, $\asqrt$ and $t^{1/2}$ be as in Corollary
\ref{Psi, Gamma well-defined}, and let $V$ be a positive-energy module for the
Neveu-Schwarz algebra.  Then the following identity holds in
$(\mbox{\em End} \; (\bigwedge_* \otimes_{\mathbb{C}} V))^0 [[t^{1/2}]]$.
\begin{multline}\label{positive energy gamma equation}
e^{- \sum_{j \in \Z} (A_j L(j) + M_{j - \frac{1}{2}}
G(j - \frac{1}{2}))} \cdot (t^{-\frac{1}{2}}\asqrt)^{-2L(0)} \cdot\\
\cdot e^{- \sum_{j \in \Z} (\mathcal{B}_j L(-j) + \mathcal{N}_{j - \frac{1}{2}} 
G(- j + \frac{1}{2}))} \cdot 
\end{multline}
\begin{multline*}
= \; e^{\sum_{j \in \Z} (\Psi_{-j}(t^{-\frac{1}{2}} \asqrt, A, M, B, N) L(-j) +
\Psi_{- j + \frac{1}{2}}(t^{-\frac{1}{2}} \asqrt, A, M, B, N) G(- j +
\frac{1}{2}))} \cdot \\
\cdot e^{\sum_{j \in \Z} (\Psi_j(t^{-\frac{1}{2}}\asqrt, A, M, B, N) L(j) +
\Psi_{j - \frac{1}{2}}(t^{-\frac{1}{2}} \asqrt, A, M, B, N) G(j -
\frac{1}{2}))}  \cdot\\
\cdot e^{2\Psi_0(t^{-\frac{1}{2}} \asqrt, A, M, B, N) L(0)} \cdot (t^{-\frac{1}{2}}
\asqrt)^{-2L(0)} \cdot e^{\Gamma(t^{-\frac{1}{2}} \asqrt, A,M,B,N)c} . 
\end{multline*}  
\end{prop}

\begin{proof} Since $V$ is a positive energy module for
$\mathfrak{ns}$ and by Corollary \ref{Psi, Gamma well-defined}
\[\Psi_j(t^{-\frac{1}{2}} \asqrt, A, M, B, N), \; \;
\Gamma(t^{-\frac{1}{2}} \asqrt, A, M, B, N) \in \mbox{$\bigwedge_*$}
[[t^\frac{1}{2}]],\]
for $j \in \frac{1}{2} \mathbb{Z}$, both the right-hand side and 
the left-hand side of (\ref{positive energy gamma equation}) are 
well-defined elements in $((\mbox{End} \; (\bigwedge_*
\otimes_{\mathbb{C}} V)) [[t^{1/2}]])^0$.  Then  by Corollary \ref{normal
order in End}, they are equal. 
\end{proof}

\begin{prop} \label{where first Theta's actually live}
The formal series $\Theta^{(1)}_j (t^{ - 1/2} \alpha_0^{1/2}, \mathcal{A}, 
\mathcal{M}, (x,\varphi))$, for $j \in \frac{1}{2} \mathbb{N}$, are in 
$\mathbb{C} [x, \varphi][\mathcal{A}][\mathcal{M}][\alpha_0^{-1/2}]
[[t^{1/2}]]$. \footnote{There is a misprint in the analogous nonsuper case of
Proposition \ref{where first Theta's actually live} given in \cite{H book}.
The series $\Theta^{(1)}_j(\mathcal{A},t^{-1},\alpha_0, y)$ of Lemma
4.3.4 in \cite{H book} are elements in $\mathbb{C}[y][\mathcal{A}]
[\alpha_0^{-1}][[t]]$, not $\mathbb{C}[y][\mathcal{A}][[t]]$ as stated.  }
\end{prop}

\begin{proof} The proof is analogous to the proof of
Proposition \ref{where Psi's actually live} except that we use
Corollary \ref{first Theta identity in End} instead of Corollary \ref{normal
order in End}. 
\end{proof}

We have the following immediate corollary of Propositions \ref{where first
Theta's actually live}.

\begin{cor} \label{first Theta cor}
Let $\asqrt \in (\bigwedge_*^0)^\times$, and let $(A,M) \in \bigwedge_*^\infty$.  
The series 
\[\Theta^{(1)}_j (t^{ - \frac{1}{2}} \asqrt, A,M, (x,\varphi)),\] 
for $j \in \frac{1}{2}\mathbb{N}$, are well defined and belong to $\bigwedge_* 
[x, \varphi][[t^{1/2}]]$. \footnote{There is a misprint in the analogous 
nonsuper case of Corollary \ref{first Theta cor} in \cite{H book}.  Corollary 
4.3.7 in \cite{H book} should state that the $\Theta_j(A,t^{-1}a_0,y)$, for 
$j \in \mathbb{N}$, belong to $\mathbb{C}[y][[t]]$, i.e.,
$a_0 \in \mathbb{C}^\times$ and $A_j \in \mathbb{C}$, for $j \in \mathbb{Z}_+$, with
$A = \{A_j\}_{j \in \mathbb{Z}_+}$ should be substituted for the formal variables
$\alpha_0$ and $\mathcal{A}_j$ in $\Theta_j(\mathcal{A},t^{-1}\alpha_0,y)$ in the
corollary.}
\end{cor}

Let $(\tilde{x}(t^{1/2}),\tilde{\varphi} (t^{1/2})) = 
H^{(1)}_{ t^{ - 1/2}\asqrt, A,M}(x,\varphi)$ where
$H_{\alpha_0^{1/2},\mathcal{A}, \mathcal{M}}^{(1)}  (x,\varphi)$ is given
by (\ref{formalzero}). 

\begin{prop}\label{last first Theta Prop} 
Let $\asqrt \in (\bigwedge_*^0)^\times$; let $(A,M) \in \bigwedge_*^\infty$; 
and let $V$ be a positive-energy module for the $N=1$  Neveu-Schwarz algebra.  Then 
in 
\[((\mbox{\em End} \; (\mbox{$\bigwedge_*$} \otimes_{\mathbb{C}} V))
[[t^{\frac{1}{2}}]] [[x,x^{-1}]][\varphi])^0, \] the following identity
holds for $\Theta^{(1)}_j =
\Theta^{(1)}_j (t^{-1/2}\asqrt, A,M,(x,\varphi))$. 
\begin{multline}\label{last first Theta equation}
\exp \Biggl(\! - \! \! \sum_{m = -1}^{\infty} \sum_{j \in \Z}
\binom{j+1}{m+1} t^j \asqrt^{-2j} x^{j - m}  \\
\biggl(\!  \Bigl( A_j + 2\left(\frac{j-m}{j+1} \right)
t^{-\frac{1}{2}} \asqrt x^{-1} \varphi M_{j - \frac{1}{2}}
\Bigr)  L(m) \\
+ \; x^{-1} \Bigl( \left(\frac{j-m}{j+1} \right) t^{-\frac{1}{2}}
\asqrt M_{j - \frac{1}{2}} + \varphi \frac{(j-m)}{2} A_j
\Bigr)  G(m + \frac{1}{2}) \biggl) \Biggr) 
\end{multline}
\begin{multline*}
= \; \exp \left( (\tilde{x}(t^\frac{1}{2}) - t \asqrt^{-2}x) L(-1) +
(\tilde{\varphi}(t^\frac{1}{2}) - t^{\frac{1}{2}}\asqrt^{-1} \varphi) 
G(-\frac{1}{2}) \right) \cdot \\
\exp \Biggl( \! - \! \sum_{j \in \Z} \Bigl( \Theta^{(1)}_j L(j) +
\Theta^{(1)}_{j - \frac{1}{2}} G(j - \frac{1}{2}) \Bigr) \! \Biggr) \! \cdot \exp
\left( - 2\Theta^{(1)}_0 L(0) \right) . 
\end{multline*}
\end{prop}

\begin{proof}  Since $V$ is a positive energy module for
$\mathfrak{ns}$ and by Corollary \ref{first Theta cor}
\[\Theta^{(1)}_j (t^{-\frac{1}{2}}\asqrt, A,M, (x,\varphi)) \in 
\mbox{$\bigwedge_*$} [x^{-1}, \varphi] [[t^{\frac{1}{2}}]], \]
for $j \in \frac{1}{2} \mathbb{N}$, both the right-hand side and the 
left-hand side of (\ref{last first Theta equation}) are well-defined 
elements in $((\mbox{End} \; (\bigwedge_* \otimes_{\mathbb{C}} V))
[[t^{1/2}]] [[x, x^{-1}]][\varphi])^0$.   Then by Corollary \ref{first
Theta identity in End}, they are equal. 
\end{proof}

\begin{prop} \label{where second Theta's actually live}
The formal series $\Theta^{(2)}_j (\{t^k \mathcal{B}_k, t^{k - 1/2}
\mathcal{N}_{k - 1/2} \}_{k \in \Z}, (x,\varphi))$, for $j \in
\frac{1}{2} \mathbb{N}$, are in $\mathbb{C} [x^{-1}, \varphi]
[\mathcal{B}][\mathcal{N}][[t^{1/2}]]$. 
\end{prop}

\begin{proof} The proof is analogous to the proof of
Proposition \ref{where Psi's actually live} except that we use
Corollary \ref{second Theta identity in End} instead of Corollary \ref{normal
order in End}.  
\end{proof}

We have the following immediate corollary of Proposition 
\ref{where second Theta's actually live}.

\begin{cor} \label{second Theta cor}
Let $(B,N) \in \bigwedge_*^\infty$.  The series 
\[\Theta^{(2)}_j (\{t^k B_k, t^{k - \frac{1}{2}} N_{k - \frac{1}{2}} \}_{k
\in \Z}, (x,\varphi))\] 
are well defined and belong to $\bigwedge_* [x^{-1}, \varphi]
[[t^{1/2}]]$. 
\end{cor}

Let $(\tilde{x}(t^{1/2}),\tilde{\varphi} (t^{1/2})) = (H_{ \{t^j B_j, 
t^{j - 1/2} N_{j - 1/2} \}_{j \in \Z}}^{(2)})^{-1} \circ I (x,\varphi)$ 
where $(H_{\mathcal{B},\mathcal{N}}^{(2)})^{-1} \circ I (x,\varphi)$ is 
given by (\ref{anotherinfty}). 

\begin{prop}\label{last second Theta Prop} 
Let $(B,N) \in \bigwedge_*^\infty$, and let $V$ be a 
positive-energy module for the $N=1$ Neveu-Schwarz algebra.  In 
$((\mbox{\em End} \; (\bigwedge_* \otimes_{\mathbb{C}} V))
[[t^{1/2}]][[x,x^{-1}]][\varphi])^0$, the  following identity holds for
$\Theta^{(2)}_j = \Theta^{(2)}_j (\{t^k  B_k, t^{k - 1/2} N_{k - 1/2}
\}_{k \in \Z}, (x,\varphi))$. 
\begin{multline}\label{last second Theta equation}
\exp \biggl(\sum_{m = -1}^{\infty} \sum_{j \in \Z} \binom{-j + 1}{m +
1} x^{-j - m} \biggl( \Bigl( t^j B_j +  2\varphi t^{j - \frac{1}{2}} N_{j
- \frac{1}{2}} \Bigr) L(m) \\
+ \; \Bigl( t^{j - \frac{1}{2}}  N_{j - \frac{1}{2}} + \varphi x^{-1}
\frac{(-j-m)}{2} t^j B_j \Bigr) G(m + \frac{1}{2}) \biggr) \biggr)
\end{multline} 
\begin{multline*}
= \; \exp \left( (\tilde{x}(t^\frac{1}{2}) - x) L(-1) +
(\tilde{\varphi}(t^\frac{1}{2}) - \varphi) G(-\frac{1}{2}) \right)
\cdot \\  
\exp \Biggl( \! - \! \sum_{j \in \Z} \Bigl( \Theta^{(2)}_j L(j) +
\Theta^{(2)}_{j - \frac{1}{2}} G(j - \frac{1}{2}) \Bigr) \! \Biggr) \! \cdot \exp
\left( - 2\Theta^{(2)}_0 L(0) \right) . 
\end{multline*} 
\end{prop}

\begin{proof}  Since $V$ is a positive energy module for
$\mathfrak{ns}$ and by Corollary \ref{second Theta cor}
\[\Theta^{(2)}_j (\{t^k B_k, t^{k - \frac{1}{2}} N_{k - \frac{1}{2}} \}_{k
\in \Z}, (x,\varphi)) \in \mbox{$\bigwedge_*$} [x^{-1}, \varphi]
[[t^{\frac{1}{2}}]], \]
for $j \in \frac{1}{2} \mathbb{N}$, both the right-hand side and 
the left-hand side of (\ref{last second Theta equation}) are 
well-defined elements in $((\mbox{End} \; (\bigwedge_*
\otimes_{\mathbb{C}} V)) [[t^{1/2}]] [[x, x^{-1}]] [\varphi])^0$.  Then
by Corollary \ref{second Theta identity in End},  they are equal. 
\end{proof}

\chapter{An analytic study of the sewing operation}

In this chapter, we study the moduli space of $N=1$ superspheres with 
tubes and the sewing operation using the formal calculus developed in 
Chapter 3.  We give a reformulation of the moduli space and a more 
explicit description of the sewing operation on this space including
a solution to the sewing equation and normalization and boundary
conditions.  We prove the analyticity and convergence of certain 
infinite series arising {}from the sewing operation.  We define a 
bracket operation on the $N=1$ supermeromorphic tangent space of the 
moduli space of superspheres with tubes at the identity element and 
show that this gives a representation of the $N=1$ Neveu-Schwarz 
algebra with central charge zero.

This chapter is organized as follows. In Section 4.1, using the 
characterization of local superconformal coordinates in terms of 
exponentials of certain infinite sums of superderivations proved in  
Chapter 3, we show that a canonical supersphere with tubes can be 
identified with certain data concerning the punctures and the 
coefficients of the infinite sums of superderivations appearing in  
the expressions for the local coordinates.  Thus we can identify the
moduli space of superspheres with tubes with the set of this data.  
In Section 4.2, we introduce an action of the symmetric group on $n$ 
letters on the moduli space of superspheres with $1 + n$ tubes.  This 
action is used in \cite{B thesis} and is a property needed to show 
that the moduli space with the sewing operation is a partial operad; 
see Remark \ref{operad remark}. 

In Section 4.3, we define supermeromorphic functions on the moduli
space.  We show in \cite{B thesis}, that such functions include 
the rational function counterparts to correlation functions for an 
$N=1$ Neveu-Schwarz vertex operator superalgebra  with odd formal 
variables \cite{B vosas}.  Using this notion of supermeromorphic 
function, we define the supermeromorphic tangent space of the moduli 
space.  We show that any supermeromorphic tangent vector can be 
expressed as an infinite sum of supermeromorphic tangent vectors 
tangent to the coordinate system of the moduli space defined by the 
characterization of the moduli space obtained in Section 4.1.

In Section 4.4, we show that any element of the moduli space can be
obtained {}from three types of elements with one, two and three tubes,
respectively.  We study some simple sewings with such elements, and
show that the partial monoid of the moduli space of $1+1$ tubes has 
two subgroups which are of interest.  Then we study a certain linear 
functional in the supermeromorphic tangent space at the identity of 
the moduli space of superspheres with $1 + 1$ tubes.  The results of 
this section are used in \cite{B thesis} and are in fact geometric 
versions of certain axioms for an $N = 1$ vertex operator superalgebra
\cite{B vosas}, namely: associativity; the fact that the $N=1$ 
vertex operator associated with the special element $\tau$ is equal 
to the formal series with coefficients satisfying the Neveu-Schwarz 
algebra relations; and the $G(-1/2)$-derivative property.

In Section 4.5, we introduce generalized superspheres with tubes which
are formal analogues of superspheres with tubes in that, in general,
the local coordinates are not necessarily convergent power series.
We define two subsets of the moduli space of generalized superspheres
with $1+1$ tubes and show that they both have group structures which
are isomorphic to the structures discussed in Chapter 3.

In Section 4.6, we give a more detailed description of the sewing
operation by giving a more explicit formula for obtaining a canonical
supersphere with tubes {}from the sewing of two canonical superspheres
with tubes.  This involves giving an explicit description of the 
transformation {}from the sewn superspheres to a canonical supersphere.  
This transformation consists of the uniformizing function taking
the sewn superspheres to a super-Riemann sphere and a superprojective
transformation taking this super-Riemann sphere to a canonical
supersphere with tubes.  The uniformizing component of this 
transformation satisfies the sewing equation, normalization 
conditions and boundary conditions introduced in Chapter 2 and  
studied formally in Chapter 3.  We show that the expansions of the two
halves $F^{(1)}$ and $F^{(2)}$ of this uniformizing function as formal 
series in terms of the local coordinate charts correspond to the formal
series $\bar{F}^{(1)} \circ s_{(z_i,\theta_i)}$ and $\bar{F}^{(2)}$
with $\bar{F}^{(1)}$ and $\bar{F}^{(2)}$ given explicitly by Theorem
\ref{uniformization}, equations (\ref{F1}) and (\ref{F2}) and 
Proposition \ref{normal order}.  Thus in order to prove the
analyticity and convergence of this expansion, we need only prove the 
analyticity and convergence of the series $\Psi_j(t^{-1/2} \asqrt, A, 
M, B, N)$, for $j \in \frac{1}{2} \mathbb{Z}$, arising algebraically
{}from the sewing equation.  

Following \cite{H book}, we prove the analyticity of the
series $\Psi_j(t^{-1/2} \asqrt, A, M, B, N)$, for $j \in \frac{1}{2} 
\mathbb{Z}$, by deforming one of the superspheres being sewn and then 
using the Fischer-Grauert Theorem in the deformation theory of complex 
manifolds to show the moduli space is locally trivial under this 
deformation.  However this theorem only holds for finite-dimensional, 
compact complex manifolds, e.g., the body of a supersphere, whereas 
super-Riemann surfaces over $\bigwedge_\infty$ are 
infinite-dimensional, non-compact complex manifolds.  In fact they are 
fiber bundles over the body with fiber isomorphic to 
$(\bigwedge_\infty)_S$.  Therefore, we define two families of global 
sections of our supersphere which cover the supersphere.  Each section 
is complex analytically isomorphic to the body, i.e., to a Riemann 
sphere.  We then use the Fischer-Grauert Theorem for each section 
which allows us to prove that the uniformizing function on each 
section of the deformed supersphere is analytic in the deformation 
variable $t^{1/2}$.  Since the two families of sections completely 
cover the supersphere, this proves the analyticity of the uniformizing 
function, which in turn allows us to prove the analyticity and 
convergence of the series $\Psi_j(t^{-1/2}\asqrt, A, M, B, N)$, for 
$j \in \frac{1}{2} \mathbb{Z}$ and $t^{1/2} \in \mathbb{C}$ with 
$|t^{1/2}| \leq 1$.  This result shows that when two canonical 
superspheres can be sewn, the uniformizing function studied 
algebraically in Chapter 3 converges and is equal to the geometric 
uniformizing transformation taking the resulting sewn superspheres to a 
super-Riemann sphere.  This proof of Proposition \ref{use of FG} is 
different and more concise in its use of the underlying geometry than 
the proof originally given in \cite{B thesis}. 

In Section 4.7, we define a bracket operation on the supermeromorphic
tangent space of the moduli space of superspheres with tubes at the
identity element and show that this gives a representation of the $N=1$
Neveu-Schwarz algebra with central charge zero.  In proving the $N=1$
Neveu-Schwarz algebra relations for this bracket operation, we use the
information about the lower order terms in the solution to the sewing
equation given formally in Theorem \ref{uniformization} and
Proposition \ref{normal order}.  In defining the bracket relations we 
see that the definition is not as straightforward as on might initially 
expect it to be.  This is related to the fact that the diffeomorphism
group of the circle does not have a Lie group complexification, i.e,
on the nonsuper level of the Virasoro Lie algebra, there is no 
corresponding complex Lie group (cf. \cite{L}).

\section[A reformulation of the moduli space]{A reformulation of the 
moduli space of superspheres with tubes}

Let
\begin{multline*}
\mathcal{H}  =  \bigl\{ (A,M) \in \mbox{$\bigwedge_\infty^\infty$} 
\; | \; \tilde{E}(A,M) \hspace{.1cm} \mbox{is an absolutely convergent power series} \bigr. \\
\bigl. \mbox{ in some neighborhood of } 0 \bigr\},
\end{multline*} 
and for $n \in \Z$, let
\[ SM^{n - 1} = \bigl\{ \bigl((z_1, \theta_1),...,(z_{n-1}, \theta_{n-1})\bigr) \; | \;
(z_i, \theta_i) \in \mbox{$\bigwedge_\infty^\times$}, \; (z_i)_B \neq (z_j)_B , \;
\mbox{for} \; i \neq j \bigr\} . \] 
Note that for $n=1$, the set $SM^0$ has exactly one element.

For any supersphere with $1 + n$ tubes, we can write the power series 
expansion of the local coordinate at the $i$-th puncture $(z_i,\theta_i)$, 
for $i = 1,...,n$, as 
\begin{eqnarray*}
H_i (w, \rho) &=& \exp
\Biggl(\! - \! \sum_{j \in \Z} \Bigl( A^{(i)}_j L_j(x,\varphi)
+  M^{(i)}_{j - \frac{1}{2}} G_{j -\frac{1}{2}} \Bigr) \! \Biggr) \cdot\\
& & \hspace{1.4in} \left. \cdot (\asqrt^{(i)})^{-2L_0(x,\varphi)} \cdot  (x, 
\varphi) \right|_{(x,\varphi) = (w - z_i - \rho \theta_i, \rho - \theta_i)}\\
&=& \hat{E}(\asqrt^{(i)}, A^{(i)}, M^{(i)}) (w - z_i - \rho \theta_i,
\rho - \theta_i) 
\end{eqnarray*}
for $(\asqrt^{(i)}, A^{(i)}, M^{(i)}) \in (\bigwedge_\infty^0)^\times 
\times \mathcal{H}$, and we can write the power series expansion of the
local coordinate at $\infty$ as  
\begin{eqnarray*}
H_0 (w, \rho) &=& \exp \Biggl(\sum_{j \in \Z} \Bigl( A^{(0)}_j
L_{-j}(w,\rho) +  M^{(0)}_{j - \frac{1}{2}} G_{-j + \frac{1}{2}}(w,\rho)
\Bigr) \! \Biggr) \cdot \Bigl(\frac{1}{w}, \frac{i \rho}{w}\Bigr)  \\ 
&=& \tilde{E}(A^{(0)}, -iM^{(0)}) \Bigl(\frac{1}{w}, \frac{i \rho}{w}\Bigr)
\end{eqnarray*}
for $(A^{(0)},-iM^{(0)}) \in \mathcal{H}$.  But $(A^{(0)},-iM^{(0)}) \in 
\mathcal{H}$ if and only if $(A^{(0)},M^{(0)}) \in \mathcal{H}$.

Thus we have the following proposition.

\begin{prop}\label{moduli2}
The moduli space of superspheres with $1 + n$ tubes, for $n \in \Z$, 
can be identified with the set   
\begin{equation}
SK(n) = SM^{n-1} \times \mathcal{H} \times \bigl((\mbox{$\bigwedge_\infty^0$})^\times 
\times \mathcal{H}\bigr)^n .
\end{equation}
\end{prop}

For superspheres with one tube, we have:
 
\begin{prop}\label{moduli1}
The moduli space of superspheres with one tube can be identified 
with the set
\[SK(0) = \bigl\{(A,M) \in \mathcal{H} \; | \; (A_1, M_{\frac{1}{2}}) =
(0,0) \bigr\} .\] 
\end{prop}

By Propositions \ref{moduli2} and \ref{moduli1}, we can identify
$SK(n)$ with the moduli space of superspheres with $1 + n$ tubes
for $n \in \mathbb{N}$, and the set 
\[SK = \bigcup_{n \in \mathbb{N}} SK(n) \]
can be identified with the moduli space of superspheres with tubes. 
The actual elements of $SK$ give the data for a canonical supersphere
representative of a given equivalence class of superspheres with
tubes modulo superconformal equivalence. {}From now on it will be 
convenient to refer to $SK$ as the moduli space of superspheres with 
tubes.  Any element of $SK(n)$, for $n \in \Z$, can be written as 
\[\bigl((z_1, \theta_1),...,(z_{n-1}, \theta_{n-1}); (A^{(0)}, M^{(0)}),
(\asqrt^{(1)}, A^{(1)}, M^{(1)}), ...,(\asqrt^{(n)}, A^{(n)}, M^{(n)} 
)\bigr)  \]
where $(z_1, \theta_1),...,(z_{n-1}, \theta_{n-1}) \in SM^{n-1}$,
$(A^{(0)}, M^{(0)}) \in \mathcal{H}$, and 
\[(\asqrt^{(1)}, A^{(1)}, M^{(1)}),...,(\asqrt^{(n)}, A^{(n)}, 
M^{(n)}) \in (\mbox{$\bigwedge_\infty^0$})^\times 
\times \mathcal{H} . \]
Thus for an element $Q \in SK$, we can think of $Q$ as consisting of the
above data, or as being a canonical supersphere with tubes corresponding
to that data.  

In Chapter 2 we introduced the sewing operation for
superspheres with tubes.  This gave rise to the sewing 
equation (\ref{sewing equation}), normalization conditions 
(\ref{normalization 1}) -- (\ref{normalization 3}) and 
boundary conditions (\ref{boundary conditions trivial infinity 1}) --
(\ref{boundary conditions trivial puncture 2}) derived
{}from taking two canonical superspheres with tubes, sewing them 
together and looking at the uniformizing function which maps the 
resulting supersphere to a super-Riemann sphere.  This resulting 
super-Riemann sphere is then superconformal to a canonical 
supersphere, i.e., an element in $SK$, via a superprojective 
transformation. 
 
Thus for any $n \in \Z$, $m \in \mathbb{N}$, and any positive integer 
$i \leq n$, the sewing operation for superspheres with tubes defined
in Chapter 2 induces an  operation $ _i\infty_0 : SK(n) \times 
SK(m) \rightarrow SK(n + m -  1)$.  We will still call it the 
{\it sewing operation}, and as in \cite{H book}, we will use the 
notation $_i\infty_0$ first introduced by Vafa \cite{V} to denote 
the sewing, although now we use it to denote the sewing operation on 
$SK$.  Recall that given two elements $Q_1$ and $Q_2$ of 
$SK$, it may not be possible to sew the $i$-th puncture of $Q_1$ with 
the $0$-th puncture of $Q_2$ without rescaling the local coordinate 
map of $Q_1$ at its $i$-th puncture.

The element of $\mathcal{H}$ with all components equal to $0$ will be 
denoted by $(\mathbf{0},\mathbf{0})$ or just $\mathbf{0}$.  Note that 
in terms of the chart $(U_\sou, \sou)$ of $S\hat{\mathbb{C}}$, the 
local coordinate chart corresponding to $(1,\mathbf{0},\mathbf{0}) 
\in (\bigwedge_\infty^0)^\times \times \mathcal{H}$ is the identity 
map on $\bigwedge_\infty$ if the puncture is at $0$, is the shift 
$s_{(z_i,\theta_i)}(w,\rho) =  (w - z_i -\rho\theta_i,\rho - \theta_i)$ 
if the puncture is at $(z_i,\theta_i)$, and is $I(w,\rho) = (1/w,
i\rho/w)$ if the puncture is at $\infty$.  We will sometimes refer to 
such coordinates as {\it standard local coordinates}.

Using the definition of sewing {}from Chapter 2, we see that $SK$ has a 
unit $e$ under the sewing operation 
\begin{equation} 
e = (\mathbf{0}, (1, \mathbf{0})) \in SK(1)
\end{equation}
in the sense that for $Q \in SK(n)$ and $0<i\leq n$, the
$i$-th puncture of $Q$ can always be sewn with the $0$-th puncture of
$e$, the first puncture of $e$ can always be sewn with the $0$-th
puncture of $Q$, and we have
\[ Q \; _i\infty_0 \; e = Q, \qquad e \; _1\infty_0 \; Q = Q . \]

{}From the geometry of sewing, the following associativity of the sewing
operation is obvious.

\begin{prop}\label{sewing associativity}
Let $l \in \Z$ and $m, n \in \mathbb{N}$ such that $l + m -1 \in \Z$, and
let $Q_1 \in SK(l)$, $Q_2 \in SK(m)$, $Q_3 \in SK(n)$, and $i,j \in \Z$
such that $1 \leq i \leq l$, and $1 \leq j \leq l + m -1$.  The iterated
sewings $(Q_1 \; _i\infty_0 \; Q_2) \; _j\infty_0 \; Q_3$ exist if and
only if one of the following holds:

(i) $j<i$ and the sewings $(Q_1 \; _j\infty_0 \; Q_3) \; _{i + n -
1}\infty_0 \; Q_2$ exist, in which case
\[(Q_1 \; _i\infty_0 \; Q_2) \; _j\infty_0 \; Q_3 = (Q_1 \; _j\infty_0
\; Q_3) \; _{i + n - 1}\infty_0 \; Q_2 ;\] 

(ii) $j \geq i + m$ and the sewings $(Q_1 \; _{j - m + 1}\infty_0
\; Q_3) \; _i\infty_0 \; Q_2$ exist, in which case
\[(Q_1 \; _i\infty_0 \; Q_2) \; _j\infty_0 \; Q_3 = (Q_1 \; _{j - m +
1}\infty_0 \; Q_3) \; _i\infty_0 \; Q_2 ;\] 

or

(iii) $i \leq j < i + m$ and the sewings $Q_1 \; _i\infty_0 \;
(Q_2 \; _{j - i + 1}\infty_0 \; Q_3)$ exist, in which case 
\[(Q_1 \; _i\infty_0 \; Q_2) \; _j\infty_0 \; Q_3 = Q_1 \; _i\infty_0
\; (Q_2 \; _{j - i + 1}\infty_0 \; Q_3) .\]
\end{prop}  

\section[An action of the symmetric groups on the moduli space]{An action of the symmetric 
group $S_n$ on the moduli space $SK(n)$}

Let $S_n$ be the group of permutations on $n$ letters, for $n \in
\Z$.   There is a natural (right) action of $S_{n - 1}$ on $SK(n)$ defined
by  permuting the ordering on the first $n - 1$ positively oriented 
punctures and their local coordinates.  More explicitly, for $\sigma 
\in S_{n-1}$, and $Q \in SK(n)$ given by
\[Q = \bigl((z_1, \theta_1),...,(z_{n-1}, \theta_{n-1}); (A^{(0)},M^{(0)}), 
(\asqrt^{(1)}, A^{(1)}, M^{(1)}),..., (\asqrt^{(n)}, A^{(n)}, M^{(n)})\bigr), \]
we define
\begin{multline*}
Q^\sigma = \bigl((z_{\sigma^{-1}(1)}, \theta_{\sigma^{-1}(1)}), ...,
(z_{\sigma^{-1}(n-1)}, \theta_{\sigma^{-1}(n-1)}); (A^{(0)}, M^{(0)}),\\
(\asqrt^{(\sigma^{-1}(1))}, A^{(\sigma^{-1}(1))}, M^{(\sigma^{-1}(1))}),
...,(\asqrt^{(\sigma^{-1}(n-1))}, A^{(\sigma^{-1}(n-1))},
M^{(\sigma^{-1}(n-1))})\\ 
(\asqrt^{(n)}, A^{(n)}, M^{(n)})\bigr) .  
\end{multline*}

To extend this to a right action of $S_n$ on $SK(n)$, we first note that
$S_n$ is generated by the symmetric group on the first $n-1$ letters
$S_{n - 1}$ and the transposition $(n-1 \; n)$.  We can let $(n-1 \;
n)$ act on $SK(n)$ by permuting the $(n- 1)$-st and $n$-th punctures
and their local coordinates for a canonical supersphere with $1 +
n$ tubes but the resulting supersphere with $1+n$ tubes is not
canonical.  To obtain the superconformally equivalent canonical
supersphere, we have to translate the new $n$-th puncture to 0.  This
translation will not change the local coordinates at positively
oriented punctures but will change the local coordinates at the
negatively oriented puncture $\infty$.  This translation is given by
\begin{eqnarray*}
T_\sou : \mbox{$\bigwedge_\infty$} & \longrightarrow & \mbox{$\bigwedge_\infty$} \\
(w,\rho) & \mapsto & (w - z_{n-1} - \rho\theta_{n-1}, \rho - \theta_{n-1}), 
\end{eqnarray*}
and thus by (\ref{T1}), we have $T :
S\hat{\mathbb{C}}
\rightarrow S\hat{\mathbb{C}}$ given by
\begin{eqnarray*}
T(p) = \left\{
  \begin{array}{ll} 
      \sou^{-1} \circ T_\sou \circ \sou (p) & \mbox{if $p \in
           U_\sou$}, \\  
      \nor^{-1} \circ T_\nor \circ \nor (p) & \mbox{if $p \in
           U_\nor \smallsetminus \nor^{-1}( \{(\frac{1}{z_{n-1}})_B \} \times
             (\bigwedge_\infty)_S )$ }  
\end{array} \right.
\end{eqnarray*}
where 
\[T_\nor = \left(\frac{w}{1 - wz_{n-1}} + \frac{i \rho \theta_{n-1}w}{(1 -
wz_{n-1})^2}, \frac{i\theta_{n-1}w}{1 - wz_{n-1}} + \frac{\rho}{1 -
wz_{n-1}} \right) . \]
The new local coordinate at infinity can be written as $\tilde{E}(\tilde{A}^{(0)}, -i
\tilde{M}^{(0)})  (1/w,i\rho/w)$, and it is determined by the old
local coordinate at infinity and the superprojective transformation $T$ via
\[\tilde{E}(\tilde{A}^{(0)}, -i \tilde{M}^{(0)}) 
\Bigl(\frac{1}{w},\frac{i\rho}{w}\Bigr) = \tilde{E}(A^{(0)}, -i M^{(0)}) \circ 
I \circ T_\sou^{-1} (w,\rho) .\]
Using Proposition \ref{Switch2}, we have
\begin{eqnarray}
& & \hspace{-.2in} \tilde{E}(\tilde{A}^{(0)}, -i \tilde{M}^{(0)}) 
\Bigl(\frac{1}{w},\frac{i\rho}{w}\Bigr) \label{Sn action}\\
\hspace{.3in} &=& \exp \Biggl(\sum_{j \in \Z} \left(\tilde{A}^{(0)}_j
L_{-j}(w,\rho) + \tilde{M}^{(0)}_{j - \frac{1}{2}} G_{-j + \frac{1}{2}}
(w,\rho) \right) \!\Biggr) \! \cdot \! \Bigl(\frac{1}{w},\frac{i\rho}{w}\Bigr)  \nonumber\\
&=&  \exp \Biggl(\sum_{j \in \Z} \left(A^{(0)}_j L_{-j}(x,\varphi) + 
M^{(0)}_{j - \frac{1}{2}} G_{-j + \frac{1}{2}}(x,\varphi)
\right) \! \Biggr) \! \cdot \nonumber \\
& & \biggl. \hspace{1.9in} \cdot \Bigl(\frac{1}{x},\frac{i\varphi}{x}\Bigr) 
\biggr|_{(x,\varphi) = (w + z_{n-1} + \rho \theta_{n-1}, \rho + \theta_{n-1})} \nonumber \\
&=& \exp \left(z_{n-1} \frac{\partial}{\partial w} + \theta_{n-1}
\Bigl(\frac{\partial}{\partial \rho} - \rho \frac{\partial}{\partial
w}\Bigr) \right) \cdot \nonumber\\
& & \hspace{.6in} \cdot \exp \Biggl(\sum_{j \in \Z} \left(A^{(0)}_j L_{-j}
(w,\rho) + M^{(0)}_{j - \frac{1}{2}} G_{-j + \frac{1}{2}}(w,\rho)
\right)\Biggr) \! \cdot \! \Bigl(\frac{1}{w},\frac{i\rho}{w}\Bigr)  \nonumber\\
&=& \exp \left(-z_{n-1} L_{-1}(w,\rho) - \theta_{n-1}
G_{-\frac{1}{2}}(w,\rho) \right) \cdot \nonumber\\
& & \hspace{.6in} \cdot \exp \Biggl(\sum_{j \in \Z} \left(A^{(0)}_j L_{-j}(w,\rho) + 
M^{(0)}_{j - \frac{1}{2}} G_{-j + \frac{1}{2}}(w,\rho)
\right)\Biggr) \! \cdot \! \Bigl(\frac{1}{w},\frac{i\rho}{w}\Bigr) . \nonumber
\end{eqnarray} 
Then 
\begin{eqnarray*}
& & \hspace{-.4in} Q^{(n-1 \; n)} \\
&=& \! \! \Bigl( \infty,(z_1, \theta_1),...,(z_{n-2}, \theta_{n-2}), 0, (z_{n-1},
\theta_{n-1}); (A^{(0)},M^{(0)}), (\asqrt^{(1)}, A^{(1)},
M^{(1)}), \Bigr.\\
& & \hspace{.15in} \Bigl. ..., (\asqrt^{(n-2)}, A^{(n-2)}, M^{(n-2)}), (\asqrt^{(n)}, A^{(n)},
M^{(n)}), (\asqrt^{(n-1)}, A^{(n-1)}, M^{(n-1)}) \Bigr) \\ 
&=& \! \! \Bigl(\infty, (z_1 - z_{n-1} - \theta_1\theta_{n-1}, \theta_1 -
\theta_{n-1}),(z_2 - z_{n-1} - \theta_2\theta_{n-1}, \theta_2 -
\theta_{n-1}),..., \Bigr.\\
& & \quad (z_{n-2} - z_{n-1} - \theta_{n-2}\theta_{n-1}, \theta_{n-2} -
\theta_{n-1}), (-z_{n-1},-\theta_{n-1}),0 ; (\tilde{A}^{(0)},
\tilde{M}^{(0)}), \\ 
& & \quad \Bigl. (\asqrt^{(1)}, A^{(1)}, M^{(1)}), ..., (\asqrt^{(n-2)}, A^{(n-2)}, M^{(n-2)}),
(\asqrt^{(n)}, A^{(n)}, M^{(n)}), \\
& & \hspace{2.3in} (\asqrt^{(n-1)}, A^{(n-1)}, M^{(n-1)}) \Bigr) \in SK(n) . 
\end{eqnarray*}     
Thus we have an action of $S_n$ on $SK(n)$.

\begin{rema}\label{operad remark} 
The moduli space of $N=1$ superspheres with tubes, $SK$, along with the sewing 
operation and the action of the symmetric group defined above, is an example
of a partial operad (cf. \cite{M}, \cite{HL2}, \cite{HL3}, \cite{H book}).
\end{rema}

\section[Supermeromorphic tangent spaces of $SK$]{Supermeromorphic superfunctions on $SK$ and
supermeromorphic tangent spaces of $SK$}

A {\it supermeromorphic superfunction on $SK(n)$}, for $n \in \Z$, is a
superfunction $F : SK(n) \rightarrow \bigwedge_\infty$ of the form
\begin{eqnarray}
\hspace{.4in} F(Q) \! \! &=& \! \! F\bigl((z_1, \theta_1),...,(z_{n-1}, 
\theta_{n-1}); (A^{(0)},M^{(0)}), (\asqrt^{(1)}, A^{(1)}, M^{(1)}), ..., 
\label{meromorphic} \bigr.\\
& & \bigl. \hspace{2.6in}(\asqrt^{(n)}, A^{(n)}, M^{(n)})\bigr) \nonumber \\
&=& \! \! F_0\bigl( (z_1,\theta_1),...,(z_{n-1}, \theta_{n-1}); (A^{(0)},M^{(0)}), 
(\asqrt^{(1)}, A^{(1)}, M^{(1)}),..., \bigr.\nonumber \\
& & \bigl. (\asqrt^{(n)},A^{(n)}, M^{(n)})\bigr) \times \biggl(\prod^{n-1}_{i
= 1} z_i^{-s_i} \prod_{1 \leq i<j \leq n-1} (z_i - z_j - \theta_i \theta_j)^{-s_{ij}}\biggr)
\nonumber 
\end{eqnarray} 
where $s_i$ and $s_{ij}$ are nonnegative integers and 
\[F_0\bigl((z_1,\theta_1),...,(z_{n-1}, \theta_{n-1}); (A^{(0)},M^{(0)}),
(\asqrt^{(1)}, A^{(1)}, M^{(1)}),...,(\asqrt^{(n)},A^{(n)}, M^{(n)})\bigr)\]  
is a polynomial in the $z_i$'s, $\theta_i$'s, $\asqrt^{(i)}$'s,
$(\asqrt^{(i)})^{-1}$'s, $A^{(i)}_j$'s, and $M^{(i)}_{j - 1/2}$'s.
For $n=0$ a {\it supermeromorphic superfunction on $SK(0)$} is a
polynomial in the components of elements of $SK(0)$, i.e., a polynomial
in the $A^{(0)}_j$'s, and $M^{(0)}_{j - 1/2}$'s.  For $F$ of the form
(\ref{meromorphic}), we say that $F$ has a pole of order $s_{ij}$ at 
$(z_i,\theta_i) = (z_j,\theta_j)$.

Since $(z_i - z_j - \theta_i \theta_j)^{- s_{ij}} = (z_i - z_j)^{-
s_{ij}} + s_{ij} \theta_i \theta_j (z_i - z_j)^{- s_{ij} - 1}$, we
can expand (\ref{meromorphic}) to be of the form 
\begin{multline}\label{meromorphic2}
\tilde{F}_0\bigl( (z_1,\theta_1),...,(z_{n-1}, \theta_{n-1}); (A^{(0)},M^{(0)}), 
(\asqrt^{(1)}, A^{(1)}, M^{(1)}),..., \bigr.  \\
\bigl. (\asqrt^{(n)},A^{(n)}, M^{(n)})\bigr) \times
\biggl(\prod^{n-1}_{i = 1} z_i^{-s_i} \prod_{1 \leq i<j \leq n-1} 
(z_i - z_j)^{-s_{ij} -1}\biggr)
\end{multline}
where $\tilde{F_0}$ is a polynomial in the $z_i$'s, $\theta_i$'s, 
$\asqrt^{(i)}$'s, $(\asqrt^{(i)})^{-1}$'s, $A^{(i)}_j$'s, and 
$M^{(i)}_{j - 1/2}$'s.  Also note that since the meaning of the 
expression $(z_i - z_j - \theta_i \theta_j)^{- s_{ij}}$ is the 
expansion about the body $((z_i)_B - (z_j)_B)^{- s_{ij}}$ in powers 
of the soul, in general this expression has an infinite number of 
negative powers of $((z_i)_B - (z_j)_B)$.  However, if we let 
$\{\zeta_1, \zeta_2,...\}$ be the fixed basis for the underlying 
vector space of $\bigwedge_\infty$, we see that each $\zeta_{i_1} 
\zeta_{i_2} \cdots \zeta_{i_{2l}}$ term  in $(z_i - z_j - \theta_i 
\theta_j)^{- s_{ij}}$ has a finite number of negative powers of 
$((z_i)_B - (z_j)_B)$, and $(z_i - z_j - \theta_i \theta_j)^{s_{ij}} 
(z_i - z_j - \theta_i \theta_j)^{- s_{ij}} = 1$.  

The set of all canonical supermeromorphic superfunctions on $SK(n)$ is
denoted by $SD(n)$.  Note that $SD(n)$ has a natural $\mathbb{Z}_2$-grading
since any supermeromorphic superfunction $F$ can be decomposed into an
even component $F^0$ and an odd component $F^1$ where $F^0(Q) \in
\bigwedge_\infty^0$, for all $Q \in SK$ and $F^1(Q) \in
\bigwedge_\infty^1$ for all $Q \in SK$. 

An {\it even supermeromorphic tangent vector} of $SK$ at 
a point $Q \in SK$ is a linear map 
\begin{equation}
\mathcal{L}_Q : SD(n) \longrightarrow \mbox{$\bigwedge_\infty$}
\end{equation}
such that for $F_1, F_2 \in SD(n)$
\begin{equation}
\mathcal{L}_Q (F_1 F_2) = \mathcal{L}_Q F_1 \cdot F_2(Q) +  F_1(Q) \cdot
\mathcal{L}_Q F_2 ,   
\end{equation}
and an {\it odd supermeromorphic tangent vector} of $SK$ at a point $Q
\in SK$ is a linear map 
\begin{equation}
\mathcal{L}_Q : SD(n) \longrightarrow \mbox{$\bigwedge_\infty$}
\end{equation}
such that for $F_1$ of homogeneous sign in $SD(n)$,
\begin{equation}
\mathcal{L}_Q (F_1 F_2) = \mathcal{L}_Q F_1 \cdot F_2(Q) + (-1)^{ \eta(F_1)}
F_1(Q) \cdot \mathcal{L}_Q F_2 .   
\end{equation}
The set of all supermeromorphic tangent vectors at $Q$ is the {\it
supermeromorphic tangent space of $SK$ at $Q$} and is denoted $T_QSK$
(or $T_QSK(n)$ when $Q \in SK(n)$).  Note that $T_QSK$ is naturally
$\mathbb{Z}_2$-graded, and thus we can define the sign function $\eta$
on $T_QSK$ (see Section 2.1).  It is obvious that 
\begin{eqnarray*}
\Bigl. \frac{\partial}{\partial z_i}\Bigr|_Q, \quad
\Bigl. \frac{\partial}{\partial \theta_i}\Bigr|_Q, && \mbox{for $i = 1,...,n-1$,}\\
\biggl.\frac{\partial}{\partial \asqrt^{(i)}}\biggr|_Q, && \mbox{for $i = 1,...,n$, and}\\
\biggl. \frac{\partial}{\partial A_j^{(i)}}\biggr|_Q,  \quad 
\biggl.\frac{\partial}{\partial M_{j - \frac{1}{2}}^{(i)}}\biggr|_Q,  && \mbox{for $i = 0,...,n$
and $j \in \mathbb{Z}_+$}, 
\end{eqnarray*}
are in $T_QSK(n)$.  Because any supermeromorphic function on $SK(n)$ depends on
only finitely many of the variables $z_i$, $\theta_i$, for $i = 1,...,n-1$, 
$\asqrt^{(i)}$, for $i = 1,...,n$, and $A^{(i)}_j$, and $M^{(i)}_{j - 1/2}$,
for $i = 0,...,n$ and $j \in \mathbb{Z}_+$, the proof of the following 
proposition is the same as that for the tangent spaces of a finite-dimensional 
manifold. 

\begin{prop}
Any supermeromorphic tangent vector $\mathcal{L}_Q$ can be expressed as 
\begin{multline}\label{tangent}
\mathcal{L}_Q = \sum_{i = 1}^{n-1} \biggl( \Bigl. \alpha_i
\frac{\partial}{\partial z_i} \Bigr|_Q + \beta_i \Bigl.
\frac{\partial}{\partial \theta_i} \Bigr|_Q \biggr) + \sum_{i =
1}^n \delta_i \biggl. \frac{\partial}{\partial \asqrt^{(i)}}
\biggr|_Q \\
+ \; \sum_{i = 0}^n \sum_{j \in \Z}
\Biggl(\gamma_i^j \biggl. \frac{\partial}{\partial A_j^{(i)}} \biggr|_Q
+ \nu_i^j \biggl. \frac{\partial}{\partial M_{j - \frac{1}{2}}^{(i)}}
\biggr|_Q \Biggr) 
\end{multline}  
where
\[\alpha_i = \mathcal{L}_Q z_i, \quad \beta_i = \mathcal{L}_Q \theta_i,
\quad \delta_i = \mathcal{L}_Q \asqrt^{(i)}, \quad \gamma_i^j = \mathcal{L}_Q
A_j^{(i)}, \quad \nu_i^j = \mathcal{L}_Q M_{j - \frac{1}{2}}^{(i)} . \]
If $\eta(\mathcal{L}_Q) = 0$, then $\alpha_i, \delta_i, \gamma_i^j \in
\bigwedge_\infty^0$ and $\beta_i, \nu_i^j \in \bigwedge_\infty^1$.  If
$\eta(\mathcal{L}_Q) = 1$, then $\alpha_i, \delta_i, \gamma_i^j \in 
\bigwedge_\infty^1$ and $\beta_i, \nu_i^j \in \bigwedge_\infty^0$.  
\end{prop}

Note that in (\ref{tangent}) we have an infinite sum.  This is
well defined since when acting on any element of $SD(n)$ it becomes a
finite sum.

\section[Sewings for superspheres with one, two, and three tubes]{The sewing 
operation and superspheres with one, two, and three tubes}

\begin{prop}\label{how to get any supersphere}
Any element $Q \in SK$ can be obtained by sewing the following types of
elements of $SK$ \\
(i) $(\mathbf{0}) \in SK(0)$, \\
(ii) $ \bigl((A^{(0)},M^{(0)}), (\asqrt^{(1)}, A^{(1)},M^{(1)})\bigr) \in SK(1)$, \\
(iii) $\bigl((z, \theta);\mathbf{0}, (1, \mathbf{0}), (1, \mathbf{0})\bigr) \in
SK(2)$. \\   
\end{prop}

\begin{proof} The proof for $Q \in SK(n)$ is by induction on $n$.  If 
$Q = (A^{(0)},M^{(0)}) \in SK(0)$, then
\[Q = \bigl((A^{(0)},M^{(0)}), (\asqrt^{(1)}, A^{(1)},M^{(1)})\bigr) \; _1\infty_0
\; (\mathbf{0})\]
is a type (ii) supersphere sewn with a type (i) supersphere.

If $Q \in SK(1)$, then it is type (ii).

If $Q\in SK(n)$ for $n \geq 2$ is given by
\begin{multline*}
Q = \bigl((z_1, \theta_1),...,(z_{n-1}, \theta_{n-1}); (A^{(0)},M^{(0)}),
(\asqrt^{(1)}, A^{(1)}, M^{(1)}),..., \bigr. \\
\bigl. (\asqrt^{(n)}, A^{(n)}, M^{(n)})\bigr) , 
\end{multline*}
then
\begin{multline*}
Q = \bigl((A^{(0)},M^{(0)}), (1,\mathbf{0})\bigr) \; _1\infty_0 \;
\biggl(\biggl(\cdots \biggl( \Bigl(\bigl((z_1, \theta_1),...,(z_{n-1}, \theta_{n-1}); 
\mathbf{0},(1,\mathbf{0}),...,  \\
 (1,\mathbf{0})\bigr) \; _1\infty_0 \; 
\bigl(\mathbf{0}, (\asqrt^{(1)}, A^{(1)}, M^{(1)}) \bigr) \Bigr) \; 
_2\infty_0 \bigl(\mathbf{0}, (\asqrt^{(2)}, A^{(2)}, M^{(2)}) \bigr) \biggr) 
\; _3\infty_0  \\
 \cdots \biggr) \; _n\infty_0 \; \bigl(\mathbf{0},
(\asqrt^{(n)}, A^{(n)}, M^{(n)}) \bigr)\biggr) 
\end{multline*}
which consists of $n + 1$ type (ii) superspheres and the supersphere  
\[Q'_n = \bigl((z_1, \theta_1),...,(z_{n-1}, \theta_{n-1});
\mathbf{0},(1,\mbox{\bf 0}),...,(1,\mathbf{0})\bigr)\] 
sewn together.  Thus we need only prove that $Q'_n$ can be obtained
{}from type (i), (ii), and (iii) superspheres for $n>2$.  Assume this is
true for $k<n$.  Let $\gamma$ be a closed Jordan curve on the body of
$Q'_n$, denoted $(Q'_n)_B$, such that 0 and $(z_{n-1})_B$ are in one 
connected component of $(Q'_n)_B \smallsetminus \gamma$ (which we will 
call the interior of $\gamma$) and $(z_1)_B,...,(z_{n-2})_B$ are in 
the other connected component of $(Q'_n)_B \smallsetminus \gamma$ 
(which we will call the exterior of $\gamma$).  By the Riemann mapping 
theorem there exists a conformal map $f_B$ {}from the interior of 
$\gamma$ to the open unit disc in $\mathbb{C}$, and we can require that 
$f_B (0) = 0$.  Then by expanding $f_B(z_B)$ about 0 and choosing
$\asqrt \in \mathbb{C}^\times$ such that $\asqrt^2 = f_B'(0)$, we can
let $H(x,\varphi)$ be the unique formally superconformal power series 
such that $\varphi H(x,\varphi) = \varphi f_B(x)$ is of the form 
(\ref{Ehat}) with even coefficient of $\varphi$ equal to $\asqrt$.  
Then $H^{-1}(x,\varphi)$ and $s_{(z_{n-1},0)} \circ H^{-1}(x + f_B (z_{n-1}),
\varphi)$ also satisfy (\ref{Ehat}) for $s_{(z_{n-1},0)} (x,\varphi) = (x - z_{n-1}, 
\varphi)$.

Let 
\begin{eqnarray*}
(\asqrt, A, \mathbf{0}) &=& \hat{E}^{-1}(H(x,\varphi)) , \\
(\asqrt', A', \mathbf{0}) &=& \hat{E}^{-1}(H^{-1}(x,\varphi)) = \bigl( \asqrt^{-1},
\{-\asqrt^{-2j} A_j, \}_{j \in \mathbb{Z}_+}, \mathbf{0} \bigr) \\ 
(\bsqrt, B, \mathbf{0}) &=& \hat{E}^{-1}\bigl(s_{(z_{n-1},0)}  \circ H^{-1}(x + f_B 
(z_{n-1}),\varphi)\bigr) .
\end{eqnarray*}
Then {}from the definition of the sewing operation, we have
\begin{eqnarray*}
Q_n' \! &=& \! \! \bigl((z_1, \theta_1),...,(z_{n-2}, \theta_{n-2});
\mathbf{0},(1,\mathbf{0}),...,(1,\mathbf{0}),(\asqrt, A, \mathbf{0}) 
\bigr) \; _{n-1}\infty_0 \; \\
& & \hspace{1.8in} \bigl(H(z_{n-1}, \theta_{n-1}); \mathbf{0},
(\bsqrt, B, \mathbf{0}), (\asqrt', A', \mathbf{0}) \bigr) \\
&=& \! \! \Bigl( \bigl((z_1, \theta_1),...,(z_{n-2}, \theta_{n-2});
\mathbf{0}, (1,\mathbf{0}),...,(1,\mathbf{0}) \bigr) \; _{n-1}\infty_0 \;
\bigl(\mathbf{0},(\asqrt, A, \mathbf{0}) \bigr)\Bigr)  \\ 
& & \hspace{1.3in} \; _{n-1}\infty_0 \; \bigl(H (z_{n-1}, \theta_{n-1}); \mathbf{0},
(\bsqrt, B, \mathbf{0}), (\asqrt', A', \mathbf{0}) \bigr) 
\end{eqnarray*}
which is obtained {}from sewing $Q_{n-1}'$, a supersphere of type (ii),
and a supersphere in $SK(2)$.  By our inductive assumption, the result
follows.  \end{proof} 

It is clear {}from the definition of sewing that $SK(1)$ is a partial monoid.
In the following proposition, we give some subgroups of $SK(1)$.

\begin{prop}\label{t and s composition}
Let $s,t \in \bigwedge_\infty^0$, and $(A,M) \in \bigwedge_\infty^\infty$.
Assume 
\[s(A,M), t(A,M) \in \mathcal{H}.\]  
Then $(s + t)(A,M) \in \mathcal{H}$, both
\[\bigl(\mathbf{0}, (1,t(A,M)\bigr) \; _1\infty_0 \; \bigl(\mathbf{0}, (1,s(A,M)\bigr)\]
and
\[ \bigl(s(A,M),(1, \mbox{\bf 0})\bigr) \; _1\infty_0 \; \bigl(t(A,M),(1, \mathbf{0})\bigr) \]
exist, and we have
\begin{eqnarray}
\bigl(\mathbf{0}, (1,(s + t)(A,M))\bigr) \! &=& \! \bigl(\mathbf{0}, (1,t(A,M)\bigr) \;
_1\infty_0 \; \bigl(\mathbf{0}, (1,s(A,M)\bigr), \label{sequencesew1} \\
\bigl((s + t)(A,M), (1, \mathbf{0})\bigr) \! &=& \! \bigl(s(A,M),(1, \mbox{\bf
0})\bigr) \; _1\infty_0 \; \bigl(t(A,M),(1, \mathbf{0})\bigr) . \label{sequencesew2} 
\end{eqnarray}
In particular, for $(A,M) \in \mathcal{H}$, the sets 
\begin{equation}\label{first subgroup}
\bigl\{\bigl(\mathbf{0},(1,t(A,M)\bigr) \; | \; t \in \mbox{$\bigwedge_\infty^0$},
\;  t(A,M) \in \mathcal{H} \bigr\}
\end{equation} 
and
\begin{equation}\label{second subgroup}
\bigl\{\bigl(t(A,M), (1,\mathbf{0})\bigr) \; | \; t \in \mbox{$\bigwedge_\infty^0$}, \; 
t(A,M) \in \mathcal{H} \bigr\}
\end{equation} 
are subgroups of $SK(1)$.  In addition we have the subgroups given by
taking $A = \mathbf{0}$ or $M = \mathbf{0}$ in (\ref{first subgroup})
and (\ref{second subgroup}). 
\end{prop}

\begin{proof} Recall that in Chapter 3, we defined a group operation on 
infinite sequences in a superalgebra $R$; see 
(\ref{sequencecompositiondef}).  Letting $R = \bigwedge_\infty$, and 
for $(A,M), (\overline{A}, \overline{M}) \in \bigwedge_\infty^\infty$ 
this operation is given by
\begin{equation}\label{group operation}
(A,M) \circ (\overline{A}, \overline{M}) = \tilde{E}^{-1}\bigl(\tilde{E} 
(\overline{A}, \overline{M}) \circ \tilde{E}(A,M) (x,\varphi)\bigr) .
\end{equation}
 
Let
\[H_s(x, \varphi) = \tilde{E}(sA,sM), \quad \mbox{and} \quad H_t(x,
\varphi) = \tilde{E}(tA,tM) . \]
Since $s(A,M)$, and $t(A,M)$ are in $\mathcal{H}$, the superfunctions
$H_s(w, \rho)$ and $H_t(w, \rho)$ are convergent in a neighborhood of
$0 \in \bigwedge_\infty$.  Hence $H_s \circ H_t (w, \rho)$ and $H_t
\circ H_s (w, \rho)$ are convergent in a neighborhood of $0$.  Thus
by (\ref{group operation}), we see that $s(A,M) \circ
t(A,M)$ and $t(A,M) \circ s(A,M)$ are in $\mathcal{H}$.  By Proposition
\ref{sequencecompositionprop}, 
\[(s + t)(A,M) = s(A,M) \circ t(A,M) = t(A,M) \circ s(A,M) .\]  
Thus $(s + t)(A,M) \in \mathcal{H}$. {}From the definition of the sewing
operation, we see that
\[\bigl(\mathbf{0}, (1,s(A,M) \circ t(A,M))\bigr) = \bigl(\mathbf{0},
(1,t(A,M)\bigr) \; _1\infty_0 \; \bigl(\mathbf{0}, (1,s(A,M)\bigr), \]
and 
\[\bigl(s(A,M) \circ t(A,M) , (1, \mathbf{0})\bigr) = \bigl(s(A,M),(1, \mbox{\bf
0})\bigr) \; _1\infty_0 \; \bigl(t(A,M),(1, \mathbf{0})\bigr) . \]
This then gives (\ref{sequencesew1}) and (\ref{sequencesew2}).
\end{proof}

\begin{prop}\label{for comm and assoc}
Let $z_1,z_2 \in \bigwedge_\infty^0$ such that $|(z_1)_B| > |(z_2)_B| >
|(z_1)_B - (z_2)_B| > 0$.  Then
\begin{eqnarray} 
& & \hspace{-.4in} \bigl((z_1,\theta_1), (z_2,\theta_2); 
\mathbf{0}, (1,\mathbf{0}), (1,\mathbf{0}), (1,\mathbf{0})\bigr)  = \nonumber \\
\hspace{.3in} &=& \! \! \bigl((z_2,\theta_2); \mathbf{0}, (1,\mathbf{0}), 
(1,\mathbf{0})\bigr) \; _1\infty_0 \; \bigl((z_1 - z_2 -\theta_1 \theta_2 , \theta_1 - 
\theta_2);  \bigr.\label{associativity1}\\  
& & \bigl.\hspace{3in} \mathbf{0}, (1,\mathbf{0}), (1,\mathbf{0})\bigr) \nonumber \\
&=& \! \! \bigl((z_1,\theta_1); \mathbf{0}, (1,\mathbf{0}), (1,\mathbf{0})\bigr) \;
_2\infty_0 \; \bigl((z_2, \theta_2); \mathbf{0}, (1,\mathbf{0}),
(1,\mathbf{0})\bigr) .\label{associativity2}
\end{eqnarray}
\end{prop}

\begin{proof} Because each supersphere involved in the proposition
has standard local coordinates at each puncture, the sewings depend only 
on the bodies of the superspheres.  The sewing of the bodies is exactly
the sewing given in Proposition 2.5.3 in \cite{H book}.  Thus as in 
Proposition 2.5.3 of \cite{H book} there exist $r_1,...,r_9 \in 
\mathbb{R}_+$ which allow the above two sewings including the soul
portions of the superspheres.  Then (\ref{associativity1}) and 
(\ref{associativity2}) follow {}from the definition of the sewing operation 
on $SK(2)$. \end{proof}  

\begin{rema} Proposition \ref{for comm and assoc} is used in 
\cite{B thesis} to prove the associativity property for an
$N=1$ Neveu-Schwarz vertex operator superalgebra with odd formal 
variables obtained {}from an $N=1$ supergeometric vertex operator 
superalgebra.  Thus Proposition \ref{for comm and assoc} can be thought of 
as a geometric version of this algebraic relation.  
\end{rema}

Let $(a,m) \in \bigwedge_\infty^\infty$, and $k,l \in \Z$, and define
\begin{multline*}
\bigl(A(a,k), M(m,l - 1/2)\bigr) =  \Bigl( \bigl\{A_j \; | \; 
A_k = a, \; A_j = 0, \; \mbox{for} \; j \neq k \bigr\}_{j \in \Z} , \Bigr. \\
\Bigl. \bigl\{M_{j - \frac{1}{2}} \; | \; M_{l - \frac{1}{2}} = m, \; 
M_{j - \frac{1}{2}} = 0, \; \mbox{for} \; j \neq l \bigr\}_{j \in \Z} \Bigr) 
\end{multline*} 
which is an element of $\bigwedge_\infty^\infty$.

\begin{prop}\label{for derivatives}
For $(z,\theta), (z_0, \theta_0) \in \bigwedge_\infty$ such that
$0 < |(z_0)_B| < |z_B|$, we have 
\begin{multline*}
\bigl((z_0 + z +  \theta_0 \theta, \theta + \theta_0); \mathbf{0},
(1,\mathbf{0}), (1,\mathbf{0})\bigr) = \\
\bigl((z, \theta); \mathbf{0}, (1,\mathbf{0}), (1,\mathbf{0}) \bigr) \; 
_1\infty_0 \; \bigl((A(-z_0,1), M(-\theta_0,1/2)), (1,\mathbf{0}) 
\bigr) .   
\end{multline*}
\end{prop} 

\begin{proof} Note that $(A(-z_0,1), M(-\theta_0,1/2)) \in \mathcal{H}^{(0)}$,
and thus 
\[((A(-z_0,1),M(-\theta_0,1/2)), (1,\mathbf{0})) \in SK(1) .\]  
In fact, $((A(-z_0,1),M(-\theta_0,1/2)), (1,\mathbf{0}))$ represents the 
equivalence class of superspheres with $1+1$ tubes with canonical supersphere 
having local coordinates (in terms of the coordinate chart about zero of the 
underlying super-Riemann sphere) at $\infty$ given by 
\[(w,\rho) \mapsto \Bigl(\frac{1}{w+z_0 +\rho \theta_0}, \frac{i(\rho + 
\theta_0)}{w+z_0 +\rho \theta_0} \Bigr) = \Bigl(\frac{1}{w + z_0} - 
\frac{\rho \theta_0}{(w + z_0)^2}, \frac{i(\rho + \theta_0)}{w+z_0} \Bigr),\] 
and at $0$ given by $(w,\rho) \mapsto (w,\rho)$.  This canonical supersphere  
with tubes is superconformally equivalent to the super-Riemann sphere with 
the negatively oriented puncture still at $\infty$, the positively oriented 
puncture at $(z_0,\theta_0)$, and with standard local coordinates at these 
punctures, i.e., the local coordinate at $\infty$ is given by $(w,\rho) 
\mapsto  (1/w,i\rho /w)$ and the local coordinate at $(z_0,\theta_0)$ is 
given by the superconformal shift $s_{(z_0,\theta_0)} : (w,\rho) \mapsto  
(w - z_0 - \rho \theta_0, \rho - \theta_0)$.  Call this non-canonical  
supersphere $C$. Then choosing $r, r_1, r_2 \in \mathbb{R}_+$ such that 
$0 < |(z_0)_B| < r_2 < r < r_1 < |z_B|$, we can sew the $0$-th puncture of 
$C$ to the $1$-st puncture of $((z, \theta); \mathbf{0}, (1,\mathbf{0}), 
(1,\mathbf{0}))$ via the sewing defined in Chapter 2.  The resulting 
supersphere $C'$ is given by 
\[C' = \Bigl((S\hat{\mathbb{C}} \smallsetminus \sou^{-1}(\bar{\mathcal{B}}_{(z,\theta)}^{r_2}))
\sqcup  (S\hat{\mathbb{C}} \smallsetminus \nor^{-1} (\bar{\mathcal{B}}_0^{1/r_1})) \Bigr)/
\sim\]  
where $p \sim q$ if 
\begin{eqnarray*}
\sou(p) = (w_1,\rho_1) &\in& (\bar{\mathcal{B}}_{(z,\theta)}^{r_1} \smallsetminus 
\bar{\mathcal{B}}_{(z,\theta)}^{r_2}) \subset (\mbox{$\bigwedge_\infty$} 
\smallsetminus \bar{\mathcal{B}}_{(z,\theta)}^{r_2})\\
\sou(q) = (w_2,\rho_2) &\in& (\bar{\mathcal{B}}_\infty^{r_2} \smallsetminus 
\bar{\mathcal{B}}_\infty^{r_1}) \subset (\mbox{$\bigwedge_\infty$} 
\smallsetminus \bar{\mathcal{B}}_\infty^{r_1}),
\end{eqnarray*}
i.e.,
\[\nor(q) = \Bigl(\frac{1}{w_2},\frac{i\rho_2}{w_2} \Bigr) \in 
(\bar{\mathcal{B}}_0^{1/r_2} \smallsetminus \bar{\mathcal{B}}_0^{1/r_1}) \subset
(\mbox{$\bigwedge_\infty$} \smallsetminus \bar{\mathcal{B}}_0^{1/r_1}),\]
and $(w_1 - z - \rho_1 \theta, \rho_1 - \theta) = (w_2,\rho_2)$.  Then the
uniformizing function which sends $C'$ to a canonical super-Riemann sphere
must satisfy 
\begin{eqnarray*}
F(p) = \left\{
  \begin{array}{ll}  
    F^{(1)}(p) & \mbox{for $p \in  (S\hat{\mathbb{C}} \smallsetminus
\sou^{-1}(\bar{\mathcal{B}}_{(z,\theta)}^{r_2}))$}, \\ 
    F^{(2)}(p) & \mbox{for $p \in  (S\hat{\mathbb{C}}  
\smallsetminus \nor^{-1} (\bar{\mathcal{B}}_0^{1/r_1}))$},
\end{array} \right.
\end{eqnarray*}
where $F^{(1)}_\sou (w_1,\rho_1) = F^{(2)}_\sou (w_2,\rho_2)$ if $\sou^{-1}
(w_1, \rho_1) \sim \sou^{-1}(w_2,\rho_2)$, i.e., $F^{(1)}_\sou (w,\rho) = 
F^{(2)}_\sou (w - z - \rho \theta, \rho - \theta)$.  In addition, we would 
like $F^{(1)}$ to send $\infty$ to $\infty$, send $0$ to $0$, and keep the 
even coefficient of $w^{-2}$ of the local coordinate at $\infty$ equal to 
zero in order for the resulting supersphere to be canonical.  Thus 
$F^{(1)}_\sou (w,\rho) = (w,\rho)$ and $F^{(2)}_\sou(w,\rho) = (w + z + \rho 
\theta, \rho + \theta)$ is a solution to the uniformizing function and is 
the unique solution sending $C'$ to a canonical supersphere. Therefore the 
canonical supersphere representative of the equivalence class of the sewing 
\[\bigl((z, \theta); \mathbf{0}, (1,\mathbf{0}), (1,\mathbf{0})) \; 
_1\infty_0 \; ((A(-z_0,1), M(-\theta_0,1/2)), (1,\mathbf{0})\bigr)\]
is given by
\[\bigl(S\hat{\mathbb{C}}; \infty, F^{(2)}_\sou(z_0,\theta_0), 0; I, 
s_{(z_0, \theta_0)} \circ (F^{(2)}_\sou)^{-1}, id_{\mbox{\tiny{$\bigwedge_\infty$}}}
\bigr), \]
but this is the canonical supersphere represented by 
\[\bigl((z_0 + z + \theta_0 \theta, \theta + \theta_0); \mathbf{0}, (1,\mathbf{0}),
(1,\mathbf{0})\bigr) \] 
as desired. \end{proof} 

Let $\epsilon \in \bigwedge_\infty^1$. {}From the definition of $SK(0)$,
we see that $(\mathbf{0}, M(\epsilon,3/2)) \in SK(0)$.  Let
$F$ be any element of $SD(1)$, and $(z, \theta) \in \bigwedge_\infty$.
We define a linear functional on $SD(1)$ by 
\begin{equation}\label{define G}
\mathcal{G}_e (z, \theta) F = \Bigl.  \frac{d}{d \epsilon} F
\Bigl( \bigl((z,\theta) ; \mathbf{0}, (1,\mathbf{0}), (1,\mathbf{0})\bigr)
\; _1\infty_0 \; \bigl(\mathbf{0}, M(\epsilon, 3/2) \bigr) \Bigr) 
\Bigr|_{\epsilon = 0} . 
\end{equation}   

\begin{prop}\label{the linear functional G}
The linear functional $\mathcal{G}_e (z, \theta)$ is in $T_e SK(1)$, and 
\begin{multline}\label{tau} 
\mathcal{G}_e (z, \theta) = \sum_{k = 0}^1 \sum_{j \in \Z}
z^{-(2k - 1)j - 2 + k} \Biggl. \frac{\partial}{\partial M_{j -
\frac{1}{2}}^{(k)}} \Biggr|_e  \\
 + \; 2 \theta \Biggl(z^{-2} \Biggl. \frac{1}{2}
\frac{\partial}{\partial \asqrt^{(1)}} \Biggr|_e + \sum_{k = 0}^1 
\sum_{j \in \Z} z^{-(2k - 1)j - 2} \Biggl. \frac{\partial}{\partial
A_j^{(k)}} \Biggr|_e \Biggr) .  
\end{multline}
\end{prop}

\begin{proof} If $\mathcal{G}_e (z, \theta)$ is given by
(\ref{tau}), then $\mathcal{G}_e (z, \theta) \in T_e SK(1)$.  Thus, the
only thing we need to prove is (\ref{tau}).  The canonical supersphere
$(\mathbf{0}, M(\epsilon, 3/2))$ is the super-Riemann sphere
with one puncture at $\infty$ and local coordinate at $\infty$ given
by 
\begin{eqnarray}
& & \hspace{-1.4in} \exp \biggl( - \epsilon w^{-1} \biggl(
\frac{\partial}{\partial \rho} - \rho \frac{\partial}{\partial w}
\biggr) \biggr) \cdot \biggl( \frac{1}{w}, \frac{i \rho}{w} \biggr)
\label{local tau coordinate} \\ 
\hspace{1.5in} &=& \left. \exp\biggl( - i\epsilon x^2 \biggl( \frac{\partial}{\partial
\varphi} - \varphi \frac{\partial}{\partial x} \biggr) \biggr) \cdot
(x, \varphi) \right|_{(x, \varphi) = \left(\frac{1}{w}, \frac{i
\rho}{w} \right)} \nonumber \\ 
&=& \biggl( \frac{1}{w} - \frac{\epsilon \rho}{w^3} , \frac{i \rho}{w}
- \frac{i \epsilon}{w^2} \biggr) .\nonumber
\end{eqnarray} 

Define $A_j^{(0)} (\epsilon), A_j^{(1)} (\epsilon) \in \bigwedge_\infty^0$, 
and $M_{j - \frac{1}{2}}^{(0)} (\epsilon),M_{j - \frac{1}{2}} ^{(1)} (\epsilon)
\in \bigwedge_\infty^1$, for $j \in \mathbb{Z}_+$, and $\asqrt^{(1)}(\epsilon)
\in (\bigwedge_\infty^0)^\times$ by
\begin{equation}\label{find tau at infinity}
\left( \frac{1}{w} + \frac{2 \epsilon \theta}{wz (w - z)} +
\frac{\rho \epsilon}{wz (w - z)} , \frac{-i \epsilon}{w (w - z)} -
\frac{i \epsilon}{wz} \hspace{1.4in} \right. 
\end{equation}
\[ \hspace{2in} \left. + \; \rho \biggl( \frac{i}{w} + \frac{i\epsilon
\theta}{w (w - z)^2} + \frac{2i \epsilon \theta}{wz (w - z)} \biggr)
\right) \]
\[ = \; \exp \Biggl(\sum_{j \in \Z} \biggl( A_j^{(0)} (\epsilon)
L_{-j}(w,\rho)  +  M_{j - \frac{1}{2}} ^{(0)} (\epsilon) G_{-j + \frac{1}{2}}
(w,\rho) \biggr) \! \Biggr) \cdot \biggl(\frac{1}{w}, \frac{i \rho}{w} \biggr) ,\] 
and
\begin{multline}\label{find tau at zero}
\biggl( w - \frac{2 \epsilon \theta}{w - z} - \frac{2 \epsilon
\theta}{z} - \rho \biggl( \frac{\epsilon}{w - z} + \frac{\epsilon}{z}
\biggr) , \frac{-\epsilon}{w - z} - \frac{\epsilon}{z} + \rho \biggl( 1
+ \frac{\epsilon \theta}{(w - z)^2} \biggr) \biggr)  \\
= \; \exp \Biggl( \! - \! \sum_{j \in \Z} \biggl( A_j^{(1)} (\epsilon) L_j(w,\rho)
+ M_{j - \frac{1}{2}}^{(1)} (\epsilon) G_{j - \frac{1}{2}}(w,\rho) \biggr) \! \Biggr)
\cdot \\
\cdot (\asqrt^{(1)}(\epsilon))^{- 2L_0(w,\rho) } \cdot (w, \rho) .
\end{multline}
By definition of the sewing operation, 
\begin{multline}\label{tau sewing}
\bigl((z, \theta); \mathbf{0}, (1,\mathbf{0}), (1,\mathbf{0})\bigr) \;
_1\infty_0 \; \bigl(\mathbf{0}, M(\epsilon, 3/2)\bigr) \\
= \; \bigl((A^{(0)}(\epsilon), M^{(0)}(\epsilon)), (\asqrt^{(1)}(\epsilon),
A^{(1)}(\epsilon), M^{(1)}(\epsilon)) \bigr) .
\end{multline}

Let $F_0, F_j^{(0)}, F_j^{(1)} \in SD^0(0)$, and $F_{j -
1/2}^{(0)}, F_{j - 1/2}^{(1)} \in SD^1(0)$, for $j \in
\Z$, be given by 
\begin{eqnarray*}
F_0 ((A^{(0)},M^{(0)}), (\asqrt^{(1)}, A^{(1)}, M^{(1)}) )  &=& \asqrt^{(1)}
, \\ 
F_j^{(0)} ((A^{(0)},M^{(0)}), (\asqrt^{(1)}, A^{(1)}, M^{(1)}) )  &=&
A_j^{(0)} ,\\ 
F_j^{(1)} ((A^{(0)},M^{(0)}), (\asqrt^{(1)}, A^{(1)}, M^{(1)}) )  &=&
A_j^{(1)} ,\\ 
F_{j - \frac{1}{2}}^{(0)} ((A^{(0)},M^{(0)}), (\asqrt^{(1)}, A^{(1)},
M^{(1)}) )  &=& M_{j - \frac{1}{2}}^{(0)} ,\\ 
F_{j - \frac{1}{2}}^{(1)} ((A^{(0)},M^{(0)}), (\asqrt^{(1)}, A^{(1)},
M^{(1)}) )  &=& M_{j - \frac{1}{2}}^{(1)} .\\
\end{eqnarray*} 
Then
\begin{eqnarray*}
\left. \frac{d}{d \epsilon} F_0 ((A^{(0)},M^{(0)}),
(\asqrt^{(1)}, A^{(1)}, M^{(1)}) )\right|_{\epsilon = 0}
&=& \left. \frac{d}{d \epsilon} \asqrt^{(1)} (\epsilon) \right|_{\epsilon
= 0} ,\\ 
\left. \frac{d}{d \epsilon}F_j^{(0)} ((A^{(0)},M^{(0)}),
(\asqrt^{(1)}, A^{(1)}, M^{(1)}) ) \right|_{\epsilon = 0}
&=& \left. \frac{d}{d \epsilon} A_j^{(0)} (\epsilon) \right|_{\epsilon
= 0} ,\\ 
\left. \frac{d}{d \epsilon}  F_j^{(1)} ((A^{(0)},M^{(0)}),
(\asqrt^{(1)}, A^{(1)}, M^{(1)}) ) \right|_{\epsilon = 0}
&=& \left. \frac{d}{d \epsilon} A_j^{(1)} (\epsilon) \right|_{\epsilon
= 0} ,\\ 
\left. \frac{d}{d \epsilon}  F_{j - \frac{1}{2}}^{(0)}
((A^{(0)},M^{(0)}), (\asqrt^{(1)}, A^{(1)}, M^{(1)}))
\right|_{\epsilon = 0} &=& \left. \frac{d}{d \epsilon} M_{j -
\frac{1}{2}}^{(0)} (\epsilon) \right|_{\epsilon = 0} ,\\ 
\left. \frac{d}{d \epsilon}  F_{j - \frac{1}{2}}^{(1)}
((A^{(0)},M^{(0)}), (\asqrt^{(1)}, A^{(1)}, M^{(1)}))
\right|_{\epsilon = 0} &=& \left. \frac{d}{d \epsilon} M_{j -
\frac{1}{2}}^{(1)} (\epsilon) \right|_{\epsilon = 0} .\\ 
\end{eqnarray*}
{}From (\ref{find tau at infinity}) and (\ref{find tau at zero}), we
have 
\begin{eqnarray*}
\left. \frac{d}{d \epsilon} \asqrt^{(1)} (\epsilon) \right|_{\epsilon =
0} &=&  \theta z^{-2} ,\\ 
\left. \frac{d}{d \epsilon} A_j^{(0)} (\epsilon) \right|_{\epsilon =
0} &=&  2 \theta z^{j - 2} ,\\ 
\left. \frac{d}{d \epsilon} A_j^{(1)} (\epsilon) \right|_{\epsilon =
0} &=&  2 \theta z^{- j - 2} ,\\  
\left. \frac{d}{d \epsilon} M_{j - \frac{1}{2}}^{(0)} (\epsilon)
\right|_{\epsilon = 0} &=&  z^{j - 2} ,\\
\left. \frac{d}{d \epsilon} M_{j - \frac{1}{2}}^{(1)} (\epsilon)
\right|_{\epsilon = 0} &=&  z^{- j - 1}, \\ 
\end{eqnarray*}
for $j \in \Z$. Using (\ref{define G}), (\ref{tau sewing}) and the
above formulas, we obtain (\ref{tau}).  \end{proof} 

\begin{rema} Propositions \ref{for derivatives} and \ref{the
linear functional G} are used in \cite{B thesis} to prove certain axioms
for an $N=1$ Neveu-Schwarz vertex operator superalgebra constructed {}from 
an $N=1$ supergeometric vertex operator superalgebra.  The axioms proved 
using these propositions are: that the vertex operator with formal 
variables associated with the Neveu-Schwarz element $\tau$ is equal to 
the formal series with coefficients satisfying the Neveu-Schwarz algebra
relations; and the $G(-1/2)$-derivative property.  One can thus
derive some motivation for looking at the functional 
$\mathcal{G}_e (z, \theta)$ by noting that the local coordinate 
(\ref{local tau coordinate}) which is given by 
\[\exp \Bigl(\epsilon G_{-\frac{1}{2}}(w,\rho) \Bigr) \cdot 
\biggl( \frac{1}{w}, \frac{i \rho}{w} \biggr) \]
is related to the algebraic aspects of $N=1$ superconformal field theory
via the fact that in an $N=1$ Neveu-Schwarz vertex operator superalgebra 
\cite{B vosas}, the Neveu-Schwarz element $\tau$ is given by $G(-1/2) 
\mathbf{1}$ where $\mathbf{1}$ is the vacuum element and $G(-1/2)$ is 
given by the representation of the $N=1$ Neveu-Schwarz algebra present
in the $N=1$ Neveu-Schwarz vertex operator superalgebra.  In the 
correspondence between the algebra and geometry, the supersphere 
represented by $(\mathbf{0}) \in SK(0)$ corresponds to the vacuum 
element $\mathbf{1}$ in the vertex operator superalgebra, and 
$(\mathbf{0}, M(\epsilon, 3/2))$ corresponds to $e^{\epsilon G(-1/2)}
\mathbf{1}$.  Therefore taking the partial derivative with respect to
$\epsilon$ and setting $\epsilon = 0$, one obtains the Neveu-Schwarz
element. Then the linear functional $\mathcal{G}_e (z, \theta)$ acting on 
the supermeromorphic function corresponding to the one-point 
correlation function gives the correlation function of the vertex 
operator associated to the Neveu-Schwarz element $\tau$, denoted 
$\langle v', Y(\tau,(z,\theta))v \rangle$; see \cite{B thesis}, 
\cite{B vosas}.
\end{rema}

\section{Generalized superspheres with tubes}

{}From Proposition \ref{t and s composition}, we know that for any
$(A,M) \in \mathcal{H}$, the sets
\begin{eqnarray}
\bigl\{ (\mathbf{0}, (1,t(A,M)) \; | \; t \in \mbox{$\bigwedge_\infty^0$},
\; t(A,M) \in \mathcal{H} \bigr\}\\
\bigl\{ (t(A,M), (1,\mathbf{0})) \; | \; t \in \mbox{$\bigwedge_\infty^0$},
\; t(A,M) \in \mathcal{H} \bigr\}
\end{eqnarray}
are subgroups of $SK(1)$, as are the corresponding sets if we restrict 
to $A = \mathbf{0}$ or $M = \mathbf{0}$.  If we consider
\begin{eqnarray}
\bigl\{ (\mathbf{0}, (1,t(A,M)) \; | \; t \in \mbox{$\bigwedge_\infty^0$},
\; (A,M) \in \mbox{$\bigwedge_\infty^\infty$}\bigr\} \label{first set}\\
\bigl\{ (t(A,M), (1,\mathbf{0})) \; | \; t \in \mbox{$\bigwedge_\infty^0$},
\; (A,M) \in \mbox{$\bigwedge_\infty^\infty$}\bigr\} \label{second set}
\end{eqnarray}
then these sets are $(1,0)$-dimensional super Lie groups over 
$\bigwedge_\infty$, i.e., they are abstract groups which are also 
$(1,0)$-dimensional superanalytic supermanifolds with superanalytic 
group multiplication and inverse mappings, (cf. \cite{Ro1}). But elements 
of (\ref{first set}) and (\ref{second set}) are in general not in $SK(1)$, 
even if $(A,M) \in \mathcal{H}$.  Therefore, in order to study these sets, 
we introduce generalized superspheres with tubes.

Let 
\[\tilde{SK}(n) = SM^{n-1} \times \mbox{$\bigwedge^\infty_\infty$} \times
((\mbox{$\bigwedge_\infty^0$})^\times \times \mbox{$\bigwedge^\infty_\infty$})^n , \]
for $n \in \Z$, and 
\[\tilde{SK}(0) = \bigl\{ \left. (A,M) \in \mbox{$\bigwedge^\infty_\infty$} \; 
\right| \;
(A_1,M_{\frac{1}{2}}) = (0,0) \bigr\} .\]
Elements of $\tilde{SK} = \bigcup_{n \in \mathbb{N}}\tilde{SK}(n)$ are
called {\it generalized} (or {\it formal}) {\it superspheres with
tubes}.  Note that $SK \subset \tilde{SK}$. {}From the definition of
supermeromorphic superfunction on $SK$, the following generalization
is obvious:
\begin{prop}
Any supermeromorphic superfunction on $SK(n)$, $n \in \mathbb{N}$, can be
extended to a superfunction on $\tilde{SK}(n)$ having the same form.
\end{prop}

Let $Q_1 \in \tilde{SK}(n)$, for $n \in \mathbb{Z}_+$, be given by 
\[((z_1, \theta_1),...,(z_{n-1}, \theta_{n-1}); (A^{(0)},M^{(0)}),
(\asqrt^{(1)}, A^{(1)}, M^{(1)}),...,(\asqrt^{(n)}, A^{(n)}, M^{(n)})), \]
and let $Q_2 = (\mathbf{0},(\bsqrt^{(1)}, B^{(1)}, N^{(1)}))$ be an element of
$\tilde{SK}(1)$.  Let
\[(B^{(1)}(\asqrt^{(i)}), N^{(1)}(\asqrt^{(i)})) = \left\{(\asqrt^{(i)})^{2j}
B^{(1)}_j, (\asqrt^{(i)})^{2j - 1} N^{(1)}_{j - \frac{1}{2}}
\right\}_{j \in \Z}  \; \in \mbox{$\bigwedge_\infty^\infty$}.\] 
We define
\begin{multline*}
Q_1 \; _i\infty_0 \; Q_2 = \Bigl((z_1, \theta_1),...,(z_{n-1},
\theta_{n-1}); (A^{(0)},M^{(0)}), (\asqrt^{(1)}, A^{(1)}, M^{(1)}),..., \Bigr. \\
(\asqrt^{(i-1)}, A^{(i-1)}, M^{(i-1)}), (\asqrt^{(i)}\bsqrt^{(1)},(A^{(i)},
M^{(i)})\circ (B^{(1)}(\asqrt^{(i)}), N^{(1)}(\asqrt^{(i)}))), \\
\Bigl. (\asqrt^{(i+1)}, A^{(i+1)}, M^{(i+1)}),...,(\asqrt^{(n)}, A^{(n)}, 
M^{(n)})\Bigr) \; \in
\tilde{SK}(n) .
\end{multline*}

Similarly, for $Q_3 = ((B^{(0)},N^{(0)}),(1,\mathbf{0})) \in
\tilde{SK}(1)$, we define 
\begin{multline*}
Q_3 \; _1\infty_0 \; Q_1 = ((z_1, \theta_1),...,(z_{n-1},
\theta_{n-1}); (A^{(0)},M^{(0)}) \circ (B^{(0)},N^{(0)}), \\
(\asqrt^{(1)}, A^{(1)}, M^{(1)}),...,(\asqrt^{(n)}, A^{(n)}, M^{(n)})) \; \in
\tilde{SK}(n) .
\end{multline*}

\begin{prop}
The subsets
\[\tilde{SK}^0(1) = \{ (\mathbf{0},(\asqrt^{(1)}, A^{(1)}, M^{(1)})) \in
\tilde{SK}(1) \} \]
and
\[\tilde{SK}^\infty(1) = \{ ((A^{(0)}, M^{(0)}),(1,\mathbf{0})) \in
\tilde{SK}(1) \} \]
with the sewing operation $\; _1\infty_0 \;$ as defined above are
groups.  If $(A,M) \in \bigwedge^\infty_\infty$ and $\asqrt \in
(\bigwedge^0_\infty)^\times$, then
\begin{eqnarray*}
t &\mapsto& (\mathbf{0}, (\asqrt^t, \mathbf{0}))\\
t &\mapsto& (\mathbf{0}, (1, t(A,M)))
\end{eqnarray*}
give homomorphisms {}from the additive group $\bigwedge_\infty^0$ to
$\tilde{SK}^0(1)$, and 
\[t \mapsto (t(A,M),(1,\mathbf{0}))\]
gives a homomorphism {}from $\bigwedge_\infty^0$ to $\tilde{SK}^\infty(1)$.
In addition, we have the subgroups obtained by setting $\asqrt^{(1)} = 1$,
$A^{(1)} = \mathbf{0}$, $M^{(1)} = \mathbf{0}$, $A^{(0)} = \mathbf{0}$, 
or $M^{(0)} = \mathbf{0}$.
\end{prop}

\begin{proof} {}From the definition of the sewing
operation defined above, $\tilde{SK}^0(1)$ and $\tilde{SK}^\infty(1)$
are closed under $\; _1\infty_0 \;$.  The identity for both
$\tilde{SK}^0(1)$ and $\tilde{SK}^\infty(1)$ is 
\[e = (\mathbf{0},(1,\mathbf{0})) . \]
In $\tilde{SK}^0(1)$, the inverse of $(\mathbf{0},(\asqrt^{(1)}, A^{(1)},
M^{(1)}))$ is 
\[(\mathbf{0},((\asqrt^{(1)})^{-1}, -A^{(1)}((\asqrt^{(1)})^{-1}),
-M^{(1)}((\asqrt^{(1)})^{-1})), \]
and in $\tilde{SK}^\infty(1)$, the inverse of $((A^{(0)}, M^{(0)}), 
(1, \mathbf{0}))$ is $((-A^{(0)}, -M^{(0)}), (1,\mathbf{0}))$.   

{}From the definition of the sewing operation, Proposition 
\ref{sequencecompositionprop} and equation (\ref{composition with a}), we have
\begin{multline*}
\left(\mathbf{0}, (\asqrt^{t_1}, t_1(A,M))) \; _1\infty_0 \; (\mathbf{0},
(\asqrt^{t_2}, t_2(A,M)) \right) = \\
\left(\mathbf{0}, (\asqrt^{t_1}\asqrt^{t_2}, t_1(A,M) \circ t_2(A(\asqrt^{t_1}),
M(\asqrt^{t_1})))\right) . 
\end{multline*}
Therefore
\[\left(\mathbf{0}, (\asqrt^{t_1}, \mathbf{0})) \; _1\infty_0 \; (\mathbf{0},
(\asqrt^{t_2}, \mathbf{0}) \right) = 
\left(\mathbf{0}, (\asqrt^{t_1 + t_2}, \mathbf{0})\right) \]
and 
\[\left(\mathbf{0}, (1, t_1(A,M))) \; _1\infty_0 \; (\mathbf{0},
(1, t_2(A,M)) \right) =  \left(\mathbf{0}, (1, (t_1 + t_2)(A,M))\right). \]
Thus $t \mapsto (\mathbf{0}, (\asqrt^t, \mathbf{0}))$ and $t \mapsto
(\mathbf{0}, (1, t(A,M))$ are homomorphisms {}from $\bigwedge_\infty^0$ to
$\tilde{SK}^0(1)$.  In addition,
\begin{eqnarray*}
\left(t_1(A,M),(1,\mathbf{0})) \; _1\infty_0 \; (t_2(A,M),(1,\mathbf{0})\right) 
&=& \left(t_1(A,M) \circ t_2(A,M),(1,\mathbf{0})\right) \\
&=& \left((t_1 + t_2)(A,M), (1,\mathbf{0})\right) .
\end{eqnarray*}
Thus $t \mapsto (t(A,M),(1,\mathbf{0}))$ is a homomorphism {}from 
$\bigwedge_\infty^0$ to $\tilde{SK}^\infty(1)$. \end{proof}

\begin{rema} The group $\tilde{SK}^0(1)$ is isomorphic to the
group $(\bigwedge_\infty^0)^\times \times \bigwedge_\infty^\infty$
discussed in Remark \ref{group remark} with $R = \bigwedge_\infty$.  The 
group $\tilde{SK}^\infty(1)$ and the subgroup  of $\tilde{SK}^0(1)$
given by $\{(\mathbf{0},(1, A^{(1)},M^{(1)})) \in \tilde{SK}^0(1) \}$
are both isomorphic to the subgroup $\bigwedge_\infty^\infty$ of
$(\bigwedge_\infty^0)^\times \times \bigwedge_\infty^\infty$.  The subset
$(\bigwedge_\infty^0)^\times \times \mathcal{H}$ is a subgroup of
$(\bigwedge_\infty^0)^\times \times \bigwedge_\infty^\infty$ and is
isomorphic to the subgroup of all elements of $\tilde{SK}^0(1)$ with
$(A^{(1)},M^{(1)}) \in \mathcal{H}$.  This subgroup is also isomorphic to
the group of local superconformal transformations vanishing at 0.  
Similarly, the subset $\mathcal{H}$ is a subgroup of 
$\bigwedge_\infty^\infty$ isomorphic to the subgroup of all elements of
$\tilde{SK}^\infty(1)$ with $(A^{(0)}, M^{(0)}) \in \mathcal{H}$.  This
subgroup is isomorphic to the group of local superconformal 
transformations vanishing at $\infty$.
\end{rema}

\section[The sewing formulas and convergence]{The sewing formulas and 
the convergence of associated series via the Fischer-Grauert Theorem}

For $m \in \Z$, let $Q_1 \in SK(m)$ be given by
\[\bigl((z_1, \theta_1),...,(z_{m-1}, \theta_{m-1});
(A^{(0)},M^{(0)}), (\asqrt^{(1)}, A^{(1)}, M^{(1)}),...,(\asqrt^{(m)},
A^{(m)}, M^{(m)})\bigr) \]
and for $n \in \mathbb{N}$, let $Q_2 \in SK(n)$ be given by  
\[\bigl((z_1', \theta_1'),...,(z_{n-1}', \theta_{n-1}');
(B^{(0)},N^{(0)}), (\bsqrt^{(1)}, B^{(1)}, N^{(1)}),...,(\bsqrt^{(n)},
B^{(n)}, N^{(n)})\bigr) .\]
For convenience, we will sometimes denote the puncture at $0$ of $Q_1$
by $(z_m, \theta_m)$. 

Let
\begin{eqnarray*}
H^{(1)} (w, \rho)\!  &=& \! \exp \Biggl(\!  - \! \sum_{j \in \Z} \biggl( A^{(i)}_j
L_j(w,\rho)  +  M^{(i)}_{j - \frac{1}{2}} G_{j - \frac{1}{2}}(w,\rho)
\biggr) \! \Biggr) \cdot \\
& & \hspace{2.4in} \cdot (\asqrt^{(i)})^{-2L_0(w,\rho)}\cdot (w, \rho ) \\
&=& \hat{E}(\asqrt^{(i)}, A^{(i)}, M^{(i)}) (w,\rho), \\
H^{(2)} (w, \rho) \! &=& \! \exp \Biggl(\sum_{j \in \Z} \biggl( B^{(0)}_j
L_{-j}(w,\rho) +  N^{(0)}_{j - \frac{1}{2}} G_{-j + \frac{1}{2}}(w,\rho)
\biggr) \! \Biggr) \cdot \Bigl(\frac{1}{w}, \frac{i \rho}{w} \Bigr) \\
&=& \tilde{E}(B^{(0)},-i N^{(0)}) \Bigl(\frac{1}{w}, \frac{i \rho}{w} \Bigr). 
\end{eqnarray*}
Then the local coordinate vanishing at the $i$-th puncture of the canonical
supersphere represented by $Q_1$ is $H^{(1)} \circ s_{(z_i,\theta_i)}(w,\rho)$,
and the local coordinate vanishing at the puncture at $\infty$ of the 
canonical supersphere represented by $Q_1$ is $H^{(2)}(w,\rho)$.

Recall the standard coordinate atlas for the super-Riemann sphere, 
$S\hat{\mathbb{C}}$, given by $\{(U_\sou,\sou),(U_\nor,\nor)\}$ with coordinate
transition given by $I = \sou \circ \nor^{-1}$.

\begin{prop}\label{actual sewing}
The i-th tube of $Q_1$ can be sewn with the 0-th tube of $Q_2$ if and
only if there exist $r_1, r_2 \in \mathbb{R}_+$, with $r_1 > r_2$ such 
that the series $(H^{(1)})^{-1}(w,\rho)$ and $(H^{(2)})^{-1}(w,\rho)$
are absolutely convergent and single-valued in $\mathcal{B}_0^{r_1}$ and
$\mathcal{B}_0^{1/r_2} \smallsetminus (\{0 \} \times (\bigwedge_\infty)_S)$,
respectively, 
\[(-z_i, -\theta_i), (z_k - z_i - \theta_k \theta_i, \theta_k -
\theta_i) \notin (H^{(1)})^{-1} (\mathcal{B}_0^{r_1}), \]
for $k = 1,..., m - 1$, $k \neq i$, and
\[0, (z_l', \theta_l') \notin (H^{(2)})^{-1}(\mathcal{B}_0^{1/r_2}
\smallsetminus (\{0 \} \times (\mbox{$\bigwedge_\infty$})_S)),\] 
for $l = 1,...,n-1$.  Moreover, in this case, there exist unique 
bijective superconformal functions $F^{(1)} (w,\rho)$ and $F^{(2)}
(w,\rho)$ defined on 
\[S\hat{\mathbb{C}} \smallsetminus \sou^{-1} \circ s_{(z_i,\theta_i)}^{-1} 
\circ (H^{(1)})^{-1} (\bar{\mathcal{B}}_0^{r_2})\] 
and 
\[U_\sou \smallsetminus \sou^{-1} \circ (H^{(2)})^{-1}
(\bar{\mathcal{B}}_0^{1/r_1} \smallsetminus (\{0\} \times 
\mbox{$(\bigwedge_\infty)_S$})) \subset S\hat{\mathbb{C}},\]  
respectively, satisfying the normalization conditions:
\begin{eqnarray}
F^{(1)}_\nor  (0) &=& 0, \label{normalize 1} \\
\lim_{w \rightarrow \infty} \frac{\partial}{\partial \rho} (F^{(1)}_\sou)^1
(w,\rho) &=& 1 \label{normalize 2} \\
F^{(2)}_\sou (0) &=& 0 \label{normalize 3} 
\end{eqnarray}
and such that in $(H^{(1)} \circ s_{(z_i, \theta_i)})^{-1} (\mathcal{B}_0^{r_1} \smallsetminus
\bar{\mathcal{B}}_0^{r_2})$, we have \footnote{There is a misprint in the analogous
nonsuper case given in \cite{H book}.  The domain of the sewing equation in Proposition
3.4.1 in \cite{H book} is $(f^{(1)})^{-1}(B^{r_1} \smallsetminus \bar{B}^{r_2}) + z_i$,
not $(f^{(2)})^{-1}(B^{1/r_2} \smallsetminus \bar{B}^{1/r_1})$, as stated.
Corrected, this corresponds to the projection onto the body of $s_{(z_i,
\theta_i)}^{-1} \circ (H^{(1)})^{-1} (\mathcal{B}_0^{r_1}
\smallsetminus
\bar{\mathcal{B}}_0^{r_2})$.}
\begin{equation}\label{F relation}
F^{(1)}_\sou  (w, \rho) = F^{(2)}_\sou  \circ (H^{(2)})^{-1} \circ I \circ H^{(1)}
\circ s_{(z_i,\theta_i)} (w, \rho) . 
\end{equation}
Furthermore, if we let 
\[H_0^{(1)}(w,\rho) = \tilde{E}(A^{(0)}, -iM^{(0)}) \Bigl(\frac{1}{w},
\frac{i\rho}{w}\Bigr) ,\]
denote the local coordinate of the puncture vanishing at $\infty$ of the
canonical supersphere represented by $Q_1$; let
\[H_k^{(1)}(w,\rho) = \hat{E}(\asqrt^{(k)}, A^{(k)}, M^{(k)}) (w, \rho), \]
so that the local coordinate vanishing at the $k$-th puncture of the 
canonical supersphere represented by $Q_1$ is given by $H_k^{(1)} \circ 
s_{(z_k,\theta_k)}$, for $k = 1,...,m$, $k \neq i$; and let 
\[H_l^{(2)}(w,\rho) = \hat{E}(\bsqrt^{(l)}, B^{(l)}, N^{(l)}) (w, \rho), \] 
so that the local coordinate vanishing at the $l$-th puncture of the 
canonical supersphere represented by $Q_2$ is given by $H_l^{(2)} \circ 
s_{(z_l',\theta_l')}$, for $l = 1,...,n$; then we have the following: \footnote{The
analogous nonsuper case to the sewing formulas (1)--(5) are given in 
\cite{H book} on pp.78--79 and pp.189--191.  However, the functions listed 
as ``local coordinate maps vanishing at the punctures" are actually the 
local coordinates normalized to vanish at zero.  That is, if we let
$g_k$, for $k=1,...,m+n-1$, denote the local coordinates vanishing at 
the puncture $z_k$ of the resulting sewn canonical sphere, in each 
of the cases (1)--(5) in \cite{H book}, the functions listed are
$g_0(1/x)$, $g_k(x + z_k)$, for $k=1,...,m+n-2$, and $g_{m+n-1}(x)$ for
$x$ a complex variable.  However, the actual formulas given for 
$Q_1 \; _i\infty_0  \; Q_2 \in K(m+n-1)$ in Appendix A of \cite{H book}
are correct.  }

(1) When $i = m$, and $n > 0$, the punctures of the canonical supersphere
represented by $Q_1 \; _m\infty_0 \; Q_2$ are 
\[\infty, \; F^{(1)}_\sou  (z_1, \theta_1),..., F^{(1)}_\sou (z_{m-1},
\theta_{m-1}), \;  F^{(2)}_\sou (z_1', \theta_1'),...,
 F^{(2)}_\sou (z_{n-1}', \theta_{n-1}'), \; 0;\]
and the local coordinates vanishing at these punctures are
\begin{gather*}
H_0^{(1)} \circ (F^{(1)}_\sou)^{-1} (w, \rho),\\    
H_1^{(1)} \circ s_{(z_1,\theta_1)} \circ (F^{(1)}_\sou)^{-1} (w,\rho), ...,
H_{m-1}^{(1)} \circ s_{(z_{m-1},\theta_{m-1})} \circ (F^{(1)}_\sou)^{-1} (w,\rho), \\
H_1^{(2)} \circ s_{(z_1',\theta_1')} \circ (F^{(2)}_\sou)^{-1} (w,\rho), ...,
H_{n-1}^{(2)} \circ s_{(z_{n-1}',\theta_{n-1}')} \circ  (F^{(2)}_\sou)^{-1}
(w,\rho) \\
H_{n}^{(2)} \circ  (F^{(2)}_\sou)^{-1}(w,\rho) ,
\end{gather*}  
respectively.  Thus
\begin{multline*}
Q_1 \; _m\infty_0 \; Q_2  \\
= \; \Bigl(F^{(1)}_\sou  (z_1, \theta_1),..., F^{(1)}_\sou (z_{m-1}, \theta_{m-1}), 
F^{(2)}_\sou (z_1', \theta_1'),..., F^{(2)}_\sou (z_{n-1}', \theta_{n-1}'); \\ 
\bigl((\tilde{E}^{-1})^0 (H_0^{(1)} \circ 
(F^{(1)}_\sou)^{-1} \circ I^{-1}  (w, \rho)), i(\tilde{E}^{-1})^1 (H_0^{(1)} \circ 
(F^{(1)}_\sou)^{-1} \circ I^{-1}  (w,\rho)\bigr),\\    
\hat{E}^{-1}(H_1^{(1)} \circ s_{(z_1,\theta_1)} \circ (F^{(1)}_\sou)^{-1} \circ 
s^{-1}_{F^{(1)}_\sou  (z_1, \theta_1)}(w,\rho)), ...,\\
\hat{E}^{-1}(H_{m-1}^{(1)} \circ s_{(z_{m-1},\theta_{m-1})} \circ 
(F^{(1)}_\sou)^{-1} \circ s^{-1}_{F^{(1)}_\sou  (z_{m-1}, \theta_{m-1})}(w,\rho)), \\   
\hat{E}^{-1}(H_1^{(2)} \circ s_{(z_1',\theta_1')} \circ 
(F^{(2)}_\sou)^{-1} \circ s^{-1}_{F^{(2)}_\sou  (z_1',\theta_1')}(w,\rho)), ...,\\
\hat{E}^{-1}(H_{n-1}^{(2)} \circ s_{(z_{n-1}',\theta_{n-1}')} \circ  
(F^{(2)}_\sou)^{-1} \circ s^{-1}_{F^{(2)}_\sou  (z_{n-1}', \theta_{n-1}')}(w,\rho)), \\
\hat{E}^{-1}(H_n^{(2)} \circ (F^{(2)}_\sou)^{-1} (w,\rho))\Bigr) ,
\end{multline*}

(2)  When $i = m = 1$ and $n = 0$, the canonical supersphere represented
by $Q_1 \; _1\infty_0 \; Q_2$ has only the one puncture at $\infty$ with 
local coordinate given by
\begin{equation}\label{coordinate with only one}
H_0^{(1)} \circ (F^{(1)}_\sou)^{-1} \circ s_{(a',m')}^{-1} (w,\rho),
\end{equation}
where $(a',m') \in \bigwedge_\infty$ is the unique element such that
$Q_1 \; _1\infty_0 \; Q_2$ represents a canonical supersphere in $SK(0)$,
i.e., such that the expansion of (\ref{coordinate with only one}) has
even coefficient of $w^{-2}$ and odd coefficient of $w^{-1}$ equal to zero.
Thus
\begin{multline*}
Q_1 \; _1\infty_0 \; Q_2 = \bigl( (\tilde{E}^{-1})^0 (H_0^{(1)} \circ 
(F^{(1)}_\sou)^{-1} \circ s_{(a',m')}^{-1} \circ I^{-1} (w,\rho)),  \\
i(\tilde{E}^{-1})^1 (H_0^{(1)} \circ
(F^{(1)}_\sou)^{-1} \circ s_{(a',m')}^{-1} \circ I^{-1} (w,\rho)) \bigr). 
\end{multline*}

(3)  When $i = m > 1$, and $n = 0$, and writing $F^{(1)}_\sou (z_{m - 1},
\theta_{m-1}) = p \in \bigwedge_\infty$, the punctures of the canonical 
supersphere represented by $Q_1 \; _i\infty_0 \; Q_2$ are
\[\infty, \; s_p \circ F^{(1)}_\sou (z_1, \theta_1),..., s_p \circ F^{(1)}_\sou 
(z_{m-2}, \theta_{m-2}), \; 0; \]
and the local coordinates vanishing at these punctures are
\begin{gather*}
H_0^{(1)} \circ (F^{(1)}_\sou)^{-1} \circ  s_p^{-1} (w, \rho), \\
H_1^{(1)} \circ s_{(z_1,\theta_1)} \circ (F^{(1)}_\sou)^{-1} \circ 
s_p^{-1} (w,\rho), ..., H_{m-2}^{(1)} \circ s_{(z_{m-2},\theta_{m-2})} 
\circ (F^{(1)}_\sou)^{-1}\circ s_p^{-1} (w,\rho)\\ 
H_{m-1}^{(1)} \circ s_{(z_{m-1},\theta_{m-1})} \circ (F^{(1)}_\sou)^{-1} \circ
s_p^{-1} (w,\rho),
\end{gather*}   
respectively. Thus 
\begin{multline*}
Q_1 \; _i\infty_0 \; Q_2 = \Bigl( s_p \circ F^{(1)}_\sou (z_1, \theta_1),..., 
s_p \circ F^{(1)}_\sou (z_{m-2}, \theta_{m-2});  \\
\bigl((\tilde{E}^{-1})^0 (H_0^{(1)} \circ (F^{(1)}_\sou)^{-1} \circ s_p^{-1} 
\circ I^{-1} (w, \rho)), i(\tilde{E}^{-1})^1 (H_0^{(1)} \circ (F^{(1)}_\sou)^{-1} 
\circ s_p^{-1} \circ I^{-1} (w,\rho)\bigr), \\  
\hat{E}^{-1}(H_1^{(1)} \circ s_{(z_1,\theta_1)} \circ (F^{(1)}_\sou)^{-1} \circ 
s^{-1}_{F^{(1)}_\sou (z_1, \theta_1)} (w,\rho)), ...,\\
\hat{E}^{-1}(H_{m-2}^{(1)} \circ s_{(z_{m-2},\theta_{m-2})} \circ 
(F^{(1)}_\sou)^{-1} \circ s^{-1}_{F^{(1)}_\sou 
(z_{m-2}, \theta_{m-2})} (w,\rho)),\\
\hat{E}^{-1}(H_{m-1}^{(1)} \circ s_{(z_{m-1},\theta_{m-1})} \circ 
(F^{(1)}_\sou)^{-1} \circ s_p^{-1} (w,\rho)) \Bigr) 
\end{multline*}   

(4) When $i < m$ and $n \neq 0$, writing $F^{(1)}_\sou (0) = p \in 
\bigwedge_\infty$, the punctures of the canonical supersphere represented 
by $Q_1 \; _i\infty_0 \; Q_2$ are
\begin{multline*}
\infty, \; s_p \circ F^{(1)}_\sou (z_1,\theta_1),..., s_p \circ 
F^{(1)}_\sou (z_{i-1}, \theta_{i-1}), \; s_p \circ F^{(2)}_\sou (z_1',\theta_1'), 
...,\\
s_p \circ F^{(2)}_\sou (z_{n-1}', \theta_{n-1}'), 
- F^{(1)}_\sou (0), \; s_p \circ F^{(1)}_\sou (z_{i+1}, \theta_{i+1}), ...,\\
s_p \circ F^{(1)}_\sou (z_{m-1}, \theta_{m-1}), \; 0 ; 
\end{multline*}
and the local coordinates vanishing at these punctures are
\begin{gather*}
H_0^{(1)} \circ (F^{(1)}_\sou)^{-1} \circ s_p^{-1} (w, \rho), \\
H_1^{(1)} \circ s_{(z_1,\theta_1)} \circ (F^{(1)}_\sou)^{-1} \circ 
s_p^{-1} (w,\rho),..., H_{i-1}^{(1)} \circ s_{(z_{i-1},\theta_{i-1})} 
\circ (F^{(1)}_\sou)^{-1} \circ s_p^{-1} (w,\rho), \\    
H_1^{(2)} \circ s_{(z_1',\theta_1')} \circ (F^{(2)}_\sou)^{-1} \circ
s_p^{-1} (w,\rho), ...,  H_{n-1}^{(2)} \circ s_{(z_{n-1}',\theta_{n-1}')} 
\circ  (F^{(2)}_\sou)^{-1} \circ s_p^{-1} (w,\rho),\\  
H_n^{(2)} \circ (F^{(2)}_\sou)^{-1} \circ s_p^{-1} (w,\rho), \\ 
H_{i+1}^{(1)} \circ s_{(z_{i+1},\theta_{i+1})} \circ (F^{(1)}_\sou)^{-1} 
\circ s_p^{-1} (w,\rho),..., \hspace{2.2in} \\
\hspace{2.2in} H_{m-1}^{(1)} \circ s_{(z_{m-1},\theta_{m-1})} 
\circ (F^{(1)}_\sou)^{-1} \circ s_p^{-1} (w,\rho), \\ 
H_m^{(1)} \circ (F^{(1)}_\sou)^{-1}\circ s_p^{-1} (w,\rho),
\end{gather*}
respectively.  Thus
\begin{multline*}
Q_1 \; _i\infty_0 \; Q_2 = \Bigl( s_p \circ F^{(1)}_\sou (z_1,
\theta_1),..., s_p \circ F^{(1)}_\sou (z_{i-1}, \theta_{i-1}),
s_p \circ F^{(2)}_\sou (z_1', \theta_1'), ...,\Bigr. \\
s_p \circ F^{(2)}_\sou (z_{n-1}', \theta_{n-1}'),  -
F^{(1)}_\sou (0), 
s_p \circ F^{(1)}_\sou (z_{i+1}, \theta_{i+1}), ...,\\
s_p \circ F^{(1)}_\sou (z_{m-1}, 
\theta_{m-1}) ; \; \; \bigl((\tilde{E}^{-1})^0 (H_0^{(1)} \circ (F^{(1)}_\sou)^{-1} \circ
s_p^{-1} \circ I^{-1}(w, \rho)), \\
i(\tilde{E}^{-1})^1 (H_0^{(1)} \circ (F^{(1)}_\sou)^{-1}
\circ s_p^{-1} \circ I^{-1} (w,\rho)\bigr), \\   
\hat{E}^{-1}(H_1^{(1)} \circ s_{(z_1,\theta_1)} \circ (F^{(1)}_\sou)^{-1} \circ
s^{-1}_{F^{(1)}_\sou (z_1,\theta_1)}(w,\rho)), ..., \\
\hat{E}^{-1}(H_{i-1}^{(1)} \circ s_{(z_{i-1},\theta_{i-1})} \circ 
(F^{(1)}_\sou)^{-1} \circ s^{-1}_{F^{(1)}_\sou 
(z_{i-1},\theta_{i-1})} (w,\rho)), \\   
\hat{E}^{-1}(H_1^{(2)} \circ s_{(z_1',\theta_1')} \circ (F^{(2)}_\sou)^{-1} 
\circ s^{-1}_{F^{(2)}_\sou (z_1',\theta_1')}(w,\rho)), ..., \\
\hat{E}^{-1}(H_{n-1}^{(2)} \circ s_{(z_{n-1}',\theta_{n-1}')} \circ  
(F^{(2)}_\sou)^{-1} \circ s^{-1}_{F^{(2)}_\sou 
(z_{n-1}',\theta_{n-1}')} (w,\rho)),\\ 
\hat{E}^{-1}(H_n^{(2)} \circ (F^{(2)}_\sou)^{-1}(w,\rho)), \\
\hat{E}^{-1} (H_{i+1}^{(1)} \circ s_{(z_{i+1},\theta_{i+1})} \circ
 (F^{(1)}_\sou)^{-1} \circ s^{-1}_{ F^{(1)}_\sou 
(z_{i+1},\theta_{i+1})} (w,\rho)),..., \\ 
\hat{E}^{-1}(H_{m-1}^{(1)} \circ s_{(z_{m-1},\theta_{m-1})} \circ 
(F^{(1)}_\sou)^{-1} \circ s^{-1}_{F^{(1)}_\sou
(z_{m-1},\theta_{m-1})}(w,\rho)),\\ 
\hat{E}^{-1}(H_m^{(1)} \circ (F^{(1)}_\sou)^{-1} \circ s_p^{-1} (w,\rho)) \Bigr) .
\end{multline*}

(5) When $i < m$ and $n = 0$, writing $F^{(1)}_\sou (0) = p \in 
\bigwedge_\infty$, the punctures of the canonical supersphere represented 
by $Q_1 \; _i\infty_0 \; Q_2$ are
\begin{multline*}
\infty, \; s_p \circ F^{(1)}_\sou (z_1,\theta_1),..., s_p \circ 
F^{(1)}_\sou (z_{i-1}, \theta_{i-1}), \; s_p \circ F^{(1)}_\sou (z_{i+1}, 
\theta_{i+1}), ..., \\
s_p \circ F^{(1)}_\sou (z_{m-1}, \theta_{m-1}), \; 0 ; 
\end{multline*}
and the local coordinates vanishing at these punctures are
\begin{gather*}
H_0^{(1)} \circ (F^{(1)}_\sou)^{-1} \circ s_p^{-1} (w, \rho), \\
H_1^{(1)} \circ s_{(z_1,\theta_1)} \circ (F^{(1)}_\sou)^{-1} \circ 
s_p^{-1} (w,\rho),..., H_{i-1}^{(1)} \circ s_{(z_{i-1},\theta_{i-1})} 
\circ (F^{(1)}_\sou)^{-1} \circ s_p^{-1} (w,\rho), \\    
H_{i+1}^{(1)} \circ s_{(z_{i+1},\theta_{i+1})} \circ (F^{(1)}_\sou)^{-1} 
\circ s_p^{-1} (w,\rho),..., \hspace{2.2in} \\
\hspace{2.2in} H_{m-1}^{(1)} \circ s_{(z_{m-1},\theta_{m-1})} 
\circ (F^{(1)}_\sou)^{-1} \circ s_p^{-1} (w,\rho), \\ 
H_m^{(1)} \circ (F^{(1)}_\sou)^{-1}\circ s_p^{-1} (w,\rho),
\end{gather*}
respectively.  Thus
\begin{multline*}
Q_1 \; _i\infty_0 \; Q_2 = \Bigl( s_p \circ F^{(1)}_\sou (z_1,
\theta_1),..., s_p \circ F^{(1)}_\sou (z_{i-1}, \theta_{i-1}),
s_p \circ F^{(1)}_\sou (z_{i+1}, \theta_{i+1}), ...,\\
s_p \circ F^{(1)}_\sou (z_{m-1}, 
\theta_{m-1}) ; \; \; \bigl((\tilde{E}^{-1})^0 (H_0^{(1)} \circ (F^{(1)}_\sou)^{-1} \circ
s_p^{-1} \circ I^{-1}(w, \rho)), \\
i(\tilde{E}^{-1})^1 (H_0^{(1)} \circ (F^{(1)}_\sou)^{-1}
\circ s_p^{-1} \circ I^{-1} (w,\rho)\bigr), \\   
\hat{E}^{-1}(H_1^{(1)} \circ s_{(z_1,\theta_1)} \circ (F^{(1)}_\sou)^{-1} \circ
s^{-1}_{F^{(1)}_\sou (z_1,\theta_1)}(w,\rho)), ..., \\
\hat{E}^{-1}(H_{i-1}^{(1)} \circ s_{(z_{i-1},\theta_{i-1})} \circ 
(F^{(1)}_\sou)^{-1} \circ s^{-1}_{F^{(1)}_\sou 
(z_{i-1},\theta_{i-1})} (w,\rho)), \\   
\hat{E}^{-1} (H_{i+1}^{(1)} \circ s_{(z_{i+1},\theta_{i+1})} \circ
 (F^{(1)}_\sou)^{-1} \circ s^{-1}_{ F^{(1)}_\sou 
(z_{i+1},\theta_{i+1})} (w,\rho)),..., \\ 
\hat{E}^{-1}(H_{m-1}^{(1)} \circ s_{(z_{m-1},\theta_{m-1})} \circ 
(F^{(1)}_\sou)^{-1} \circ s^{-1}_{F^{(1)}_\sou
(z_{m-1},\theta_{m-1})}(w,\rho)),\\ 
\hat{E}^{-1}(H_m^{(1)} \circ (F^{(1)}_\sou)^{-1} \circ s_p^{-1} (w,\rho)) \Bigr) .
\end{multline*}
\end{prop}

\begin{proof}  If the $i$-th tube of $Q_1$ can be sewn
with the 0-th tube of $Q_2$, then by definition there exist $r_1, r_2
\in \mathbb{R}_+$ satisfying $r_1 > r_2$ and DeWitt open neighborhoods
$U_i^{(1)}$ and $U_0^{(2)}$ of $(z_i, \theta_i) \in Q_1$ and $(\infty,
0) \in Q_2$, respectively, such that: $(z_i, \theta_i)$ and $(\infty,
0)$ are the only punctures in $U_i^{(1)}$ and $U_0^{(2)}$,
respectively; $H^{(1)} \circ s_{(z_i, \theta_i)} (w,\rho)$ and
$H^{(2)}(w,\rho)$ are convergent in $U_i^{(1)}$ and $U_0^{(2)}$,
respectively; and $\bar{\mathcal{B}}_0^{r_1} \subset H^{(1)} \circ
s_{(z_i, \theta_i)}(U_i^{(1)})$ and $\bar{\mathcal{B}}_0^{1/r_2} 
\subset H^{(2)} (U_0^{(2)})$.  The positive numbers $r_1$ and $r_2$
satisfy all the properties needed.  Conversely, if there exist
$r_1, r_2 \in \mathbb{R}_+$, with $r_1 > r_2$ such that
$(H^{(1)})^{-1}(w,\rho)$ and $(H^{(2)})^{-1}(w,\rho)$ are convergent
and single-valued in $\mathcal{B}_0^{r_1}$ and 
$\mathcal{B}_0^{1/r_2} \smallsetminus (\{0 \} \times (\bigwedge_\infty)_S)$,
respectively, 
\[(-z_i, -\theta_i), (z_k - z_i - \theta_k \theta_i, \theta_k -
\theta_i) \notin (H^{(1)})^{-1} (\mathcal{B}_0^{r_1}), \]
for $k = 1,...,m-1$, $k \neq i$, and 
\[0, (z_l', \theta_l') \notin (H^{(2)})^{-1}(\mathcal{B}_0^{1/r_2}
\smallsetminus (\{0 \} \times (\mbox{$\bigwedge_\infty$})_S)),\]
for $l = 1,...,n-1$, then we can choose $r$ to be any positive number 
between $r_1$ and $r_2$, and we can let
\[U_i^{(1)} = s^{-1}_{(z_i,\theta_i)} \circ (H^{(1)})^{-1} (\mathcal{B}_0^{r_1}) \]
and \[U_0^{(2)} = (H^{(2)})^{-1}(\mathcal{B}_0^{1/r_2}) .\]
By assumption, $(z_i, \theta_i)$ and $(\infty, 0)$ are the only
punctures in $U_i^{(1)}$ and $U_0^{(2)}$, respectively.  In addition,
\[\bar{\mathcal{B}}_0^r \subset \mathcal{B}_0^{r_1} = H^{(1)} \circ s_{(z_i,
\theta_i)}(U_i^{(1)}) \]
and 
\[\bar{\mathcal{B}}_0^{1/r} \subset \mathcal{B}_0^{1/r_2} = H^{(2)}
(U_0^{(2)}) . \]
Thus the $i$-th tube of $Q_1$ can be sewn with the 0-th tube of $Q_2$.
This finishes the first statement of the proposition. 

Now assume that the $m$-th tube of $Q_1$ can be sewn with the 0-th
tube of $Q_2$ and $n > 0$.  By the definition of the sewing operation,
the canonical supersphere $Q_1 \; _m\infty_0 \; Q_2$ 
is in the same superconformal equivalence class as the supersphere
with tubes
\begin{equation}\label{non-canonical supersphere}
\bigl(S;p_0,p_1,...,p_{m + n -1}; (U_0, \Omega_0),...,(U_{m + n -1},
\Omega_{m + n - 1})\bigr) 
\end{equation}
where the supersphere $S$ is given by the local coordinate system
\[\{(W^{(1)}, R^{(1)}), (W^{(2)}, R^{(2)})\} \]
with the coordinate transition function given by
\[R^{(2)} \circ (R^{(1)})^{-1} (w,\rho) = (H^{(2)})^{-1} \circ I \circ 
H^{(1)} \circ s_{(z_i,\theta_i)} (w,\rho) \]  
for $(w,\rho) \in R^{(1)}(W^{(1)} \cap W^{(2)})$ such that
\begin{eqnarray*}
R^{(1)} (W^{(1)}) &=& S\hat{\mathbb{C}} \smallsetminus s_{(z_i,\theta_i)}^{-1}
\circ (H^{(1)})^{-1} (\bar{\mathcal{B}}_0^{r_2}) \\
R^{(2)} (W^{(2)}) &=& \mbox{$\bigwedge_\infty$} \smallsetminus
(H^{(2)})^{-1}(\bar{\mathcal{B}}_0^{1/r_1}) ;
\end{eqnarray*}
$p_0,...,p_{m-1} \in W^{(1)}$, with coordinates $\infty, (z_1,
\theta_1),...,(z_{m-1},\theta_{m-1})$, respectively;
$p_m,...,p_{m + n -1} \in W^{(2)}$,
with coordinates $(z_1', \theta_1'),...,(z_{n-1}', \theta_{n-1}'), 0$,
respectively; in terms of the coordinate $(w,\rho) = R^{(1)}(p)$
for $p \in  W^{(1)}$,
\[\Omega_k(p) = H_k^{(1)} \circ s_{(z_k,\theta_k)} (w,\rho), \; \; 
\mbox{for} \; \; k = 0,..., m-1 ;\]
and in terms of the coordinate $(w,\rho) = R^{(2)}(p)$ for $p \in
W^{(2)}$, 
\[\Omega_{m-1 + l} (p) = H_l^{(2)} \circ s_{(z_l',\theta_l')} (w,\rho), \; \; 
\mbox{for} \; \; l = 1,..., n-1 .\]
By Proposition \ref{canonicalcriteria}, there exists a
superconformal equivalence $F$ {}from this supersphere with tubes to some
canonical supersphere with tubes 
\begin{multline}\label{canonical supersphere}
\bigl(S\hat{\mathbb{C}}; \nor^{-1}(0), \sou^{-1}(x_1,\varphi_1),....., \sou^{-1}(x_{m
+ n - 2}, \varphi_{m + n - 2}), \sou^{-1}(0);  \\
(\sou^{-1}(\mathcal{B}^{r_0'}_\infty) \cup \nor^{-1}( \{0 \} \times
(\mbox{$\bigwedge_\infty$})_S ), H_0 ), (\sou^{-1}(\mathcal{B}^{r_1'}_{x_1}), H_1
\circ \sou),...,\\ 
(\sou^{-1}(\mathcal{B}^{r_{m + n - 2}'}_{x_{m + n -
2}}), H_{m + n - 2} \circ \sou), (\sou^{-1}(\mathcal{B}^{r_{m + n -
1}'}_0), H_{m + n -1} \circ \sou) \bigr) .
\end{multline}  

Let 
\[ F^{(1)}_\sou (w,\rho) = F \circ (R^{(1)})^{-1}(w,\rho) \] 
for $(w,\rho) \in R^{(1)} (W^{(1)})$, and let 
\[ F^{(2)}_\sou (w,\rho) = F \circ (R^{(2)})^{-1}(w,\rho) \]
for $(w,\rho) \in R^{(2)} (W^{(2)})$.  Then in $W^{(1)} \cap W^{(2)}$,
we have
\begin{equation}\label{FR equation}
F^{(1)}_\sou \circ R^{(1)} (p) = F^{(2)}_\sou \circ R^{(2)}(p) . 
\end{equation}
It is clear that we can choose $F^{(1)}$ and $F^{(2)}$ to satisfy
(\ref{normalize 1}) -- (\ref{normalize 3}).  On
the other hand, for $p \in W^{(1)} \cap W^{(2)}$,
\begin{equation}\label{R equation}
R^{(2)} (p) = (H^{(2)})^{-1} \circ I \circ 
H^{(1)} \circ s_{(z_i,\theta_i)} \circ R^{(1)}  (p) . 
\end{equation}
By (\ref{FR equation}) and (\ref{R equation}), $F^{(1)}$ and $F^{(2)}$
satisfy (\ref{F relation}).  The formula in case (1) follows
immediately {}from the relation between canonical superspheres and
elements of $SK(m + n - 1)$.  This finishes the proof for case (1) 
$i = m$ and $n > 0$ up to uniqueness of $F^{(1)}$ and $F^{(2)}$.
Cases (2) -- (5) are proved similarly.  

To prove that $F^{(1)}$ and $F^{(2)}$ are unique, assume that there
exists another pair of functions $\tilde{F}^{(1)}$ and
$\tilde{F}^{(2)}$ which also satisfy (\ref{normalize 1}) -- 
(\ref{F relation}).  If we define
\[\tilde{F} (p) = \left\{ \begin{array}{lll}
                  \tilde{F}^{(1)}_\sou \circ R^{(1)} (p) & \mbox{for} & p
                   \in W^{(1)}, \\  
                  \tilde{F}^{(2)}_\sou \circ R^{(2)} (p) & \mbox{for} & p
                   \in W^{(2)}, \end{array} \right. \]
then $\tilde{F}$ is a superconformal equivalence {}from the supersphere
with tubes (\ref{non-canonical supersphere}) to a canonical
supersphere with tubes.  Thus $\tilde{F} \circ F^{-1}$ is a 
superconformal equivalence {}from one canonical supersphere with tubes 
to another canonical supersphere with tubes.  By Corollary 
\ref{bijection}, we  see that $\tilde{F} \circ F^{-1}$ must be the 
identity map, or equivalently $\tilde{F} = F$.  Therefore $\tilde{F}^{(1)}
=F^{(1)}$ and $\tilde{F}^{(2)} = F^{(2)}$.  This proves the uniqueness of 
$F^{(1)}$ and $F^{(2)}$.  \end{proof} 

Note that  
\[s_{(z,\theta)}(w,\rho) = \exp \bigl(zL_{-1}(w,\rho) + \theta G_{-\frac{1}{2}}
(w,\rho) \bigr) \cdot (w,\rho) .\]
Thus using the formal solution to the sewing equation given in Chapter 3 to
express the components of the uniformizing function formally as
\begin{eqnarray*}
F_\sou^{(1)} (x,\varphi) &=& \left. \bar{F}^{(1)} \circ s_{(z_i,\theta_i)} (x,\varphi)
\right|_{(\alpha_0^{1/2}, \mathcal{A}, \mathcal{M}, \mathcal{B}, \mathcal{N}) = 
(\asqrt^{(i)},A^{(i)},M^{(i)},B^{(0)},N^{(0)})}\\
F_\sou^{(2)} (x,\varphi) &=& \left. \bar{F}^{(2)} (x,\varphi)
\right|_{(\alpha_0^{1/2}, \mathcal{A}, \mathcal{M}, \mathcal{B}, \mathcal{N}) = 
(\asqrt^{(i)},A^{(i)},M^{(i)},B^{(0)},N^{(0)})} ,
\end{eqnarray*}
we see that formally all the terms in the sewing formulas given in Proposition
\ref{actual sewing} can be expressed in terms of infinitesimal superconformal
transformations.  In order to use the formal solution as an actual analytic
and geometric solution, we need to show that the power series expansions
of $F_\sou^{(1)}$ and $F_\sou^{(2)}$ converge and are equal to $\bar{F}^{(1)} 
\circ s_{(z_i,\theta_i)}$ and $\bar{F}^{(2)}$, respectively, evaluated at 
$(\alpha_0^{1/2}, \mathcal{A}, \mathcal{M}, \mathcal{B}, \mathcal{N}) =
(\asqrt^{(i)},A^{(i)},M^{(i)},B^{(0)},N^{(0)})$.  The absolute convergence of 
these series is equivalent to the absolute convergence of the $x$ and 
$\varphi$ coefficients of
\begin{multline*}
\bar{F}^{(1)}_{Q_1, Q_2, t^{1/2}} (x,\varphi) \; = \; \exp \Biggl(\sum_{j \in
\Z} \biggl( \Psi_{-j}(t^{-\frac{1}{2}} \asqrt^{(i)}, A^{(i)}, M^{(i)},
B^{(0)}, N^{(0)}) L_{-j}(x,\varphi) \biggr. \Biggr. \\
\Biggl. \biggl. + \; \Psi_{-j+ \frac{1}{2}}
(t^{-\frac{1}{2}} \asqrt^{(i)}, A^{(i)}, M^{(i)}, B^{(0)}, N^{(0)}) 
G_{-j + \frac{1}{2}}(x,\varphi)
\biggr) \! \Biggr) \cdot (x,\varphi)
\end{multline*}
and 
\begin{multline*}
\bar{F}^{(2)}_{Q_1, Q_2, t^{1/2}} (x,\varphi) = 
\exp \left( - \Psi_0(t^{-\frac{1}{2}} \asqrt^{(i)}, A^{(i)}, M^{(i)},
B^{(0)}, N^{(0)}) 2L_0(x,\varphi) \right) \cdot \\ 
\cdot (\asqrt^{(i)})^{ 2L_0 (x,\varphi)} \cdot \exp \Biggl(\! - \! 
\sum_{j \in \Z} \biggl( \Psi_j(t^{-\frac{1}{2}} \asqrt^{(i)}, A^{(i)}, 
M^{(i)}, B^{(0)}, N^{(0)}) L_j(x,\varphi) \biggr. \Biggr. \\ 
\Biggl. \biggl. + \; \Psi_{j-\frac{1}{2}}
(t^{-\frac{1}{2}} \asqrt^{(i)}, A^{(i)}, M^{(i)}, B^{(0)}, N^{(0)}) 
G_{j - \frac{1}{2}}(x,\varphi) \biggr) \! \Biggr) \cdot (x,\varphi) 
\end{multline*} 
as power series in the complex variable $t^{1/2}$ for $|t^{1/2}| \leq 1$.

Define
\begin{multline*}
Q_1(t^\frac{1}{2}) = \bigl((z_1, \theta_1),...,(z_{m-1}, \theta_{m-1}); 
(A^{(0)},M^{(0)}),...,(\asqrt^{(i-1)}, A^{(i-1)},
M^{(i-1)}), \\
(t^{-\frac{1}{2}}\asqrt^{(i)}, A^{(i)}, M^{(i)}),(\asqrt^{(i+1)}, A^{(i+1)},
M^{(i+1)}),..., (\asqrt^{(m)}, A^{(m)}, M^{(m)})\bigr) ,
\end{multline*}
and apply Proposition \ref{actual sewing} to $Q_2$ sewn into the $i$-th
puncture of $Q_1(t^{1/2})$ for $0 < |t^{1/2}| \leq 1$.  Denote 
the two functions $F^{(1)}$ and $F^{(2)}$ giving the canonical supersphere 
in this case by $F^{(1)}_{t^{1/2}}$ and $F^{(2)}_{t^{1/2}}$, respectively.  
If we can prove that the expansion coefficients of $\sou \circ 
F^{(1)}_{t^{1/2}} \circ \sou^{-1} = F^{(1)}_{t^{1/2}, \sou} $ and 
$\sou \circ F^{(2)}_{t^{1/2}} \circ \sou^{-1} = F^{(2)}_{t^{1/2},\sou}$ as 
analytic functions in $w$ and $\rho$ are analytic in $t^{1/2}$ for 
$|t^{1/2}| \leq 1$, and their expansions as Laurent series in $t^{1/2}$ 
are equal to the corresponding $x$ and $\varphi$ coefficients of 
$\bar{F}^{(1)}_{Q_1, Q_2, t^{1/2}} \circ s_{(z_i,\theta_i)} (x,\varphi)$ and 
$\bar{F}^{(2)}_{Q_1, Q_2, t^{1/2}} (x,\varphi)$, respectively, then these 
$x$ and $\varphi$ coefficients of $\bar{F}^{(1)}_{Q_1, Q_2,t^{1/2}} \circ 
s_{(z_i,\theta_i)} (x,\varphi)$ and $\bar{F}^{(2)}_{Q_1, Q_2, t^{1/2}}(x,
\varphi)$ are absolutely convergent for $|t^{1/2}| \leq 1$.  In particular, 
taking $t^{1/2}=1$, we would have that each $x$ and $\varphi$ coefficient of
$\bar{F}^{(1)}_{Q_1, Q_2} \circ s_{(z_i,\theta_i)}(x,\varphi)$ and 
$\bar{F}^{(2)}_{Q_1, Q_2}(x,\varphi)$ is absolutely convergent to
the corresponding $w$ and $\rho$ coefficients of $F^{(1)}_\sou$ and 
$F^{(2)}_\sou$.

\begin{rema}
Since by Theorem \ref{uniformization}, $\bar{F}^{(1)}_{Q_1, Q_2} 
(x,\varphi)$ and $\bar{F}^{(2)}_{Q_1, Q_2}(x,\varphi)$ depend
algebraically on 
\begin{multline*}
H^{(1)}_{\asqrt^{1/2}, A, M} (x, \varphi) \\
= \; \exp \Biggl( \! - \! \sum_{j \in \Z} \left( A_j L_j(x,\varphi)
 + M_{j - \frac{1}{2}} G_{j - \frac{1}{2}} (x,\varphi)
\right) \! \Biggr) \cdot (\asqrt^\frac{1}{2})^{-2L_0(x,\varphi)}
\cdot (x, \varphi) 
\end{multline*}
and
\[H^{(2)}_{ B, N} (x, \varphi) =  
\exp \Biggl( \sum_{j \in \Z} \left( B_j L_{-j}(x,\varphi) + 
N_{j - \frac{1}{2}}  G_{-j + \frac{1}{2}}(x,\varphi) 
\right) \! \Biggr) \cdot \Bigl(\frac{1}{x}, \frac{i \varphi}{x}\Bigr), \]
and are determined uniquely by the formal sewing equation and formal
boundary conditions, the convergence described above and proved below
implies that $F^{(1)}$ and $F^{(2)}$ depend algebraically on the local 
coordinate maps $H^{(1)}(w,\rho)$ and $H^{(2)}(w - z_i -\rho \theta_i, 
\rho - \theta_i)$ in the sense that the formal series obtained {}from 
the expansions of $F^{(1)}$ and $F^{(2)}$ depend algebraically on the 
formal series obtained {}from the expansions of $H^{(1)}(w,\rho)$ and 
$H^{(2)}(w - z_i -\rho \theta_i, \rho - \theta_i)$.  In addition,
we see that $F^{(1)}$ and $F^{(2)}$ are determined directly and
explicitly by the sewing equation, normalization conditions, boundary
conditions and the algebraic dependency on $H^{(1)}(w,\rho)$ and 
$H^{(2)}(w - z_i -\rho \theta_i,\rho - \theta_i)$.  Thus a consequence
of the convergence result below is the answers to the questions posed 
at the end of Section 2.6 regarding the dependency of the uniformizing
function on the local coordinates.
\end{rema}

We follow the argument in \cite{H book} using a theorem proved by 
Fischer and Grauert \cite{FG} in the deformation theory of complex 
manifolds to prove the analyticity in $t^{1/2}$.  However, the 
supersphere $Q_1(t^{1/2}) \; _i\infty_0 \; Q_2$ is of course an  
infinite-dimensional non-compact complex manifold, and the  
Fischer-Grauert Theorem only applies to finite-dimensional compact 
manifolds, e.g., the body of this supersphere.  But we do have
that $Q_1(t^{1/2}) \; _i\infty_0 \; Q_2$ is a fiber bundle over the 
body, with fiber isomorphic to $(\bigwedge_\infty)_S$, and we can 
define two families of global sections for each $(w_S,\rho) \in 
(\bigwedge_\infty)_S$ such that each section is complex analytically 
isomorphic to the Riemann sphere.  We can then apply the 
Fischer-Grauert Theorem to each section and analyze the behavior in 
$t^{1/2}$.   This will allow us to prove certain analyticity 
properties of the uniformizing function $F_{t^{1/2}}$ for 
$0 < |t^{1/2}| \leq 1$, on each section.  The sections that we 
will define completely cover the fiber bundle that the canonical
supersphere represented by $Q_1(t^{1/2}) \; _i\infty_0 \; Q_2$ defines 
in such a way as to give the analyticity of the uniformizing function 
$F_{t^{1/2}}$ for $0 < |t^{1/2}| \leq 1$ on the entire supersphere.  
We will then analyze the nature of the singularity at $t^{1/2} = 0$ of 
the components of the uniformizing function $F^{(1)}_{t^{1/2},\sou} (w,
\rho)$ and $F^{(2)}_{t^{1/2},\sou} (w,\rho)$ and expand each function 
as a Laurent series in $w$, $\rho$ and $t^{1/2}$.  We then show that 
the formal series corresponding to these expansions with $(w,\rho)$  
replaced by $(x,\varphi)$ are the same as the formal series 
$\bar{F}^{(1)}_{Q_1,Q_2, t^{1/2}} \circ s_{(z_i,\theta_i)} (x,\varphi)$ 
and $\bar{F}^{(2)}_{Q_1, Q_2, t^{1/2}} (x,\varphi)$, respectively.  
This proves the analyticity and convergence of the series $\Psi_j 
(t^{-\frac{1}{2}}a, A, M, B, N)$, for $j \in \frac{1}{2} \mathbb{Z}$ 
and $|t^\frac{1}{2}| \leq 1$.  
 
Now we present a brief description of the Fischer-Grauert Theorem
following \cite{H book}.  Let $\mathcal{D}$ be a connected complex 
manifold, and let $\mathfrak{M} = \{M_t \; | \; t \in \mathcal{D} \}$ 
be a family of compact complex manifolds parameterized by $t \in 
\mathcal{D}$.  We say that $M_t$ {\it depends holomorphically} (or {\it 
complex analytically}) {\it on $t$} and that $\mathfrak{M} = \{M_t \; | \; 
t \in \mathcal{D} \}$ {\it forms a complex analytic family} if there is 
a complex manifold $\mathbf{M}$ and a holomorphic map $\bar{\omega}$ {}from $\mathbf{M}$ onto $\mathcal{D}$ such that

(i) $\bar{\omega}^{-1} (t) = M_t$ for each $t \in \mathcal{D}$, and

(ii) the rank of the Jacobian of $\bar{\omega}$ is equal to the
complex dimension of $\mathcal{D}$ at each point of {\bf M}.

If each point $t \in \mathcal{D}$ has a neighborhood $\Delta_t$ such
that $\bar{\omega}^{-1} (\Delta_t)$ is complex analytically isomorphic
to $M_t \times \Delta_t$ and the diagram
\[\begin{array}{ccc}
\bar{\omega}^{-1} (\Delta_t)  & \stackrel{\sim}{\longrightarrow} & 
M_t \times \Delta_t \\
\Big\downarrow\vcenter{%
     \rlap{$\bar{\omega}$}} & &
\Big\downarrow\vcenter{%
     \rlap{proj.}}  \\ 
\Delta_t & \longrightarrow & \Delta_t
\end{array} \]
is commutative, then we say that $\mathfrak{M} = \{M_t \; | \; t \in \mathcal{D}
\}$ is {\it locally trivial} (complex analytically). 

\begin{thm}\label{Fischer-Grauert}(Fischer-Grauert)
If for all $t \in \mathcal{D}$, the $M_t$'s are complex analytically
isomorphic, then $\mathfrak{M}$ is locally trivial.
\end{thm}

The proof of this theorem can be found in \cite{FG}.

Next we construct two sets of global sections for the super-Riemann sphere
$S\hat{\mathbb{C}}$ viewed as a fiber bundle over the Riemann sphere
$\mathbb{C}\cup \infty = S\hat{\mathbb{C}}_B$.  Recall that 
$S\hat{\mathbb{C}}$ is given by the superconformal coordinate atlas 
$\{(U_\sou,\sou), (U_\nor, \nor)\}$ with coordinate transition function 
$\sou \circ \nor^{-1} = I$, and the Riemann sphere can then be represented
by the coordinate atlas $\{((U_\sou)_B,\sou_B), ((U_\nor)_B, \nor_B)\}$. 
For $(w_S,\rho) \in (\bigwedge_\infty)_S$ define the global section 
$\sigma_{(w_S,\rho)} : S\hat{\mathbb{C}}_B \rightarrow S\hat{\mathbb{C}}$ 
by 
\begin{equation}
\sigma_{(w_S,\rho)} (p_B)  =  \left\{\begin{array}{ll} 
       \sou^{-1}(\sou_B(p_B) + w_S, \rho) 
           &\mbox{for $p_B \in (U_\sou)_B$} ,\\
      \nor^{-1} (0)  
         & \mbox{for $p_B \notin (U_\sou)_B$}.\\ 
\end{array}  \right.
\end{equation}
Similarly, for $(w_S,\rho) \in (\bigwedge_\infty)_S$ define the global section 
$\tau_{(w_S,\rho)} : S\hat{\mathbb{C}}_B \rightarrow S\hat{\mathbb{C}}$ 
by 
\begin{equation}
\tau_{(w_S,\rho)} (p_B)  =  \left\{\begin{array}{ll} 
       \nor^{-1}(\nor_B(p_B) + w_S, \rho) 
           &\mbox{for $p_B \in (U_\nor)_B$} ,\\
      \sou^{-1} (0)  
         & \mbox{for $p_B \notin (U_\nor)_B$}.\\ 
\end{array}  \right.
\end{equation}
Then for $(w_S,\rho) \in (\bigwedge_\infty)_S$, the spaces $\sigma_{(w_S,\rho)} 
(S\hat{\mathbb{C}}_B)$ and $\tau_{(w_S,\rho)} (S\hat{\mathbb{C}}_B)$ each 
define Riemann spheres.  Note that if $M$ is any genus zero superconformal 
manifold with uniformizing function $F: M \longrightarrow S\hat{\mathbb{C}}$, 
for each $(w_S,\rho) \in (\bigwedge_\infty)_S$, we have the global sections 
of the fiber bundle $M$ given by $F^{-1} \circ \sigma_{(w_S,\rho)} \circ F_B:
M_B \longrightarrow M$ and $F^{-1} \circ \tau_{(w_S,\rho)} \circ F_B:
M_B \longrightarrow M$.

For $Q_1 \in SK(m)$ and $Q_2 \in SK(n)$ as in Proposition \ref{actual
sewing}, and
\begin{multline*}
Q_1(t^\frac{1}{2}) = \bigl( (z_1, \theta_1),...,(z_{m-1}, \theta_{m-1});
(A^{(0)}, M^{(0)}), ...,(\asqrt^{(i-1)}, A^{(i-1)}, M^{(i-1)}),  \\
(t^{-\frac{1}{2}}\asqrt^{(i)}, A^{(i)}, M^{(i)}),
(\asqrt^{(i+1)}, A^{(i+1)}, M^{(i+1)}), ..., (a^{(m)}, A^{(m)}, M^{(m)}) \bigr), 
\end{multline*}  
define
\[H_{t^{1/2}}^{(1)} (x,\varphi) = \hat{E} (t^{-\frac{1}{2}} \asqrt^{(i)},
A^{(i)}, M^{(i)}) (x,\varphi). \] 
Since the $i$-th tube of $Q_1$ can be sewn with the 0-th tube of $Q_2$, 
it is easy to see {}from Proposition \ref{actual sewing} that the $i$-th 
tube of $Q_1(t^{1/2})$ can be also be sewn with the $0$-th tube of $Q_2$ 
when $0<|t^{1/2}|<r$ for some $r >1$.  By Proposition \ref{actual sewing}, 
for any such $t^{1/2}$ there exist $r_1(t^{1/2})$, $r_2(t^{1/2}) \in 
\mathbb{R}_+$ satisfying $r_1(t^{1/2}) > r_2(t^{1/2})$ and bijective 
superconformal functions $F_{t^{1/2}}^{(1)}$ and $F_{t^{1/2}}^{(2)}$ on 
$S\hat{\mathbb{C}} \smallsetminus \sou^{-1} \circ s_{(z_i, \theta_i)}^{-1} 
\circ (H_{t^{1/2}}^{(1)})^{-1} (\bar{\mathcal{B}}_0^{r_2(t^{1/2})})$, and 
$U_\sou \smallsetminus \sou^{-1} \circ (H_{t^{1/2}}^{(2)})^{-1} 
(\bar{\mathcal{B}}_0^{1/r_1(t^{1/2})} \smallsetminus (\{0\} \times
(\bigwedge_\infty)_S))$, respectively, such that the conclusion of 
Proposition \ref{actual sewing} holds if we replace $Q_1$, $H^{(1)}$, 
$r_1$, $r_2$, $F^{(1)}$, and $F^{(2)}$ by $Q_1(t^{1/2})$, 
$H_{t^{1/2}}^{(1)}$, $r_1(t^{1/2})$, $r_2(t^{1/2})$, $F_{t^{1/2}}^{(1)}$, 
and $F_{t^{1/2}}^{(2)}$, respectively.  Denote $r_1(1)$ and $r_2(1)$ by 
$r_1$ and $r_2$, respectively. {}From the definition of 
$H_{t^{1/2}}^{(1)}$, it is clear that we can choose $r_1(t^{1/2}) = r_1$ 
and $r_2(t^{1/2}) = r_2$ for all $t^{1/2}$.     

\begin{prop}\label{use of FG}
There exist $r \in \mathbb{R}$, $r > 1$ such that the superconformal
superfunctions $F^{(1)}_{t^{1/2}}$ and $F^{(2)}_{t^{1/2}}$ are
analytic in $t^{1/2}$ for $0< |t^{1/2}| < r$.  Furthermore, the 
singularity at $t^{1/2} = 0$ is a removable singularity of 
$F^{(1)}_{t^{1/2},\sou} (w,\rho)$ and $F^{(2)}_{t^{1/2},\sou} (w,\rho)$, 
and writing 
\[F^{(2)}_{t^{1/2},\sou} (w,\rho) = \left((F^{(2)}_{t^{1/2},\sou} (w,\rho))^0,
(F^{(2)}_{t^{1/2},\sou} (w,\rho))^1\right) \in \mbox{$\bigwedge_\infty^0 \oplus 
\bigwedge_\infty^1$}\]
we have that $t^{1/2} = 0$ is a second-order zero of 
$(F^{(2)}_{t^{1/2},\sou} (w,\rho))^0$ and a first-order zero of 
$(F^{(2)}_{t^{1/2},\sou} (w,\rho))^1$.    
\end{prop}

\begin{proof}  $Q_1(t^{1/2}) \; _i\infty_0 \; Q_2$
is represented by a supersphere with tubes $M_{t^{1/2}}$ given by the
local coordinate system  
\[\{(W_{t^{1/2}}^{(1)}, R_{t^{1/2}}^{(1)}),
(W_{t^{1/2}}^{(2)}, R_{t^{1/2}}^{(2)})\} \] 
with the coordinate transition function given by
\begin{equation}\label{coordinate transition}
R_{t^{1/2}}^{(2)} \circ (R_{t^{1/2}}^{(1)})^{-1} (w, \rho)
= (H^{(2)})^{-1} \circ I \circ H_{t^{1/2}}^{(1)} \circ
s_{(z_i,\theta_i)} (w, \rho)   
\end{equation}
for $(w,\rho) \in R_{t^{1/2}}^{(1)} (W_{t^{1/2}}^{(1)}
\cap W_{t^{1/2}}^{(2)})$ such that 
\[R_{t^{1/2}}^{(1)} (W_{t^{1/2}}^{(1)}) = S\hat{\mathbb{C}} \smallsetminus
\sou^{-1} \circ s_{(z_i,\theta_i)}^{-1} \circ (H_{t^{1/2}}^{(1)})^{-1}
(\bar{\mathcal{B}}_0^{r_2})\] 
and
\[R_{t^{1/2}}^{(2)} (W_{t^{1/2}}^{(2)}) = U_\sou \smallsetminus
\sou^{-1} \circ (H^{(2)})^{-1} (\bar{\mathcal{B}}_0^{1/r_1} \smallsetminus 
(\{0\} \times\mbox{$(\bigwedge_\infty)_S)$}) = R^{(2)} (W^{(2)}) . \]
We know that for each $t^{1/2}$ the coordinate transition
function $(H^{(2)})^{-1} \circ I \circ H_{t^{1/2}}^{(1)} \circ
s_{(z_i,\theta_i)} (w,\rho)$ is superconformal. 

{}From the uniformizing function at $t^{1/2} = 1$, i.e., for $F=F_1 : 
M_1 \longrightarrow S\hat{\mathbb{C}}$, we have the global sections 
of $M_1$ given by $F_1^{-1} \circ \sigma_{(w_S,\rho)} \circ (F_1)_B$ 
and $F_1^{-1} \circ \tau_{(w_S,\rho)} \circ (F_1)_B$, for each $(w_S,
\rho) \in (\bigwedge_\infty)_S$.  Letting $t^{1/2}$ vary, we obtain 
global sections we denote by 
\begin{eqnarray*}
\sigma_{(w_S,\rho),t^{1/2}} : (M_{t^{1/2}})_B &\rightarrow& 
M_{t^{1/2}} \\
\tau_{(w_S,\rho),t^{1/2}} : (M_{t^{1/2}})_B &\rightarrow& 
M_{t^{1/2}},
\end{eqnarray*}
respectively.  Then defining 
\begin{eqnarray*}
M_{t^{1/2}}^{\sigma(w_S,\rho)} &=& \sigma_{(w_S,\rho),t^{1/2}} ((M_{t^{1/2}})_B)\\ 
M_{t^{1/2}}^{\tau(w_S,\rho)} &=& \tau_{(w_S,\rho),t^{1/2}} ((M_{t^{1/2}})_B), 
\end{eqnarray*}
we have that $M_{t^{1/2}}^{\sigma(w_S,\rho)}$ and 
$M_{t^{1/2}}^{\tau(w_S,\rho)}$ are genus zero compact complex manifolds.  

Let
\[ W^{(k)}_{t^{1/2}, \sigma(w_S,\rho)} = (R^{(k)}_{t^{1/2}})^{-1} \circ
\sigma_{(w_S,\rho),t^{1/2}}(W^{(k)}_{t^{1/2},B}) \]
for $k = 1,2$, and let $R^{(k)}_{t^{1/2},\sigma(w_S,\rho)} : W^{(k)}_{t^{1/2},
\sigma(w_S,\rho)} \rightarrow \bigwedge_\infty$ be the restriction of 
$R^{(k)}_{t^{1/2}}$ to $W^{(k)}_{t^{1/2},\sigma(w_S,\rho)}$, for $k=1,2$.
Then 
\[\left\{\left( W^{(1)}_{t^{1/2}, \sigma(w_S,\rho)}, R^{(1)}_{t^{1/2},
\sigma(w_S,\rho)}\right), \left(W^{(2)}_{t^{1/2}, \sigma(w_S,\rho)}, 
R^{(2)}_{t^{1/2},\sigma(w_S,\rho)}\right) \right\} \]
is a coordinate system for $M_{t^{1/2}}^{\sigma(w_S,\rho)}$, with coordinate 
transition function given by 
\[ R^{(2)}_{t^{1/2},\sigma(w_S,\rho)} \circ (R^{(1)}_{t^{1/2},
\sigma(w_S,\rho)})^{-1} (w,\rho) = (H^{(2)})^{-1} \circ I \circ 
H_{t^{1/2}}^{(1)} \circ s_{(z_i,\theta_i)} (w, \rho) \]
for $(w,\rho) \in R^{(1)}_{t^{1/2}, \sigma(w_S,\rho)} (W^{(1)}_{t^{1/2}, 
\sigma(w_S,\rho)} \cap W^{(2)}_{t^{1/2}, \sigma(w_S,\rho)})$.  It is clear that 
$M_{t^{1/2}}^{\sigma(w_S,\rho)}$ is complex analytically isomorphic to
$(M_{t^{1/2}})_B$.

Let $\mathcal{D}_r = \{t^{1/2} \in \mathbb{C} \; | \; 0<|t^{1/2}| < r, 
\; -\pi/2 < \mbox{arg} \; t^{1/2} \leq \pi/2 \}$ for $r \in \mathbb{R}_+$
and let $\tilde{\mathcal{D}}_r = \{t^{1/2} \in \mathbb{C} \; | \; 
0<|t^{1/2}| < r, \; \pi/2 < \mbox{arg} \; t^{1/2} \leq 3\pi/2 \}$. 
Consider $\mathfrak{M}^{\sigma(w_S,\rho)} = \{M_{t^{1/2}}^{\sigma(w_S,\rho)} 
\; | \; t^{1/2} \in \mathcal{D}_r \}$ and 
$\tilde{\mathfrak{M}}^{\sigma(w_S,\rho)} = \{M_{t^{1/2}}^{\sigma(w_S,\rho)} 
\; | \; t^{1/2} \in \tilde{\mathcal{D}}_r\}$.  In \cite{H book}, Huang 
uses the Fischer-Grauert Theorem on the body component 
$\mathfrak{M}^{\sigma(0)}$ to prove that it is locally trivial in $t$, for 
$0 < |t| < r$.  We  will follow this argument in the more general case of
$\mathfrak{M}^{\sigma(w_S,\rho)}$ and 
$\tilde{\mathfrak{M}}^{\sigma(w_S,\rho)}$ to prove that they are both 
locally trivial in $t^{1/2}$ for all $(w_S,\rho) \in (\bigwedge_\infty)_S$. 
The reason we split the domain of $t^{1/2}$ into $\mathcal{D}_r$ and
$\tilde{\mathcal{D}}_r$ is that we will need certain bijective properties
to form a manifold {}from the family $\mathfrak{M}^{\sigma(w_S,\rho)}$, and
thus have to restrict the domain of $t^{1/2}$ in order to insure that there
is no double cover in the $t^{1/2}$ coordinate.  We will then use the
local triviality properties of $\mathfrak{M}^{\sigma(w_S,\rho)}$ and
$\tilde{\mathfrak{M}}^{\sigma(w_S,\rho)}$ to show that $F^{(1)}_{t^{1/2}}$  
and $F^{(2)}_{t^{1/2}}$ restricted to the sections 
$M_{t^{1/2}}^{\sigma(w_S,\rho)}$ are analytic in $t^{1/2}$ for $0< |t^{1/2}| 
< r$.  Using the family of sections $\tau_{(w_S,\rho)}$ instead of
$\sigma_{(w_S,\rho)}$, the analogous argument for 
$\mathfrak{M}^{\tau(w_S,\rho)}$ and $\tilde{\mathfrak{M}}^{\tau(w_S,\rho)}$
shows that they are also locally trivial in $t^{1/2}$.  Then the analogous
argument can be used to prove that $F^{(1)}_{t^{1/2}}$ and 
$F^{(2)}_{t^{1/2}}$ restricted to the sections $M_{t^{1/2}}^{\tau(w_S,\rho)}$ 
are analytic in $t^{1/2}$ for $0< |t^{1/2}|  < r$.  This will allow us to
conclude that $F^{(1)}_{t^{1/2}}$ and $F^{(2)}_{t^{1/2}}$ themselves are 
analytic in $t^{1/2}$ for $0< |t^{1/2}|  < r$.

Let $\mathbf{M}^{\sigma(w_S,\rho)} = \bigcup_{t^{1/2} \in \mathcal{D}_r} 
M_{t^{1/2}}^{\sigma(w_S,\rho)}$.  Let    
\begin{eqnarray*}
U^{(1)}_{\sigma(w_S,\rho)} &=& \bigcup_{t^\frac{1}{2} \in \mathcal{D}_r}
W_{t^{1/2},\sigma(w_S,\rho)}^{(1)},\\
U^{(2)}_{\sigma(w_S,\rho)} &=& \bigcup_{t^\frac{1}{2} \in \mathcal{D}_r}
W_{t^{1/2},\sigma(w_S,\rho)}^{(2)}.
\end{eqnarray*}  
Then $\mathbf{M}^{\sigma(w_S,\rho)} =  U^{(1)}_{\sigma(w_S,\rho)} \cup
U^{(2)}_{\sigma(w_S,\rho)}$.  Let   
\[V^{(1)}_{\sigma(w_S,\rho)}  = \bigl\{  (w_B,t^\frac{1}{2}) \in 
\mathbb{C}^2 \; | \; t^\frac{1}{2} \in \mathcal{D}_r, \; w_B \in \pi_B 
\circ R_{t^{1/2},\sigma(w_S,\rho)}^{(1)} (W_{t^{1/2},
\sigma(w_S,\rho)}^{(1)}) \bigr\}, \] 
and 
\[V^{(2)}_{\sigma(w_S,\rho)}  = \bigl\{  (w_B,t^\frac{1}{2}) \in 
\mathbb{C}^2 \; | \; t^\frac{1}{2} \in \mathcal{D}_r, \;  w_B \in  
\pi_B \circ R_{t^{1/2}, \sigma(w_S,\rho)}^{(2)} (W_{t^{1/2},
\sigma(w_S,\rho)}^{(2)}) \bigr\} .\]
Obviously $V^{(k)}_{\sigma(w_S,\rho)}$ is an open set in $\mathbb{C}^2$ 
for $k=1,2$.  We define 
\begin{eqnarray*}
\beta^{(k)}_{\sigma(w_S,\rho)} : U^{(k)}_{\sigma(w_S,\rho)} 
&\longrightarrow& V^{(k)}_{\sigma(w_S,\rho)} \\
p &\mapsto& (\pi_B \circ R_{t^{1/2}, \sigma(w_S,\rho)}^{(k)} (p), 
t^\frac{1}{2})
\end{eqnarray*}
for $p \in W_{t^{1/2},\sigma(w_S,\rho)}^{(k)}$.  It is clear that the 
$\beta^{(k)}_{\sigma(w_S,\rho)}$ are bijections.  Thus 
\[\bigl\{(U^{(k)}_{\sigma(w_S,\rho)}, \beta^{(k)}_{\sigma(w_S,\rho)})
\bigr\}_{k = 1,2}\] 
is a local coordinate system for $\mathbf{M}^{\sigma(w_S,\rho)}$ with 
the coordinate transition function given by   
\[\beta^{(2)}_{\sigma(w_S,\rho)} \circ (\beta^{(1)}_{\sigma(w_S,\rho)})^{-1} 
(w_B,t^\frac{1}{2}) = \Bigl( R_{t^{1/2},\sigma(w_S,\rho)}^{(2)} \circ 
(R_{t^{1/2},\sigma(w_S,\rho)}^{(1)})^{-1} (w_B), t^\frac{1}{2} \Bigr)\]
for $w_B \in \beta^{(2)}_{\sigma(w_S,\rho)} (U^{(1)}_{\sigma(w_S,\rho)} 
\cap U^{(2)}_{\sigma(w_S,\rho)})$.

This gives a complex manifold structure to $\mathbf{M}^{\sigma(w_S,\rho)}$.  
It is clear that the projection {}from $\mathbf{M}^{\sigma(w_S,\rho)}$ to 
$\mathcal{D}_r$ is complex analytic and the rank of the Jacobian is one.  
Thus $\mathfrak{M}^{\sigma(w_S,\rho)}  = \{M_{t^{1/2}}^{\sigma(w_S,\rho)} \; 
| \; t^{1/2} \in \mathcal{D}_r \}$ is a complex analytic family.  Since 
all the $M_{t^{1/2}}^{\sigma(w_S,\rho)}$ are complex analytically isomorphic 
to $S\hat{\mathbb{C}}_B = \mathbb{C} \cup \{\infty\}$, they are all complex 
analytically isomorphic to each other.  By the Fischer-Grauert Theorem,
$\mathfrak{M}^{\sigma(w_S,\rho)} $ is locally trivial for each $(w_S,\rho)
\in (\Lambda_\infty)_S$.  Thus given a section $\sigma_{t^{1/2}}(w_S,\rho)$,
for every $t_0^{1/2} \in \mathcal{D}_r$, there exists a neighborhood of 
$t_0^{1/2}$, denoted $\Delta_{t_0^{1/2}}^{\sigma(w_S,\rho)} \subset 
\mathcal{D}_r$, such that complex analytically
\[\bigcup_{t^{1/2} \in \Delta_{t_0^{1/2}}^{\sigma(w_S,\rho)}} 
M_{t^{1/2}}^{\sigma(w_S,\rho)} \cong M_{t_0^{1/2}}^{\sigma(w_S,\rho)} \times 
\Delta_{t_0^{1/2}}^{\sigma(w_S,\rho)} .\]
Thus there exist complex analytic isomorphisms   
\[\gamma_{t^{1/2}}^{\sigma(w_S,\rho)} : M_{t^{1/2}}^{\sigma(w_S,\rho)} \longrightarrow
M_{t_0^{1/2}}^{\sigma(w_S,\rho)}\]   
for any $t^{1/2} \in \Delta_{t_0^{1/2}}^{\sigma(w_S,\rho)}$ such that if 
we use the local coordinates of $M_{t^{1/2}}^{\sigma(w_S,\rho)}$ and 
$M_{t_0^{1/2}}^{\sigma(w_S,\rho)}$ to express 
$\gamma_{t^{1/2}}^{\sigma(w_S,\rho)}$, then
$\gamma_{t^{1/2}}^{\sigma(w_S,\rho)}$ is also  analytic in $t^{1/2}$.  

For every $M_{t^{1/2}}^{\sigma(w_S,\rho)}$, we already have a superconformal
isomorphism $F_{t^{1/2}}^{\sigma(w_S,\rho)}$ {}from $M_{t^{1/2}}^{\sigma(w_S,\rho)}$ to
$\sigma_{(w_S,\rho)}(S\hat{\mathbb{C}}_B)$ given by
\begin{eqnarray*}
F_{t^{1/2}}^{\sigma(w_S,\rho)} (p) = \left\{\begin{array}{ll} 
       F_{t^{1/2}}^{(1)}  (p) &
          \mbox{for $p \in W_{t^{1/2},\sigma(w_S,\rho)}^{(1)}$} ,\\ 
\\
      F_{t^{1/2}}^{(2)} (p) &
          \mbox{for $p \in W_{t^{1/2},\sigma(w_S,\rho)}^{(2)}$}.\\ 
\end{array}  \right.
\end{eqnarray*}
For $a \in (\bigwedge_\infty^0)^\times$, define the superprojective
transformation $T_a : S\hat{\mathbb{C}} \longrightarrow  S\hat{\mathbb{C}}$
by (\ref{T1}) and 
\begin{eqnarray*}
(T_a)_\sou : \mbox{$\bigwedge_\infty$} &\longrightarrow&  \mbox{$\bigwedge_\infty$}  \\
(w,\rho) &\mapsto& (a^2w,a\rho) .
\end{eqnarray*}
Choose a point $(w_0)_B \in \mathbb{C}^\times$.  We
can always find some function $\alpha^{\sigma(w_S,\rho)} (t^{1/2}) \neq 0$ {}from $\Delta_{t_0^{1/2}}^{\sigma(w_S,\rho)}$ to $\mathbb{C}^\times$ such that 
in terms of the coordinate atlas $\{(U_\sou,\sou),(U_\nor,\nor)\}$ of 
$S\hat{\mathbb{C}}$ restricted to the section 
$\sigma(w_S,\rho)(S\hat{\mathbb{C}}_B)$, the function
\begin{multline*}
\sou \circ T_{\alpha^{\sigma(w_S,\rho)} (t_0^{1/2})} \circ 
F_{t_0^{1/2}}^{\sigma(w_S,\rho)} \circ \gamma_{t^{1/2}}^{\sigma(w_S,\rho)} \\
\circ (F_{t^{1/2}}^{\sigma(w_S,\rho)})^{-1} \circ 
T_{\alpha^{\sigma(w_S,\rho)}(t^{1/2})}^{-1} \circ \sou^{-1} ((w_0)_B + w_S,\rho) 
\end{multline*}  
is analytic in $t^{1/2}$ for $t^{1/2} \in \Delta_{t_0^{1/2}}^{\sigma(w_S,\rho)}$.  
Note that 
\[(T_{\alpha^{\sigma(w_S,\rho)}(t_0^{1/2})}\circ F_{t_0^{1/2}}^{\sigma(w_S,\rho)}) \circ
\gamma_{t^{1/2}}^{\sigma(w_S,\rho)} \circ (T_{\alpha^{\sigma(w_S,\rho)} (t^{1/2})} \circ
F_{t^{1/2}}^{\sigma(w_S,\rho)})^{-1}\] 
is a family of analytic isomorphisms {}from $\sigma(w_S,\rho)(S\hat{\mathbb{C}}_B) 
\equiv S\hat{\mathbb{C}}_B$ to itself.  Hence it must be a family of linear 
fractional transformations depending on $t^{1/2}$.  Any linear fractional 
transformation is determined by its  values on three complex variables.  
Furthermore, it is clear that if the values at these three points depend 
analytically on $t^{1/2}$, then the value at any point depends analytically 
on $t^{1/2}$.  Consider the three points $\sou^{-1}(0 + w_S,\rho)$, 
$\nor^{-1}(0)$ and $\sou^{-1} (w_0 + w_S,\rho) \in S\hat{\mathbb{C}}$.  Since 
\begin{multline*}
T_{\alpha^{\sigma(w_S,\rho)} (t_0^{1/2})} \circ F_{t_0^{1/2}}^{\sigma(w_S,\rho)}  
\circ
\gamma_{t^{1/2}}^{\sigma(w_S,\rho)}  \circ  (T_{\alpha^{\sigma(w_S,\rho)}  (t^{1/2})} 
\circ F_{t^{1/2}}^{\sigma(w_S,\rho)} )^{-1} \circ
\sou^{-1} (0) \\  = T_{\alpha^{\sigma(w_S,\rho)}  (t_0^{1/2})} \circ
F_{t_0^{1/2}}^{\sigma(w_S,\rho)}  \circ
\gamma_{t^{1/2}}^{\sigma(w_S,\rho)}  
\circ (R^{(2)}_{t^{1/2},\sigma(w_S,\rho) })^{-1} (0) ,
\end{multline*} 
is analytic in $t^{1/2}$, and the value of any superfunction at $(w_B+w_s,\rho)$
is the Taylor expansion about the body $w_B$, we have that 
\[T_{\alpha^{\sigma(w_S,\rho)}  (t_0^{1/2})} \circ F_{t_0^{1/2}}^{\sigma(w_S,\rho)} 
\circ \gamma_{t^{1/2}}^{\sigma(w_S,\rho)} \circ  (T_{\alpha^{\sigma(w_S,\rho)}  
(t^{1/2})} \circ F_{t^{1/2}}^{\sigma(w_S,\rho)})^{-1} 
\circ \sou^{-1} (0+w_S,\rho)\] 
is analytic in $t^{1/2}$.  Furthermore,
\begin{multline*}
T_{\alpha^{\sigma(w_S,\rho)} (t_0^{1/2})} \circ F_{t_0^{1/2}}^{\sigma(w_S,\rho)}  
\circ \gamma_{t^{1/2}}^{\sigma(w_S,\rho)} \circ (T_{\alpha^{\sigma(w_S,\rho)} 
(t^{1/2})} \circ F_{t^{1/2}}^{\sigma(w_S,\rho)} )^{-1} \circ
\nor^{-1}(0) \\ 
= T_{\alpha^{\sigma(w_S,\rho)}  (t_0^{1/2})} \circ F_{t_0^{1/2}}^{\sigma(w_S,\rho)}  
\circ \gamma_{t^{1/2}}^{\sigma(w_S,\rho)} \circ (R^{(1)}_{t^{1/2}, 
\sigma(w_S,\rho)})^{-1} (0) ,
\end{multline*}  
is analytic in $t^{1/2}$, and by our choice of $\alpha^{\sigma(w_S,\rho)}$,
\begin{multline*}
T_{\alpha^{\sigma(w_S,\rho)}  (t_0^{1/2})} \circ F_{t_0^{1/2}}^{\sigma(w_S,\rho)}  
\circ \gamma_{t^{1/2}}^{\sigma(w_S,\rho)} \circ (T_{\alpha^{\sigma(w_S,\rho)}  
(t^{1/2})} \circ F_{t^{1/2}}^{\sigma(w_S,\rho)} )^{-1} \\
\circ \sou^{-1} ((w_0)_B + w_S, \rho)
\end{multline*} 
is analytic in $t^{1/2}$.  Thus 
\[(T_{\alpha^{\sigma(w_S,\rho)}  (t_0^{1/2})} \circ
F_{t_0^{1/2}}^{\sigma(w_S,\rho)} ) \circ \gamma_{t^{1/2}}^{\sigma(w_S,\rho)}  
\circ (T_{\alpha^{\sigma(w_S,\rho)} (t^{1/2})} 
\circ F_{t^{1/2}}^{\sigma(w_S,\rho)} )^{-1}\] 
is analytic in $t^{1/2}$.  This implies that $T_{\alpha^{\sigma(w_S,\rho)}  
(t^{1/2})} \circ F_{t^{1/2}}^{\sigma(w_S,\rho)}$ is analytic in
$t^{1/2}$.   By the definition of $F_{t^{1/2}}^{\sigma(w_S,\rho)}$ and the 
normalization conditions that it satisfies, we have 
\begin{multline*}
\rho(\sou \circ T_{\alpha^{\sigma(w_S,\rho)} (t^{1/2})} \circ 
(F_{t^{1/2}}^{\sigma(w_S,\rho)})^{(1)}
\circ  (R^{(1)}_{t^{1/2},\sigma(w_S,\rho)})^{-1} (w, \rho))^0 \\
= \; \rho (\alpha^{\sigma(w_S,\rho)}(t^\frac{1}{2}) w + \; \mbox{terms of lower order in $w$}).
\end{multline*} 
Since $T_{\alpha^{\sigma(w_S,\rho)} (t^{1/2})} \circ 
F_{t^{1/2}}^{\sigma(w_S,\rho)}$ is analytic in $t^{1/2}$, its coefficients 
in $w$ and $\rho$ are also analytic in $t^{1/2}$.  Hence 
$\alpha^{\sigma(w_S,\rho)}  (t^{1/2})$ is analytic in $t^{1/2}$.  We conclude 
that $F_{t^{1/2}}^{\sigma(w_S,\rho)}$ is analytic in $t^{1/2}$.   Since $t_0^{1/2}$ 
is an arbitrary complex number in $\mathcal{D}_r$, we have that
$F_{t^{1/2}}^{\sigma(w_S,\rho)}$ is analytic in $t^{1/2}$ for any $t^{1/2} \in 
\mathcal{D}_r$. 

Following the same argument above with $\mathcal{D}_r$ replaced by
$\tilde{\mathcal{D}}_r$, we conclude that $F_{t^{1/2}}^{\sigma(w_S,\rho)}$ is 
analytic in $t^{1/2}$ for $0< |t^{1/2}| < r$.  Then following a similar
argument using the sections $\tau(w_S,\rho)$, we can prove that $F_{t^{1/2}}$
restricted to $M_{t^{1/2}}^{\tau(w_S,\rho)}$ is analytic in $t^{1/2}$ 
for $0< |t^{1/2}| < r$.  Then since
\[M_{t^{1/2}} = \bigcup_{(w_S,\rho) \in \mbox{$(\bigwedge_\infty)_S$}}
(M_{t^{1/2}}^{\sigma(w_S,\rho)}
\cup M_{t^{1/2}}^{\tau(w_S,\rho)}), \]
we have that $F_{t^{1/2}}$ is analytic in $t^{1/2}$ for every point 
$p \in M_{t^{1/2}}$.  We conclude that $F_{t^{1/2}}$ is analytic in 
$t^{1/2}$ for $0< |t^{1/2}| < r$, and thus $F_{t^{1/2}}^{(1)}$ and
$F_{t^{1/2}}^{(2)}$ are analytic in $t^{1/2}$ for $0< |t^{1/2}| < r$.

We now prove the second statement of the proposition.  Let
$(H^{(2)})^0(w,\rho)$ and $(H^{(2)})^1(w,\rho)$ be the even and odd
superfunction components of $H^{(2)}(w,\rho)$, and let
\[H_{t^{1/2}}^{(2)} (w,\rho) = \left( t^{-1} (H^{(2)})^0(t^{-1}w, 
t^{-\frac{1}{2}} \rho), t^{-\frac{1}{2}} (H^{(2)})^1(t^{-1}w, 
t^{-\frac{1}{2}}\rho) \right) \]   
for  $t^{1/2} \neq 0$.  For $t^{1/2} = 0$, we define
\[H_0^{(2)} (w,\rho) = \Bigl( \frac{1}{w}, \frac{i\rho}{w} \Bigr) .\]
Then $H_{t^{1/2}}^{(2)} (w,\rho)$ is analytic in $t^{1/2}$ for
$|t^{1/2}| \leq 1$.

Let \footnote{There is a misprint in the analogous nonsuper case to
the definition of $Q_2(t)$ given in \cite{H book}.  On p.85, in 
\cite{H book}, the coordinate at $\infty$ of the sphere $Q_2(t)$ should
be $\tilde{E}^{-1}(f_t^{(2)}(1/x))$, not $\tilde{E}(f_t^{(2)}(1/x))$
as stated.}
\begin{multline*}
Q_2(t^\frac{1}{2}) = \bigl((z_1', \theta_1'),...,(z_{n-1}',
\theta_{n-1}'); \bigl((\tilde{E}^{-1})^0(H_{t^{1/2}}^{(2)} \circ I^{-1} (w,\rho)), \\
i(\tilde{E}^{-1})^1(H_{t^{1/2}}^{(2)} \circ I^{-1} (w,\rho))\bigr), 
(b^{(1)}, B^{(1)}, N^{(1)}),...,(b^{(n)}, B^{(n)}, N^{(n)}) \bigr) 
\end{multline*} 
It is easy to see that the $i$-th tube of $Q_1$ can be sewn with the
0-th tube of $Q_2(t^{1/2})$ when $|t^{1/2}| \leq 1$.  By
Proposition \ref{actual sewing}, there exist
$\tilde{F}_{t^{1/2}}^{(1)}$ and $\tilde{F}_{t^{1/2}}^{(2)}$ such 
that
\begin{eqnarray*}
\tilde{F}_{t^{1/2},\nor}^{(1)} (0) &=& (0),\\
\lim_{w \rightarrow \infty} \frac{\partial}{\partial \rho}
(\tilde{F}_{t^{1/2},\sou}^{(1)})^1(w,\rho) &=& 1 ,\\
\tilde{F}_{t^{1/2},\sou}^{(2)} (0) &=& (0),
\end{eqnarray*}
and in some region
\[\tilde{F}_{t^{1/2},\sou}^{(1)} (w,\rho) = 
\tilde{F}_{t^{1/2},\sou}^{(2)} \circ (H_{t^{1/2}}^{(2)})^{-1} \circ 
I \circ H^{(1)} \circ s_{(z_i,\theta_i)} (w,\rho) \]  
holds for $|t^{1/2}| \leq 1$.  Moreover $\tilde{F}_{t^{1/2}}^{(1)}$ 
and $\tilde{F}_{t^{1/2}}^{(2)}$ are unique.  Using the same method as 
that used in the proof of the analyticity of $F_{t^{1/2}}^{(1)}$ and 
$F_{t^{1/2}}^{(2)}$, we can show that $\tilde{F}_{t^{1/2}}^{(1)}$ and 
$\tilde{F}_{t^{1/2}}^{(2)}$ are analytic in $t^{1/2}$ for $|t^{1/2}| 
\leq 1$.  But {}from   
\[F_{t^{1/2},\sou}^{(1)} (w,\rho) = F_{t^{1/2},\sou}^{(2)} \circ 
(H^{(2)})^{-1} \circ I \circ H_{t^{1/2}}^{(1)} \circ s_{(z_i,
\theta_i)} (w,\rho) \]  
and
\begin{multline*}
H_{t^{1/2}}^{(1)} \circ s_{(z_i,\theta_i)} (w,\rho)  \\
= \; \left( t^{-1} (H^{(1)})^0(w - z_i - \rho \theta_i,\rho
- \theta_i),  t^{-\frac{1}{2}} (H^{(1)})^1(w - z_i - \rho
\theta_i,\rho - \theta_i) \right) ,
\end{multline*} 
we see that $F_{t^{1/2},\sou}^{(1)}(w,\rho)$ and 
$F_{t^{1/2},\sou}^{(2)}(t^{-1}w, t^{-{1/2}}\rho)$ satisfy the above 
equation for $\tilde{F}_{t^{1/2},\sou}^{(1)}$ and $\tilde{F}_{t^{1/2},
\sou}^{(2)}$.  By uniqueness, 
\begin{eqnarray*}
\tilde{F}_{t^{1/2},\sou}^{(1)} (w,\rho) &=& F_{t^{1/2},\sou}^{(1)} 
(w,\rho)\\ 
\tilde{F}_{t^{1/2},\sou}^{(2)} (w,\rho) &=& F_{t^{1/2},\sou}^{(2)}
(t^{-1}w, t^{-\frac{1}{2}}\rho) ,
\end{eqnarray*} 
i.e., $F_{t^{1/2},\sou}^{(2)} (w,\rho) = \tilde{F}_{t^{1/2},\sou}^{(2)}
(t w, t^{1/2} \rho)$.  Since $\tilde{F}_{t^{1/2}}^{(1)}$ and
$\tilde{F}_{t^{1/2}}^{(2)}$ are analytic in $t^{1/2}$ not only 
when $0 < |t^{1/2}| \leq 1$, but also when $t^{1/2} = 0$, it 
must be that $t^{1/2} = 0$ is a removable singularity of 
$F_{t^{1/2}}^{(1)}$ and $F_{t^{1/2}}^{(2)}$.  Finally, since 
$\rho (\tilde{F}_{t^{1/2}}^{(2)}(w,\rho))^0 = \rho \cdot O(w)$, and 
$\rho (\tilde{F}_{t^{1/2}}^{(2)}(w,\rho))^1 = \rho \cdot O(w)$, 
and for any $t^{1/2}$ with $|t^{1/2}| \leq 1$, the superfunction 
$\tilde{F}_{t^{1/2}}^{(2)}(w,\rho)$ is nonzero for $(w,\rho) \neq 0$, we 
see that $t^{1/2} = 0$ is a second-order zero of $(F_{t^{1/2}}^{(2)}(w,\rho))^0$ 
and a first-order zero of $(F_{t^{1/2}}^{(2)}(w,\rho))^1$. \end{proof}

{}From Proposition \ref{normal order}, we have 
the formal series
\[\Psi_j(t^\frac{1}{2}) (t^{-\frac{1}{2}} \asqrt^{(i)}, A^{(i)}, M^{(i)},
B^{(0)}, N^{(0)}) , \] 
in $\bigwedge_\infty [[t^{1/2}]]$ for $j \in \frac{1}{2} \mathbb{Z}$, 
$(A^{(i)}, M^{(i)}), (B^{(0)}, N^{(0)})  \in \bigwedge_\infty^\infty$,
and $\asqrt^{(i)} \in (\bigwedge_\infty^0)^\times$.  Let
$\tilde{\Psi}_j(t^{1/2}) = \tilde{\Psi}_j(t^{1/2})
(t^{-{1/2}}\asqrt^{(i)}, A^{(i)}, M^{(i)}, B^{(0)}, N^{(0)})$, for $j \in
\frac{1}{2} \mathbb{Z}$, be defined by 
\[\Bigl(1, \bigl\{\tilde{\Psi}_{-j}(t^\frac{1}{2}), -i\tilde{\Psi}_{-j +
\frac{1}{2}}(t^\frac{1}{2}) \bigr\}_{j \in \Z} \Bigr) = \hat{E}^{-1}\bigl(I \circ 
F_{t^{1/2},\sou}^{(1)} \circ s_{(z_i,\theta_i)}^{-1} \circ I^{-1} (x,\varphi)\bigr) \] 
\[\Bigl(t^{-\frac{1}{2}}\asqrt^{(i)} \exp(\tilde{\Psi}_0(t^\frac{1}{2})),
\bigl\{\tilde{\Psi}_j(t^\frac{1}{2}), \tilde{\Psi}_{j -
\frac{1}{2}}(t^\frac{1}{2}) \bigr\}_{j \in \Z} \Bigr) =
\hat{E}^{-1}\bigl(F_{t^{1/2},\sou}^{(2)} (x,\varphi)\bigr) .\] 
Then we have the following proposition.

\begin{prop}\label{Psi converges}
Let $Q_1 \in SK(m)$ and $Q_2 \in SK(n)$ be given by
\begin{multline*}
Q_1 = \bigl((z_1, \theta_1),...,(z_{m-1}, \theta_{m-1});
(A^{(0)},M^{(0)}), (\asqrt^{(1)}, A^{(1)}, M^{(1)}), \\
...,(\asqrt^{(m)}, A^{(m)},M^{(m)})\bigr) 
\end{multline*} 
and
\begin{multline*}
Q_2 = \bigl((z_1', \theta_1'),...,(z_{n-1}', \theta_{n-1}'); (B^{(0)},
N^{(0)}), (\bsqrt^{(1)}, B^{(1)}, N^{(1)}),\\
...,(\bsqrt^{(n)}, B^{(n)}, N^{(n)})\bigr)
\end{multline*} 
for $m \in \Z$ and $n \in \mathbb{N}$.  If $Q_1 \; _i\infty_0 \; Q_2$
exists, then the series $\Psi_j(t^{1/2})$, for $j \in \frac{1}{2}
\mathbb{Z}$, are convergent when $|t^{1/2}| \leq 1$, and the values 
of these convergent series are equal to $\tilde{\Psi}_j(t^{1/2})$. 
\end{prop}

\begin{proof}  By Proposition \ref{use of FG}, $F_{t^{1/2}}^{(1)}$
and $F_{t^{1/2}}^{(2)}$ are analytic in $t^{1/2}$ for $|t^{1/2}| 
\leq 1$, and therefore $F_{t^{1/2},\sou}^{(1)} \circ 
s_{(z_i,\theta_i)}^{-1} (w,\rho)$ is also analytic in $t^{1/2}$ for 
$|t^{1/2}| \leq 1$.   Thus $F_{t^{1/2},\sou}^{(1)} \circ 
s_{(z_i,\theta_i)}^{-1} (w,\rho)$ and $F_{t^{1/2},\sou}^{(2)}$ can 
be expanded as power series in $t^{1/2}$.  The functions 
$\tilde{\Psi}_j (t^{1/2})$ as polynomials in the $w$ and $\rho$ 
coefficients of $F_{t^{1/2},\sou}^{(1)} \circ s_{(z_i,\theta_i)}^{-1} 
(w,\rho)$ and $F_{t^{1/2},\sou}^{(2)}(w,\rho)$ can also be expanded 
as power series  in $t^{1/2}$.  Since  $F_{t^{1/2},\sou}^{(1)} \circ 
s_{(z_i,\theta_i)}^{-1} (w,\rho)$ and $F_{t^{1/2},\sou}^{(2)}(w,\rho)$
satisfy the sewing equation
\[F_{t^{1/2},\sou}^{(1)} (w,\rho) = F_{t^{1/2},\sou}^{(2)} \circ 
(H_{t^{1/2}}^{(2)})^{-1} \circ I \circ H_{t^{1/2}}^{(1)} \circ 
s_{(z_i,\theta_i)} (w,\rho) \]
or equivalently
\[F_{t^{1/2},\sou}^{(1)} \circ s_{(z_i,\theta_i)}^{-1} (w,\rho) = 
F_{t^{1/2},\sou}^{(2)} \circ (H_{t^{1/2}}^{(2)})^{-1} \circ I \circ 
H_{t^{1/2}}^{(1)} (w,\rho) \]
and the obvious boundary conditions, the formal series 
$F_{t^{1/2},\sou}^{(1)} \circ s_{(z_i,\theta_i)}^{-1} (x,\varphi)$ and
$F_{t^{1/2},\sou}^{(2)}(x,\varphi)$ in $\bigwedge_\infty [[x,x^{-1}]]
[\varphi][[t^{1/2}]]$ corresponding to the expansions of $F_{t^{1/2},
\sou}^{(1)} \circ s_{(z_i,\theta_i)}^{-1} (w,\rho)$ and $F_{t^{1/2},
\sou}^{(2)}(w,\rho)$, respectively, satisfy the equation
\begin{equation}\label{another sewing equation}
F_{t^{1/2},\sou}^{(1)} \circ s_{(z_i,\theta_i)}^{-1} (x,\varphi) =  
F_{t^{1/2},\sou}^{(2)} \circ (H_{t^{1/2}}^{(2)})^{-1} \circ I \circ  
H_{t^{1/2}}^{(1)} (x,\varphi)
\end{equation} 
in $\bigwedge_\infty [[x,x^{-1}]][\varphi][[t^{1/2}]]$ and the 
corresponding formal boundary conditions.  Note that the coefficients
of the right-hand side of (\ref{another sewing equation}) are, in general,
infinite sums.  Thus $F_{t^{1/2},\sou}^{(1)} \circ s_{(z_i,
\theta_i)}^{-1} (x,\varphi)$ and $F_{t^{1/2},\sou}^{(2)}(x,\varphi)$ 
satisfying (\ref{another sewing equation}) means that the coefficients of 
the right-hand side are absolutely convergent to the coefficients of the
left-hand side.

Note that equation (\ref{another sewing equation}) and the corresponding
formal boundary conditions can be obtained {}from the formal sewing
equation and formal boundary conditions in Theorem \ref{uniformization}
by substituting $A_j$, $M_{j - 1/2}$, $B_j$, $N_{j - 1/2}$, and $t^{-1/2} 
\asqrt^{(i)}$ for $\mathcal{A}_j$, $\mathcal{M}_{j - 1/2}$, 
$\mathcal{B}_j$, $\mathcal{N}_{j - 1/2}$, and $\alpha_0^{1/2}$, 
respectively, for $j \in \Z$.  Since the solution of the formal sewing
equation and the formal boundary conditions in Theorem \ref{uniformization}
is unique, the solution $F_{t^{1/2},\sou}^{(1)} \circ s_{(z_i,
\theta_i)}^{-1} (x,\varphi)$ and $F_{t^{1/2},\sou}^{(2)}(x,\varphi)$ to
formula (\ref{another sewing equation}) and the corresponding boundary
conditions can be obtained by substituting $A_j$, $M_{j - 1/2}$, $B_j$, 
$N_{j - 1/2}$, and $t^{-1/2} \asqrt^{(i)}$ for $\mathcal{A}_j$, 
$\mathcal{M}_{j - 1/2}$, $\mathcal{B}_j$, $\mathcal{N}_{j - 1/2}$, and 
$\alpha_0^{1/2}$, respectively, for $j \in \Z$, into the solution of
the formal sewing equation and the formal boundary conditions given
in Theorem \ref{uniformization}.  Thus we have \footnote{There are a couple
of misprints in the analogous nonsuper version to equation (\ref{label
for footnote}) and the setting for nonsuper version of the proof of 
Proposition \ref{Psi converges} given in \cite{H book}.  In the 
proof of Proposition 3.4.5 in \cite{H book}, in the last paragraph on
p.87, in both instances it should be stated that equation (3.4.2) and 
the corresponding formal boundary conditions can be obtained {}from the 
formal sewing equation and formal boundary conditions in Theorem 2.2.4 by 
substituting $A_j$, $B_j$, $j \in \mathbb{Z}_+$ and $t^{-1}a_0^{(i)}$,
not $A_j$, $t^jB_j$, $j \in \mathbb{Z}_+$ and $t^{-j}a_0^{(i)}$ as stated, 
for $\mathcal{A}_j$, $\mathcal{B}_j$, $j \in \mathbb{Z}_+$ and $\alpha_0$, 
respectively.  In addition, the first displayed equation on p.88, should
read $1/F_t^{(1)}(x^{-1} + z_i) = \hat{E}(1,\{\Psi_{-j}(t)\}_{j \in \mathbb{Z}_+})$,
not $F_t^{(1)}(x + z_i) = \hat{E}(1,\{\Psi_{-j}(t)\}_{j \in \mathbb{Z}_+})$
as stated. }
\begin{multline}\label{label for footnote} 
\hat{E}\Bigl(1, \bigl\{\tilde{\Psi}_{-j}(t^\frac{1}{2}), -i\tilde{\Psi}_{-j +
\frac{1}{2}}(t^\frac{1}{2}) \bigr\}_{j \in \Z} \Bigr) (x,\varphi) \\
= \; I \circ F_{t^{1/2},\sou}^{(1)} \circ s_{(z_i,\theta_i)}^{-1} 
\circ I^{-1} (x,\varphi) 
\end{multline}
and
\begin{equation}
\hat{E} \Bigl(t^{-\frac{1}{2}}\asqrt^{(i)} \exp(\tilde{\Psi}_0(t^\frac{1}{2})),
\bigl\{\tilde{\Psi}_j(t^\frac{1}{2}), \tilde{\Psi}_{j -
\frac{1}{2}}(t^\frac{1}{2}) \bigr\}_{j \in \Z} \Bigr) = F_{t^{1/2},\sou}^{(1)} 
(x,\varphi) .
\end{equation} 
By the definition of $\tilde{\Psi}_j (t^{1/2})$, for $j \in \frac{1}{2}
\mathbb{Z}$, the expansion of $\tilde{\Psi}_j (t^{1/2})$ is equal to
$\Psi_j(t^{1/2})$ for $j \in \frac{1}{2}\mathbb{Z}$, i.e., $\Psi_j(t^{1/2})$ 
is convergent to $\tilde{\Psi}_j(t^{1/2})$ for $|t^{1/2}| \leq1$. 
\end{proof}

\section[A Neveu-Schwarz algebra structure on the 
tangent space]{An $N=1$ Neveu-Schwarz algebra
structure of central charge zero on the supermeromorphic tangent 
space of $SK(1)$ at its identity}

The super-moduli space $SK(1)$ of superspheres with $1 +1$ tubes with 
the sewing operation is a partial monoid.  But in general, elements 
of $SK(1)$ do not have inverses.  In Lie theory, one uses invariant 
vector fields of a Lie group to define the Lie algebra.  We will 
define a bracket operation on a subspace $\hat{T}_e SK(1)$ of the 
supermeromorphic tangent space $T_e SK(1)$ of $SK(1)$ at the identity 
$e$.  We will then show that $\hat{T}_e SK(1)$ with this bracket operation 
is the Neveu-Schwarz algebra with central charge zero.  

Let $\alpha_0$, $\mathcal{A}_j^{(0)}$, 
$\mathcal{A}_j^{(1)}$, $\beta_0$, $\mathcal{B}_j^{(0)}$, and 
$\mathcal{B}_j^{(1)}$ be even formal variables, for $j \in \Z$; and 
let $\mathcal{M}_{j- 1/2}^{(0)}$, $\mathcal{M}_{j- 1/2}^{(1)}$, 
$\mathcal{N}_{j-1/2}^{(0)}$,  and $\mathcal{N}_{j-1/2}^{(1)}$ be odd 
formal variables, for $j \in \Z$.  For convenience, we will let use
the notation
\begin{eqnarray*}
\mathcal{A}\mathcal{M}^{(01)} &=& ((\mathcal{A}^{(0)},
\mathcal{M}^{(0)}), (\alpha_0^{1/2}, \mathcal{A}^{(1)}, 
\mathcal{M}^{(1)}))\\
\mathcal{B}\mathcal{N}^{(01)} &=& ((\mathcal{B}^{(0)}, 
\mathcal{N}^{(0)}), (\beta_0^{1/2},\mathcal{B}^{(1)}, 
\mathcal{N}^{(1)})) .
\end{eqnarray*}   

Define 
\begin{multline*}
\lambda_0^\frac{1}{2} (\mathcal{A}\mathcal{M}^{(01)},
\mathcal{B}\mathcal{N}^{(01)}), \Phi^{(0)}_j (\mathcal{A}\mathcal{M}^{(01)},
\mathcal{B}\mathcal{N}^{(01)}), \Phi^{(1)}_j (\mathcal{A}\mathcal{M}^{(01)},
\mathcal{B}\mathcal{N}^{(01)}) \\
\in \mathbb{C}[\mathcal{A},\mathcal{M},
\mathcal{B},\mathcal{N}][[\alpha_0^{\frac{1}{2}},\beta_0^{\frac{1}{2}},
\alpha_0^{-\frac{1}{2}},\beta_0^{-\frac{1}{2}}]] ,
\end{multline*}
for $j\in\frac{1}{2}\Z$, by
\begin{multline*}
\Bigl( \left\{ \Phi^{(0)}_j (\mathcal{A}\mathcal{M}^{(01)},
\mathcal{B}\mathcal{N}^{(01)}), -i \Phi_{j - \frac{1}{2}}^{(0)}
(\mathcal{A}\mathcal{M}^{(01)}, \mathcal{B}\mathcal{N}^{(01)}) 
\right\}_{j \in \Z} \Bigr)\\ 
= \; \tilde{E}^{-1} \Bigl( \tilde{E} \left(\mathcal{A}^{(0)}, -i \mathcal{M}^{(0)}
\right) \circ  \tilde{E} \Bigl( \left\{-\Psi_{-j} (\alpha_0^\frac{1}{2}, 
\mathcal{A}^{(1)}, \mathcal{M}^{(1)},\mathcal{B}^{(0)}, \mathcal{N}^{(0)}), 
\right. \Bigr. \Bigr. \\
\Bigl. \Bigl. \left. i \Psi_{- j + \frac{1}{2}}(\alpha_0^\frac{1}{2},
\mathcal{A}^{(1)}, \mathcal{M}^{(1)}, \mathcal{B}^{(0)}, \mathcal{N}^{(0)} )
\right\}_{j \in \Z} \Bigr) (x,\varphi) \Bigr) 
\end{multline*} 
and 
\begin{multline*}
\Bigl( \lambda_0^\frac{1}{2} (\mathcal{A}\mathcal{M}^{(01)},
\mathcal{B}\mathcal{N}^{(01)}), \left\{\Phi_j^{(1)}(\mathcal{A}\mathcal{M}^{(01)},
\mathcal{B}\mathcal{N}^{(01)}), \Phi_{j- \frac{1}{2}}^{(1)}(\mathcal{A}
\mathcal{M}^{(01)}, \mathcal{B} \mathcal{N}^{(01)}) \right\}_{j \in \Z}
\Bigr)\\
= \; \hat{E}^{-1} \Bigl( \hat{E} \bigl(\beta_0^\frac{1}{2}, \mathcal{B}^{(1)},
\mathcal{N}^{(1)} \bigr) \circ \hat{E} \Bigl(\alpha_0^\frac{1}{2} \exp(-\Psi_0
(\alpha_0^\frac{1}{2}, \mathcal{A}^{(1)}, \mathcal{M}^{(1)}, \mathcal{B}^{(0)},
\mathcal{N}^{(0)})), \Bigr. \Bigr.\\
\Bigl. \Bigl. \left\{-\Psi_j(\alpha_0^\frac{1}{2}, \mathcal{A}^{(1)}, 
\mathcal{M}^{(1)}, \mathcal{B}^{(0)}, \mathcal{N}^{(0)}), - \Psi_{j - 
\frac{1}{2}} (\alpha_0^\frac{1}{2}, \mathcal{A}^{(1)}, \mathcal{M}^{(1)}, 
\mathcal{B}^{(0)}, \mathcal{N}^{(0)})  \right\}_{j \in \Z} \Bigr) \\
(x,\varphi) \Bigr). 
\end{multline*}
Formally, these series give the coordinates of the supersphere 
\[\mathcal{A}\mathcal{M}^{(01)} \; _1\infty_0 \; \mathcal{B}\mathcal{N}^{(01)} \in 
SK(2)\] 
at infinity and zero, respectively.  For convenience, we denote
the resulting supersphere  
\begin{multline*}
\Bigl( \Bigl( \left\{ \Phi^{(0)}_j (\mathcal{A}\mathcal{M}^{(01)}, \mathcal{B}
\mathcal{N}^{(01)}), \Phi_{j - \frac{1}{2}}^{(0)}(\mathcal{A}
\mathcal{M}^{(01)}, \mathcal{B}\mathcal{N}^{(01)}) \right\}_{j \in \Z} \Bigr),
\Bigr. \\
\Bigl. \Bigl( \lambda_0^\frac{1}{2} (\mathcal{A}\mathcal{M}^{(01)}, \mathcal{B}
\mathcal{N}^{(01)}),  \left\{\Phi_j^{(1)}(\mathcal{A}\mathcal{M}^{(01)}, \mathcal{B}
\mathcal{N}^{(01)}), \Phi_{j- \frac{1}{2}}^{(1)}(\mathcal{A}\mathcal{M}^{(01)},
\mathcal{B} \mathcal{N}^{(01)}) \right\}_{j \in \Z} \Bigr) \Bigr) 
\end{multline*}
by
\[\Phi^{(01)} (\mathcal{A}\mathcal{M}^{(01)}, \mathcal{B}\mathcal{N}^{(01)} ) . \]
Writing 
\[\Phi^{(k)}_j = \Phi^{(k)}_j (\mathcal{A}\mathcal{M}^{(01)},\mathcal{B}
\mathcal{N}^{(01)}), \] 
for $k=0,1$ and $j \in \frac{1}{2}\Z$, 
\[\lambda_0^{\frac{1}{2}} = \lambda_0^{\frac{1}{2}}(\mathcal{A}\mathcal{M}^{(01)}, 
\mathcal{B} \mathcal{N}^{(01)}), \]
and 
\[\Psi_j = \Psi_j(\alpha_0^{\frac{1}{2}}, \mathcal{A}^{(1)}, \mathcal{M}^{(1)},
\mathcal{B}^{(0)}, \mathcal{N}^{(0)}), \]
for $j \in \frac{1}{2}\mathbb{Z}$, we can give the local coordinates of $\Phi^{(01)}
(\mathcal{A}\mathcal{M}^{(01)}, \mathcal{B}\mathcal{N}^{(01)})$ at infinity and 
zero explicitly by
\begin{multline*}
\exp\Bigl(\sum_{j \in \Z} \bigl(\Phi^{(0)}_j L_{-j}(x,\varphi) + 
\Phi_{j - \frac{1}{2}}^{(0)} G_{-j + \frac{1}{2}} (x,\varphi) \bigr)\Bigr) \cdot
\Bigl(\frac{1}{x},\frac{i\varphi}{x}\Bigr)\\  
= \; \exp\Bigl(  - \! \sum_{j \in \Z} \bigl(\Psi_{-j}L_{-j} (x,\varphi)  
+ \Psi_{- j + \frac{1}{2}} G_{-j + \frac{1}{2}}(x,\varphi) \Bigr) \cdot\\
\exp\Bigl( \sum_{j \in \Z} \bigl( \mathcal{A}_j^{(0)} L_{-j}(x,\varphi) + 
\mathcal{M}_{j - \frac{1}{2}}^{(0)} G_{-j + \frac{1}{2}}(x,\varphi) \bigr) \Bigr) \cdot
\Bigl(\frac{1}{x},\frac{i\varphi}{x}\Bigr)
\end{multline*}
for the local coordinate at infinity, and
\begin{multline}\label{explicitly at zero for misprint}
\exp \Bigl( - \! \sum_{j \in \Z} \bigl( \Phi_j^{(1)}L_j(x,\varphi) +  \Phi_{j-
\frac{1}{2}}^{(1)} G_{j - \frac{1}{2}}(x,\varphi) \bigr) \Bigr) 
\cdot ( \lambda_0^\frac{1}{2})^{-L_0(x,\varphi)} \cdot (x,\varphi)\\
= \; \exp \Bigl( \sum_{j\in\Z} \bigl( \Psi_j L_j(x,\varphi) + \Psi_{j - 
\frac{1}{2}} G_{j - \frac{1}{2}}(x,\varphi) \bigr) \Bigr)
\cdot (\alpha_0^\frac{1}{2})^{-L_0(x,\varphi)}
\cdot \exp(\Psi_0 L_0(x,\varphi)) \cdot\\
\exp \Bigl( -\sum_{j \in \Z} \bigl(\mathcal{B}_j^{(1)} L_j(x,\varphi) + 
\mathcal{N}_{j - \frac{1}{2}}^{(1)} G_{j - \frac{1}{2}}(x,\varphi) \bigr)\Bigr) \cdot 
(\beta_0^\frac{1}{2})^{-L_0(x,\varphi)} \cdot (x,\varphi) ,
\end{multline} 
for the local coordinate at zero. \footnote{There is a misprint to 
the analogous nonsuper case to equation (\ref{explicitly at zero for misprint}) 
given in \cite{H book}.  On p.90 of \cite{H book}, the right-hand side of the
equation giving the coordinate at zero which is equal to $\lambda_0(\mathcal{A}^{(01)},
\mathcal{B}^{(01)})e_{\Lambda^{(1)}(\mathcal{A}^{(01)},\mathcal{B}^{(01)})}
(x)x$ should be acting on $x$.}
It is obvious that 
\begin{equation} \label{Phi boundary conditions} 
\Phi^{(01)} (\mathcal{A}\mathcal{M}^{(01)}, e ) = \mathcal{A}\mathcal{M}^{(01)}
, \quad \mbox{and} \quad \Phi^{(01)} (e , \mathcal{B}\mathcal{N}^{(01)} ) = 
\mathcal{B}\mathcal{N}^{(01)}  . 
\end{equation}

Let $\hat{T}_e SK(1)$ be the subspace of $T_e SK(1)$ (the
supermeromorphic tangent space of $SK(1)$ at $e$) consisting of all
finite linear combinations of $- \left. \frac{\partial}{\partial a}
\right|_e$, $- \Bigl. \frac{\partial}{\partial A_j^{(0)}}\Bigr|_e$,
$- \Bigl. \frac{\partial}{\partial M_{j - 1/2}^{(0)}}
\Bigr|_e$, $- \Bigl. \frac{\partial}{\partial A_j^{(1)}}\Bigr|_e$,
and $\Bigl. - \frac{\partial}{\partial M_{j - 1/2}^{(1)}}\Bigr|_e$, 
for $j \in \Z$.  Let $F \in SD(0)$.  We will use the notation 
\[F(\Phi^{(01)} (\mathcal{A}\mathcal{M}^{(01)}, \mathcal{B}
\mathcal{N}^{(01)})) \Bigr|_e  = F(\Phi^{(01)} (\mathcal{A}\mathcal{M}^{(01)}, 
\mathcal{B}\mathcal{N}^{(01)})) \Bigr|_{\begin{scriptsize}\begin{array}{l}
\mathcal{A}\mathcal{M}^{(01)} = \\
\mathcal{B}\mathcal{N}^{(01)} = e  \end{array} \end{scriptsize}}. \] 
We define a bracket operation on $\hat{T}_e SK(1)$ as 
follows: For $j,k \in \mathbb{Z}_+$, let
\begin{multline*}
\biggl[\biggl. \frac{\partial}{\partial A_j^{(0)}}\biggr|_e,
\biggl. \frac{\partial}{\partial A_k^{(0)}}\biggr|_e  \biggr] F \\
= \biggl. \biggl( \! \biggl( \frac{\partial}{\partial \mathcal{A}_j^{(0)}}
\frac{\partial}{\partial \mathcal{B}_k^{(0)}} - \frac{\partial}{\partial
\mathcal{A}_k^{(0)}} \frac{\partial}{\partial \mathcal{B}_j^{(0)}} \biggr)
F(\Phi^{(01)} (\mathcal{A}\mathcal{M}^{(01)}, \mathcal{B}\mathcal{N}^{(01)}))
\biggr) \biggr|_e  ; 
\end{multline*}  
\begin{multline*}
\biggl[\biggl. \frac{\partial}{\partial a}\biggr|_e,
\biggl. \frac{\partial}{\partial A_j^{(0)}}\biggr|_e  \biggr] F\\
= \biggl. \biggl( \! \biggl( \frac{\partial}{\partial \alpha_0^\frac{1}{2}}
\frac{\partial}{\partial \mathcal{B}_j^{(0)}} - \frac{\partial}{\partial
\mathcal{A}_j^{(0)}} \frac{\partial}{\partial \alpha_0^\frac{1}{2}}
\biggr) F(\Phi^{(01)} (\mathcal{A}\mathcal{M}^{(01)}, 
\mathcal{B}\mathcal{N}^{(01)})) \biggr) \biggr|_e ;
\end{multline*}
when $j \geq k$, let
\begin{multline*}
\biggl[\biggl. \frac{\partial}{\partial A_j^{(1)}}\biggr|_e,
\biggl. \frac{\partial}{\partial A_k^{(0)}}\biggr|_e  \biggr] F\\
= \biggl. \biggl( \! \biggl( \frac{\partial}{\partial \mathcal{A}_j^{(1)}}
\frac{\partial}{\partial \mathcal{B}_k^{(0)}} - \frac{\partial}{\partial
\mathcal{B}_k^{(0)}} \frac{\partial}{\partial \mathcal{B}_j^{(1)}} \biggr)
F(\Phi^{(01)} (\mathcal{A}\mathcal{M}^{(01)}, \mathcal{B}\mathcal{N}^{(01)}))
\biggr) \biggr|_e;
\end{multline*}
when $j<k$, let 
\begin{multline*}
\biggl[\biggl. \frac{\partial}{\partial A_j^{(1)}}\biggr|_e,
\biggl. \frac{\partial}{\partial A_k^{(0)}}\biggr|_e  \biggr] F\\
= \biggl. \biggl( \! \biggl( \frac{\partial}{\partial \mathcal{A}_j^{(1)}}
\frac{\partial}{\partial \mathcal{B}_k^{(0)}} - \frac{\partial}{\partial
\mathcal{A}_k^{(0)}} \frac{\partial}{\partial \mathcal{A}_j^{(1)}} \biggr)
F(\Phi^{(01)} (\mathcal{A}\mathcal{M}^{(01)}, \mathcal{B}\mathcal{N}^{(01)}))
\biggr) \biggr|_e ; 
\end{multline*}
\begin{multline*}
\biggl[\Bigl. \frac{\partial}{\partial a}\Bigr|_e,
\biggl. \frac{\partial}{\partial A_j^{(1)}}\biggr|_e  \biggr] F\\
=\biggl. \biggl( \! \biggl( \frac{\partial}{\partial \alpha_0^\frac{1}{2}}
\frac{\partial}{\partial \mathcal{B}_j^{(1)}} - \frac{\partial}{\partial
\mathcal{A}_j^{(1)}} \frac{\partial}{\partial \beta_0^\frac{1}{2}}
\biggr) F(\Phi^{(01)} (\mathcal{A}\mathcal{M}^{(01)}, \mathcal{B}
\mathcal{N}^{(01)})) \biggr)
\biggr|_e ; 
\end{multline*}
\begin{multline*}
\biggl[\biggl. \frac{\partial}{\partial A_j^{(1)}}\biggr|_e,
\biggl. \frac{\partial}{\partial A_k^{(1)}}\biggr|_e  \biggr] F \\
= \biggl. \biggl( \! \biggl( \frac{\partial}{\partial \mathcal{A}_j^{(1)}}
\frac{\partial}{\partial \mathcal{B}_k^{(1)}} - \frac{\partial}{\partial
\mathcal{A}_k^{(1)}} \frac{\partial}{\partial \mathcal{B}_j^{(1)}} \biggr)
F(\Phi^{(01)} (\mathcal{A}\mathcal{M}^{(01)}, \mathcal{B}\mathcal{N}^{(01)}))
\biggr) \biggr|_e ;
\end{multline*}
\begin{multline*}
\Biggl[\Biggl. \frac{\partial}{\partial M_{j - \frac{1}{2}}^{(0)}}
\Biggr|_e, \Biggl. \frac{\partial}{\partial M_{k - \frac{1}{2}}^{(0)}}
\Biggr|_e  \Biggr] F \\
= \Biggl. \Biggl( \! \Biggl( \frac{\partial}{\partial \mathcal{M}_{j -
\frac{1}{2}}^{(0)}} \frac{\partial}{\partial \mathcal{N}_{k -
\frac{1}{2}}^{(0)}} + \frac{\partial}{\partial \mathcal{M}_{k -
\frac{1}{2}}^{(0)}} \frac{\partial}{\partial \mathcal{N}_{j -
\frac{1}{2}}^{(0)}} \Biggr) 
F(\Phi^{(01)} (\mathcal{A}\mathcal{M}^{(01)},
\mathcal{B}\mathcal{N}^{(01)})) \! \Biggr) \Biggr|_e ; 
\end{multline*} 
\begin{multline*}
\Biggl[\Biggl. \frac{\partial}{\partial M_{j - \frac{1}{2}}^{(1)}}
\Biggr|_e, \Biggl. \frac{\partial}{\partial M_{k - \frac{1}{2}}^{(0)}}
\Biggr|_e  \Biggr] F\\
= \Biggl. \Biggl( \! \Biggl( \frac{\partial}{\partial \mathcal{M}_{j -
\frac{1}{2}}^{(1)}} \frac{\partial}{\partial \mathcal{N}_{k -
\frac{1}{2}}^{(0)}} + \frac{\partial}{\partial \mathcal{M}_{k -
\frac{1}{2}}^{(0)}} \frac{\partial}{\partial \mathcal{M}_{j -
\frac{1}{2}}^{(1)}} \Biggr) F(\Phi^{(01)} (\mathcal{A}\mathcal{M}^{(01)},
\mathcal{B}\mathcal{N}^{(01)})) \! \Biggr) \Biggr|_e ;
\end{multline*}
\begin{multline*}
\Biggl[\Biggl. \frac{\partial}{\partial M_{j - \frac{1}{2}}^{(1)}}
\Biggr|_e, \Biggl. \frac{\partial}{\partial M_{k - \frac{1}{2}}^{(1)}}
\Biggr|_e  \Biggr] F\\
= \Biggl. \Biggl( \! \Biggl( \frac{\partial}{\partial \mathcal{M}_{j -
\frac{1}{2}}^{(1)}} \frac{\partial}{\partial \mathcal{N}_{k -
\frac{1}{2}}^{(1)}} + \frac{\partial}{\partial \mathcal{M}_{k -
\frac{1}{2}}^{(1)}} \frac{\partial}{\partial \mathcal{N}_{j -
\frac{1}{2}}^{(1)}} \Biggr) F(\Phi^{(01)} (\mathcal{A}\mathcal{M}^{(01)},
\mathcal{B}\mathcal{N}^{(01)})) \! \Biggr) \Biggr|_e ; 
\end{multline*}
\begin{multline*}
\Biggl[ \biggl. \frac{\partial}{\partial A_j^{(0)}}\biggr|_e,
\Biggl. \frac{\partial}{\partial M_{k - \frac{1}{2}}^{(0)}}\Biggr|_e
\Biggr] F \\
= \Biggl. \Biggl( \! \Biggl( \frac{\partial}{\partial \mathcal{A}_j^{(0)}}
\frac{\partial}{\partial \mathcal{N}_{k - \frac{1}{2}}^{(0)}} -
\frac{\partial}{\partial \mathcal{M}_{k- \frac{1}{2}}^{(0)}}
\frac{\partial}{\partial \mathcal{B}_j^{(0)}} \Biggr) F(\Phi^{(01)}
(\mathcal{A}\mathcal{M}^{(01)}, \mathcal{B}\mathcal{N}^{(01)})) \! \Biggr)
\Biggr|_e  ;
\end{multline*}
\begin{multline*}
\Biggl[\Bigl. \frac{\partial}{\partial a}\Bigr|_e,
\Biggl. \frac{\partial}{\partial M_{j - \frac{1}{2}}^{(0)}}\Biggr|_e 
\Biggr] F \\
= \Biggl. \Biggl( \! \Biggl( \frac{\partial}{\partial \alpha_0^\frac{1}{2}}
\frac{\partial}{\partial \mathcal{N}_{j - \frac{1}{2}}^{(0)}} -
\frac{\partial}{\partial \mathcal{M}_{j - \frac{1}{2}}^{(0)}}
\frac{\partial}{\partial \alpha_0^\frac{1}{2}} \Biggr) F(\Phi^{(01)}
(\mathcal{A}\mathcal{M}^{(01)}, \mathcal{B}\mathcal{N}^{(01)})) \! \Biggr)
\Biggr|_e ;
\end{multline*}
when $j \geq k$, let
\begin{multline*}
\Biggl[ \biggl. \frac{\partial}{\partial A_j^{(1)}}\biggr|_e,
\Biggl. \frac{\partial}{\partial M_{k - \frac{1}{2}}^{(0)}} \Biggr|_e
\Biggr] F \\
= \Biggl. \Biggl( \! \Biggl( \frac{\partial}{\partial \mathcal{A}_j^{(1)}}
\frac{\partial}{\partial \mathcal{N}_{k - \frac{1}{2}}^{(0)}} -
\frac{\partial}{\partial \mathcal{N}_{k - \frac{1}{2}}^{(0)}}
\frac{\partial}{\partial \mathcal{B}_j^{(1)}} \Biggr) F(\Phi^{(01)}
(\mathcal{A}\mathcal{M}^{(01)}, \mathcal{B}\mathcal{N}^{(01)})) \! \Biggr)
\Biggr|_e ;
\end{multline*}
when $j<k$, let
\begin{multline*}
\Biggl[\biggl. \frac{\partial}{\partial A_j^{(1)}}\biggr|_e,
\Biggl. \frac{\partial}{\partial M_{k - \frac{1}{2}}^{(0)}}\Biggr|_e
\Biggr] F\\
= \Biggl. \Biggl( \! \Biggl( \frac{\partial}{\partial \mathcal{A}_j^{(1)}}
\frac{\partial}{\partial \mathcal{N}_{k - \frac{1}{2}}^{(0)}} -
\frac{\partial}{\partial \mathcal{M}_{k - \frac{1}{2}}^{(0)}}
\frac{\partial}{\partial \mathcal{A}_j^{(1)}} \Biggr) F(\Phi^{(01)}
(\mathcal{A}\mathcal{M}^{(01)}, \mathcal{B}\mathcal{N}^{(01)}))\!\Biggr)
\Biggr|_e ;
\end{multline*} 
when $j \geq k$, let
\begin{multline*}
\Biggl[\Biggl. \frac{\partial}{\partial M_{j - \frac{1}{2}}^{(1)}}
\Biggr|_e, \biggl. \frac{\partial}{\partial A_k^{(0)}}\biggr|_e
\Biggr] F \\
= \Biggl. \Biggl( \! \Biggl( \frac{\partial}{\partial \mathcal{M}_{j -
\frac{1}{2}}^{(1)}} \frac{\partial}{\partial \mathcal{B}_k^{(0)}} -
\frac{\partial}{\partial \mathcal{B}_k^{(0)}} \frac{\partial}{\partial
\mathcal{N}_{j - \frac{1}{2}}^{(1)}} \Biggl) F(\Phi^{(01)} (\mathcal{A}
\mathcal{M}^{(01)}, \mathcal{B}\mathcal{N}^{(01)})) \! \Biggr) 
\Biggr|_e  ;
\end{multline*} 
when $j<k$, let 
\begin{multline*}
\Biggl[\Biggl. \frac{\partial}{\partial M_{j - \frac{1}{2}}^{(1)}}
\Biggr|_e, \biggl. \frac{\partial}{\partial A_k^{(0)}}\biggr|_e
\Biggr] F \\
= \Biggl. \Biggl( \! \Biggl( \frac{\partial}{\partial \mathcal{M}_{j -
\frac{1}{2}}^{(1)}} \frac{\partial}{\partial \mathcal{B}_k^{(0)}} -
\frac{\partial}{\partial \mathcal{A}_k^{(0)}} \frac{\partial}{\partial
\mathcal{M}_{j - \frac{1}{2}}^{(1)}} \Biggr) F(\Phi^{(01)} (\mathcal{A}
\mathcal{M}^{(01)}, \mathcal{B}\mathcal{N}^{(01)})) \!  \Biggr) 
\Biggr|_e  ; 
\end{multline*}
\begin{multline*}
\Biggl[\Bigl. \frac{\partial}{\partial a}\Bigr|_e,
\Biggl. \frac{\partial}{\partial M_{j - \frac{1}{2}}^{(1)}} \Biggr|_e
\Biggr] F \\
= \Biggl. \Biggl( \! \Biggl( \frac{\partial}{\partial \alpha_0^\frac{1}{2}}
\frac{\partial}{\partial \mathcal{N}_{j - \frac{1}{2}}^{(1)}} -
\frac{\partial}{\partial \mathcal{M}_{j - \frac{1}{2}}^{(1)}}
\frac{\partial}{\partial \beta_0^\frac{1}{2}} \Biggr) F(\Phi^{(01)}
(\mathcal{A}\mathcal{M}^{(01)}, \mathcal{B}\mathcal{N}^{(01)})) \! \Biggr)
\Biggr|_e ;
\end{multline*}  
\begin{multline*}
\Biggl[\biggl. \frac{\partial}{\partial A_j^{(1)}}\biggr|_e,
\Biggl. \frac{\partial}{\partial M_{k - \frac{1}{2}}^{(1)}} \Biggr|_e
\Biggr] F \\
= \Biggl. \Biggl( \! \Biggl( \frac{\partial}{\partial \mathcal{A}_j^{(1)}}
\frac{\partial}{\partial \mathcal{N}_{k - \frac{1}{2}}^{(1)}} -
\frac{\partial}{\partial \mathcal{M}_{k - \frac{1}{2}}^{(1)}}
\frac{\partial}{\partial \mathcal{B}_j^{(1)}} \Biggr) F(\Phi^{(01)}
(\mathcal{A}\mathcal{M}^{(01)}, \mathcal{B}\mathcal{N}^{(01)})) \! \Biggr)
\Biggr|_e  .
\end{multline*}

\begin{prop}\label{NS bracket}
The vector space $\hat{T}_e SK(1)$ with the bracket operation defined
above is the Neveu-Schwarz algebra with central charge zero.  The
basis is given by 
\begin{eqnarray}
\mathcal{L} (j) &=& \biggl. - \frac{\partial}{\partial A_{-j}^{(0)}}
\biggr|_e, \quad \mbox{for} \; \; - j \in \Z ,\\
\mathcal{L} (j) &=& \biggl. - \frac{\partial}{\partial A_j^{(1)}}
\biggr|_e, \quad \mbox{for} \; \; j \in \Z ,\\
\mathcal{L} (0) &=& \Bigl. - \frac{1}{2} \frac{\partial}{\partial a}
\Bigr|_e, \\ 
\mathcal{G} (j + \frac{1}{2}) &=& \Biggl. - \frac{\partial}{\partial M_{-j
- \frac{1}{2}}^{(0)}} \Biggr|_e, \quad \mbox{for} \; \; - j \in \Z ,\\
\mathcal{G} (j - \frac{1}{2}) &=& \Biggl. - \frac{\partial}{\partial M_{j
- \frac{1}{2}}^{(1)}} \Biggr|_e, \quad \mbox{for} \; \; j \in \Z .
\end{eqnarray}
\end{prop}

\begin{proof} We prove the bracket formula for
\[ \Bigl[ \mathcal{G} (j + \frac{1}{2}), \mathcal{G} (k - \frac{1}{2})
\Bigr] = 2 \mathcal{L} ( j + k) \]
for $- j \in \Z$, $k \in \Z$, and $-j < k$.  The proofs for the other
cases are similar.  To simplify notation, we will write, for example, 
\[\biggl. \frac{\partial}{\partial \Phi_k^{(0)}} F(\Phi^{(01)})
\biggr|_e \]
instead of 
\[\biggl. \frac{\partial}{\partial \Phi_k^{(0)}( \mathcal{AM}^{(01)},
\mathcal{BN}^{(01)})} F(\Phi^{(01)}(\mathcal{AM}^{(01)}, \mathcal{BN}^{(01)}))
\biggr|_{\begin{scriptsize}\begin{array}{l}
\mathcal{AM}^{(01)} = \\
\mathcal{BN}^{(01)} = e \end{array}\end{scriptsize}} .\] 
For $- j \in \Z$ and $k \in \Z$, 
\begin{eqnarray}\label{calculate brackets} 
& & \hspace{-.8in} \Bigl[ \mathcal{G} (j + \frac{1}{2}), \mathcal{G} (k - \frac{1}{2}) 
\Bigr] F  \nonumber \\
&=& \! \! \Biggl[ \Biggl. \frac{\partial}{\partial M_{-j-\frac{1}{2}}^{(0)}} 
\Biggr|_e, \Biggl. \frac{\partial}{\partial M_{k-\frac{1}{2}}^{(1)}}
\Biggr|_e \Biggr] F \nonumber \\
&=& \! \! \Biggl. \Biggl( \frac{\partial}{\partial \mathcal{M}_{k -
\frac{1}{2}}^{(1)}} \frac{\partial}{\partial \mathcal{N}_{-j -
\frac{1}{2}}^{(0)}} + \frac{\partial}{\partial \mathcal{M}_{-j -
\frac{1}{2}}^{(0)}} \frac{\partial}{\partial \mathcal{M}_{k -
\frac{1}{2}}^{(1)}} \Biggr) F(\Phi^{(01)}) \Biggr|_e \nonumber \\
&=& \! \! \! \sum_{n \in \Z} \Biggl( \Biggl. \biggl. \frac{\partial}{\partial
\mathcal{M}_{k - \frac{1}{2}}^{(1)}} \frac{\partial}{\partial 
\mathcal{N}_{-j - \frac{1}{2}}^{(0)}} \Phi_n^{(0)} \Biggr|_e
\frac{\partial}{\partial \Phi_n^{(0)}} F(\Phi^{(01)}) \biggr|_e
\Biggr. - \Biggl. \frac{\partial}{\partial \mathcal{N}_{-j -
\frac{1}{2}}^{(0)}} \Phi_n^{(0)} \Biggr|_e \cdot \\
& & \hspace{-.4in} \cdot \Biggl. \frac{\partial}{\partial \mathcal{M}_{k - \frac{1}{2}}^{(1)}}
\frac{\partial}{\partial \Phi_n^{(0)}} F(\Phi^{(01)}) \Biggr|_e \! 
 + \Biggl. \Biggl. \frac{\partial}{\partial \mathcal{M}_{k - \frac{1}{2}}^{(1)}}
\frac{\partial}{\partial \mathcal{N}_{-j -
\frac{1}{2}}^{(0)}} \Phi_{n -\frac{1}{2}}^{(0)} \Biggr|_e 
\frac{\partial}{\partial \Phi_{n - \frac{1}{2}}^{(0)}} F(\Phi^{(01)})
\Biggr|_e \nonumber \\
&  &\hspace{1.1in} \Biggl. - \; \Biggl. \frac{\partial}{\partial \mathcal{N}_{-j -
\frac{1}{2}}^{(0)}} \Phi_{n -\frac{1}{2}}^{(0)} \Biggr|_e
\Biggl. \frac{\partial}{\partial \mathcal{M}_{k - \frac{1}{2}}^{(1)}}
\frac{\partial}{\partial \Phi_{n - \frac{1}{2}}^{(0)}} F(\Phi^{(01)})
\Biggr|_e \Biggr) \nonumber 
\end{eqnarray}
\begin{multline*}
\qquad \quad \Biggl. \biggl. + \; \frac{\partial}{\partial \mathcal{M}_{k -
\frac{1}{2}}^{(1)}} \frac{\partial}{\partial \mathcal{N}_{-j -
\frac{1}{2}}^{(0)}} \lambda_0^\frac{1}{2} \Biggr|_e
\frac{\partial}{\partial \lambda_0^\frac{1}{2}} F(\Phi^{(01)})
\biggr|_e   \\
- \; \Biggl. \frac{\partial}{\partial \mathcal{N}_{-j -
\frac{1}{2}}^{(0)}} \lambda_0^\frac{1}{2} \Biggr|_e  \Biggl. 
\frac{\partial}{\partial \mathcal{M}_{k - \frac{1}{2}}^{(1)}}
\frac{\partial}{\partial \lambda_0^\frac{1}{2}} F(\Phi^{(01)})
\Biggr|_e 
\end{multline*}
\begin{multline*}
\qquad \quad + \sum_{n \in \Z} \Biggl( \biggl. \Biggl.
\frac{\partial}{\partial \mathcal{M}_{k - \frac{1}{2}}^{(1)}}
\frac{\partial}{\partial \mathcal{N}_{-j - \frac{1}{2}}^{(0)}}
\Phi_n^{(1)} \Biggr|_e \frac{\partial}{\partial \Phi_n^{(1)}}
F(\Phi^{(01)}) \biggr|_e \Biggr.  - \Biggl. \frac{\partial}{\partial 
\mathcal{N}_{-j - \frac{1}{2}}^{(0)}} \Phi_n^{(1)} \Biggr|_e \cdot \\
\cdot \Biggl. \frac{\partial}{\partial \mathcal{M}_{k - \frac{1}{2}}^{(1)}}
\frac{\partial}{\partial \Phi_n^{(1)}} F(\Phi^{(01)}) \Biggr|_e
+ \Biggl. \Biggl. \frac{\partial}{\partial \mathcal{M}_{k - \frac{1}{2}}^{(1)}}
\frac{\partial}{\partial \mathcal{N}_{-j -
\frac{1}{2}}^{(0)}} \Phi_{n -\frac{1}{2}}^{(1)} \Biggr|_e 
\frac{\partial}{\partial \Phi_{n - \frac{1}{2}}^{(1)}} F(\Phi^{(01)})
\Biggr|_e \\
\Biggl. - \; \Biggl. \frac{\partial}{\partial 
\mathcal{N}_{-j - \frac{1}{2}}^{(0)}} \Phi_{n -\frac{1}{2}}^{(1)} \Biggr|_e
\Biggl. \frac{\partial}{\partial \mathcal{M}_{k - \frac{1}{2}}^{(1)}}
\frac{\partial}{\partial \Phi_{n - \frac{1}{2}}^{(1)}} F(\Phi^{(01)})
\Biggr|_e \Biggr) 
\end{multline*}
\begin{multline*}
\qquad \quad + \sum_{n \in \Z} \Biggl( \Biggl. \biggl.
\frac{\partial}{\partial \mathcal{M}_{-j - \frac{1}{2}}^{(0)}}
\frac{\partial}{\partial \mathcal{M}_{k - \frac{1}{2}}^{(1)}}
\Phi_n^{(0)} \Biggr|_e \frac{\partial}{\partial \Phi_n^{(0)}}
F(\Phi^{(01)}) \biggr|_e \Biggr.  - \Biggl. \frac{\partial}{\partial 
\mathcal{M}_{k - \frac{1}{2}}^{(1)}} \Phi_n^{(0)} \Biggr|_e \cdot \\
\cdot \Biggl. \frac{\partial}{\partial \mathcal{M}_{-j - \frac{1}{2}}^{(0)}} 
\frac{\partial}{\partial \Phi_n^{(0)}} F(\Phi^{(01)}) \Biggr|_e
+ \Biggl. \Biggl. \frac{\partial}{\partial \mathcal{M}_{-j - \frac{1}{2}}^{(0)}}
\frac{\partial}{\partial \mathcal{M}_{k -
\frac{1}{2}}^{(1)}}  \Phi_{n -\frac{1}{2}}^{(0)} \Biggr|_e  
\frac{\partial}{\partial \Phi_{n - \frac{1}{2}}^{(0)}} F(\Phi^{(01)})
\Biggr|_e \\
\Biggl. - \; \Biggl. \frac{\partial}{\partial \mathcal{M}_{k -
\frac{1}{2}}^{(1)}} \Phi_{n -\frac{1}{2}}^{(0)} \Biggr|_e
\Biggl. \frac{\partial}{\partial \mathcal{M}_{-j - \frac{1}{2}}^{(0)}}
\frac{\partial}{\partial \Phi_{n - \frac{1}{2}}^{(0)}} F(\Phi^{(01)}) 
\Biggr|_e \Biggr) 
\end{multline*}
\begin{multline*}
\qquad \quad \Biggl. \biggl. + \; \frac{\partial}{\partial \mathcal{M}_{-j -
\frac{1}{2}}^{(0)}} \frac{\partial}{\partial \mathcal{M}_{k -
\frac{1}{2}}^{(1)}} \lambda_0^\frac{1}{2} \Biggr|_e
\frac{\partial}{\partial \lambda_0^\frac{1}{2}} F(\Phi^{(01)})
\biggr|_e \\
 - \; \Biggl. \frac{\partial}{\partial \mathcal{M}_{k -
\frac{1}{2}}^{(1)}} \lambda_0^\frac{1}{2} \Biggr|_e \Biggl. 
\frac{\partial}{\partial \mathcal{M}_{-j - \frac{1}{2}}^{(0)}}
\frac{\partial}{\partial \lambda_0^\frac{1}{2}} F(\Phi^{(01)})
\Biggr|_e 
\end{multline*} 
\begin{multline*}
\qquad \quad + \sum_{n \in \Z} \Biggl( \Biggl. \biggl.
\frac{\partial}{\partial \mathcal{M}_{-j - \frac{1}{2}}^{(0)}}
\frac{\partial}{\partial \mathcal{M}_{k - \frac{1}{2}}^{(1)}}
\Phi_n^{(1)} \Biggr|_e \frac{\partial}{\partial \Phi_n^{(1)}}
F(\Phi^{(01)}) \biggr|_e \Biggr. - \Biggl. \frac{\partial}{\partial 
\mathcal{M}_{k - \frac{1}{2}}^{(1)}} \Phi_n^{(1)} \Biggr|_e \cdot \\
\cdot \Biggl. \frac{\partial}{\partial \mathcal{M}_{-j -\frac{1}{2}}^{(0)}}
\frac{\partial}{\partial \Phi_n^{(1)}} F(\Phi^{(01)}) \Biggr|_e 
+  \Biggl. \Biggl. \frac{\partial}{\partial \mathcal{M}_{-j - \frac{1}{2}}^{(0)}}
\frac{\partial}{\partial \mathcal{M}_{k - \frac{1}{2}}^{(1)}} 
\Phi_{n -\frac{1}{2}}^{(1)} \Biggr|_e \frac{\partial}{\partial \Phi_{n - 
\frac{1}{2}}^{(1)}} F(\Phi^{(01)}) \Biggr|_e \\
\Biggl. - \; \Biggl. \frac{\partial}{\partial \mathcal{M}_{k -
\frac{1}{2}}^{(1)}} \Phi_{n -\frac{1}{2}}^{(1)} \Biggr|_e
\Biggl. \frac{\partial}{\partial \mathcal{M}_{-j - \frac{1}{2}}^{(0)}}
\frac{\partial}{\partial \Phi_{n - \frac{1}{2}}^{(1)}} F(\Phi^{(01)}) 
\Biggr|_e \Biggr) . 
\end{multline*}
Using \ref{Phi boundary conditions}, we have
\begin{eqnarray*}
\begin{array}{rclcrcl}
\Biggl. \frac{\partial}{\partial \mathcal{N}_{-j - \frac{1}{2}}^{(0)}}
\Phi_n^{(0)} \biggr|_e &=& 0, & \qquad & \Biggl. \frac{\partial}{\partial
\mathcal{M}_{k - \frac{1}{2}}^{(1)}} \Phi_n^{(0)} \biggr|_e &=& 0, \\
& & \\
\biggl. \frac{\partial}{\partial \mathcal{N}_{-j - \frac{1}{2}}^{(0)}}
\Phi_{n - \frac{1}{2}}^{(0)} \biggr|_e &=& \delta_{n,-j} , & \qquad &
\biggl. \frac{\partial}{\partial \mathcal{M}_{k - \frac{1}{2}}^{(1)}}
\Phi_{n - \frac{1}{2}}^{(0)} \biggr|_e &=& 0, \\ 
& & \\
\biggl. \frac{\partial}{\partial \mathcal{N}_{-j - \frac{1}{2}}^{(0)}}
\lambda_0^\frac{1}{2} \biggr|_e &=& 0, & \qquad &
\biggl. \frac{\partial}{\partial \mathcal{M}_{k - \frac{1}{2}}^{(1)}}
\lambda_0^\frac{1}{2} \biggr|_e &=& 0, \\ 
& & \\
\biggl. \frac{\partial}{\partial \mathcal{N}_{-j - \frac{1}{2}}^{(0)}}
\Phi_n^{(1)} \biggr|_e &=& 0, & \qquad &
\biggl. \frac{\partial}{\partial \mathcal{M}_{k - \frac{1}{2}}^{(1)}}
\Phi_n^{(1)} \biggr|_e &=& 0, \\ 
& & \\
\biggl. \frac{\partial}{\partial \mathcal{N}_{-j - \frac{1}{2}}^{(0)}}
\Phi_{n - \frac{1}{2}}^{(1)} \biggr|_e &=& 0 , & \qquad &
\biggl. \frac{\partial}{\partial \mathcal{M}_{k - \frac{1}{2}}^{(1)}}
\Phi_{n - \frac{1}{2}}^{(1)} \biggr|_e &=& \delta_{n,k} . \\ 
\end{array}
\end{eqnarray*}
Also, we have
\begin{eqnarray*}
\lefteqn{\Biggl. \frac{\partial}{\partial \mathcal{M}_{k -
\frac{1}{2}}^{(1)}} \frac{\partial}{\partial \Phi_{n -
\frac{1}{2}}^{(0)}} F(\Phi^{(01)}) \Biggr|_e}\\
&=& \sum_{m \in \Z} \Biggl( \Biggl. \frac{\partial}{\partial \mathcal{M}_{k
- \frac{1}{2}}^{(1)}} \Phi_m^{(0)} \Biggr|_e
\Biggl. \frac{\partial}{\partial \Phi_m^{(0)}}   
\frac{\partial}{\partial \Phi_{n - \frac{1}{2}}^{(0)}} F(\Phi^{(01)})
\Biggr|_e \Biggr. \\
& & \hspace{1.6in} + \; \Biggl. \Biggl. \frac{\partial}{\partial \mathcal{M}_{k -
\frac{1}{2}}^{(1)}} \Phi_{m - \frac{1}{2}}^{(0)} \Biggr|_e
\Biggl. \frac{\partial}{\partial \Phi_{m - \frac{1}{2}}^{(0)}}
\frac{\partial}{\partial \Phi_{n - \frac{1}{2}}^{(0)}} F(\Phi^{(01)}) 
\Biggr|_e \Biggr) \\ 
& & + \; \Biggl. \frac{\partial}{\partial \mathcal{M}_{k
- \frac{1}{2}}^{(1)}} \lambda_0^\frac{1}{2} \Biggr|_e
\Biggl. \frac{\partial}{\partial \lambda_0^\frac{1}{2}}  
\frac{\partial}{\partial \Phi_{n - \frac{1}{2}}^{(0)}} F(\Phi^{(01)})
\Biggr|_e \\
& & + \sum_{m \in \Z} \Biggl( \Biggl. \frac{\partial}{\partial
\mathcal{M}_{k - \frac{1}{2}}^{(1)}} \Phi_m^{(1)} \Biggr|_e
\Biggl. \frac{\partial}{\partial \Phi_m^{(1)}}  
\frac{\partial}{\partial \Phi_{n - \frac{1}{2}}^{(0)}} F(\Phi^{(01)})
\Biggr|_e \Biggr. \\
& & \hspace{1.6in} + \; \Biggl. \Biggl. \frac{\partial}{\partial \mathcal{M}_{k -
\frac{1}{2}}^{(1)}} \Phi_{m - \frac{1}{2}}^{(1)} \Biggr|_e
\Biggl. \frac{\partial}{\partial \Phi_{m - \frac{1}{2}}^{(1)}}
\frac{\partial}{\partial \Phi_{n - \frac{1}{2}}^{(0)}} F(\Phi^{(01)})
\Biggr|_e \Biggr) \\ 
&=& \Biggl. \frac{\partial}{\partial \Phi_{k - \frac{1}{2}}^{(1)}}
\frac{\partial}{\partial \Phi_{n - \frac{1}{2}}^{(0)}} F(\Phi^{(01)})
\Biggr|_e ,
\end{eqnarray*}
and similarly
\[\Biggl. \frac{\partial}{\partial \mathcal{M}_{n - \frac{1}{2}}^{(0)}}
\frac{\partial}{\partial \Phi_{k - \frac{1}{2}}^{(1)}} F(\Phi^{(01)}) 
\Biggr|_e = \Biggl. \frac{\partial}{\partial \Phi_{n -
\frac{1}{2}}^{(0)}} \frac{\partial}{\partial \Phi_{k -
\frac{1}{2}}^{(1)}} F(\Phi^{(01)} \Biggr|_e  .\]

For simplicity we will write $\Psi_n  = \Psi_n (\alpha_0^{1/2},
\mathcal{A}^{(1)}, \mathcal{M}^{(1)}, \mathcal{B}^{(0)}, \mathcal{N}^{(0)})$, for
$n \in \frac{1}{2} \mathbb{Z}$.  Let 
\begin{eqnarray*}
H_{\mathcal{A}^{(0)}, \mathcal{M}^{(0)}}^{(1)} (x, \varphi) &=& \tilde{E}
(\mathcal{A}^{(0)}, i\mathcal{M}^{(0)}) \circ I(x,\varphi) \\
H_{\beta_0^{1/2}, \mathcal{B}^{(1)}, \mathcal{N}^{(1)}}^{(2)}
(x,\varphi) &=& \hat{E} (\beta_0^\frac{1}{2}, \mathcal{B}^{(1)}, 
\mathcal{N}^{(1)}) (x,\varphi) \\  
F_{\alpha_0^{1/2}, \mathcal{A}^{(1)}, \mathcal{M}^{(1)}, \mathcal{B}^{(0)},
\mathcal{N}^{(0)}}^{(1)} (x,\varphi) &=& I \circ \tilde{E}(\{ \Psi_{-n},
i \Psi_{-n + \frac{1}{2}} \}_{n \in \Z} ) \circ I (x, \varphi) \\  
F_{\alpha_0^{1/2}, \mathcal{A}^{(1)}, \mathcal{M}^{(1)}, \mathcal{B}^{(0)},
\mathcal{N}^{(0)}}^{(2)} (x,\varphi) &=& \hat{E} ( \exp (\Psi_0 )
\alpha_0^{-\frac{1}{2}}, \{\exp(2n \Psi_0) \alpha_0^{-n} \Psi_n, \\ 
& & \qquad \exp(2(n - \frac{1}{2}) \Psi_0 ) \alpha_0^{-n +
\frac{1}{2}} \Psi_{n - \frac{1}{2}} \}_{n \in \Z} ) (x,\varphi)  
\end{eqnarray*}
Then {}from the definition of $\Phi^{(01)} (\mathcal{AM}^{(01)}, 
\mathcal{BN}^{(01)})$, we have
\begin{eqnarray*}
\varphi \Phi_n^{(0)} &=& \varphi \mbox{Res}_x x^n (H_{\mathcal{A}^{(0)}, 
\mathcal{M}^{(0)}}^{(1)} \circ (F_{\alpha_0^{1/2}, \mathcal{A}^{(1)}, \mathcal{M}^{(1)},
\mathcal{B}^{(0)}, \mathcal{N}^{(0)}}^{(1)})^{-1} (x,\varphi))^0 \\  
\varphi \Phi_{n - \frac{1}{2}}^{(0)} &=& i \varphi \mbox{Res}_x x^{n
-1} (H_{\mathcal{A}^{(0)}, \mathcal{M}^{(0)}}^{(1)} \circ (F_{\alpha_0^{1/2},
\mathcal{A}^{(1)}, \mathcal{M}^{(1)}, \mathcal{B}^{(0)}, 
\mathcal{N}^{(0)}}^{(1)})^{-1} (x,\varphi))^1 \\ 
\varphi \lambda_0^\frac{1}{2} &=& \varphi \mbox{Res}_x x^{-2}
(H_{\beta_0^{1/2}, \mathcal{B}^{(1)}, \mathcal{N}^{(1)}}^{(2)} \circ
(F_{\alpha_0^{1/2}, \mathcal{A}^{(1)}, \mathcal{M}^{(1)}, \mathcal{B}^{(0)},
\mathcal{N}^{(0)}}^{(2)})^{-1} (x,\varphi))^0 \\
\varphi \Phi_n^{(1)} &=& \varphi \mbox{Res}_x x^{-n-2}
(H_{\beta_0^{1/2}, \mathcal{B}^{(1)}, \mathcal{N}^{(1)}}^{(2)} \circ 
(F_{\alpha_0^{1/2}, \mathcal{A}^{(1)}, \mathcal{M}^{(1)}, \mathcal{B}^{(0)},
\mathcal{N}^{(0)}}^{(2)})^{-1} (x,\varphi))^0 \\
\varphi \Phi_{n - \frac{1}{2}}^{(1)} &=& \varphi \mbox{Res}_x x^{-n-1}
(H_{\beta_0^{1/2}, \mathcal{B}^{(1)}, \mathcal{N}^{(1)}}^{(2)} \circ 
(F_{\alpha_0^{1/2}, \mathcal{A}^{(1)}, \mathcal{M}^{(1)}, \mathcal{B}^{(0)},
\mathcal{N}^{(0)}}^{(2)})^{-1} (x,\varphi))^1 . 
\end{eqnarray*}
Thus using (\ref{specific F1}), (\ref{specific F2}), (\ref{psi
condition 1}), (\ref{psi condition 2}), and (\ref{psi condition 3}),
we have 
\begin{eqnarray*}
\begin{array}{rclcrcl}
\biggl. \frac{\partial}{\partial \mathcal{M}_{k - \frac{1}{2}}^{(1)}}
\frac{\partial}{\partial \mathcal{N}_{-j - \frac{1}{2}}^{(0)}}
\Phi_n^{(0)} \biggr|_e \! &=& \! - 2 \delta_{- n, j + k} & \! \quad & \biggl. 
\frac{\partial}{\partial \mathcal{M}_{-j - \frac{1}{2}}^{(0)}}
\frac{\partial}{\partial \mathcal{M}_{k - \frac{1}{2}}^{(1)}}
\Phi_n^{(0)} \biggr|_e \! &=& \! 0 \\
& & \\
\biggl. \frac{\partial}{\partial \mathcal{M}_{k - \frac{1}{2}}^{(1)}}
\frac{\partial}{\partial \mathcal{N}_{-j - \frac{1}{2}}^{(0)}}
\Phi_{n - \frac{1}{2}}^{(0)} \biggr|_e \! &=& \! 0 & \! \quad & \biggl. 
\frac{\partial}{\partial \mathcal{M}_{-j - \frac{1}{2}}^{(0)}} 
\frac{\partial}{\partial \mathcal{M}_{k - \frac{1}{2}}^{(1)}} \Phi_{n -
\frac{1}{2}}^{(0)} \biggr|_e \! &=& \! 0\\
& & \\
\biggl. \frac{\partial}{\partial \mathcal{M}_{k - \frac{1}{2}}^{(1)}}
\frac{\partial}{\partial \mathcal{N}_{-j - \frac{1}{2}}^{(0)}}
\lambda_0^\frac{1}{2} \biggr|_e \! &=& \! - 2 \delta_{0, j + k} & \! \quad &
\biggl. \frac{\partial}{\partial \mathcal{M}_{-j - \frac{1}{2}}^{(0)}} 
\frac{\partial}{\partial \mathcal{M}_{k - \frac{1}{2}}^{(1)}}
\lambda_0^\frac{1}{2} \biggr|_e \! &=& \! 0 \\
& & \\
\biggl. \frac{\partial}{\partial \mathcal{M}_{k - \frac{1}{2}}^{(1)}}
\frac{\partial}{\partial \mathcal{N}_{-j - \frac{1}{2}}^{(0)}}
\Phi_n^{(1)} \biggr|_e \! &=& \! - 2 \delta_{n, j + k} & \! \quad & \biggl. 
\frac{\partial}{\partial \mathcal{M}_{-j - \frac{1}{2}}^{(0)}} 
\frac{\partial}{\partial \mathcal{M}_{k - \frac{1}{2}}^{(1)}}
\Phi_n^{(1)} \biggr|_e \! &=& \! 0 \\
& & \\
\biggl. \frac{\partial}{\partial \mathcal{M}_{k - \frac{1}{2}}^{(1)}}
\frac{\partial}{\partial \mathcal{N}_{-j - \frac{1}{2}}^{(0)}}
\Phi_{n - \frac{1}{2}}^{(1)} \biggr|_e \! &=& \! 0 & \! \quad & \biggl. 
\frac{\partial}{\partial \mathcal{M}_{-j - \frac{1}{2}}^{(0)}} 
\frac{\partial}{\partial \mathcal{M}_{k - \frac{1}{2}}^{(1)}} \Phi_{n -
\frac{1}{2}}^{(1)} \biggr|_e \! &=& \! 0 . \\
\end{array}
\end{eqnarray*}
Substituting these calculations into (\ref{calculate brackets}), we
have
\begin{eqnarray*}
\lefteqn{\Bigl[ \mathcal{G} (j + \frac{1}{2}), \mathcal{G} (k - \frac{1}{2}) 
\Bigr] F} \\
&=& \sum_{n \in \Z} \Biggl( - 2 \delta_{-n, j + k} \biggl. 
\frac{\partial}{\partial \Phi_n^{(0)}} F(\Phi^{(01)}) \biggr|_e +
\delta_{n,-j} \Biggl. \frac{\partial}{\partial
\Phi_{k-\frac{1}{2}}^{(1)}} \frac{\partial}{\partial \Phi_{n -
\frac{1}{2}}^{(0)}} F(\Phi^{(01)}) \Biggr|_e  \Biggr) \\
& &  - \; 2 \delta_{0, j + k} \left. \frac{\partial}{\partial
\lambda_0^\frac{1}{2}} F(\Phi^{(01)}) \right|_e \\
& &  + \sum_{n \in \Z} \Biggl( - 2 \delta_{n, j + k}
\biggl. \frac{\partial}{\partial \Phi_n^{(1)}} F(\Phi^{(01)}) \biggr|_e
+ \delta_{n,k} \Biggl. \frac{\partial}{\partial \Phi_{-j -
\frac{1}{2}}^{(0)}} \frac{\partial}{\partial \Phi_{n
-\frac{1}{2}}^{(1)}} F(\Phi^{(01)}) \Biggr|_e \Biggr) \\
&=& \left\{ \begin{array}{lll}
- 2 \Bigl. \frac{\partial}{\partial \Phi_{-j - k}^{(0)}} F(\Phi^{(01)})
\Bigr|_e & \mbox{if} & k < -j \\
-2 \Bigl. \frac{\partial}{\partial \lambda_0^\frac{1}{2}} F(\Phi^{(01)})
\Bigr|_e & \mbox{if} & k = -j \\
-2 \Bigl. \frac{\partial}{\partial \Phi_{j + k}^{(1)}} F(\Phi^{(01)})
\Bigr|_e  & \mbox{if} & k > -j 
\end{array} \right. \\
& & \\
&=& 2 \mathcal{L}(j + k) .
\end{eqnarray*}
\end{proof}

\backmatter

\end{document}